\documentclass[11pt,a4]{article}
\usepackage[T2A]{fontenc}
\usepackage[cp866]{inputenc}
\usepackage[mathscr]{eucal}
\usepackage{amssymb}
\usepackage{amsmath}
\usepackage{graphicx}
\usepackage[russian, english]{babel}  % затем убрать
\usepackage{array}  % добавили новый пакет для таблицы

\renewcommand{\leq}{\leqslant}
\renewcommand{\geq}{\geqslant}

\newcommand{\R}{{\mathbb R}}

\newcommand{\K}{{\mathbb K}}
\renewcommand{\P}{{\mathbb P}}

\renewcommand{\k}{\rule{0.7em}{0.7em}}
\usepackage{fancyhdr}
\begin{document}
\sloppy

{\normalsize

\thispagestyle{empty}

\mbox{}
\\[-2.25ex]
\centerline{
{\large
\bf
PROJECTIVE ATLAS OF TRAJECTORIES  
}
}
\\[0.75ex]
\centerline{
{\large
\bf
OF DIFFERENTIAL SYSTEMS 
}
}
\\[2.5ex]
\centerline{
\bf 
V.N. Gorbuzov
}
\\[2ex]
\centerline{
\it 
Yanka Kupala Grodno State University
{\rm(}Ozeshko 22, Grodno, 230023, Belarus{\rm)}
}
\\[1.5ex]
\centerline{
E-mail: gorbuzov@grsu.by
}
\\[5.5ex]
\centerline{{\large\bf Abstract}}
\\[1ex]
\indent
Topological bases of behaviour of trajectories for autonomous differential systems of the second order 
on the projective phase plane are stated. 
By means of Poincar\'{e}'s circles the projective atlas of trajectories is constructed. 
Differential connections between trajectories of projectively conjugated differential systems are established. 
The behaviour of trajectories in an neighbourhood of infinitely remote straight line of the projective phase plane and 
property of closed trajectories on the projective phase plane is investigated. 
Examples of full qualitative research of trajectories of differential systems on the projective phase plane are given.
\\[1.5ex]
\\[1.25ex]
%34A26 - Geometric methods in differential equations
%34C05 - Location of integral curves, singular points, limit cycles
\\[7.5ex]
\centerline{{\large\bf Contents}}
\\[1.5ex]
{\bf  Introduction}                   \dotfill\ 2
\\[1ex]
{\bf \S 1. Poincar\'{e}'s sphere}
                                                 \dotfill \ 2
\\[0.5ex]
\mbox{}\hspace{1.35em}
1. Map of plane to Poincar\'{e}'s sphere 
                                                 \dotfill \ 2
\\[0.5ex]
\mbox{}\hspace{1.35em}
2. Atlas of  Poincar\'{e}'s sphere 
                                                 \dotfill \ 4
\\[0.5ex]
\mbox{}\hspace{1.35em}
3. Poincar\'{e}'s circle 
                                                 \dotfill \ 5
\\[0.5ex]
\mbox{}\hspace{1.35em}
4. Maps of Poincar\'{e}
                                                 \dotfill \ 6
\\[0.5ex]
\mbox{}\hspace{1.35em}
5. Atlas of projective circles for projective plane
                                                 \dotfill \ 9
\\[1ex]
\noindent
{\bf \S 2. 
Transformations of Poincar\'{e} for differential systems
}
                                                 \dotfill \ 10
\\[0.75ex]
\mbox{}\hspace{1.35em}
6. Projectively reduced systems
                                                 \dotfill \ 10
\\[0.5ex]
\mbox{}\hspace{1.35em}
7. Projective type of differential system 
                                                 \dotfill \ 11
\\[0.5ex]
\mbox{}\hspace{1.35em}
8. Projective atlas of trajectories for differential systems
                                                 \dotfill \ 16
\\[1ex]
\noindent
{\bf \S 3. Trajectories on Poincar\'{e}'s sphere
}
                                                 \dotfill \ 19
\\[0.75ex]
\mbox{}\hspace{1.35em}
9.\ \, Trajectories on projective phase plane
                                                 \dotfill \ 19
\\[0.5ex]
\mbox{}\hspace{1.35em}
10. Trajectories of the first projectively reduced system
                                                 \dotfill \ 23
\\[0.5ex]
\mbox{}\hspace{1.35em}
11. Trajectories of the second projectively reduced system
                                                 \dotfill \ 25
\\[0.5ex]
\mbox{}\hspace{1.35em}
12.  Linear and open limit cycles
                                                 \dotfill \ 26
\\[0.5ex]
\mbox{}\hspace{1.35em}
13. Symmetry of the phase field of directions
                                                 \dotfill \ 33
\\[0.5ex]
\mbox{}\hspace{1.35em}
14. Sets of projectively nonsingular and projectively singular systems
                                                 \dotfill \ 35
\\[0.5ex]
\mbox{}\hspace{1.35em}
15. Topological equivalence of differential systems 
\\[0.25ex]
\mbox{}\hspace{3.1em}
on projective circle and on projective sphere  
                                                 \dotfill \ 36
\\[0.5ex]
\mbox{}\hspace{1.35em}
16. Examples of global qualitative research of trajectories
\\[0.25ex]
\mbox{}\hspace{3.1em}
for differential systems on projective phase plane
                                                 \dotfill \ 38
\\[1ex]
{\bf References}
                                              \dotfill \ 60

\newpage

\mbox{}
\\[-1.75ex]
\centerline{\large\bf  Introduction}
\\[1.5ex]
\indent
Research object is ordinary autonomous polynomial differential system of the second order
\\[2ex]
\mbox{}\hfill                                   %(D)
$
\displaystyle 
\dfrac{dx}{dt} =
\sum \limits_{k=0}^{n}\, X_k^{}(x,y)\equiv
 X(x,y), 
\quad 
\dfrac{dy}{dt} =
\sum \limits_{k=0}^{n}\,Y_k^{}(x,y)\equiv
Y(x,y),
$
\hfill (D)
\\[2.25ex]
where $X_k^{}$ and $Y_k^{}$ 
\vspace{0.5ex}
are homogeneous polynomials of degree $k,\ k=0,1,\ldots,n,$ on variables $x$ and $y$
such that 
\vspace{0.35ex}
$|X_n^{}(x,y)|+|Y_n^{}(x,y)|\not\equiv0$ on $\R^2,$ and polynomials $X$ and $Y$ are relatively prime, 
\vspace{0.25ex}
i.e. they have no the common divisors which are distinct from real numbers.

For the purpose of study of behaviour of trajectories of system (D) 
H. Poincar\'{e} along with equilibrium states lying in the final part of the phase plane 
also investigated infinitely removed equilibrium states [1, p. 23 -- 31]. 
For this purpose the phase plane $(x,y)$ has been replenished by infinitely removed points 
and trajectories of system (D) has been projected on the sphere with their subsequent image on the circle [1, p. 84 -- 91]. 
Thereby the beginning of the global qualitative theory of autonomous ordinary differential systems of the second order 
was given [2 -- 4].

In this paper taking into account researches [5] the behaviour of trajectories of the polynomial differential system (D) 
on the projective phase plane $\R\P(x,y)$ is considered. 
We recognise that Poincar\'{e}'s sphere (two-dimensional unit sphere with the identified antipodal points [6, p. 749]) 
is diffeomorphic to the projective plane, 
and the projective plane is two-dimensional variety which is analytically described by means 
of three local rectangular Cartesian coordinate systems 
[7, p. 96; 8, p. 421 -- 423].
\\[5ex]
\centerline{
{\bf\large \S\;\!1. Poincar\'{e}'s sphere}}
\\[2ex]
\centerline{
{\bf  1. 
Map of plane to Poincar\'{e}'s sphere
}
}
\\[1.5ex]
\indent
Let's introduce three-dimensional rectangular Cartesian coordinate system 
\vspace{0.15ex}
$O ^ {\star} x ^ {\star} y ^ {\star} z ^ {\star}, $ 
combined with the right rectangular Cartesian coordinate system $Oxy, $ 
\vspace{0.15ex}
meeting conditions: 
the straight line $OO ^ {\star} $ is orthogonal to plane $Oxy,$ 
\vspace{0.15ex}
the length of the segment $OO ^ {\star} $ is equal to one unit of the scale of the system of coordinate $Oxy;$ 
\vspace{0.15ex}
the axis $O ^ {\star} x ^ {\star} $ is codirected with the axis $Ox,$ 
the axis $O ^ {\star} y ^ {\star} $ is codirected with the axis $Oy,$ 
\vspace{0.15ex}
and the axis $O ^ {\star} z ^ {\star} $ is directed so that the system of coordinate $O ^ {\star} x ^ {\star} y ^ {\star} z ^ {\star}$ will be right; 
\vspace{0.15ex}
a scale in the system of coordinate $O ^ {\star} x ^ {\star} y ^ {\star} z ^ {\star} $ same, as well as in the system of coordinate $Oxy.$
\vspace{0.25ex}
We will construct the sphere with the centre $O^ {\star}$ of unit radius: 
\\[2.25ex]
\mbox{}\hfill                                   %(1.1)
$
S^2=\bigl\{
(x^{\star},y^{\star},z^{\star})
\colon
x^{\star}{}^{\,^{\scriptstyle 2}}+
y^{\star}{}^{\,^{\scriptstyle 2}}+
z^{\star}{}^{\,^{\scriptstyle 2}}=1\bigr\}.
$
\hfill (1.1)
\\[2.25ex]
\noindent
Points $N (0,0,1) $ and $S (0,0, {}-1) $ 
\vspace{0.25ex}
are according to northern and southern poles of this sphere. 
\vspace{0.25ex}
Thus the south pole $\!S(0,0, {}-1)\!$ coincides with the beginning $\!O(0,0)\!$ of the system of coordinate $Oxy. $
\vspace{0.25ex}
The equation $\!z ^ {\star}\! = -\;\!1\!$ is the equation in the system of coordinate 
$\!O ^ {\star} x ^ {\star} y ^ {\star} z ^ {\star}\!$ of the plane $\!Oxy.\!$
The plane $Oxy$ 
\vspace{0.35ex}
concerns by sphere (1.1) in the south pole $S (0,0, {}-1). $

On the plane $Oxy $ 
\vspace{0.15ex}
arbitrarily we will choose the point $M $ and we will spend the ray with the beginning 
$M$ through the centre $O ^ {\star} $ of the sphere (1.1). 
\vspace{0.15ex}
The ray $MO ^ {\star} $ intersects the sphere in two antipodal (diametrically opposite) points 
$M ^ {\star} $ and $M _ {\star} ^ {} $ (Fig. 1.1).
\vspace{0.15ex}
Thereby to each point of the plane $Oxy $ there correspond two antipodal points of the sphere (1.1).
\vspace{0.15ex}
And on the contrary to any pair of antipodal points of the sphere (1.1), 
\vspace{0.15ex}
except for only points lying on the equator of the sphere 
(circles of the sphere lying in the coordinate plane $O ^ {\star} x ^ {\star} y ^ {\star}), $ one point of the plane $Oxy$ is compared. 
\vspace{0.15ex}

To spread correspondence to all sphere (1.1), we will arrive as follows.
Each straight line $l _ {_ 0}, $ passing through the point $O $ 
and lying in the coordinate plane $Oxy, $ we will supplement 
by the infinitely removed point $L,$ 
lying on its <<extremities>>. 
Thus to the different straight lines passing through the point $O $ and lying in the plane $Oxy, $ 
there correspond different infinitely removed points, and to a bundle of parallel straight lines there corresponds 
one infinitely removed point.
\vspace{0.25ex}
Straight line $l, $ continuously continued 
by the infinitely removed point $L,$ we will designate through $l_{_L}$ or $ \overline {l}. $
\vspace{0.25ex}
The plane $(x, y), $ added by the infinitely removed points of all straight lines lying on it, 
\vspace{0.35ex}
is [9, p. 247] the projective plane which we will designate as $ \R\P^2$ or $ \R \P (x, y). $
\vspace{0.25ex}
The set of infinitely removed points of the projective plane is [9, p. 245] infinitely removed straight line of this projective plane. 
\vspace{0.25ex}
The projective plane $ \R \P (x, y) $ is disjunctive association of the plane $ (x, y) $ and infinitely removed straight line.
\\[3.75ex]
\mbox{}\hfill
{\unitlength=1mm
\begin{picture}(77,50)
\put(0,0){\includegraphics[width=96.7mm,height=49.84mm]{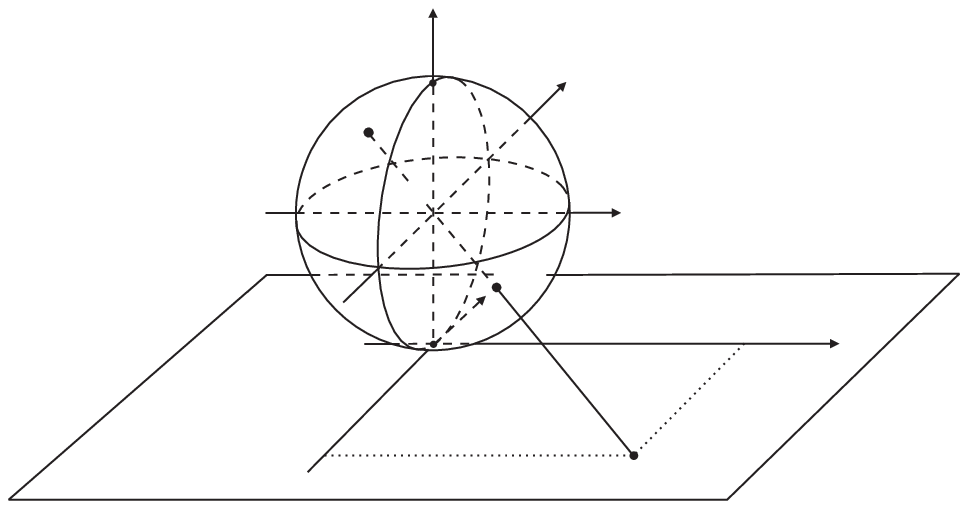}}
 
\put(44.9,13.7){\makebox(0,0)[cc]{\scriptsize $O$}}
\put(41.5,17.7){\makebox(0,0)[cc]{\scriptsize $S$}}
\put(75,17.5){\makebox(0,0)[cc]{\scriptsize $x$}}
\put(83,17.5){\makebox(0,0)[cc]{\scriptsize $x$}}
\put(30,5){\makebox(0,0)[cc]{\scriptsize $y$}}
\put(49,19){\makebox(0,0)[cc]{\scriptsize $y$}}
%\put(41.9,3.6){\makebox(0,0)[cc]{\scriptsize $(x,y)$}}
\put(66,4){\makebox(0,0)[cc]{\scriptsize $M$}}

\put(40.9,30.7){\makebox(0,0)[cc]{\scriptsize $O^{\star}$}}
\put(52,22.8){\makebox(0,0)[cc]{\scriptsize $M^{\star}$}}
\put(34.2,36){\makebox(0,0)[cc]{\scriptsize $M_{\star}^{}$}}
%\put(32.5,27){\makebox(0,0)[cc]{\scriptsize $L^{\star}_{\scriptscriptstyle-}$}}
%\put(52.8,31.3){\makebox(0,0)[cc]{\scriptsize $L^{\star}_{\scriptscriptstyle+}$}}
%\put(9,2.7){\makebox(0,0)[cc]{\scriptsize $L^{}_{\scriptscriptstyle-}$}}
%\put(65,19.7){\makebox(0,0)[cc]{\scriptsize $L^{}_{\scriptscriptstyle+}$}}
%\put(17,10){\makebox(0,0)[cc]{\scriptsize $l^{}_{\scriptscriptstyle 0}$}}
\put(61.6,27.3){\makebox(0,0)[cc]{\scriptsize $x^{\star}$}}
\put(58,30.5){\makebox(0,0)[cc]{\scriptsize $1$}}
\put(27,30.6){\makebox(0,0)[cc]{\scriptsize ${\scriptstyle-}1$}}
\put(57.3,40){\makebox(0,0)[cc]{\scriptsize $y^{\star}$}}
\put(49.6,37.5){\makebox(0,0)[cc]{\scriptsize $1$}}
\put(35.3,25.3){\makebox(0,0)[cc]{\scriptsize ${\scriptstyle-}1$}}

\put(41.3,44.3){\makebox(0,0)[cc]{\scriptsize $N$}}
\put(44.5,44.3){\makebox(0,0)[cc]{\scriptsize $1$}}
\put(45.5,49.5){\makebox(0,0)[cc]{\scriptsize $z^{\star}$}}

\put(48,-5){\makebox(0,0)[cc]{\rm Fig. 1.1}}
\end{picture}}
\hfill\mbox{}
\\[9.75ex]
\mbox{}\hfill
{\unitlength=1mm
\begin{picture}(78.9,50)
\put(3,0){\includegraphics[width=75.9mm,height=49.84mm]{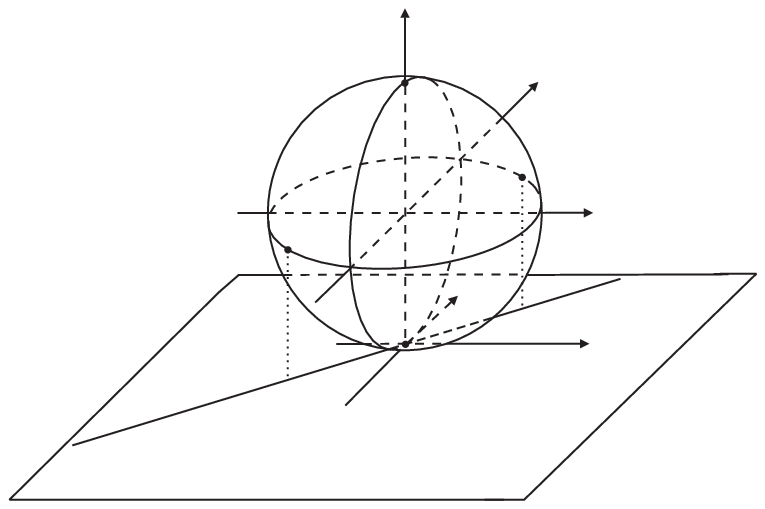}}
 
\put(44.9,13.7){\makebox(0,0)[cc]{\scriptsize $O$}}
\put(41.7,17.7){\makebox(0,0)[cc]{\scriptsize $S$}}
\put(60,14){\makebox(0,0)[cc]{\scriptsize $x$}}
\put(49,19){\makebox(0,0)[cc]{\scriptsize $y$}}
%\put(41.9,3.6){\makebox(0,0)[cc]{\scriptsize $(x,y)$}}
%\put(66,4){\makebox(0,0)[cc]{\scriptsize $M$}}

\put(40.9,30.7){\makebox(0,0)[cc]{\scriptsize $O^{\star}$}}
%\put(52,22.8){\makebox(0,0)[cc]{\scriptsize $M^{\star}$}}
%\put(37.6,38.9){\makebox(0,0)[cc]{\scriptsize $M_{\star}^{}$}}
\put(32.5,27){\makebox(0,0)[cc]{\scriptsize $L_{\star}$}}
\put(52.8,31.3){\makebox(0,0)[cc]{\scriptsize $L^{\star}$}}
\put(12,4){\makebox(0,0)[cc]{\scriptsize $L$}}
\put(65,20){\makebox(0,0)[cc]{\scriptsize $L$}}
\put(23,12){\makebox(0,0)[cc]{\scriptsize $l_{{}_0}$}}
\put(61.6,27.3){\makebox(0,0)[cc]{\scriptsize $x^{\star}$}}
\put(58,30.5){\makebox(0,0)[cc]{\scriptsize $1$}}
\put(27,30.6){\makebox(0,0)[cc]{\scriptsize ${\scriptstyle-}1$}}
\put(57.3,40){\makebox(0,0)[cc]{\scriptsize $y^{\star}$}}
\put(49.6,37.5){\makebox(0,0)[cc]{\scriptsize $1$}}
\put(35.3,25.3){\makebox(0,0)[cc]{\scriptsize ${\scriptstyle-}1$}}

\put(41.3,44.3){\makebox(0,0)[cc]{\scriptsize $N$}}
\put(44.5,44.3){\makebox(0,0)[cc]{\scriptsize $1$}}
\put(45.5,49.5){\makebox(0,0)[cc]{\scriptsize $z^{\star}$}}

\put(39,-5){\makebox(0,0)[cc]{\rm Fig. 1.2}}
\end{picture}}
\hfill\mbox{}
\\[7ex]
\indent
On the coordinate plane $Oxy $ arbitrarily we will choose the straight line 
\vspace{0.25ex}
$l _ {_ 0} \, $ passing through the origin of coordinate $O.$ 
Let's construct the plane $ \Pi _ {l _ {_ 0}} ^ {},$ 
\vspace{0.25ex}
passing through the axis $O ^ {\star} z ^ {\star} $ and the straight line $l _ {_ 0}. $
\vspace{0.35ex}
The plane $ \Pi _ {l _ {_ 0}} ^ {} $ intersects equator of the sphere (1.1) in two antipodal points: 
$L ^ {\star} $ and $L _ {\star} ^ {}$ (Fig. 1.2).
\vspace{0.25ex}
To the infinitely removed point $L $ of the straight line $ \overline {l}$ 
we will put in correspondence two antipodal points $L ^ {\star} $ and $L _ {\star} ^ {}$ 
\vspace{0.25ex}
of the equator the sphere (1.1) 
lying on the plane $ \Pi _ {l _ {_ 0}} ^ {}. $
\vspace{0.35ex}
Then to the continued straight line $ \overline {l} $ to the projective plane $ \R \P (x, y) $ 
\vspace{0.25ex}
on the sphere (1.1) there will correspond the circle of the big radius passing through points $L ^ {\star} $ and $L _ {\star} ^ {}. $
\vspace{0.25ex}

So, binary two-valued correspondence between the projective plane $ \R \P (x, y) $ and the sphere (1.1) at which image of each point of the plane $ \R \P (x, y) $ is the set consisting of two antipodal points of the sphere (1.1) is established.
Thus, such sphere (1.1) we will name \textit{projective sphere} of the plane $Oxy $ and to designate 
\vspace {0.5ex}
$ \P {\mathbb S} (x, y). $

Projective sphere $ \P {\mathbb S} (x, y)$ 
\vspace{0.25ex}
with the identified antipodal points we name [6, p. 749] 
\textit{Poincar\'{e}'s sphere} of the plane $Oxy $ also we will speak about Poincar\'{e}'s sphere $ \P {\mathbb S} (x, y). $

Then 
\vspace{0.15ex}
the introduced binary two-valued correspondence between the projective plane $\R \P (x, y)$ 
\vspace{0.25ex}
and the projective sphere $\P {\mathbb S} (x, y)$ establishes bijective map of the projective plane 
$ \R \P (x, y) $ on Poincar\'{e}'s sphere $\P {\mathbb S} (x, y).$
\vspace{0.25ex}

Hence, Poincar\'{e}'s sphere is two-dimensional variety [7, p. 92 -- 93], homeomorphic to the projective plane.
\\[2.75ex]
\centerline{
{\bf  2. Atlas of  Poincar\'{e}'s sphere}
}
\\[1.5ex]
\indent
The bijection between Poincar\'{e}'s sphere $ \P {\mathbb S} (x, y) $ 
\vspace{0.15ex}
and the projective plane $ \R\P (x, y) $ allows to construct the atlas of maps for Poincar\'{e}'s sphere 
on the basis of the atlas of maps of the projective plane. 
\vspace{0.35ex}
For this purpose, following [7, p. 96], along with the system of coordinate $Oxy $ 
we will introduce two more flat rectangular Cartesian coordinate systems.
\vspace{0.35ex}

On the plane, that concerning sphere (1.1) in the point with coordinates 
\vspace{0.35ex}
$x ^ {\star} =1, \, y ^ {\star} =0,$ $z ^ {\star} =0, $ 
\vspace{0.35ex}
we will introduce the right rectangular Cartesian coordinate system 
$O ^ {(1)} _ {\phantom {1}} \xi \theta $ so that its beginning $O ^ {(1)} _ {\phantom {1}} $ 
\vspace{0.35ex}
coincides with the point of the tangency of the sphere (1.1), 
the axis $O ^ {(1)} _ {\phantom {1}} \xi $ codirected with the axis $O ^ {\star} y ^ {\star}$ 
\vspace{0.55ex}
(so, and with the axis $O y $ \!), the axis $O ^ {(1)} _ {\phantom {1}} \theta $
is opposite directed with the axis $O ^ {\star} z ^ {\star} $ (Fig. 2.1).
\\[4.75ex]
\mbox{}\hfill
{\unitlength=1mm
\begin{picture}(76.272,44)
\put(0,0){\includegraphics[width=81.272mm,height=44mm]{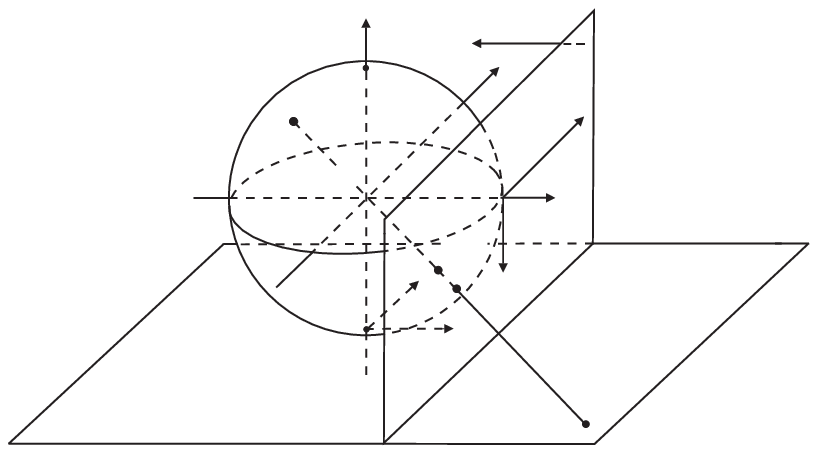}}

\put(34.5,9.7){\makebox(0,0)[cc]{\scriptsize $O$}}
%\put(36.7,13.7){\makebox(0,0)[cc]{\scriptsize $S$}}
%\put(71,13.5){\makebox(0,0)[cc]{\scriptsize $x$}}
\put(45,10.5){\makebox(0,0)[cc]{\scriptsize $x$}}
%\put(26.6,2.9){\makebox(0,0)[cc]{\scriptsize $y$}}
\put(42.7,15){\makebox(0,0)[cc]{\scriptsize $y$}}
%\put(36.9,3.6){\makebox(0,0)[cc]{\scriptsize $(x,y)$}}
\put(59.8,4.7){\makebox(0,0)[cc]{\scriptsize $M$}}
\put(46,19.5){\makebox(0,0)[cc]{\scriptsize $M^{\star}$}}
\put(50.2,15.9){\makebox(0,0)[cc]{\scriptsize $M^{(1)}_{\phantom1}$}}
\put(31,34.3){\makebox(0,0)[cc]{\scriptsize $M_{\star}$}}

%\put(19.5,27.3){\makebox(0,0)[cc]{\scriptsize $O^{(1)}_{\star}$}}
\put(52.7,30.3){\makebox(0,0)[cc]{\scriptsize $O^{(1)}_{\phantom{1}}$}}
\put(55.5,23.2){\makebox(0,0)[cc]{\scriptsize $x^{\star}$}}

\put(42,33){\makebox(0,0)[cc]{\scriptsize $O^{(2)}_{\phantom{1}}$}}
%\put(31.7,16.6){\makebox(0,0)[cc]{\scriptsize $O^{(2)}_{\star}$}}
\put(33.9,26.7){\makebox(0,0)[cc]{\scriptsize $O^{\star}$}}
%\put(14.5,22.7){\makebox(0,0)[cc]{\scriptsize $z^{\star(1)}$}}
%\put(53,26.5){\makebox(0,0)[cc]{\scriptsize $1$}}
%\put(22,26.6){\makebox(0,0)[cc]{\scriptsize $-1$}}
\put(47.4,38){\makebox(0,0)[cc]{\scriptsize $y^{\star}$}}
\put(39,42){\makebox(0,0)[cc]{\scriptsize $z^{\star}$}}

\put(48.8,42.5){\makebox(0,0)[cc]{\scriptsize $\vec{n}$}}
\put(51.8,18.8){\makebox(0,0)[cc]{\scriptsize $\theta$}}
\put(56.5,34){\makebox(0,0)[cc]{\scriptsize $\xi$}}
%\put(44.6,33.5){\makebox(0,0)[cc]{\scriptsize $1$}}
%\put(30.3,21.3){\makebox(0,0)[cc]{\scriptsize $-1$}}

%\put(36.3,40.3){\makebox(0,0)[cc]{\scriptsize $N$}}
%\put(39.5,40.3){\makebox(0,0)[cc]{\scriptsize $1$}}
%\put(34,4.5){\makebox(0,0)[cc]{\scriptsize $\theta^{\star}$}}

\put(40.6,-5){\makebox(0,0)[cc]{\rm Fig. 2.1}}
\end{picture}}
\hfill\mbox{}
\\[10.5ex]
\mbox{}\hfill
{\unitlength=1mm
\begin{picture}(87.153,44)
\put(0,0){\includegraphics[width=87.153mm,height=44mm]{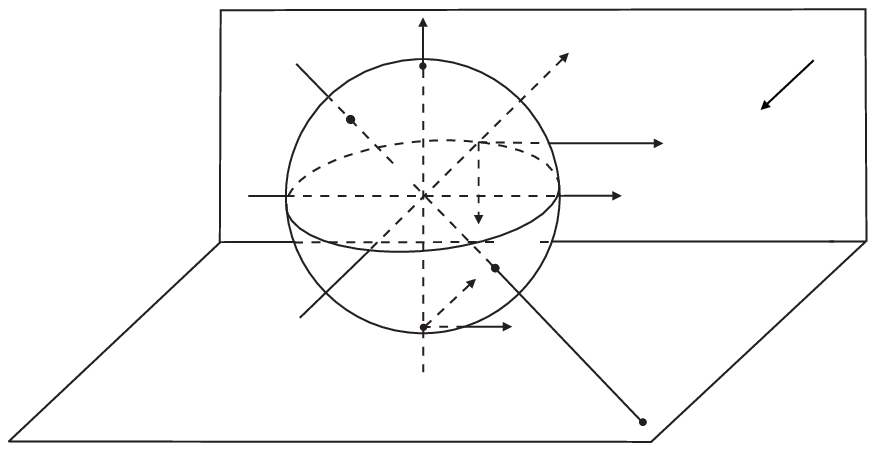}}

\put(40.6,9.7){\makebox(0,0)[cc]{\scriptsize $O$}}
%\put(45.7,13.7){\makebox(0,0)[cc]{\scriptsize $S$}}
%\put(80,13.5){\makebox(0,0)[cc]{\scriptsize $x$}}
\put(51,10.5){\makebox(0,0)[cc]{\scriptsize $x$}}
%\put(35.6,2.9){\makebox(0,0)[cc]{\scriptsize $y$}}
\put(48.3,15){\makebox(0,0)[cc]{\scriptsize $y$}}
%\put(45.9,3.6){\makebox(0,0)[cc]{\scriptsize $(x,y)$}}
\put(65.8,4.7){\makebox(0,0)[cc]{\scriptsize $M$}}
\put(51.6,19.8){\makebox(0,0)[cc]{\scriptsize $M^{\star}$}}
\put(31.5,40.5){\makebox(0,0)[cc]{\scriptsize $M^{(2)}_{\phantom1}$}}
\put(36.9,34.6){\makebox(0,0)[cc]{\scriptsize $M_{\star}$}}

\put(39.9,26.7){\makebox(0,0)[cc]{\scriptsize $O^{\star}$}}
%\put(25.5,27.3){\makebox(0,0)[cc]{\scriptsize $O^{(1)}_{\star}$}}
\put(59.7,26.7){\makebox(0,0)[cc]{\scriptsize $O^{(1)}_{\phantom{1}}$}}
\put(47.3,33){\makebox(0,0)[cc]{\scriptsize $O^{(2)}_{\phantom{1}}$}}
\put(37.7,17){\makebox(0,0)[cc]{\scriptsize $O^{(2)}_{\star}$}}
\put(60.6,23.1){\makebox(0,0)[cc]{\scriptsize $x^{\star}$}}
%\put(53,26.5){\makebox(0,0)[cc]{\scriptsize $1$}}
%\put(22,26.6){\makebox(0,0)[cc]{\scriptsize $-1$}}
%\put(31,10){\makebox(0,0)[cc]{\scriptsize $z^{\star(2)}$}}
\put(49.4,23){\makebox(0,0)[cc]{\scriptsize $\eta$}}
\put(65,28){\makebox(0,0)[cc]{\scriptsize $\zeta$}}
\put(57.8,37){\makebox(0,0)[cc]{\scriptsize $y^{\star}$}}
\put(44.6,42){\makebox(0,0)[cc]{\scriptsize $z^{\star}$}}

%\put(44.6,33.5){\makebox(0,0)[cc]{\scriptsize $1$}}
%\put(30.3,21.3){\makebox(0,0)[cc]{\scriptsize $-1$}}

%\put(36.3,40.3){\makebox(0,0)[cc]{\scriptsize $N$}}
%\put(39.5,40.3){\makebox(0,0)[cc]{\scriptsize $1$}}
%\put(45,5){\makebox(0,0)[cc]{\scriptsize $\eta^{\star}$}}
\put(76,37){\makebox(0,0)[cc]{\scriptsize $\vec{n}$}}

\put(43.5,-5){\makebox(0,0)[cc]{\rm Fig. 2.2}}
\end{picture}}
\hfill\mbox{}
\\[7ex]
\indent
On the plane, that concerning the sphere (1.1) in the point with coordinates 
\vspace{0.35ex}
$x ^ {\star}\!=0, \, y ^ {\star}\!=1,$ $z ^ {\star} =0, $ 
\vspace{0.35ex}
we will introduce the right rectangular Cartesian coordinate system 
$O ^ {(2)} _ {\phantom {1}} \eta \zeta $ so that its beginning $O ^ {(2)} _ {\phantom {1}} $ 
\vspace{0.35ex}
coincides with the point of the tangency of the sphere (1.1), 
the axis $O ^ {(2)} _ {\phantom {1}} \eta $ is opposite directed with the axis 
\vspace{0.35ex}
$O ^ {\star} z ^ {\star}, $ the axis $O ^ {(2)} _ {\phantom {1}} \zeta $ 
codirected with the axis $O ^ {\star} x ^ {\star} $ (Fig. 2.2).
Scale in the coordinate systems 
$Oxy, \ O ^ {\star} x ^ {\star} y ^ {\star} z ^ {\star}, \ O ^ {(1)} _ {\phantom {1}} \xi\theta $ 
\vspace{0.75ex}
and $O ^ {(2)} _ {\phantom {1}} \eta \zeta $ is the identical.

Let's cover the sphere of Poincar\'{e} (with the account  of identifycity of antipodal points) three hemispheres without edge
\\[0.5ex]
\mbox {} \hfill
$
\displaystyle
U_1 ^ {} =\Bigl \{
(x ^ {\star}, y ^ {\star}, z ^ {\star}) \colon \ z ^ {\star} = {}-\sqrt {
1- {x ^ {\star}} ^ {\, ^ {\scriptstyle 2}} - {y ^ {\star}} ^ {\, ^ {\scriptstyle 2}}}\,, 
\, \ \
{x ^ {\star}} ^ {\, ^ {\scriptstyle 2}} + {y ^ {\star}} ^ {\, ^ {\scriptstyle 2}} <1\Bigr \},
\hfill
$
\\[2.75ex]
\mbox {} \hfill
$
\displaystyle
U_2 ^ {} =\Bigl \{
(x ^ {\star}, y ^ {\star}, z ^ {\star}) \colon \ x ^ {\star} = \sqrt {
1 -{y ^ {\star}} ^ {\, ^ {\scriptstyle 2}} - {z ^ {\star}} ^ {\, ^ {\scriptstyle 2}}}\,, \, \ \
{y ^ {\star}} ^ {\, ^ {\scriptstyle 2}} + {z ^ {\star}} ^ {\, ^ {\scriptstyle 2}} <1\Bigr \},
\hfill
$
\\[2.75ex]
\mbox {} \hfill
$
\displaystyle
U_3 ^ {} =\Bigl \{
(x ^ {\star}, y ^ {\star}, z ^ {\star}) \colon \ y ^ {\star} = \sqrt {
1 -{x ^ {\star}} ^ {\, ^ {\scriptstyle 2}} - {z ^ {\star}} ^ {\, ^ {\scriptstyle 2}}}\,, \, \ \
{x ^ {\star}} ^ {\, ^ {\scriptstyle 2}} + {z ^ {\star}} ^ {\, ^ {\scriptstyle 2}} <1\Bigr \}.
\hfill
$
\\[2.75ex]
\indent
Also we will introduce bijective maps [7, p. 96]
\\[2.25ex]
\mbox{}\hfill                                                              % (2.1)
$
\displaystyle
\varphi_1^{}\colon (x^{\star}, y^{\star}, z^{\star})\to \
\bigl(
x(x^{\star}, y^{\star}, z^{\star}),\, y(x^{\star}, y^{\star}, z^{\star})
\bigr), 
\hfill
$
\\
\mbox{}\hfill (2.1)
\\
\mbox{}\hfill
$
x={}-\dfrac{x^{\star}}{z^{\star}}\,, 
\quad 
y={}-\dfrac{y^{\star}}{z^{\star}}
$
\ \
for all 
$
(x^{\star}, y^{\star}, z^{\star})\in U_1^{}\;\!,
\hfill
$
\\[3.5ex]
\mbox{}\hfill                                                              % (2.2)
$
\displaystyle
\varphi_2^{}\colon (x^{\star}, y^{\star}, z^{\star})\to \
\bigl(
\xi(x^{\star}, y^{\star}, z^{\star}),\, \theta(x^{\star}, y^{\star}, z^{\star})
\bigr), 
\hfill
$
\\
\mbox{}\hfill (2.2)
\\
\mbox{}\hfill
$
\xi=\dfrac{y^{\star}}{x^{\star}}\,, 
\,
\quad 
\theta={}-\dfrac{z^{\star}}{x^{\star}}
$ 
\ \ for all 
$
(x^{\star}, y^{\star}, z^{\star})\in U_2^{}\;\!,
\hfill
$
\\[3.5ex]
\mbox{}\hfill                                                              % (2.3)
$
\displaystyle
\varphi_3^{}\colon (x^{\star}, y^{\star}, z^{\star})\to \
\bigl(
\eta(x^{\star}, y^{\star}, z^{\star}),\, \zeta(x^{\star}, y^{\star}, z^{\star})
\bigr), 
\hfill
$
\\
\mbox{}\hfill (2.3)
\\
\mbox{}\hfill
$
\eta={}-\dfrac{z^{\star}}{y^{\star}}\,, 
\,\quad 
\zeta=\dfrac{x^{\star}}{y^{\star}}
$
\ \ for all 
$
(x^{\star}, y^{\star}, z^{\star})\in U_3^{}\;\!.
\hfill
$
\\[3ex]
\indent
Thus, 
\vspace{0.5ex}
three maps $(U_\tau ^ {}, \varphi_\tau ^ {}), \, \tau=1,2,3,$ of  Poincar\'{e}'s sphere $ \P {\mathbb S} (x, y)$ are constructed. 
Set of maps $ (U_\tau ^ {}, \varphi_\tau ^ {}), \, \tau=1,2,3, $ organise the atlas of maps for Poincar\'{e}'s sphere $\P {\mathbb S} (x, y). $
\vspace{0.5ex}

Let's notice that the atlas of maps of the projective sphere $ \P {\mathbb S} (x, y) $ consists of six maps.
\\[2.75ex]
\centerline{
{\bf  3. Poincar\'{e}'s circle}
}
\\[1.5ex]
\indent
Let the point $M $ be 
\vspace{0.25ex}
arranged in the final part $(x, y)$ of the projective plane $\R\P (x, y)$ and 
has coordinates $M (x, y)$ in the coordinate system $Oxy $ (Fig. 1.1).
\vspace{0.35ex}
Then this point in the space coordinate system $O ^ {\star} x ^ {\star} y ^ {\star} z ^ {\star}$ 
\vspace{0.35ex}
has coordinates $M (x, y, {}-1). $ The straight line $MO ^ {\star} $ 
in the coordinate system $O ^ {\star} x ^ {\star} y ^ {\star} z ^ {\star} $ is defined by the system of equations
\\[2.5ex]
\mbox{}\hfill                                    %(3.1)
$
\dfrac{x^{\star}}{x}=
\dfrac{y^{\star}}{y}=
\dfrac{z^{\star}}{{}-1}\,.
$
\hfill (3.1)
\\[2.5ex]
To point $M (x, y, {}-1) $ there correspond points 
\vspace{0.35ex}
$M ^ {\star} (x ^ {\star}, y ^ {\star}, z ^ {\star}) $ 
and $M _ {\star} ^ {} (x _ {\star} ^ {} \, y _ {\star} ^ {} \, z _ {\star} ^ {}), $ 
which are cross points of the straight line $MO ^ {\star} $ and the sphere (1.1). 
\vspace{0.35ex}
Therefore the coordinates $x ^ {\star}, \, y ^ {\star}, \, z ^ {\star} $ and 
$x _ {\star} ^ {}\,, \ y _ {\star} ^ {}\,, \ z _ {\star} ^ {}$ of points $M ^ {\star} $ and $M _ {\star} ^ {} $ 
are  solutions of the algebraic system 
\\[2.5ex]
\mbox{}\hfill                                   %(3.2)
$
\dfrac{x^{\star}}{x}=
\dfrac{y^{\star}}{y}=
\dfrac{z^{\star}}{{}-1}\,, \quad \ \,
x^{\star}{}^{\,^{\scriptstyle 2}}+
y^{\star}{}^{\,^{\scriptstyle 2}}+
z^{\star}{}^{\,^{\scriptstyle 2}}=1.
$
\hfill (3.2)
\\[2.5ex]
\indent
Let us consider the point 
\vspace{0.35ex}
$M ^ {\star} (x ^ {\star}, y ^ {\star}, z ^ {\star}) $ lies in the southern hemisphere $S^2 _ {\scriptscriptstyle-}.$ 
Then its $z\!$-coordinate $z ^ {\star} \in [{}-1; 0). $ 
\vspace{0.5ex}
Having resolved the system of equations (3.2) rather $x ^ {\star}, \, y ^ {\star}, \, z ^ {\star} $ at $ {}-1\leq z ^ {\star} <0, $ 
we will receive the bijective reflexion 
\\[2ex]
\mbox{}\hfill                                   %(3.3)
$
\varphi_1^{{}-1}\colon (x,y)\to\ 
\bigl(x^{\star}(x,y), y^{\star}(x,y), z^{\star}(x,y)\bigr),
\hfill                                  
$
\\[-0.75ex]
\mbox{}\hfill {\rm (3.3)}
\\[0.75ex]
\mbox{}\hfill 
$
x^{\star}=\dfrac{x}{\sqrt{1+x^2+y^2}}\,, \ \ \
y^{\star}=\dfrac{y}{\sqrt{1+x^2+y^2}}\,, \ \ \
z^{\star}={}-\dfrac{1}{\sqrt{1+x^2+y^2}}
$
\quad
for all 
$
(x,y)\in \R^2
\hfill 
$
\\[2.75ex]
of the final part $(x, y)$ 
\vspace{0.5ex}
of the projective plane $ \R\P (x, y) $ on the southern hemisphere $S^2 _ {\scriptscriptstyle-}$ 
without edge $ \partial S^2 _ {\scriptscriptstyle-} $ (without the equator of sphere  (1.1)).
\vspace{0.35ex}

The coordinate functions of map (3.3) are continuously differentiable.
The coordinate function $z ^ {\star} $ expresses through the coordinate functions $x ^ {\star}$ and $y ^ {\star} $ under the formula
\\[2.5ex]
\mbox{}\hfill                                  
$
z^{\star}(x,y)={}-\sqrt{1- x^{\star}{}^{\,^{\scriptstyle 2}}(x,y)- y^{\star}{}^{\,^{\scriptstyle 2}}(x,y) }
$
\ \, for all 
$
(x,y)\in \R^2.
\hfill
$
\\[2.25ex]
Passage from coordinates $x, \, y $ to coordinates $x ^ {\star}, \, y ^ {\star} $ has the Jacobian
\\[2.25ex]
\mbox{}\hfill                                  
$
\dfrac{{\sf D} (x^{\star},y^{\star})}
{{\sf D} (x,y)}=
\dfrac{1}{(1+x^2+y^2)^2}
\ne 0
$
\ \, for all 
$
(x,y) \in \R^2.
\hfill
$
\\[2.25ex]
Hence, the bijective map (3.3) is a diffeomorphism. 

The maps (2.1) and (3.3) are mutually inverse.
\vspace{0.35ex}
Diffeomorphicity of maps (3.3) means diffeomorphicity of maps (2.1) 
southern hemispheres without edge $S^2 _ {\scriptscriptstyle-} \backslash \partial S^2 _ {\scriptscriptstyle-} = U_1 ^ {} $ on $ \R^2. $
\vspace{0.35ex}

By means 
\vspace{0.25ex}
of the bijective map (3.3) the final part $(x, y)$ of the projective plane $ \R\P (x, y) $ 
is diffeomorphically  mapped on the southern hemisphere $S^2 _ {\scriptscriptstyle-} $ 
\vspace{0.25ex}
without edge $ \partial S ^ {2} _ {\scriptscriptstyle-},$ being the equator of the projective sphere $ \P {\mathbb S} (x, y) $ (Fig. 1.1).
To each infinitely removed point of the projective plane $\R\P (x, y)$ 
\vspace{0.25ex}
there correspond two antipodal points of equator $ \partial S ^ {2} _ {\scriptscriptstyle-} $ 
of the projective sphere $\P {\mathbb S} (x, y).$
\vspace{0.25ex}
Therefore the southern hemisphere $S^2 _ {\scriptscriptstyle-} $ 
with the identified antipodal points of equator $ \partial S ^ {2} _ {\scriptscriptstyle-} $ of the projective sphere $ \P {\mathbb S} (x, y) $
is space model of the projective plane $ \R\P (x, y). $
\\[2.5ex]
\mbox{}\hfill
{\unitlength=1mm
\begin{picture}(109,53)
\put(0,0){\includegraphics[width=109mm,height=52.46mm]{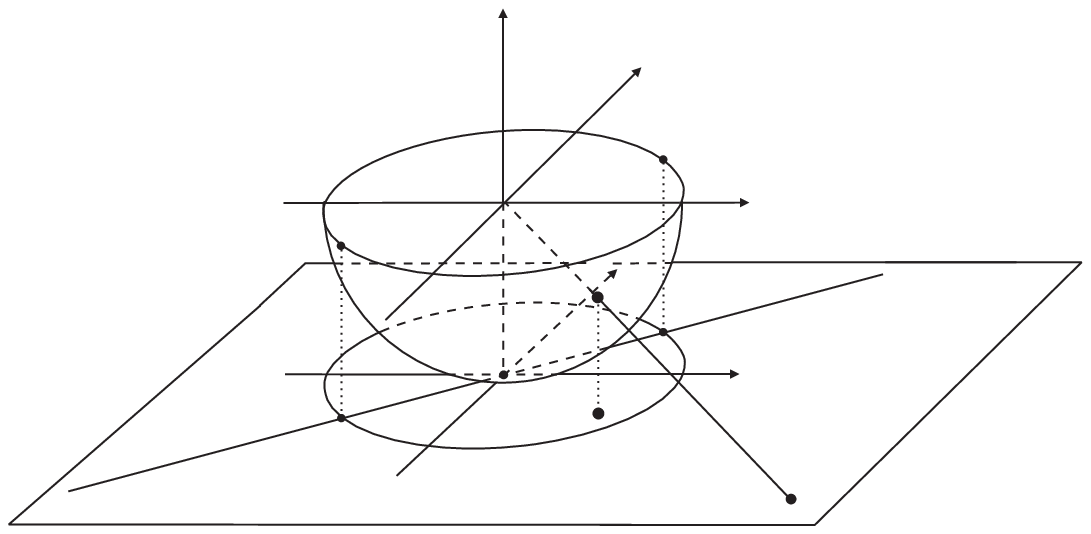}}

\put(51.9,13){\makebox(0,0)[cc]{\scriptsize $O$}}
%\put(36.7,13.7){\makebox(0,0)[cc]{\scriptsize $S$}}
%\put(71,13.5){\makebox(0,0)[cc]{\scriptsize $x$}}
\put(73,14){\makebox(0,0)[cc]{\scriptsize $x$}}
%\put(26.6,2.9){\makebox(0,0)[cc]{\scriptsize $y$}}
\put(59.8,25.7){\makebox(0,0)[cc]{\scriptsize $y$}}
%\put(36.9,3.6){\makebox(0,0)[cc]{\scriptsize $(x,y)$}}
\put(76,2.7){\makebox(0,0)[cc]{\scriptsize $M$}}

\put(48,34.4){\makebox(0,0)[cc]{\scriptsize $O^{\star}$}}
\put(62.8,23.3){\makebox(0,0)[cc]{\scriptsize $M^{\star}$}}
\put(55.2,11.2){\makebox(0,0)[cc]{\scriptsize $\varkappa(M)$}}
%\put(34,8.2){\makebox(0,0)[cc]{\scriptsize $D^{}_{\scriptscriptstyle-}$}}
%\put(71,18){\makebox(0,0)[cc]{\scriptsize $D^{}_{\scriptscriptstyle+}$}}
\put(36,29.5){\makebox(0,0)[cc]{\scriptsize $L^{}_{\star}$}}
\put(68.5,38.5){\makebox(0,0)[cc]{\scriptsize $L^{\star}$}}
\put(10,2){\makebox(0,0)[cc]{\scriptsize $L$}}
\put(88,22.7){\makebox(0,0)[cc]{\scriptsize $L$}}
\put(22,10){\makebox(0,0)[cc]{\scriptsize $l^{}_{\scriptscriptstyle 0}$}}
\put(73,31){\makebox(0,0)[cc]{\scriptsize $x^{\star}$}}
%\put(53,26.5){\makebox(0,0)[cc]{\scriptsize $1$}}
%\put(22,26.6){\makebox(0,0)[cc]{\scriptsize $-1$}}
\put(62,46.5){\makebox(0,0)[cc]{\scriptsize $y^{\star}$}}
%\put(54,19){\makebox(0,0)[cc]{\scriptsize $z$}}
%\put(59.3,34){\makebox(0,0)[cc]{\scriptsize $u$}}
%\put(44.6,33.5){\makebox(0,0)[cc]{\scriptsize $1$}}
%\put(30.3,21.3){\makebox(0,0)[cc]{\scriptsize $-1$}}

%\put(36.3,40.3){\makebox(0,0)[cc]{\scriptsize $N$}}
%\put(39.5,40.3){\makebox(0,0)[cc]{\scriptsize $1$}}
\put(53,51.5){\makebox(0,0)[cc]{\scriptsize $z^{\star}$}}

\put(54.5,-5){\makebox(0,0)[cc]{\rm Fig. 3.1}}
\end{picture}}
\hfill\mbox{}
\\[6.5ex]
\indent
The natural projection of the southern hemisphere 
\vspace{0.35ex}
$S^2 _ {\scriptscriptstyle-}$ on the plane $(x, y) $ is the circle of unit radius 
$K (x, y) = \{(x, y) \colon x^2+y^2\leq1 \} $ (Fig. 3.1).
\vspace{0.35ex}
By means of this projection diffeomorphic map $p $ of the southern hemisphere 
\vspace{0.35ex}
of Poincar\'{e} $S^2 _ {\scriptscriptstyle-} $ on the circle $K (x, y)$ is established. 
Superposition of maps $ \varkappa=p\circ \varphi_1 ^ {{}-1} $ is the diffeomorphic map
\\[2ex]
\mbox{}\hfill                                 
$
\varkappa \colon (x,y)\to \
\biggl( \dfrac{x}{\sqrt{1+x^2+y^2}}\,,\ \dfrac{y}{\sqrt{1+x^2+y^2}}\biggr)
$
\ \ 
for all 
$
(x,y) \in \R^2
\hfill  
$
\\[2.25ex]
of the final part $(x, y)$ 
\vspace{0.25ex}
of the projective plane $ \R\P (x, y) $ on the open circle $K (x, y) \backslash\partial K (x, y). $ 
To each infinitely removed point of the projective plane $ \R\P (x, y) $ 
\vspace{0.25ex}
there correspond two anti\-po\-dal points of the boundary circle 
\vspace{0.25ex}
$ \partial K (x, y) = \{(x, y) \colon x^2+y^2=1 \} $ of the circle $K (x, y). $ 
Such circle $K (x, y) $ we will name {\it projective circle} 
\vspace{0.35ex}
of plane $Oxy $ and to designate $ \P\K (x, y). $

The projective circle $\P\K (x, y)$ 
\vspace{0.25ex}
with the identified antipodal points of the boundary circle $\partial K (x, y) $ 
\vspace{0.25ex}
we name {\it Poincar\'{e}'s circle} of plane $Oxy$ and we will tell about Poincar\'{e}'s circle $\P\K (x, y). $
\vspace{0.25ex}
The circle of Poincar\'{e} $\P\K (x, y)$ is diffeomorphic to the projective plane $ \R\P (x, y), $ so, is flat compact model of this projective plane.
\\[2ex]
\centerline{
{\bf  4. Maps of Poincar\'{e}}
}
\\[1.5ex]
\indent
Let's establish connections between local coordinate systems 
\vspace{0.35ex}
$Oxy,\ O ^ {(1)} _ {\phantom {1}} \xi\theta, \ O ^ {(2)} _ {\phantom {1}} \eta\zeta $ 
of the atlas of maps for Poincar\'{e}'s sphere $\P {\mathbb S} (x, y).$
\vspace{0.35ex}

The plane $O ^ {(1)} _ {\phantom {1}} \xi \theta $ in the coordinate system $O ^ {\star} x ^ {\star} y ^ {\star} z ^ {\star}$ 
\vspace{0.35ex}
is defined by the equation $x ^ {\star} =1.$ 
Arbitrarily we will choose the point $M (x, y)$ in the final part $(x, y)$ 
\vspace{0.35ex}
of the projective plane $\R\P (x, y)$  
such that this point did not lie on the axis $Oy.$
\vspace{0.35ex}
Then (Fig. 2.1) the straight line $MO ^ {\star} $ intersects the plane $O ^ {(1)} _ {\phantom {1}} \xi \theta $ in the point $M ^ {(1)}. $ 
\vspace{0.35ex}
The straight line $MO ^ {\star}$ is defined by the system of equations (3.1), 
and the plane $O ^ {(1)} _ {\phantom {1}} \xi \theta$  is defined by the equation $x ^ {\star} =1.$ 
\vspace{1ex}

Therefore 
\vspace{0.75ex}
the point $M ^ {(1)}$ in the coordinate system $O ^ {\star} x ^ {\star} y ^ {\star} z ^ {\star}$ has the abscissa $x ^ {\star} =1.$ 
The ordinate $y ^ {\star} = \dfrac {y} {x}$ and the z-coordinate $z ^ {\star} = {}-\dfrac{1}{x}$ of the point $M ^ {(1)}$ 
\vspace{1ex}
are discovered from the system of equations (3.1) with $x^{\star}=1.$
\vspace{1ex}
The point $M _ {\phantom {1}} ^ {(1)}$ in the coordinate system $O ^ {(1)} _ {\phantom {1}} \xi\theta $ has 
the abscissa $ \xi=y ^ {\star} = \dfrac {y} {x}$ and the ordinate  $\theta={}-z^{\star}=\dfrac{1}{x}\,.$
\vspace{1ex}

Thus the point $M (x, y) $ with the abscissa $x\ne 0,$
\vspace{0.35ex}
which is lying in the final part $ (x, y) $ of the projective plane $ \R \P (x, y),$  
\vspace{0.5ex}
corresponds to the point $M _ {\phantom {1}} ^ {(1)} (\xi, \theta)$
on the coordinate plane $O ^ {(1)} _ {\phantom {1}} \xi\theta.$
\vspace{0.35ex}
The coordinates of the point $M _ {\phantom {1}} ^ {(1)} (\xi, \theta)$
through the coordinates of the point $M (x, y) $ express under the formulas
\\[1.5ex]
\mbox{}\hfill       % (4.1)
$
\xi=\dfrac{y}{x}\,,
\quad \ \ 
\theta=\dfrac{1}{x}\,.
$
\hfill (4.1)
\\[1.5ex]
\indent
Thereby the bijective map 
\\[1.75ex] 
\mbox {} \hfill      % (4.2) 
$
P ^ {(1)} _ {} \colon (x, y) \to\  
\Bigl (\;\!\dfrac {y} {x}\,,\;\! \dfrac{1}{x}\;\!\Bigr) 
$
\ for all 
$
(x , y) \in\R^2\backslash \{(x, y) \colon x=0 \}$ 
\hfill (4.2) 
\\[1.75ex] 
of the plane $Oxy, $ from which the axis $Oy$ is removed,  
\vspace{0.35ex}
on the plane $O ^ {(1)} _ {\phantom {1}} \xi\theta, $ from which the axis $O ^ {(1)} _ {\phantom {1}} \xi$ is removed, is established. 
\vspace{0.35ex}
The map (4.2) is said to be {\it the first map of Poincar\'{e}} of the plane $Oxy. $ 
The map (4.2) is a diffeomorphism with the Jacobian
\\[2ex]
\mbox{}\hfill
$
\dfrac{{\sf D}(\xi,\theta)}{{\sf D}(x,y)}=\dfrac{1}{x^3}\ne 0
$
\ for all 
$
(x,y)\in\R^2\backslash \{(x,y)\colon x=0\}.
\hfill
$
\\[2ex]
\indent
Having resolved equalities (4.1) rather $x$ and $y, $ we obtain the formulas
\\[2ex]
\mbox {} \hfill % (4.3)
$
x =\dfrac {1} {\theta}\,, 
\quad \ \
y =\dfrac {\xi} {\theta} \,.
$
\hfill (4.3)
\\[2ex]
\indent
By the formulas (4.3), 
\vspace{0.35ex}
the coordinates of the point $M (x, y)$ with $x\ne 0$ express through 
the coordinates of the point $M ^ {(1)}$ with $ \theta\ne 0. $
\vspace{0.35ex}
The formulas (4.3) are called {\it the first transformation of Poincar\'{e}} of the plane $Oxy $ [1, p. 31].
\vspace{0.35ex}

The functions (4.1) are functions of passage from the coordinates $(\xi, \theta)$ to the coordinates $(x, y),$ 
and the functions (4.3) are functions of passage from the coordinates $ (x, y) $ to the coordinates $(\xi, \theta).$

Arbitrarily we take the straight line $y=a\;\!x$ with parameter $a\in\R $ and we will transform it under the formulas (4.3).
\vspace{0.35ex}
As a result we will receive the straight line $ \xi=a $ on the plane $O ^ {(1)} _ {\phantom {1}} \xi\theta. $
\vspace{0.35ex}
Thus the infinitely removed point $L, $ which is lying on <<extremities>> of the straight line $y=ax,$ 
corresponds to the point $L _ {\phantom {1}} ^ {(1)} (a, 0), $ which is lying on the axis 
\vspace{0.35ex}
$O ^ {(1)} _ {\phantom {1}} \xi $ of the plane $O ^ {(1)} _ {\phantom {1}} \xi\theta. $
\vspace{0.35ex}
And, on the contrary, to each point $L _ {\phantom {1}} ^ {(1)} (a, 0),$ 
which is lying on the coordinate axis $O ^ {(1)} _ {\phantom {1}} \xi $ of the plane $O ^ {(1)} _ {\phantom {1}} \xi\theta, $ 
\vspace{0.35ex}
corresponds the infinitely removed point $L, $ which is lying on <<extremities>> of the straight line $y=ax. $ 

Hence we have
\vspace{0.35ex}

{\bf Property 4.1.} 
\vspace{0.35ex}
{\it
The first transformation of Poincar\'{e} $ (4.3) $ establishes diffeomorphic map of the projective plane $\R \P (x, y)$ 
\vspace{0.35ex}
without the straight line $x=0$ on the combined coordinate plane $O ^ {(1)} _ {\phantom {1}}\xi\theta.$  
\vspace{0.35ex}
Thus by image of the infinitely removed straight line of the projective plane $ \R \P (x, y),$ 
\vspace{0.35ex}
from which it is removed the point lying on <<extremities>> of the straight line $x=0, $ 
is the coordinate axis $O ^ {(1)} _ {\phantom {1}} \xi $ {\rm (}the straight line $ \theta=0 {\rm )}$  
of the plane $O ^ {(1)} _ {\phantom {1}} \xi\theta. $
}
\vspace{1ex}

The plane $O ^ {(2)} _ {\phantom {1}} \eta\zeta$
\vspace{0.5ex}
in the coordinate system $O ^ {\star} x ^ {\star} y ^ {\star} z ^ {\star}$ 
is defined by the equation $y ^ {\star}=1.$
Arbitrarily we will choose the point $M (x, y)$
\vspace{0.5ex}
in the final part $(x, y)$ of the projective plane $\R\P (x, y)$  such that this point did not lie on the axis $Ox.$ 
\vspace{0.5ex}
Then (Fig. 2.2) the straight line $MO^{\star}$ intersects the plane $O^{(2)} _ {\phantom {1}} \eta\zeta $ in the point $M^{(2)}.$
\vspace{0.5ex}
The straight line $MO^{\star}$ is defined by the system of equations (2.1), and the plane  
$O^{(2)} _ {\phantom {1}} \eta\zeta $ is defined by the equation $y^{\star} =1.$ 
\vspace{0.5ex}

Therefore the point $M _ {\phantom {1}}^{(2)}$
\vspace{0.75ex}
in the coordinate system $O ^ {\star} x ^ {\star} y ^ {\star} z ^ {\star}$ has the ordinate $y^{\star} =1.$ 
The abscissa $x^{\star} =\dfrac{x}{y}$ and the z-coordinate $z^{\star} = {}-\dfrac{1}{y}$ 
\vspace{0.75ex}
of the point $M _ {\phantom {1}}^{(2)}$ are discovered from the system of equations (2.1) with $y^{\star} =1.$
\vspace{0.75ex}
The point $M _ {\phantom {1}}^{(2)}$ in the coordinate system 
$O^{(2)} _ {\phantom {1}} \eta\zeta$ has the abscissa $\eta = {}-z ^ {\star} = \dfrac{1}{y}$ 
and the ordinate $\zeta=x^{\star}=\dfrac{x}{y}\,.$
\vspace {1ex}

Thus the point $M (x, y)$ with the ordinate $y\ne 0,$ which is
\vspace{0.35ex}
lying in the final part $(x, y)$ of the projective plane $\R\P (x, y),$ 
\vspace{0.5ex}
corresponds to the point $M _ {\phantom {1}} ^ {(2)} (\eta, \zeta)$ 
on the coordinate plane $O ^ {(2)} _ {\phantom {1}} \eta\zeta.$
\vspace{0.35ex}
The coordinates of the point $M _ {\phantom {1}} ^ {(2)} (\eta, \zeta)$ 
through the coordinates of the point $M (x, y) $ express under the formulas
\\[2.15ex]
\mbox{}\hfill       % (4.4)
$
\eta=\dfrac{1}{y}\,,
\quad \ \ 
\zeta=\dfrac{x}{y}\,.
$
\hfill (4.4)
\\[2ex]
\indent
Thereby the bijective map 
\\[2.15ex]
\mbox{}\hfill       % (4.5)
$
P^{(2)}_{}\colon (x,y)\to\  
\Bigl(\;\!\dfrac{1}{y}\,,\;\!\dfrac{x}{y}\;\!\Bigr)
$
\ for all 
$
(x,y)\in\R^2\backslash\{(x,y)\colon y=0\}
$
\hfill (4.5)
\\[2ex]
of the plane $Oxy,$ from which the axis $Ox$ is removed,  
\vspace{0.5ex}
on the plane $O ^ {(2)} _ {\phantom {1}} \eta\zeta,$ from which the axis $O^{(2)}_{\phantom{1}}\zeta$
is removed, is established.
\vspace{0.5ex}
The map (4.5) is said to be {\it the second map of Poincar\'{e}} of the plane $Oxy. $
The map (4.5) is a diffeomorphism with the Jacobian
\\[2ex]
\mbox{}\hfill
$
\dfrac{{\sf D}(\eta,\zeta)}{{\sf D}(x,y)}=\dfrac{1}{y^3}\ne 0
$
\ for all 
$
(x,y)\in\R^2\backslash \{(x,y)\colon y=0\}.
\hfill
$
\\[2ex]
\indent
Having resolved equalities (4.4) rather $x$ and $y, $ we obtain the formulas
\\[2ex]
\mbox {} \hfill % (4.6)
$
x =\dfrac {\zeta} {\eta}\,,
\quad \ \ 
y =\dfrac {1} {\eta}\,.
$
\hfill (4.6)
\\[2.15ex]
\indent
By the formulas (4.6),
\vspace{0.5ex}
the coordinates of the point $M (x, y) $ with the ordinate $y\ne 0$ express through the coordinates 
of the point $M _ {\phantom {1}} ^ {(2)} (\eta, \zeta)$ with $\eta\ne 0.$
\vspace{0.5ex}
The formulas (4.6) are called {\it the second transformation of Poincar\'{e}} of the plane $Oxy$ [1, p. 31].
\vspace{0.35ex}

The functions (4.4) 
\vspace{0.25ex}
are functions of passage from the coordinates $(\eta, \zeta)$ to the coordinates $(x, y),$ 
and the functions (4.6) are functions of passage from the coordinates $(x, y) $ to the coordinates $(\eta, \zeta).$
\vspace{0.25ex}

Arbitrarily we take the straight line $x=b\;\!y $ 
\vspace{0.35ex}
with parameter $b\in\R$ and we will transform it under the formulas (4.6).
\vspace{0.5ex}
As a result we will receive the straight line $\zeta=b$ on the plane $O ^ {(2)} _ {\phantom {1}} \eta\zeta.$
\vspace{0.5ex}
Thus the infinitely removed point $L,$ which is lying on <<extremities>> of the straight line $x=b\;\!y,$ 
corresponds to the point $L _ {\phantom {1}} ^ {(2)}(0, b),$ which is lying on the axis 
\vspace{0.5ex}
$O ^ {(2)} _ {\phantom {1}} \zeta$ of the plane $O ^ {(2)} _ {\phantom {1}} \eta\zeta.$
\vspace{0.5ex}
And, on the contrary, to each point $L _ {\phantom {1}} ^ {(2)} (0, b),$ which is 
lying on the coordinate axis $O ^ {(2)} _ {\phantom {1}} \zeta $ of the plane $O ^ {(2)} _ {\phantom {1}} \eta\zeta,$ 
\vspace{0.5ex}
corresponds the infinitely removed point $L,$ which is lying on <<extremities>> of the straight line $x=b\;\!y. $ 
\vspace{0.35ex}

Therefore we obtain
\vspace{0.5ex}

{\bf Property 4.2.} 
\vspace{0.35ex}
{\it
The second transformation of Poincar\'{e} {\rm (4.6)} establishes diffeomorphic map 
of the projective plane $\R\P (x, y)$ without the straight line $y=0$ 
\vspace{0.5ex}
on the combined coordinate plane
$O ^ {(2)} _ {\phantom {1}} \eta\zeta.$
\vspace{0.5ex}
Thus by image of infinitely removed straight line of the projective plane $\R \P (x, y),$ from which it is removed 
\vspace{0.5ex}
the point lying on <<extremities>> of the straight line $y=0, $ is the coordinate axis 
$O ^ {(2)} _ {\phantom {1}} \zeta$ {\rm(}the straight line $\eta=0)$ of the plane 
$O ^ {(2)} _ {\phantom {1}} \eta\zeta.$
}
\vspace {1.5ex}

{\bf Theorem 4.1.} 
{\it
The identity mapping 
\vspace{0.5ex}
$I\colon (x, y) \to (x, y)$ for all $(x, y) \in\R^2, $ 
the first {\rm (4.2)} and the second {\rm (4.5)} maps of Poincar\'{e} are a group of the third order}:
\\[2ex]
\mbox{}\hfill                                   %(4.7)
$
I\circ I=I,
\quad 
P^{(1)}_{}\circ I=I\circ P^{(1)}_{}=P^{(1)}_{},
\quad
P^{(2)}_{}\circ I=I\circ P^{(2)}_{}=P^{(2)}_{},
\hfill        
$
\\[-0.75ex]
\mbox{}\hfill        
\hfill (4.7)       
\\[0.75ex]
\mbox{}\hfill        
$
P^{(1)}_{}\circ P^{(2)}_{}=P^{(2)}_{}\circ P^{(1)}_{}=I,
\quad
P^{(1)}_{}\circ P^{(1)}_{}=P^{(2)}_{},
\quad
P^{(2)}_{}\circ P^{(2)}_{}=P^{(1)}_{}.
\hfill
$
\\[2ex]
\indent
{\sl Really,}  the first and the second transformations of Poincar\'{e} are mutually inverse:
\\[2ex]
\mbox{}\hfill
$
P^{(1)}_{}\circ P^{(2)}_{}=(x,y)\stackrel{\stackrel{\scriptstyle P^{(2)}_{}}
{\mbox{}}}{\rightarrow}\Bigl(\dfrac{1}{y}\,,\dfrac{x}{y}\Bigr)
\stackrel{\stackrel{\scriptstyle P^{(1)}_{}}{\mbox{}}}
{\rightarrow}\Bigl(\dfrac{x}{y}\,:\dfrac{1}{y}\,,1:
\dfrac{1}{y}\Bigr)=(x,y)=I;
\hfill
$
\\[2.25ex]
\mbox{}\hfill
$
P^{(2)}_{}\circ P^{(1)}_{}=(x,y)\stackrel{\stackrel{\scriptstyle P^{(1)}_{}}
{\mbox{}}}{\rightarrow}\Bigl(\dfrac{y}{x}\,,\dfrac{1}{x}\Bigr)
\stackrel{\stackrel{\scriptstyle P^{(2)}_{}}{\mbox{}}}
{\rightarrow}\Bigl(1:\dfrac{1}{x}\,,\dfrac{y}{x}\,:
\dfrac{1}{x}\Bigr)=(x,y)=I.
\hfill
$
\\[2ex]
\indent
Besides, 
\\[2ex]
\mbox{}\hfill        
$
P^{(1)}_{}\circ P^{(1)}_{}=(x,y)\stackrel{\stackrel{\scriptstyle P^{(1)}_{}}{\mbox{}}}
{\rightarrow}\Bigl(\dfrac{y}{x}\,,\dfrac{1}{x}\Bigr)
\stackrel{\stackrel{\scriptstyle P^{(1)}_{}}{\mbox{}}}
{\rightarrow}\Bigl(\dfrac{1}{x}:\dfrac{y}{x}\,,1:
\dfrac{y}{x}\Bigr)=\Bigl(\dfrac{1}{y}\,,\dfrac{x}{y}\Bigr)=P^{(2)}_{},
\hfill 
$
\\[2.25ex]
\mbox{}\hfill      
$
P^{(2)}_{}\circ P^{(2)}_{}=(x,y)\stackrel{\stackrel{\scriptstyle P^{(2)}_{}}{\mbox{}}}
{\rightarrow}\Bigl(\dfrac{1}{y}\,,\dfrac{x}{y}\Bigr)
\stackrel{\stackrel{\scriptstyle P^{(2)}_{}}{\mbox{}}}{\rightarrow}
\Bigl(1:\dfrac{x}{y}\,,\dfrac{1}{y}:\dfrac{x}{y}\Bigr)=
\Bigl(\dfrac{y}{x}\,,\dfrac{1}{x}\Bigr)=P^{(1)}_{}.
\hfill 
$
\\[3ex]
\indent
It is obvious that 
\\[2ex]
\mbox{}\hfill
$
P^{(1)}_{}\circ I=I\circ P^{(1)}_{}=P^{(1)}_{},
\quad \
P^{(2)}_{}\circ I=I\circ P^{(2)}_{}=P^{(2)}_{},
\quad \
I\circ I=I.\
\k
\hfill
$
\\[4.25ex]
\centerline{
{\bf  5. Atlas of projective circles for projective plane}
}
\\[1.5ex]
\indent
We take the atlas of maps $(U_\tau ^ {}, \varphi_\tau ^ {}), \ \tau=1,2,3, $ 
\vspace{0.5ex}
of Poincar\'{e}'s sphere $\P {\mathbb S} (x, y). $
On the local rectangular Cartesian coordinate systems 
\vspace{0.5ex}
$Oxy,\ O ^ {(1)} \xi\theta,\ O^{(2)} \eta\zeta$ 
we will construct (Fig.~5.1) the projective circles
$\P\K(x,y),\ \P\K(\xi,\theta),\ \P\K(\eta,\zeta).$
\vspace{0.5ex}

The group property (Theorem 4.1) of Poincar\'{e}'s maps (4.2) and (4.5) 
\vspace{0.35ex}
allows to establish the correspondences between the projective circles 
$\P\K (x, y), \ \P\K (\xi, \theta),$ and $\P\K (\eta, \zeta),$ which are shown on Fig. 5.1. 
\vspace{0.25ex}

By means of the points $M_1 ^ {}, \ M_2 ^ {}, \ M_3 ^ {}, \ M_4 ^ {} $ it is shown 
\vspace{0.25ex}
in what open coordinate quarter of the projective circles these points lie, 
\vspace{0.25ex}
when we  passage from one projective circle to another. 

The numbers $1, \ldots, 12$ 
\vspace{0.25ex}
reflect the correspondences between halfneighborhoods 
of the points lying on boundaries of coordinate quarters of these projective circles.
\\[3ex]
\mbox{}\hfill
{\unitlength=1mm
\begin{picture}(42,42)
\put(0,0){\includegraphics[width=42mm,height=42mm]{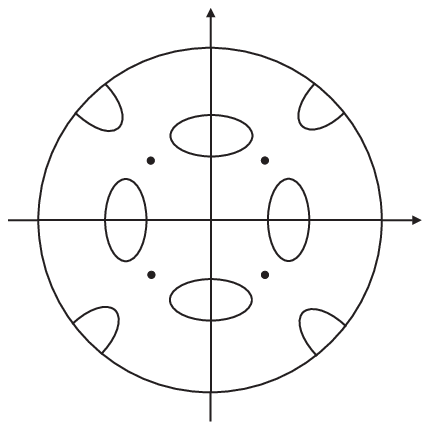}}

\put(18,41){\makebox(0,0)[cc]{ $y$}}
\put(40.2,18.2){\makebox(0,0)[cc]{ $x$}}

\put(28.5,22.3){\makebox(0,0)[cc]{\scriptsize $1$}}
\put(31.3,31.5){\makebox(0,0)[cc]{\scriptsize $2$}}
\put(22.3,29){\makebox(0,0)[cc]{\scriptsize $3$}}

\put(19,29){\makebox(0,0)[cc]{\scriptsize $4$}}
\put(9.7,31.5){\makebox(0,0)[cc]{\scriptsize $5$}}
\put(12,22.3){\makebox(0,0)[cc]{\scriptsize $6$}}

\put(12,18.6){\makebox(0,0)[cc]{\scriptsize $7$}}
\put(9.3,9.5){\makebox(0,0)[cc]{\scriptsize $8$}}
\put(19,12.3){\makebox(0,0)[cc]{\scriptsize $9$}}

\put(22.3,12.3){\makebox(0,0)[cc]{\scriptsize $10$}}
\put(31.5,9.5){\makebox(0,0)[cc]{\scriptsize $11$}}
\put(28.4,18.8){\makebox(0,0)[cc]{\scriptsize $12$}}

\put(28,28){\makebox(0,0)[cc]{\scriptsize $M_1^{}$}}
\put(14,28.8){\makebox(0,0)[cc]{\scriptsize $M_2^{}$}}
\put(12.5,13.5){\makebox(0,0)[cc]{\scriptsize $M_3^{}$}}
\put(28.8,13.5){\makebox(0,0)[cc]{\scriptsize $M_4^{}$}}

%\put(22.5,-6){\makebox(0,0)[cc]{\rm $K(x,y)$}}
\end{picture}
}
\quad
{\unitlength=1mm
\begin{picture}(42,42)
\put(0,0){\includegraphics[width=42mm,height=42mm]{r05-01.eps}}

\put(18,41){\makebox(0,0)[cc]{ $\theta$}}
\put(40.2,17.8){\makebox(0,0)[cc]{ $\xi$}}

\put(28.5,22.3){\makebox(0,0)[cc]{\scriptsize $2$}}
\put(31.3,31.5){\makebox(0,0)[cc]{\scriptsize $3$}}
\put(22.3,29){\makebox(0,0)[cc]{\scriptsize $1$}}

\put(18.8,29){\makebox(0,0)[cc]{\scriptsize $12$}}
\put(9.5,31.5){\makebox(0,0)[cc]{\scriptsize $10$}}
\put(12,22.1){\makebox(0,0)[cc]{\scriptsize $11$}}

\put(12,18.6){\makebox(0,0)[cc]{\scriptsize $5$}}
\put(9.3,9.5){\makebox(0,0)[cc]{\scriptsize $4$}}
\put(19,12.3){\makebox(0,0)[cc]{\scriptsize $6$}}

\put(22.3,12.3){\makebox(0,0)[cc]{\scriptsize $7$}}
\put(31.5,9.5){\makebox(0,0)[cc]{\scriptsize $9$}}
\put(28.4,18.8){\makebox(0,0)[cc]{\scriptsize $8$}}

\put(28,28){\makebox(0,0)[cc]{\scriptsize $M_1^{}$}}
\put(14,28.8){\makebox(0,0)[cc]{\scriptsize $M_4^{}$}}
\put(12.5,13.5){\makebox(0,0)[cc]{\scriptsize $M_2^{}$}}
\put(28.8,13.5){\makebox(0,0)[cc]{\scriptsize $M_3^{}$}}

%\put(22.5,-6){\makebox(0,0)[cc]{\rm $K(\xi,\theta)$}}
\put(21,-6){\makebox(0,0)[cc]{\rm Fig. 5.1}}
\end{picture}
}
\quad
{\unitlength=1mm
\begin{picture}(42,42)
\put(0,0){\includegraphics[width=42mm,height=42mm]{r05-01.eps}}

\put(18,41){\makebox(0,0)[cc]{ $\zeta$}}
\put(40.2,18){\makebox(0,0)[cc]{ $\eta$}}

\put(28.5,22.3){\makebox(0,0)[cc]{\scriptsize $3$}}
\put(31.3,31.5){\makebox(0,0)[cc]{\scriptsize $1$}}
\put(22.3,29){\makebox(0,0)[cc]{\scriptsize $2$}}

\put(19,29){\makebox(0,0)[cc]{\scriptsize $8$}}
\put(9.7,31.5){\makebox(0,0)[cc]{\scriptsize $7$}}
\put(12,22.3){\makebox(0,0)[cc]{\scriptsize $9$}}

\put(12,18.7){\makebox(0,0)[cc]{\scriptsize $10$}}
\put(9.3,9.7){\makebox(0,0)[cc]{\scriptsize $12$}}
\put(18.8,12.3){\makebox(0,0)[cc]{\scriptsize $11$}}

\put(22.3,12.3){\makebox(0,0)[cc]{\scriptsize $5$}}
\put(31.5,9.5){\makebox(0,0)[cc]{\scriptsize $6$}}
\put(28.4,18.8){\makebox(0,0)[cc]{\scriptsize $4$}}

\put(28,28){\makebox(0,0)[cc]{\scriptsize $M_1^{}$}}
\put(14,28.8){\makebox(0,0)[cc]{\scriptsize $M_3^{}$}}
\put(12.5,13.5){\makebox(0,0)[cc]{\scriptsize $M_4^{}$}}
\put(28.8,13.5){\makebox(0,0)[cc]{\scriptsize $M_2^{}$}}

%\put(22.5,-6){\makebox(0,0)[cc]{\rm $K(\eta,\zeta)$}}
\end{picture}
}
\hfill\mbox{}
\\[7ex]
\indent
The ordered triple $ (\P\K (x, y), \, \P\K (\xi, \theta), \, \P\K (\eta, \zeta))$ 
\vspace{0.25ex}
is said to be {\it the atlas of projective circles} of the projective plane $\R\P (x, y).$ 
\vspace{0.25ex}

Then the ordered triple $ (\P\K (\xi, \theta), \, \P\K (\eta, \zeta), \, \P\K (x, y))$ 
\vspace{0.35ex}
is the atlas of projective circles of the projective plane $\R \P (\xi, \theta),$ 
and the ordered triple $(\P\K (\eta, \zeta), \, \P\K (x, y), \, \P\K (\xi, \theta))$ 
\vspace{0.25ex}
there is the atlas of projective circles of the projective plane $ \R \P (\eta, \zeta). $
\\[4.75ex]
\centerline{
{\bf\large \S\;\!2. Transformations of Poincar\'{e} for differential systems}
}
\\[2.25ex]
\centerline{
{\bf  6. Projectively reduced systems
}
}
\\[1.5ex]
\indent
By the first transformation of Poincar\'{e} (4.3), it follows that the differential system (D) is reduced to the system
\\[2ex]
\mbox{}\hfill                             % (6.1)
$
\begin{array}{l}
\dfrac{d\;\!\xi}{dt}\, ={}-\xi\;\! \theta\;\! 
X\!\Bigl(\;\!\dfrac{1}{\theta}\,,\;\!\dfrac{\xi}{\theta}\;\!\Bigr) +\, 
\theta\;\! Y\!\Bigl(\;\!\dfrac{1}{\theta}\,,\;\!\dfrac{\xi}{\theta}\;\!\Bigr) \equiv\;\!
\Xi_{\phantom1}^{(1)}(\xi,\theta),
\\[4.25ex]
\dfrac{d\;\!\theta}{dt}={}-\theta^{\;\!2}
X\Bigl(\;\!\dfrac{1}{\theta}\,,\;\!\dfrac{\xi}{\theta}\;\!\Bigr)\equiv\;\! 
\Theta_{\phantom1}^{(1)}(\xi,\theta).
\end{array}
$
\hfill (6.1)
\\[2.5ex]
\indent
Since $X $ and $Y $ are polynomials, we see that the system (6.1) has the form
\\[2ex]
\mbox{}\hfill                               %(6.2)
$
\dfrac{d\;\!\xi}{dt}\;\!=\;\!\dfrac{\Xi(\xi,\theta)}{\theta^{m}}\,,
\qquad
\dfrac{d\;\!\theta}{dt}\;\!=\;\!\dfrac{\Theta(\xi,\theta)}{\theta^{m}}\,,
$
\hfill (6.2)
\\[2.25ex]
where $ \Xi $ and $ \Theta $ are polynomials, which are not sharing simultaneously on $\theta,$ 
and the number $m$ is a whole nonnegative.
Using the system (6.2), we can make the system
\\[2ex]
\mbox{}\hfill                                     %(6.3)
$
\dfrac{d\;\!\xi}{d\tau}=\Xi(\xi,\theta),
\qquad
\dfrac{d\;\!\theta}{d\tau}=\Theta(\xi,\theta),
$
\hfill (6.3)
\\[2.25ex]
where $\theta ^ {m}\;\!d\tau=dt,$ at which right members $ \Xi $ and $ \Theta $ are relatively prime polynomials.
\vspace{0.5ex}

By the second transformation of Poincar\'{e} (4.6), it follows that the differential system (D) is reduced to the system
\\[2ex]
\mbox{}\hfill                     %(6.4)
$
\begin{array}{l}
\dfrac{d\;\!\eta}{dt}={}-\eta^2\,Y\Bigl(\;\!\dfrac{\zeta}{\eta}\,,\;\!\dfrac{1}{\eta}\;\!\Bigr)
\equiv H_{\phantom1}^{(2)}(\eta,\zeta),
\\[4.25ex]
\dfrac{d\;\!\zeta}{dt}=\eta\,X\Bigl(\;\!\dfrac{\zeta}{\eta}\,,\;\!\dfrac{1}{\eta}\;\!\Bigr)
-\eta\;\!\zeta\,Y\Bigl(\;\!\dfrac{\zeta}{\eta}\,,\;\!\dfrac{1}{\eta}\;\!\Bigr) \equiv
Z_{\phantom1}^{(2)}(\eta,\zeta).
\end{array}
$
\hfill (6.4)
\\[2.5ex]
\indent
Since $X $ and $Y $ are polynomials, we see that the system (6.4) has the form
\\[2ex]
\mbox{}\hfill                    %(6.5)
$
\dfrac{d\eta}{dt}\;\!=\;\!\dfrac{H(\eta,\zeta)}{\eta^{m}}\,,
\qquad
\dfrac{d\zeta}{dt}\;\!=\;\!\dfrac{Z(\eta,\zeta)}{\eta^{m}}\,,
$
\hfill (6.5)
\\[2ex]
where $H $ and $Z $ are polynomials, which are not sharing simultaneously on $\eta,$ 
and the number $m$ is a whole nonnegative.
Using the system (6.5), we can make the system
\\[2ex]
\mbox{}\hfill                    %(6.6)
$
\dfrac{d\eta}{d\nu}=H(\eta,\zeta),
\qquad
\dfrac{d\zeta}{d\nu} =Z(\eta,\zeta),
$
\hfill (6.6)
\\[2.25ex]
where $\eta^{m}\;\!d\nu=dt,$ at which right members $H $ and $Z $ are relatively prime polynomials.
\vspace{0.75ex}

The autonomous polynomial differential systems (6.3) and (6.6), which are
\vspace{0.25ex}
obtained on the basis of systems (6.1) and (6.4),  
\vspace{0.25ex}
are said to be {\it the projectively reduced systems} or 
\linebreak
P-{\it\!\! reduced systems} of system (D). 
\vspace{0.25ex}
The system (6.3) is called {\it the first projectively reduced system} or (P-1)-{\it\!\! reduced system} of systems (D), 
\vspace{0.25ex}
and the system (6.6) is called  {\it the second projectively reduced system} or (P-2)-{\it\!\! reduced system} of systems (D). 
\vspace{0.25ex}
The planes $ (\xi, \theta) $ and $ (\eta, \zeta) $ are the phase planes for P-reduced systems (6.3) and (6.6), respectively.
\\[3.25ex]
\mbox{}\hfill
{\unitlength=1mm
\begin{picture}(35,26)
\put(0,0){\includegraphics[width=34.78mm,height=25.02mm]{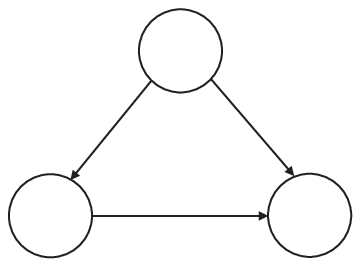}}

\put(4.6,4.3){\makebox(0,0)[cc]{\scriptsize (6.3)}}
\put(30.5,4.3){\makebox(0,0)[cc]{\scriptsize (6.6)}}
\put(17.5,20.5){\makebox(0,0)[cc]{\scriptsize (D)}}

\put(8,14.2){\makebox(0,0)[cc]{\scriptsize $P^{(1)}_{}$}}
\put(28.4,14.2){\makebox(0,0)[cc]{\scriptsize $P^{(2)}_{}$}}
\put(17.5,6.5){\makebox(0,0)[cc]{\scriptsize $P^{(1)}_{}$}}

\put(18,-6){\makebox(0,0)[cc]{\rm Fig. 6.1}}
\end{picture}}
\qquad
\
{\unitlength=1mm
\begin{picture}(35,26)
\put(0,0){\includegraphics[width=34.78mm,height=25.02mm]{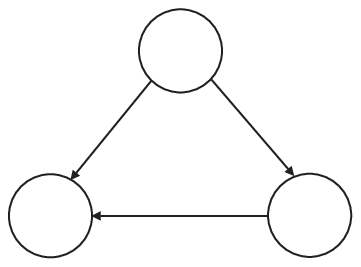}}

\put(4.6,4.3){\makebox(0,0)[cc]{\scriptsize (6.3)}}
\put(30.5,4.3){\makebox(0,0)[cc]{\scriptsize (6.6)}}
\put(17.5,20.5){\makebox(0,0)[cc]{\scriptsize (D)}}

\put(8,14.2){\makebox(0,0)[cc]{\scriptsize $P^{(1)}_{}$}}
\put(28.4,14.2){\makebox(0,0)[cc]{\scriptsize $P^{(2)}_{}$}}
\put(17.5,6.5){\makebox(0,0)[cc]{\scriptsize $P^{(2)}_{}$}}

\put(18,-6){\makebox(0,0)[cc]{\rm Fig. 6.2}}
\end{picture}}
\qquad
\
{\unitlength=1mm
\begin{picture}(61.3,8.5)
\put(0,9){\includegraphics[width=61.3mm,height=8.5mm]{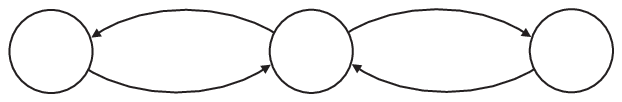}}

\put(4.6,13.3){\makebox(0,0)[cc]{\scriptsize (6.3)}}
\put(57,13.3){\makebox(0,0)[cc]{\scriptsize (6.6)}}
\put(30.6,13.3){\makebox(0,0)[cc]{\scriptsize (D)}}

\put(17.5,20){\makebox(0,0)[cc]{\scriptsize $P^{(1)}_{}$}}
\put(44,5.7){\makebox(0,0)[cc]{\scriptsize $P^{(1)}_{}$}}
\put(44,20){\makebox(0,0)[cc]{\scriptsize $P^{(2)}_{}$}}
\put(17.5,5.7){\makebox(0,0)[cc]{\scriptsize $P^{(2)}_{}$}}

\put(30.5,-6){\makebox(0,0)[cc]{\rm Fig. 6.3}}
\end{picture}}
\hfill\mbox{}
\\[7.75ex]
\indent
Usin the group property (Theorem 4.1) of Poincar\'{e}'s maps,  
\vspace{0.5ex}
we can establish the connections 
between the systems (D), (6.3), (6.6). 
\vspace{0.5ex}
The sequence scheme $P ^ {(1)} _ {} \circ P ^ {(1)} _ {} =P ^ {(2)} _ {}$ of Poincar\'{e}'s maps for 
the systems (D), (6.3), (6.6) is represented on Fig. 6.1, 
\vspace{0.5ex}
the sequence scheme $P ^ {(2)} _ {} \circ P ^ {(2)} _ {} =P ^ {(1)} _ {} $ is represented on Fig. 6.2, 
\vspace{0.5ex}
and the schemes of sequences $P ^ {(1)} _ {} \circ P ^ {(2)} _ {} =I $ and $P ^ {(2)} _ {} \circ P ^ {(1)} _ {} =I $ 
are represented on Fig. 6.3. 
\vspace{1ex}

{\bf Property 6.1.} 
\vspace {0.75ex}
{\it 
The first projectively reduced system {\rm (6.3)} by means of the first Poincar\'{e}'s transformation 
$\xi =\dfrac{1}{\zeta}\,,\  \theta=\dfrac{\eta}{\zeta}$ is reduced to the second projectively reduced system {\rm (6.6)}. 
}
\vspace{1.5ex}

{\bf Property 6.2.} 
\vspace {0.75ex}
{\it 
The first projectively  reduced system {\rm (6.3)} by means of the second Poincar\'{e}'s transformation  
$
\xi =\dfrac {y} {x}\,, \ 
\theta =\dfrac {1} {x}
$
is reduced to the system {\rm (D)}. 
}
\vspace{1.5ex}

{\bf Property 6.3.} 
\vspace {0.75ex}
{\it 
The second projectively  reduced system {\rm (6.6)} by means of the first Poincar\'{e}'s transformation  
$
\eta =\dfrac{1}{y}\,, \ 
\zeta =\dfrac{x}{y}
$
is reduced to the system {\rm (D)}. 
}
\vspace{1.5ex}

{\bf Property 6.4.} 
\vspace {0.75ex}
{\it 
The second projectively  reduced system {\rm (6.6)} by means of the second Poin\-ca\-re's transformation  
$
\eta =\dfrac{\theta}{\xi}\,, \ 
\zeta =\dfrac {1} {\xi}
$
is reduced to the first projectively reduced system {\rm (6.3)}. 
}
\\[4.5ex]
\centerline{
{\bf  7.
 Projective type of differential system
}
}
\\[1.5ex]
\indent
For the differential system (D) the form of the first and the second projectively  reduced systems (6.3) and (6.6) 
depends on, whether the polynomial
\\[2ex]
\mbox{}\hfill                           
$
W_n^{}\colon (x,y)\to\  x\;\!Y_n^{}(x, y)\;\!-\;\!y\;\!X_n^{}(x, y)
$
\ for all 
$
(x,y)\in \R^2
\hfill
$
\\[2ex]
is identical zero on the plane $\R^2$ or not. 
\vspace{0.5ex}

If $W_n ^ {} (x, y) \not\equiv 0$ on $ \R^2, $ then the first projectively  reduced system (6.3) has the form
\\[2ex]
\mbox{}\hfill                            % (7.1)
$
\begin{array}{l}
\displaystyle
\dfrac{d\;\!\xi}{d\tau}\,=\ \sum\limits_{i=0}^{n}\,\theta^{n-i}\,Y_{i}^{}(1,\xi)-
\xi\;\!\sum\limits_{i=0}^{n}\,\theta^{n - i}X_{i}^{}(1, \xi)
\,\equiv\,
\sum\limits_{i=0}^{n}\,\theta^{n - i}\,W_{i}^{}(1, \xi)
\equiv\;\!
\widetilde{\Xi}(\xi,\theta),
\\[4.75ex]
\displaystyle
\dfrac{d\;\!\theta}{d\tau}\, =
{}-\theta\;\! \sum\limits_{i=0}^{n}\,\theta^{n - i}X_{i}^{}(1, \xi)
\equiv\;\!
\widetilde{\Theta}(\xi,\theta),
\end{array}
$
\hfill (7.1)
\\[2.5ex]
where 
\vspace{0.75ex}
$\theta^{n-1}\;\!d\tau=dt,\ \, W_i^{}(x,y)=x\;\!Y_i^{}(x,y)-yX_i^{}(x,y)$ for all $(x,y)\in\R^2,\ i=0,1,\ldots,n,$
$W_n ^ {} (1, \xi) \not\equiv0$ on the field $\R,$ 
and the second projectively  reduced system (6.6) has the form
\\[2.5ex]
\mbox{}\hfill                             % (7.2)
$
\begin{array}{l}
\displaystyle
\dfrac{d\;\!\eta}{d\nu}\,=
{}-\eta\;\!\sum\limits_{i=0}^{n}\, \eta^{n - i}\, Y_{i}^{}(\zeta,1)
\equiv\;\!
\widetilde{H}(\eta,\zeta),
\\[4.75ex]
\displaystyle
\dfrac{d\;\!\zeta}{d\nu}\, =\
\sum\limits_{i = 0}^{n}\, \eta^{n-i}X_{i}^{}(\zeta,1)\;\! -\;\!
\zeta\;\!\sum\limits_{i = 0}^{n}\,\eta^{n - i}\, Y_{i}^{}(\zeta,1)
\;\!\equiv
{}-\sum\limits_{i = 0}^{n}\, \eta^{n - i}\,W_{i}^{}(\zeta,1)
\equiv\;\!
\widetilde{Z}(\eta,\zeta),
\end{array}
$
\hfill (7.2)
\\[2.75ex]
where $ \eta ^ {n-1}\,d\nu=dt, $ and $W_n ^ {} (\zeta, 1) \not\equiv 0$ on the field $\R.$
\vspace{1.35ex}

If $W_n ^ {} (x, y) \equiv 0$ on $\R^2,$ then the first projectively  reduced system (6.3) has the form 
\\[2.25ex]
\mbox{}\hfill                           % (7.3)
$
\begin{array}{l}
\displaystyle
\dfrac{d\;\!\xi}{d\tau}\,=\
\sum\limits_{j=0}^{n-1}\,\theta^{n-j-1}\,Y_{j}^{}(1, \xi)\;\!-\;\!
\xi\;\!\sum\limits_{j = 0}^{n-1}\,\theta^{n-j-1}X_{j}^{}(1, \xi)
\;\!\equiv\;\!
\sum\limits_{j = 0}^{n-1}\,\theta^{n-j-1}\,W_{j}^{}(1, \xi)
\equiv\;\!
\widehat{\Xi}(\xi,\theta),
\\[5ex]
\displaystyle
\dfrac{d\;\!\theta}{d\tau}\,=
{}-\sum\limits_{i = 0}^{n}\,\theta^{n - i}X_{i}^{}(1, \xi)
\equiv\;\!
\widehat{\Theta}(\xi,\theta),
\end{array}
$
\hfill (7.3)
\\[2ex]
where 
\vspace{0.5ex}
$ \theta ^ {n-2}\,d\tau=dt,\ X_n ^ {} (1, \xi) \not\equiv 0$ on the field $\R,$ 
and the second projectively reduced system (6.6) has the form
\\[2ex]
\mbox{}\hfill                             % (7.4)
$
\begin{array}{l}
\displaystyle
\dfrac{d\;\!\eta}{d\nu}\,=
{}-\sum\limits_{i=0}^{n}\,\eta^{n - i}\,Y_{i}^{}(\zeta,1)
\equiv \;\!
\widehat{H}(\eta,\zeta),
\\[5ex]
\displaystyle
\dfrac{d\zeta}{d\nu}\,=\
\sum\limits_{j = 0}^{n - 1}\,\eta^{n-j-1}X_{j}^{}(\zeta,1)\;\!-\;\!
\zeta\;\!\sum\limits_{j = 0}^{n - 1}\, \eta^{n-j-1}\,Y_{j}^{}(\zeta,1)\;\!
\equiv
{}-\sum\limits_{j = 0}^{n - 1}\, \eta^{n-j-1}\,W_{j}^{}(\zeta,1)
\equiv\;\!
\widehat{Z}(\eta,\zeta),
\end{array}
\!\!\!\!\!\!\!\!\!\!
$
\hfill (7.4)
\\[2.5ex]
where $\eta^{n-2}\,d\nu=dt,$ and $Y_n ^ {} (\zeta, 1) \not\equiv 0$ on the field $\R.$
\vspace{1.25ex}

If $W_n ^ {} (x, y) \not\equiv 0$ on $ \R^2,$ 
\vspace{0.5ex}
then the system (D) is called {\it projectively nonsingular} or 
\linebreak
P-{\it\!\! nonsingular}. 
\vspace{0.5ex}
Otherwise, i.e. if $W_n ^ {} (x, y) =0$ for all $(x, y)\in\R^2,$ then the system (D) is called 
{\it projectively singular} or P-{\it\! singular}. 
\vspace{0.25ex}
All differential systems (D) are divided into two classes: P-singular and P-nonsingular. 
\vspace{0.25ex}
The form of the projectively reduced systems (6.3) and (6.6) depends on membership of system  (D) to this or that class.
\vspace{1ex}

{\bf Property 7.1}.
\vspace{0.25ex}
{\it 
The first projectively  reduced system {\rm (6.3)} for the projectively nonsingular system {\rm (D)}  has the form {\rm (7.1),} 
and the second projectively  reduced system {\rm (6.6)} 
\vspace{0.25ex}
for the projectively nonsingular system {\rm (D)}  has the form} (7.2).
\vspace{1ex}

{\bf Property 7.2}.
\vspace{0.25ex}
{\it 
The first projectively  reduced system {\rm (6.3)} for the projectively singular system {\rm (D)} has the form {\rm (7.3),} 
and the second projectively  reduced system {\rm (6.6)} 
\vspace{0.25ex}
for the projectively singular system {\rm (D)} has the form} (7.4).
\vspace{1ex}

The number $n =\max \{{\rm deg} \, X, \, {\rm deg} \, Y \}$ is called {\it degree of system} (D) and is denoted by ${\rm deg}\, {\rm(D)}.$

\newpage

The degrees of the P-reduced systems (7.1), (7.2), (7.3), and (7.4) 
\vspace{0.25ex}
depends on the degree $n$ of system (D), and also from that 
the straight lines $x=0$ and $y=0$ consist or not of trajectories of system (D). 
\vspace{0.25ex}

Consider three numbers $\delta,\ \delta ^ {(1)},$ and $\delta ^ {(2)}.$ 
\vspace{0.25ex}
If the system (D) is projectively singular, then $\delta=0;$ and if the system (D) is projectively nonsingular, then $\delta=1.$
\vspace{0.35ex}

If the straight line $x=0$ does not consist of trajectories of system (D), 
\vspace{0.35ex}
then $\delta ^ {(1)} =0; $ 
and if the straight line $x=0$ consists of trajectories of system (D), then $\delta ^ {(1)} =1.$
\vspace{0.25ex}

If the straight line $y=0$ does not consist of trajectories of system  (D), 
\vspace{0.35ex}
then $\delta ^ {(2)}=0;$ 
and if the straight line $y=0$ consists of trajectories of system (D), then $\delta ^ {(2)} =1.$ 
\vspace{0.25ex}

Hence we get
\vspace{0.35ex}

{\bf Property 7.3}.
{\it 
The degree of the first projectively reduced system} (6.3) {\it is}
\\[2ex]
\mbox {} \hfill % (7.5)
$
{\rm deg} \, (6.3) =n +\delta-\delta ^ {(1)},
$
\hfill (7.5)
\\[2ex]
{\it and the degree of the second projectively reduced system} (6.6) {\it is}
\\[2ex]
\mbox{}\hfill                              %(7.6)
$
{\rm deg}\,(6.6) =n+\delta-\delta^{(2)}.
$
\hfill (7.6)
\\[2ex]
\indent
It follows from the formulas (7.5) and (7.6) that the degrees of the P-reduced systems (7.1) --- (7.4) are
\\[2ex]
\mbox{}\hfill                              %(7.7)
$
{\rm deg}\,(7.1) =n-\delta^{(1)}+1,
\quad \ \
{\rm deg}\,(7.2) =n-\delta^{(2)}+1,
\hfill
$
\\
\mbox{}\hfill (7.7)
\\
\mbox{}\hfill
$
{\rm deg}\,(7.3) =n-\delta^{(1)},
\quad \ \
{\rm deg}\,(7.4) =n-\delta^{(2)}.
\hfill
$
\\[2.5ex]
\indent
{\bf Examples}
\vspace{0.75ex}

{\bf 7.1.}
Let us consider the autonomous system
\\[2ex]
\mbox{}\hfill                           %(7.8)
$
\dfrac{dx}{dt}=a_{_0}\equiv X(x,y),
\quad
\dfrac{dy}{dt}=b_{_{0}}\equiv Y(x,y),
\qquad
|a_{_0}|+|b_{_0}|\ne 0.
$
\hfill (7.8)
\\[2.5ex]
\indent
The condition $|a_{_0}|+|b_{_0}|\ne 0$ is justified only by that 
\vspace{0.75ex}
$X_{_0}(x,y)=a_{_0},\ Y_{_0}(x,y)=b_{_0},$ and the system (D) at $n=0$ should be such that 
$|X_{_0}(x,y)|+|Y_{_0}(x,y)|\not\equiv0$ on $\R^2.$
\vspace{1ex}

The polynomial 
\vspace {0.75ex}
$W _ {_ 0} (x, y) =b _ {_ 0} x-a _ {_ 0} y\not\equiv0$ on $ \R^2$ at 
$ |a _ {_ 0} | + |b _ {_ 0} | \ne 0.$ 

The system (7.8) is a projectively nonsingular system.
\vspace {0.35ex}

The first projectively reduced system of system (7.8)  is the system
\\[2ex]
\mbox{}\hfill                           % (7.9)
$
\dfrac{d\xi}{d\tau}=b_{_0}-a_{_0}\xi
\equiv\widetilde{\Xi}(\xi,\theta),
\quad
\dfrac{d\theta}{d\tau}={}-a_{_0}\theta
\equiv\widetilde{\Theta}(\xi,\theta),
$
\ where \  $d\tau=\theta dt,\ |a_{_0}|+|b_{_0}|\ne 0;$ 
\hfill\mbox{} (7.9)
\\[2.35ex]
and the second projectively  reduced system of system (7.8) is the system
\\[2.1ex]
\mbox{}\hfill                             % (7.10)
$
\dfrac{d\eta}{d\nu}={}-b_{_0}\eta
\equiv \widetilde{H}(\eta,\zeta),
\quad
\dfrac{d\zeta}{d\nu}=a_{_0}-b_{_0}\zeta
\equiv \widetilde{Z}(\eta,\zeta),
$
\ where \ $d\nu=\eta dt,\ |a_{_0}|+|b_{_0}|\ne 0.$ 
\hfill\mbox{} (7.10)
\\[3.25ex]
\indent
{\bf 7.2.}
Consider the linear autonomous system
\\[2ex]
\mbox{}\hfill                                           %(7.11)
$
\dfrac{dx}{dt}  =  a_{_0}  +  a^{}_{1}x  +  a^{}_{2}y
\equiv X(x,y),
\quad
\dfrac{dy}{dt} =  b_{_0}  +  b_{1}^{}x  +  b_{2}^{}y
\equiv Y(x,y)
$
\hfill (7.11)
\\[2.5ex]
with coefficients such that $|a_{1}^{}| + |a_{2}^{}| +|b_1^{}| + |b_2^{}|\ne 0.$
\vspace{1ex}

The polynomial 
\vspace {1ex}
$W_1^{}(x,y)=b_1^{}x^2+(b_2^{}-a_1^{})xy-a_2^{}y^2$ for all $(x,y)\in\R^2.$

If $|a_2^{}|+|b_1^{}|+|a_1^{}-b_2^{}|\ne 0,$ then the system (7.11) is P-nonsingular; 
\vspace {0.75ex}
and if $a_2^{}=b_1^{}=0,$ $a_1^{}=b_2^{}\ne 0,$ 
then the system (7.11) is P-singular.

\newpage

Suppose 
$|a_{2}^{}|+|b_{1}^{}| +|a_{1}^{}-b_{2}^{}|\ne 0,$
\vspace{0.75ex}
i.e. the system (7.11) is P-nonsingular.

Using the first Poincar\'{e}'s transformation $x=\dfrac{1}{\theta}\,,\ y=\dfrac{\xi}{\theta}\,,$ we get
\vspace{0.75ex}
the P-nonsingular system (7.11) is reduced to the  (P-1)-reduced system
\\[2.5ex]
\mbox{}\hfill                                 %(7.12)
$
\begin{array}{l}
\dfrac{d\xi}{dt} = 
b_{1}^{} - (a_{1}^{} - b_{2}^{})\xi + b_{_0}\theta -
a_{2}^{}\xi^{2} - a_{_0}\xi\theta
\,\equiv \widetilde{\Xi} (\xi,\theta),
\\[3.75ex]
\dfrac{d\theta}{dt} = {} - a_{1}^{}\theta - a_{2}^{}\xi\theta - a_{_0}\theta^{2}
\,\equiv  \widetilde{\Theta}(\xi,\theta),
\end{array}
$
\hfill\mbox{} (7.12)
\\[2.25ex]
where $|a_{2}^{}|+|b_{1}^{}| +|a_{1}^{}-b_{2}^{}|\ne 0.$
\vspace{0.5ex}

Using the second Poincar\'{e}'s transformation
$x=\dfrac{\zeta}{\eta}\,,\ y=\dfrac{1}{\eta}\,,$ we get
\vspace{0.75ex}
the P-nonsingular system (7.11) is reduced to the (P-2)-reduced system
\\[2ex]
\mbox{}\hfill                                         %(7.13)
$
\begin{array}{l}
\dfrac{d\eta}{dt} = {} - b_{2}^{}\eta -b_{_0}\eta^{2} -b_{1}^{}\eta \zeta
\, \equiv \widetilde{H}(\eta,\zeta),
\\[3.75ex]
\dfrac{d\zeta}{dt} = a_{2}^{} +a_{_0}\eta+ (a_{1}^{} - b_{2}^{})\zeta
- b_{_0}\eta\zeta - b_{1}^{}\zeta^{2}
\,\equiv \widetilde{Z} (\eta,\zeta),
\end{array}
$
\hfill\mbox{} (7.13)
\\[2.5ex]
where $|a_{2}^{}|+|b_{1}^{}| +|a_{1}^{}-b_{2}^{}|\ne 0.$
\vspace{1.25ex}

Suppose  
\vspace{0.75ex}
$a_{2}^{} =  b_{1}^{} = 0,\ b_{2}^{}=a_{1}^{}\ne 0,$
i.e. the system (7.11) is P-singular.

Using the first Poincar\'{e}'s transformation $x=\dfrac{1}{\theta}\,,\ y=\dfrac{\xi}{\theta}\,,$ we obtain
\vspace{0.75ex}
the P-singular system (7.11) is reduced to the (P-1)-reduced system
\\[2.25ex]
\mbox{}\hfill                                      %(7.14)
$
\dfrac{d\xi}{d\tau} = b_{_0} -  a_{_0}\xi
\;\!\equiv\;\! 
\widehat{\Xi}(\xi,\theta),
\qquad
\dfrac{d\theta}{d\tau} = {}-a_{1}^{}-a_{_0}\theta
\;\!\equiv\;\! 
\widehat{\Theta}(\xi,\theta),
$
\hfill\mbox{} (7.14)
\\[2.75ex]
where $d\tau=\theta dt,$ the coefficient $a_1^{}\ne 0.$
\vspace{0.5ex}

Using the second Poincar\'{e}'s transformation
$x=\dfrac{\zeta}{\eta}\,,\ y=\dfrac{1}{\eta}\,,$ we obtain
\vspace{0.75ex}
the P-singular system (7.11) is reduced to the (P-2)-reduced system
\\[2ex]
\mbox{}\hfill                                         %(7.15)
$
\dfrac{d\eta}{d\nu} = {}-a_{1}^{} - b_{_0}\eta
\;\!\equiv\;\! \widehat{H}(\eta,\zeta),
\qquad
\dfrac{d\zeta}{d\nu} = a_{_0} - b_{_0}\zeta
\;\!\equiv\;\! \widehat{Z}(\eta,\zeta),
$
\hfill\mbox{} (7.15)
\\[2.5ex]
where $d\nu=\theta dt,$ the coefficient $a_1^{}\ne 0.$
\vspace{1.25ex}

{\bf 7.3.}
Let us consider the autonomous quadratic system
\\[2.25ex]
\mbox{}\hfill                                         %(7.16)
$
\begin{array}{l}
\dfrac{dx}{dt} = a_{_0}+a_1^{} x+a_2^{} y+a_3^{}x^2+a_4^{} xy+a_5^{} y^2
\equiv X(x,y),
\\[4ex]
\dfrac{dy}{dt} = 
b_{_0}+b_1^{} x+b_2^{} y+b_3^{}x^2+b_4^{} xy+b_5^{} y^2
\equiv Y(x,y),
\end{array}
$
\hfill\mbox{} (7.16)
\\[2.75ex]
where
$|a_3^{}| +|a_{4}^{}|+|a_{5}^{}| +|b_{3}^{}|+|b_{4}^{}| + |b_{5}^{}|\ne  0.$
\vspace{1.25ex}

The polynomial 
\vspace{0.75ex}
$W_2^{}(x,y)=b_3^{}x^3+(b_4^{}-a_3^{})x^2y+(b_5^{}-a_4^{})xy^2 -a_5^{}y^3$ for all $(x,y)\in\R^2.$ 

If 
\vspace{0.75ex}
$|a_5^{}| +|b_{3}^{}|+|a_{3}^{}-b_{4}^{}|\,+
|a_{4}^{}-b_{5}^{}|\ne  0,$ then the system (7.16) 
is P-nonsingular, and if
$a_5^{} =b_{3}^{}=0,\ a_{3}^{}=b_{4}^{},\ a_{4}^{}=b_{5}^{},\ |b_4^{}|+|b_5^{}|\ne  0,$ 
then the system (7.16) is P-singular.
\vspace{0.75ex}

Suppose 
$|a_5^{}| +|b_{3}^{}|+|a_{3}^{}-b_{4}^{}|+|a_{4}^{}-b_{5}^{}|\ne  0,$ i.e. the system (7.16) is P-nonsingular.
\vspace{0.75ex}

The first projectively  reduced system for the P-nonsingular system (7.16) is
\\[2ex]
\mbox{}\hfill                                         %(7.17)
$
\begin{array}{l}
\dfrac{d\;\!\xi}{d\tau} \ =\,  b_{3}^{} + (b_{4}^{} - a_{3}^{})\xi +
b_{1}^{}\theta + (b_{5}^{} - a_{4}^{})\xi^{2} +
(b_{2}^{} - a_{1}^{})\xi\theta\ +
\\[2.25ex]
\mbox{}\quad\  \ 
+\, b_{_0}\theta^{2} - a_{5}^{}\xi^{3} - a_{2}^{}\xi^{2}\theta - a_{_0}\xi\theta^{2}
\;\! \equiv\;\!  \widetilde{\Xi}(\xi,\theta),
\\[3ex]
\dfrac{d\;\!\theta}{d\tau}\, = {} -  a_{3}^{}\theta - a_{4}^{}\xi\theta - a_{1}^{}\theta^{2} -
a_{5}^{}\xi^{2}\theta - a_{2}^{}\xi\theta^{2} - a_{_0}\theta^{3}
\;\!\equiv\;\!  \widetilde{\Theta} (\xi,\theta),
\end{array}
$
\hfill\mbox{} (7.17)
\\[3ex]
where 
\vspace{0.5ex}
$z\,d\tau = dt,\ 
|a_5^{}| +|b_{3}^{}|+|a_{3}^{}-b_{4}^{}|+|a_{4}^{}-b_{5}^{}|\ne  0,$ 
and the second projectively  reduced system for the P-nonsingular system (7.16) is
\\[2.25ex]
\mbox{}\hfill                                         %(7.18)
$
\begin{array}{l}
\dfrac{d\;\!\eta}{d\nu}\  = {} -b_{5}^{}\eta - b_{2}^{}\eta^{2}- b_{4}^{}\eta\;\!\zeta  - b_{_0}\eta^{3}-
b_{1}^{}\eta^{2}\zeta-b_{3}^{}\eta\;\!\zeta^{2}
\;\!\equiv \, \widetilde{H}(\eta,\zeta),
\\[4ex]
\dfrac{d\zeta}{d\nu}\,  = \, a_{5}^{} + a_{2}^{}\eta
+ (a_{4}^{} - b_{5}^{})\zeta
+ a_{_0}\eta^{2}
+(a_{1}^{} - b_{2}^{})\eta\;\!\zeta  
+  (a_{3}^{} - b_{4}^{})\zeta^{2}\ -
\\[2.5ex]
\mbox{}\quad\ \
-\, b_{_0}\eta^{2}\;\!\zeta-b_{1}^{}\eta\;\!\zeta^{2}
- b_{3}^{}\zeta^{3}
\;\!\equiv \, \widetilde{Z} (\eta,\zeta),
\end{array}
$
\hfill\mbox{} (7.18)
\\[2.5ex]
where $z\,d\nu = dt,\ 
|a_5^{}| +|b_{3}^{}|+|a_{3}^{}-b_{4}^{}|+|a_{4}^{}-b_{5}^{}|\ne  0.$ 
\vspace{1ex}

Suppose $a_{5}^{} = b_{3}^{}=0,\ b_{4}^{} =  a_{3}^{},\ b_{5}^{} =  a_4^{},\ |a_{3}^{}|+|a_{4}^{}|\ne0,$ 
\vspace{0.75ex}
i.e. the system (7.16) is P-singular. 

The first projectively reduced system of the P-nonsingular system (7.16)  is
\\[1.75ex]
\mbox{}\hfill                                         %(7.19)
$
\begin{array}{l}
\dfrac{d\;\!\xi}{dt}\;\! = b_{1}^{}+(b_{2}^{}- a_{1}^{})\xi + b_{_0}\;\!\theta -
a_{2}^{}\;\!\xi^2 - a_{_0}\xi\;\!\theta
\, \equiv \, \widehat{\Xi}(\xi,\theta),
\\[3.75ex]
\dfrac{d\;\!\theta}{dt}\;\! ={}- a_{3}^{} - a_{4}^{}\;\!\xi -
a_{1}^{}\;\!\theta  -  a_{2}^{}\;\!\xi\;\!\theta  - a_{_0}\;\!\theta^2
\, \equiv\;\! \widehat{\Theta}(\xi,\theta),
\ \ 
\text{where} \,\  |a_3^{}|+|a_4^{}|\ne 0;
\end{array}
$
\hfill\mbox{} (7.19)
\\[2.25ex]
and the second projectively  reduced system of the P-nonsingular system (7.16)  is
\\[2ex]
\mbox{}\hfill                                         %(7.20)
$
\begin{array}{l}
\dfrac{d\;\!\eta}{dt}\;\! = 
{}- a_{4}^{}  - b_{2}^{}\;\!\eta - a_{3}^{}\;\!\zeta - b_{_0}\;\!\eta^2
-b_{1}^{}\;\!\eta\;\!\zeta
\;\!\equiv\;\! \widehat{H}(\eta,\zeta),
\\[3.75ex]
\dfrac{d\;\!\zeta}{dt}\;\! =\;\! a_{2}^{} + a_{_0}\;\!\eta + (a_{1}^{} - b_{2}^{})\zeta
- b_{_0}\eta\;\!\zeta - b_{1}^{}\zeta^2
\;\!\equiv\;\! 
\widehat{Z} (\eta,\zeta),
\ \ 
\text{where} \,\  |a_3^{}|+|a_4^{}|\ne 0.
\end{array}
$
\hfill\mbox{} (7.20)
\\[2.75ex]
\indent
{\bf 7.4.}
The differential system
\\[2ex]
\mbox{}\hfill                                         %(7.21)
$
\begin{array}{l}
\dfrac{dx}{dt}=a_{_0}+a_1^{}x+a_2^{}y+x(c_{_0}+c_1^{}x+c_2^{}y)\equiv X(x,y),
\\[3ex]
\dfrac{dy}{dt}=b_{_0}+b_1^{}x+b_2^{}y+y(c_{_0}+c_1^{}x+c_2^{}y)\equiv Y(x,y),
\end{array}
$
\hfill\mbox{} (7.21)
\\[2.5ex]
where real constants $a_i ^ {}, \ b_i ^ {}, \ c_i ^ {}, \ i=1,2,3,$ 
\vspace{0.5ex}
such that polynomials $X $ and $Y $ are relatively prime,  is said to be [10] {\it Jacobi's system}.
Since for the differential system (7.21) 
\\[1.25ex]
\mbox{}\hfill
$
X_2^{}(x,y)=x(c_1^{}x+c_2^{}y),
\quad 
Y_2^{}(x,y)=y(c_1^{}x+c_2^{}y)
$ 
\ for all 
$
(x,y)\in\R^2,
\hfill
$
\\[0.5ex]
we obtain
\vspace{0.35ex}

{\bf Property 7.4.} 
{\it 
The system {\rm (7.16)} is projectively singular if and only if this system is Jacobi's system.
}
\vspace{0.75ex}

{\bf 7.5.}
The differential system
\\[2.25ex]
\mbox{}\hfill                                         %(7.22)
$
\dfrac{dx}{dt}=A_{m}^{}(x,y)+x\;\!C_{n-1}^{}(x,y)\equiv X(x,y),
\ \ \ \
\dfrac{dy}{dt}=B_{m}^{}(x,y)+y\;\!C_{n-1}^{}(x,y)\equiv Y(x,y),
$
\hfill\mbox{} (7.22)
\\[2.75ex]
where 
\vspace{0.5ex}
$A_m^{},\ B_m^{},$ and $C_{n-1}^{}$ are homogeneous polynomials of degrees $m$ and $n-1,\ m\leq n-1,$
on variables $x,\ y,$  such that polynomials $X$ and $Y $ are relatively prime, 
is said to be [10] {\it Darboux's system}.
Since $m\leq n-1,$ we see that for the  differential system (7.22)  
\\[1.25ex]
\mbox{}\hfill
$
X_n^{}(x,y)=x\;\!C_{n-1}^{}(x,y),
\quad 
Y_n^{}(x,y)=y\;\!C_{n-1}^{}(x,y)
$ 
\ for all 
$
(x,y)\in\R^2,
\hfill
$
\\[0.5ex]
and we have
\vspace{0.15ex}

{\bf Property 7.5.} 
{\it 
Darboux's system {\rm (7.22)} is projectively singular.
}
\\[3.25ex]
\centerline{
{\bf  8. Projective atlas of trajectories for differential systems
}
}
\\[1ex]
\indent
The behaviour of trajectories of system (D) on the projective sphere $ \P {\mathbb S} (x, y)$ 
\vspace{0.25ex}
is defined by the behaviour of trajectories of systems (D), (6.3), (6.6) 
\vspace{0.35ex}
in final parts of projective phase planes 
$\R \P (x, y), \ \R \P (\xi, \theta), \ \R \P (\eta, \zeta),$ correspondingly. 
\vspace{0.5ex}

Using the connections of the projective circles 
\vspace{0.35ex}
$ \P\K (x, y), \ \P\K (\xi, \theta), \ \P\K (\eta, \zeta), $ specified in Fig. 5.1, we get
by the behaviour of trajectories of systems (D), (6.3), (6.6) in final parts of their projective phase planes 
we can define the motion of trajectories for each of these systems on all its projective phase plane.

Trajectories of systems (D), (6.3), (6.6) on the 
\vspace{0.35ex}
projective circles of the atlas $ (\P\K (x, y), \P\K (\xi, \theta), \P\K (\eta, \zeta))$ 
\vspace{0.35ex}
are called {\it the projective atlas} of trajectories for system (D).

Then trajectories of systems (6.3), (6.6), (D) 
\vspace{0.35ex}
on the projective circles of the atlas 
$(\P\K (\xi, \theta), \P\K (\eta, \zeta), \P\K (x, y)) $ 
\vspace{0.35ex}
are the projective atlas of trajectories for system (6.3); 
and trajectories of systems (6.6), (D), (6.3) 
\vspace{0.35ex}
on the projective circles of the atlas 
$(\P\K (\eta, \zeta), \P\K (x, y), \P\K (\xi, \theta)) $ are the projective atlas of trajectories for system (6.6).
\vspace{0.35ex}

Thus, if we constructed the projective atlas of trajectories for one of systems (D), (6.3), (6.6), then
we receive the projective atlases of trajectories for two other systems.
\vspace{0.5ex}

{\bf Examples}
%\vspace{0.35ex}

{\bf 8.1.}
\vspace{0.25ex}
On Fig. 8.1 the projective atlas of trajectories $a _ {_ 0} y-b _ {_ 0} x=C $ 
for systems (7.8) at $a _ {_ 0} =0, \, b _ {_ 0} \ne0$ is constructed.  
\vspace{0.25ex}
On Fig. 8.2 the projective atlas of trajectories $a _ {_ 0} y-b _ {_ 0} x=C $ 
for systems (7.8) at $a _ {_ 0} =b _ {_ 0} \ne 0$ is constructed. 
\vspace{0.25ex}
The direction of movement along trajectories is defined by the constant vector 
$\vec{a}(x,y)=a_{_0}\vec{\imath}+b_{_0}\vec{\jmath}$ for all $(x,y)\in \R^2.$
\\[2.75ex]
\mbox{}\hfill
{\unitlength=1mm
\begin{picture}(42,42)
\put(0,0){\includegraphics[width=42mm,height=42mm]{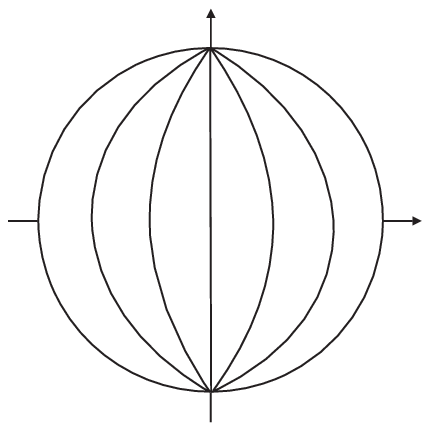}}
\put(18,41){\makebox(0,0)[cc]{ $y$}}
\put(40.2,18.2){\makebox(0,0)[cc]{ $x$}}
\end{picture}}
\quad
{\unitlength=1mm
\begin{picture}(42,42)
\put(0,0){\includegraphics[width=42mm,height=42mm]{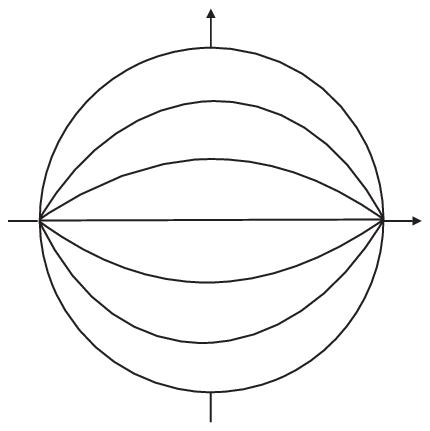}}
\put(18,41){\makebox(0,0)[cc]{ $\theta$}}
\put(40.2,17.8){\makebox(0,0)[cc]{ $\xi$}}
\put(21,-6){\makebox(0,0)[cc]{Fig. 8.1}}
\end{picture}}
\quad
{\unitlength=1mm
\begin{picture}(42,42)
\put(0,0){\includegraphics[width=42mm,height=42mm]{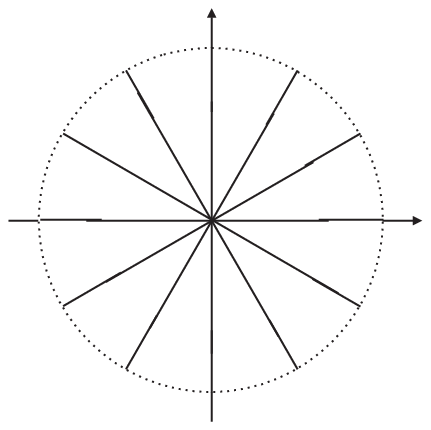}}
\put(18,41){\makebox(0,0)[cc]{ $\zeta$}}
\put(40.2,18){\makebox(0,0)[cc]{ $\eta$}}
\end{picture}}
\hfill\mbox{}
\\[7.75ex]
\mbox{}\hfill
{\unitlength=1mm
\begin{picture}(42,42)
\put(0,0){\includegraphics[width=42mm,height=42mm]{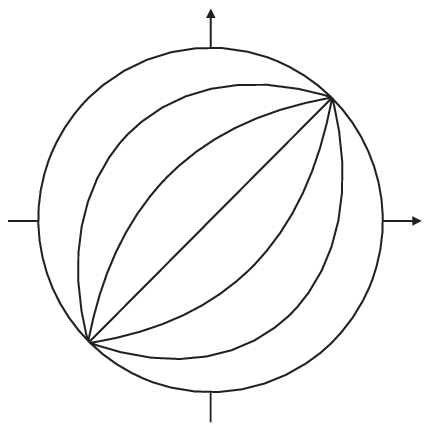}}
\put(18,41){\makebox(0,0)[cc]{ $y$}}
\put(40.2,18.2){\makebox(0,0)[cc]{ $x$}}
\end{picture}}
\quad
{\unitlength=1mm
\begin{picture}(42,42)
\put(0,0){\includegraphics[width=42mm,height=42mm]{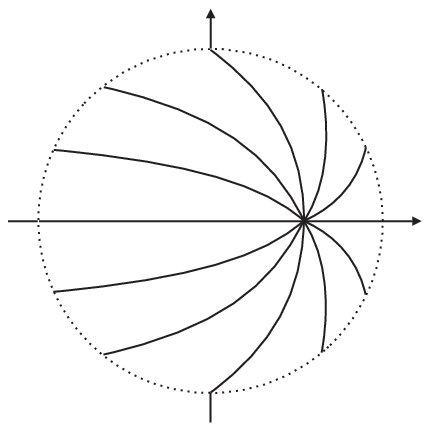}}
\put(18,41){\makebox(0,0)[cc]{ $\theta$}}
\put(40.2,17.8){\makebox(0,0)[cc]{ $\xi$}}
\put(21,-6){\makebox(0,0)[cc]{Fig. 8.2}}
\end{picture}}
\quad
{\unitlength=1mm
\begin{picture}(42,42)
\put(0,0){\includegraphics[width=42mm,height=42mm]{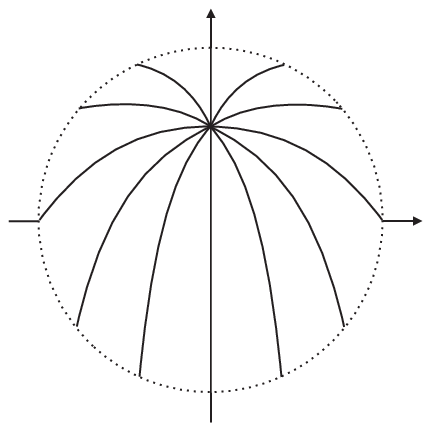}}
\put(18,41){\makebox(0,0)[cc]{ $\zeta$}}
\put(40.2,18){\makebox(0,0)[cc]{ $\eta$}}
\end{picture}}
\hfill\mbox{}
\\[7ex]
\indent
{\bf 8.2.}
On Fig. 8.3 --- 8.8 the projective atlases of trajectories for the linear autonomous system (7.11) 
in the cases specified in Table 8.1 
are constructed (correspond to the cases reduced in Theorem 2.15 from [11, p. 39]).
%\marginpar{[103]}

%\newpage

\mbox{}
\\[-3ex]
\mbox{}\hfill
{
\renewcommand{\arraystretch}{2.5} %регулировка высоты; 
\newcommand{\PreserveBackslash}[1]{\let\temp=\\#1\let\\=\temp}
\let\PBS=\PreserveBackslash
\begin{tabular}{l}
\;\;\;Table 8.1
\\
\begin{tabular}                   
              {|>{\PBS\centering\hspace{0pt}}p{3.5cm}
               |>{\PBS\centering\hspace{0pt}}p{5.3cm}
               |>{\PBS\centering\hspace{0pt}}p{3.3cm}
               |>{\PBS\centering\hspace{0pt}}p{0.7cm}|}
\hline 
System   & Family of trajectories & Equilibrium states   &  Fig. \\ 
\hline  
$
\begin{array}{l}
\dfrac{dx}{dt}=x,
\\[-0.3ex]
\dfrac{dy}{dt}=2y 
\end{array}
$
\hfill (8.1) 
& 
$C_1^{}y+C_2^{}x^2=0$
&
\!\!\!\!\!
\begin{tabular}{l}
$O$ --- node,\\[-3ex]
$O^{(1)}$ --- saddle,\\[-3ex]
$O^{(2)}$ --- node
\end{tabular}
\hfill\mbox{}
& 8.3\\ \hline  
$
\begin{array}{l}
\dfrac{dx}{dt}=x,
\\[-0.3ex] 
\dfrac{dy}{dt}={}-y 
\end{array}
$
\hfill (8.2) 
& 
$xy=C$
&
\!\!\!\!\!
\begin{tabular}{l}
$O$ --- saddle,\\[-2.75ex]
$O^{(1)}$ --- node,\\[-2.75ex]
$O^{(2)}$ --- node
\end{tabular}
\hfill\mbox{}
& 8.4\\ \hline  
$
\begin{array}{l}
\dfrac{dx}{dt}=y,
\\[-0.3ex] 
\dfrac{dy}{dt}={}-x 
\end{array}
$
\hfill (8.3) 
& 
$x^2+y^2=C$
&
$O$ --- center
\hfill\mbox{}
& 8.5\\ \hline  
$
\begin{array}{l}
\dfrac{dx}{dt}=x-y,
\\[-0.3ex] 
\dfrac{dy}{dt}=x +y
\end{array}
$
\hfill (8.4) 
& 
$(x^2+y^2)\exp\Bigl({}-2\arctan\dfrac{y}{x}\Bigr)=C$
&
$O$ --- focus
\hfill\mbox{}
& 8.6\\ \hline  
$
\begin{array}{l}
\dfrac{dx}{dt}=x+y,
\\[-0.3ex] 
\dfrac{dy}{dt}=y
\end{array}
$
\hfill (8.5) 
& 
$y\exp\Bigl({}-\dfrac{x}{y}\Bigr)=C$
&
\mbox{}\hfill
\!\!\!\!\!\!
\begin{tabular}{l}
$O$ --- degenerate
\\[-3.75ex]
node
\\[-2.75ex]
$O^{(1)}\!$ --- saddle-node
\end{tabular}
\hfill\mbox{}
& 8.7\\ \hline  
$
\begin{array}{l}
\dfrac{dx}{dt}=x,
\\[-0.3ex] 
\dfrac{dy}{dt}=y 
\end{array}
$
\hfill (8.6) 
&$C_1^{}y+C_2^{}x=0$& 
\!\!\!\!\!
\begin{tabular}{l}
$O$ --- dicritical
\\[-4ex]
node
\end{tabular}
\hfill\mbox{}
& 8.8\\ \hline  
\end{tabular}
\end{tabular}
}
\hfill\mbox{}
\\[6ex]
\mbox{}\hfill
{\unitlength=1mm
\begin{picture}(42,42)
\put(0,0){\includegraphics[width=42mm,height=42mm]{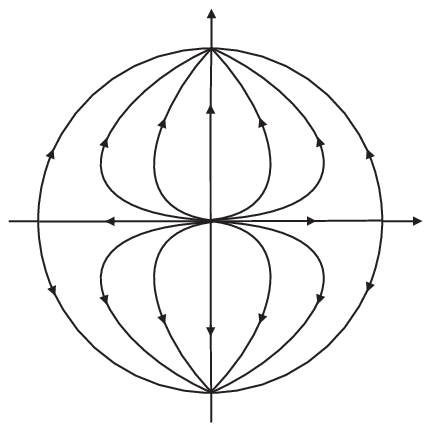}}
\put(18,41){\makebox(0,0)[cc]{ $y$}}
\put(40.2,18.5){\makebox(0,0)[cc]{ $x$}}
\end{picture}}
\quad
{\unitlength=1mm
\begin{picture}(42,42)
\put(0,0){\includegraphics[width=42mm,height=42mm]{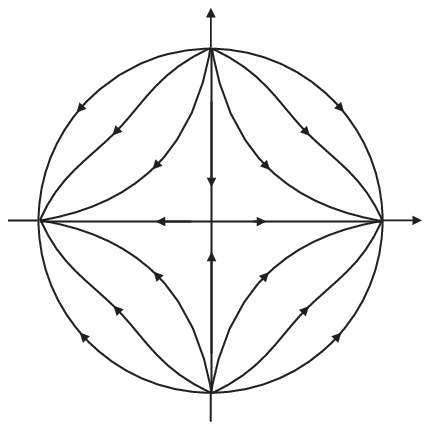}}
\put(18,41){\makebox(0,0)[cc]{ $\theta$}}
\put(40.2,18){\makebox(0,0)[cc]{ $\xi$}}
\put(21,-7){\makebox(0,0)[cc]{Fig. 8.3}}
\end{picture}}
\quad
{\unitlength=1mm
\begin{picture}(42,42)
\put(0,0){\includegraphics[width=42mm,height=42mm]{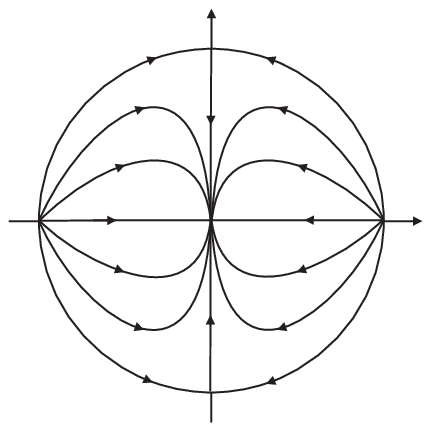}}
\put(18,41){\makebox(0,0)[cc]{ $\zeta$}}
\put(40.2,18){\makebox(0,0)[cc]{ $\eta$}}
\end{picture}}
\hfill\mbox{}
\\[10ex]
\mbox{}\hfill
{\unitlength=1mm
\begin{picture}(42,42)
\put(0,0){\includegraphics[width=42mm,height=42mm]{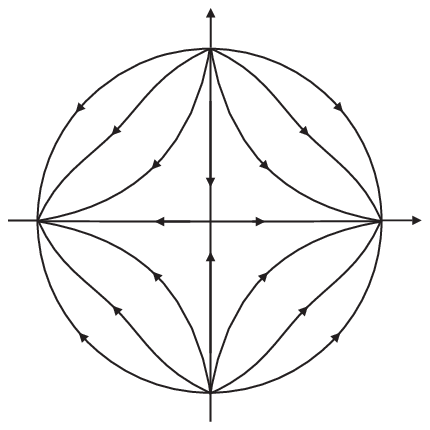}}
\put(18,41){\makebox(0,0)[cc]{ $y$}}
\put(40.2,18.5){\makebox(0,0)[cc]{ $x$}}
\end{picture}}
\quad
{\unitlength=1mm
\begin{picture}(42,42)
\put(0,0){\includegraphics[width=42mm,height=42mm]{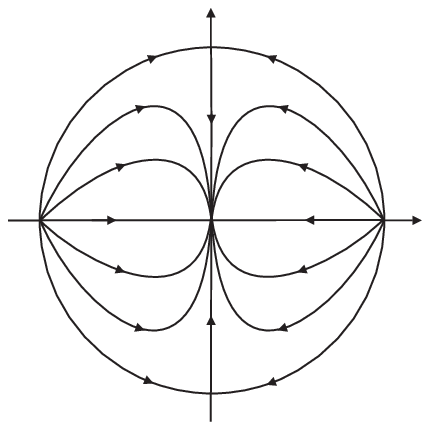}}
\put(18,41){\makebox(0,0)[cc]{ $\theta$}}
\put(40.2,18){\makebox(0,0)[cc]{ $\xi$}}
\put(21,-7){\makebox(0,0)[cc]{Fig. 8.4}}
\end{picture}}
\quad
{\unitlength=1mm
\begin{picture}(42,42)
\put(0,0){\includegraphics[width=42mm,height=42mm]{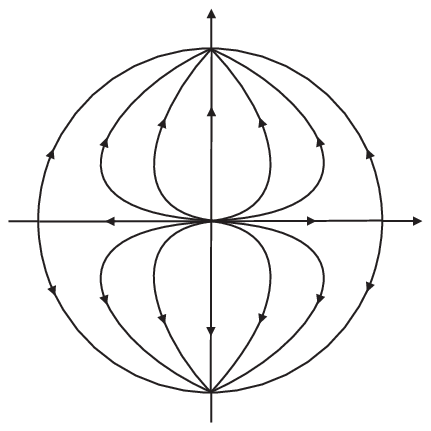}}
\put(18,41){\makebox(0,0)[cc]{ $\zeta$}}
\put(40.2,18){\makebox(0,0)[cc]{ $\eta$}}
\end{picture}}
\hfill\mbox{}
\\[13ex]
\mbox{}\hfill
{\unitlength=1mm
\begin{picture}(42,42)
\put(0,0){\includegraphics[width=42mm,height=42mm]{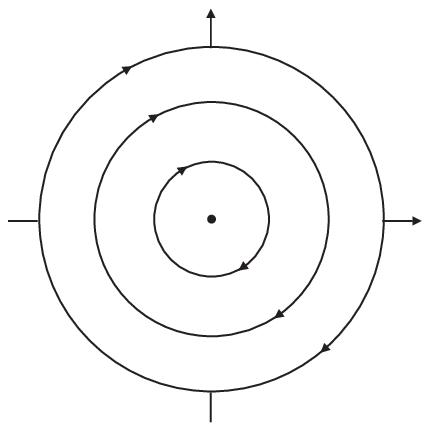}}
\put(18,41){\makebox(0,0)[cc]{ $y$}}
\put(40.2,18.2){\makebox(0,0)[cc]{ $x$}}
\end{picture}}
\quad
{\unitlength=1mm
\begin{picture}(42,42)
\put(0,0){\includegraphics[width=42mm,height=42mm]{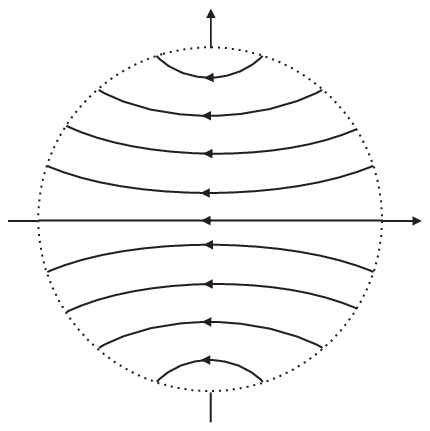}}
\put(18,41){\makebox(0,0)[cc]{ $\theta$}}
\put(40.2,17.8){\makebox(0,0)[cc]{ $\xi$}}
\put(21,-7){\makebox(0,0)[cc]{Fig. 8.5}}
\end{picture}}
\quad
{\unitlength=1mm
\begin{picture}(42,42)
\put(0,0){\includegraphics[width=42mm,height=42mm]{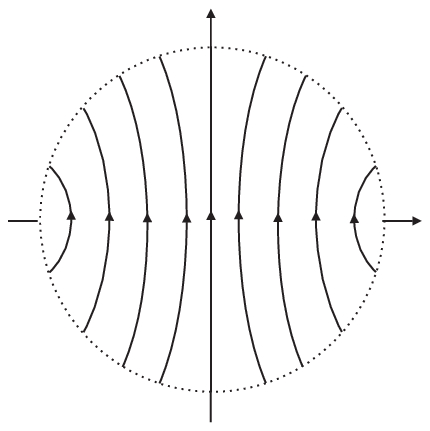}}
\put(18,41){\makebox(0,0)[cc]{ $\zeta$}}
\put(40.2,18){\makebox(0,0)[cc]{ $\eta$}}
\end{picture}}
\hfill\mbox{}
\\[13ex]
\mbox{}\hfill
{\unitlength=1mm
\begin{picture}(42,42)
\put(0,0){\includegraphics[width=42mm,height=42mm]{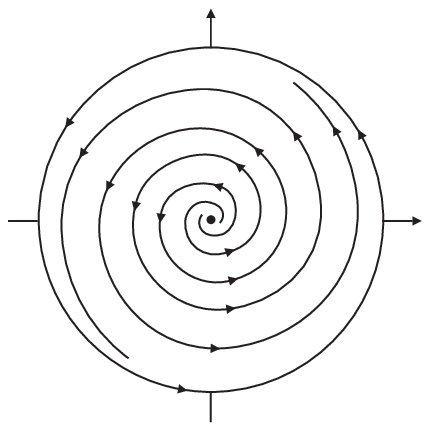}}
\put(18,41){\makebox(0,0)[cc]{ $y$}}
\put(40.2,18.2){\makebox(0,0)[cc]{ $x$}}
\end{picture}}
\quad
{\unitlength=1mm
\begin{picture}(42,42)
\put(0,0){\includegraphics[width=42mm,height=42mm]{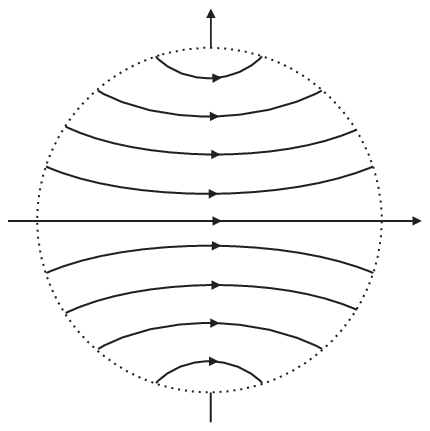}}
\put(18,41){\makebox(0,0)[cc]{ $\theta$}}
\put(40.2,17.8){\makebox(0,0)[cc]{ $\xi$}}
\put(21,-7){\makebox(0,0)[cc]{Fig. 8.6}}
\end{picture}}
\quad
{\unitlength=1mm
\begin{picture}(42,42)
\put(0,0){\includegraphics[width=42mm,height=42mm]{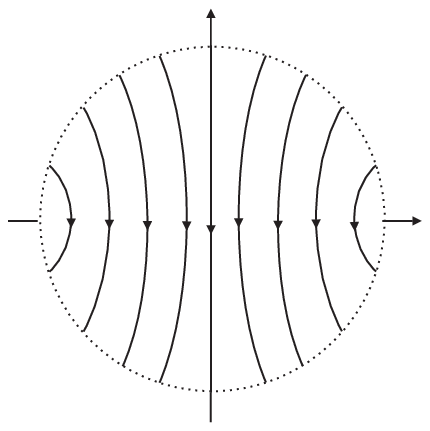}}
\put(18,41){\makebox(0,0)[cc]{ $\zeta$}}
\put(40.2,18){\makebox(0,0)[cc]{ $\eta$}}
\end{picture}}
\hfill\mbox{}
\\[13ex]
\mbox{}\hfill
{\unitlength=1mm
\begin{picture}(42,42)
\put(0,0){\includegraphics[width=42mm,height=42mm]{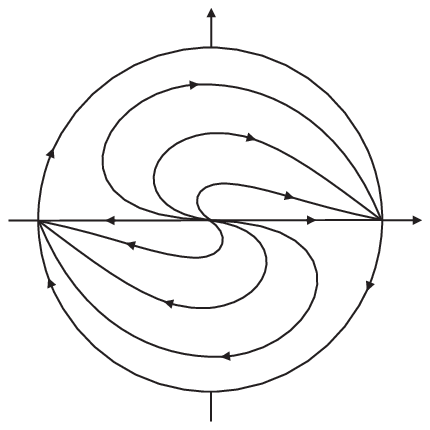}}
\put(18,41){\makebox(0,0)[cc]{ $y$}}
\put(40.2,18.2){\makebox(0,0)[cc]{ $x$}}
\end{picture}}
\quad
{\unitlength=1mm
\begin{picture}(42,42)
\put(0,0){\includegraphics[width=42mm,height=42mm]{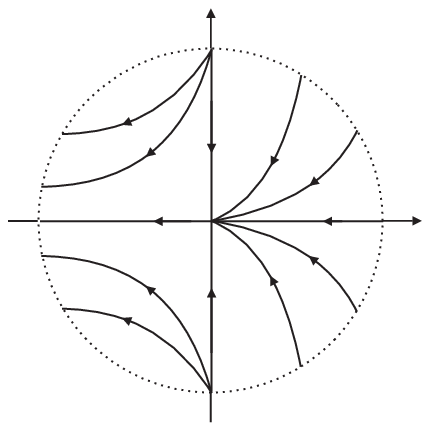}}
\put(18,41){\makebox(0,0)[cc]{ $\theta$}}
\put(40.2,17.8){\makebox(0,0)[cc]{ $\xi$}}
\put(21,-7){\makebox(0,0)[cc]{Fig. 8.7}}
\end{picture}}
\quad
{\unitlength=1mm
\begin{picture}(42,42)
\put(0,0){\includegraphics[width=42mm,height=42mm]{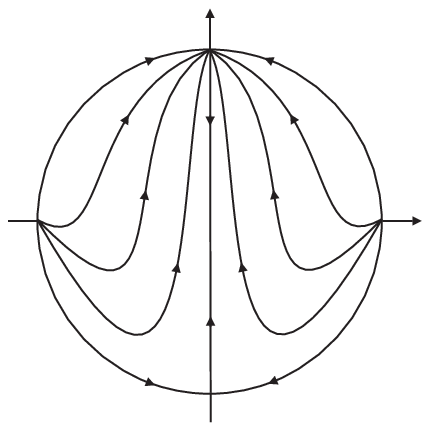}}
\put(18,41){\makebox(0,0)[cc]{ $\zeta$}}
\put(40.2,18){\makebox(0,0)[cc]{ $\eta$}}
\end{picture}}
\hfill\mbox{}
\\[13ex]
\mbox{}\hfill
{\unitlength=1mm
\begin{picture}(42,42)
\put(0,0){\includegraphics[width=42mm,height=42mm]{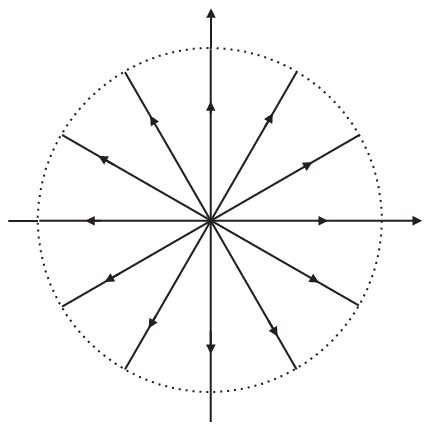}}
\put(18,41){\makebox(0,0)[cc]{ $y$}}
\put(40.2,18.2){\makebox(0,0)[cc]{ $x$}}
\end{picture}}
\quad
{\unitlength=1mm
\begin{picture}(42,42)
\put(0,0){\includegraphics[width=42mm,height=42mm]{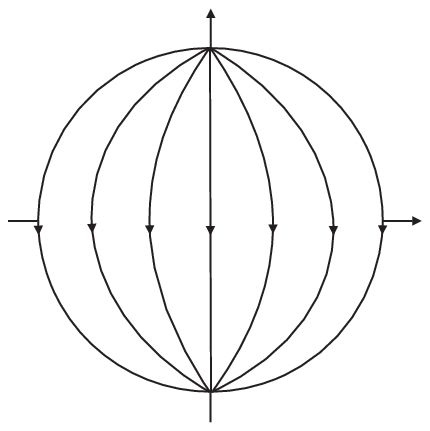}}
\put(18,41){\makebox(0,0)[cc]{ $\theta$}}
\put(40.2,17.8){\makebox(0,0)[cc]{ $\xi$}}
\put(21,-7){\makebox(0,0)[cc]{Fig. 8.8}}
\end{picture}}
\quad
{\unitlength=1mm
\begin{picture}(42,42)
\put(0,0){\includegraphics[width=42mm,height=42mm]{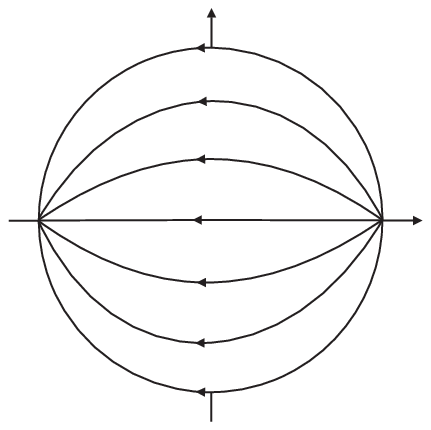}}
\put(18,41){\makebox(0,0)[cc]{ $\zeta$}}
\put(40.2,18){\makebox(0,0)[cc]{ $\eta$}}
\end{picture}}
\hfill\mbox{}
\\[7.75ex]
\indent
{\bf 8.3.}
\vspace{0.5ex}
Using the topological pictures of trajectories on the projective circle $\P\K (x, y)$ 
\linebreak
\text{[12, pp. 61 -- 65]} for the differential system $(i =\sqrt {{}-1}\,, \ q=4,$ and $q=5)$ 
\\[2.25ex]
\mbox{}\hfill        %(8.7)
$
2\,\dfrac{dx}{dt}=2y+i\;\!(x-i\;\!y)^q-i\;\!(x+i\;\!y)^q,
\quad \ 
2\,\dfrac{dy}{dt}={}-2x+(x-i\;\!y)^q+ (x+i\;\!y)^q,
$
\hfill (8.7)
\\[2.5ex]
we obtain 
\vspace{0.35ex}
on Fig. 8.9 the projective atlas of trajectories for system (8.7) at $q=4$ is constructed  
and 
on Fig. 8.10 the projective atlas of trajectories for system (8.7) at $q=5$ is constructed.  
\\[4ex]
\mbox{}\hfill
{\unitlength=1mm
\begin{picture}(42,42)
\put(0,0){\includegraphics[width=42mm,height=42mm]{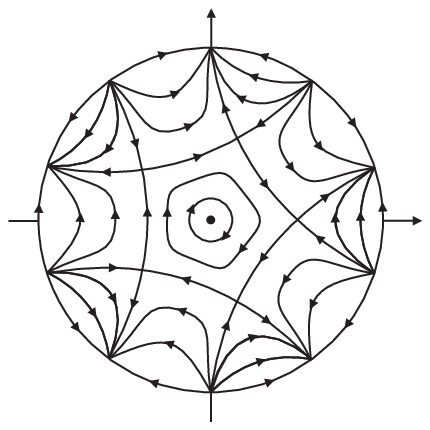}}
\put(18,41){\makebox(0,0)[cc]{ $y$}}
\put(40.2,18.2){\makebox(0,0)[cc]{ $x$}}
\end{picture}}
\quad
{\unitlength=1mm
\begin{picture}(42,42)
\put(0,0){\includegraphics[width=42mm,height=42mm]{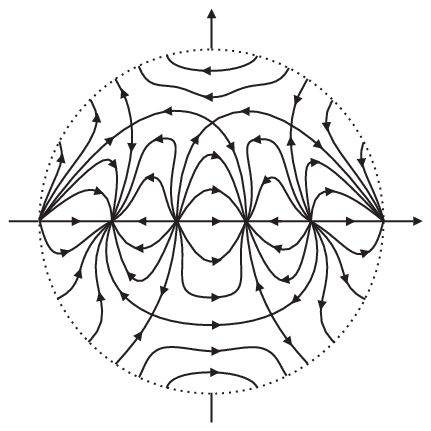}}
\put(18,41){\makebox(0,0)[cc]{ $\theta$}}
\put(40.2,17.8){\makebox(0,0)[cc]{ $\xi$}}
\put(21,-7){\makebox(0,0)[cc]{Fig. 8.9}}
\end{picture}}
\quad
{\unitlength=1mm
\begin{picture}(42,42)
\put(0,0){\includegraphics[width=42mm,height=42mm]{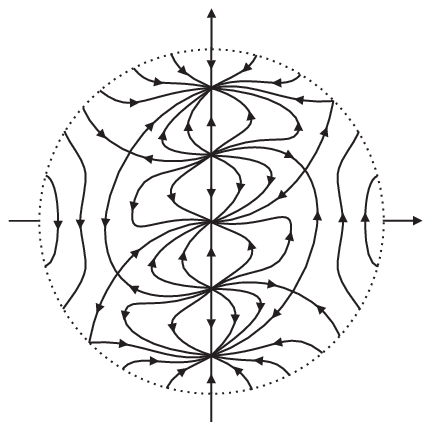}}
\put(18,41){\makebox(0,0)[cc]{ $\zeta$}}
\put(40.2,18){\makebox(0,0)[cc]{ $\eta$}}
\end{picture}}
\hfill\mbox{}
\\[10ex]
\mbox{}\hfill
{\unitlength=1mm
\begin{picture}(42,42)
\put(0,0){\includegraphics[width=42mm,height=42mm]{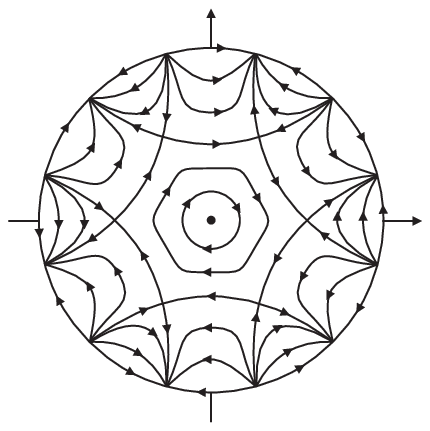}}
\put(18,41){\makebox(0,0)[cc]{ $y$}}
\put(40.2,18.2){\makebox(0,0)[cc]{ $x$}}
\end{picture}}
\quad
{\unitlength=1mm
\begin{picture}(42,42)
\put(0,0){\includegraphics[width=42mm,height=42mm]{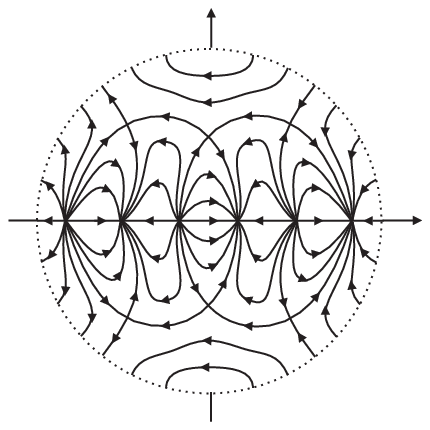}}
\put(18,41){\makebox(0,0)[cc]{ $\theta$}}
\put(40.2,17.8){\makebox(0,0)[cc]{ $\xi$}}
\put(21,-7){\makebox(0,0)[cc]{Fig. 8.10}}
\end{picture}}
\quad
{\unitlength=1mm
\begin{picture}(42,42)
\put(0,0){\includegraphics[width=42mm,height=42mm]{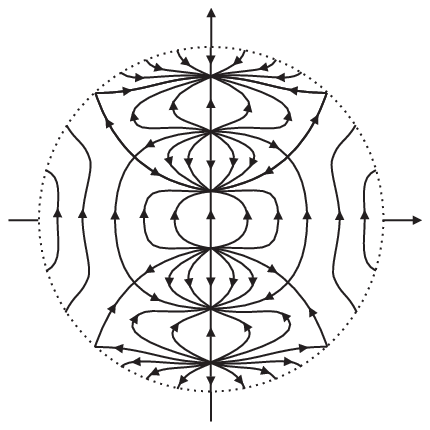}}
\put(18,41){\makebox(0,0)[cc]{ $\zeta$}}
\put(40.2,18){\makebox(0,0)[cc]{ $\eta$}}
\end{picture}}
\hfill\mbox{}
\\[9.75ex]
\centerline{
{\bf\large \S\;\!3. Trajectories on Poincar\'{e}'s sphere}
}
\\[2ex]
\centerline{
{\bf  9. Trajectories on projective phase plane
}
}
\\[1.5ex]
\indent
The qualitative behaviour of trajectories of system (D) on the projective phase plane is  defined 
by equilibrium states, limit cycles, and also by specific properties of trajectories of concrete system 
(symmetry of the phase directional field, presence of the contact points, known trajectories etc.). 
\vspace{0.25ex}

First of all we will formulate the properties connected 
by that image on the projective sphere of infinitely removed straight line of the projective phase plane 
\vspace{0.75ex}
is equator of the sphere.

{\bf Property 9.1.}
{\it 
The following statements are equivalent}\;\!:
\vspace{0.35ex}

1. 
{\it
The infinitely remote straight line of the projective phase plane $\R\P (x, y)$ consists of trajectories
of system} (D);
\vspace{0.35ex}

2. 
{\it 
The polynomial $W ^ {} _ n $ is not identical zero on $\R^2;$
} 
\vspace{0.35ex}

3. 
{\it
The system} (D) {\it is projectively nonsingular}\;\!;
\vspace{0.35ex}

4. 
{\it 
The first projectively reduced system for the system} (D) {\it  is the system} (7.1);
\vspace{0.35ex}

5. 
{\it
The straight line $ \theta=0$ 
\vspace{0.25ex}
consists of trajectories of the first projectively reduced system for the differential system} (D);
\vspace{0.35ex}

6. 
{\it
The second projectively reduced system for the system} (D) {\it  is the system} (7.2);
\vspace{0.35ex}

7. 
{\it
The straight line $\eta=0$ 
\vspace{0.25ex}
consists of trajectories of the second projectively reduced system for the differential system} (D).
\vspace{1ex}

{\bf Property 9.2.}
{\it 
The following statements are equivalent}\;\!:
\vspace{0.35ex}

1. 
{\it
The infinitely remote straight line of the projective phase plane $\R\P (x, y)$ 
does not consists of trajectories of system} (D);
\vspace{0.5ex}

2. 
{\it
The polynomial} $W_n^{}(x,y)=0$  for all $(x,y)\in \R^2;$
\vspace{0.5ex}

3. 
{\it
The system} (D) {\it is projectively singular}\;\!;
\vspace{0.35ex}

4. 
{\it 
The first projectively reduced system for the system} (D) {\it is the system} (7.3);
\vspace{0.35ex}

5.
{\it
The straight line $\theta=0$ 
\vspace{0.25ex}
does not consists of trajectories of the first projectively reduced system 
for the differential system } (D);
\vspace{0.35ex}

6. 
{\it
The second projectively reduced system for the system} (D) {\it is the system} (7.4);
\vspace{0.35ex}

7. 
{\it
The straight line $\eta=0$ 
\vspace{0.25ex}
does not consists of trajectories of the second projectively reduced system 
for the differential system } (D).
\vspace{0.75ex}

Using Properties 9.3 -- 9.6, we can establish the existence
\vspace{0.25ex}
of infinitely removed equilibrium states on the  projective phase plane of system (D).
\vspace{0.75ex}

{\bf Property 9.3.}
{\it 
The following statements are equivalent}\;\!:
\vspace{0.25ex}

1. 
{\it
The projectively nonsingular system} (D) 
\vspace{0.25ex}
{\it has infinitely remote equilibrium state,
which is lying on <<extremities>> of the straight line} $y=ax;$
\vspace{0.35ex}

2. 
{\it 
The point} $ (a, 0) $ {\it is an equilibrium state of the first projectively reduced system} (7.1);
\vspace{0.35ex}

3. 
$ \xi=a $ {\it is a solution of the equation} $W_n ^ {} (1, \xi) =0. $
\vspace{1ex}

{\bf Property 9.4.}
{\it 
The following statements are equivalent}\;\!:
\vspace{0.35ex}

1. 
{\it 
The projectively nonsingular system} (D) 
\vspace{0.25ex}
{\it has infinitely remote equilibrium state, which is lying on <<extremities>> of the straight line} $x=ay;$
\vspace{0.35ex}

2. 
{\it 
The point} $ (0, a) $ 
\vspace{0.35ex}
{\it is an equilibrium state of the second projectively reduced system} (7.2);

3. 
$\zeta=a$ {\it is a solution of the equation} $W_n^{}(\zeta,1)=0.$
\vspace{1ex}

{\bf Property 9.5.}
{\it 
The following statements are equivalent}\;\!:
\vspace{0.35ex}

1. 
{\it 
The projectively singular system} (D) 
\vspace{0.25ex}
{\it has infinitely remote equilibrium state, which is lying on <<extremities>> of the straight line} $y=ax; $
\vspace{0.35ex}

2. 
{\it 
The point} $ (a, 0) $ {\it is an equilibrium state of the first projectively reduced system} (7.3);
\vspace{0.5ex}

3. 
$ \xi=a $ {\it is a solution of the system of equations} $X_n ^ {} (1, \xi) =0,\ W _ {n-1} ^ {} (1, \xi) =0. $
\vspace{1ex}

{\bf Property 9.6.}
{\it 
The following statements are equivalent}\;\!:
\vspace{0.35ex}

1. 
{\it 
The projectively singular system} (D) 
\vspace{0.25ex}
{\it has infinitely remote equilibrium state, which is lying on <<extremities>> of the straight line} $x=ay;$
\vspace{0.35ex}

2. 
{\it 
The point} $(0, a)$ 
\vspace{0.5ex}
{\it is an equilibrium state of the second projectively reduced system} (7.4);

3. 
$\zeta=a$ {\it is a solution of the system of equations} $Y_n^{}(\zeta,1)=0,\ W_{n-1}^{}(\zeta,1)=0.$
%\vspace{0.75ex}

\newpage

By Properties 9.7 -- 9.12, we can establish the type 
\vspace{0.25ex}
of infinitely removed equilibrium states on the  projective phase plane of system (D).
\vspace{0.75ex}

{\bf Property 9.7.}
\vspace{0.15ex}
{\it
Suppose an infinitely removed equilibrium state of the projectively non\-sin\-gu\-lar system} (D) 
{\it is lying on <<extremities>> of the straight line $y=ax.$ 
\vspace{0.15ex}
Then the type of this infinitely removed equilibrium state 
\vspace{0.15ex}
{\rm(}to within a direction of movement along trajectories  adjoining it{\rm)}
\vspace{0.15ex}
is the same as the type of the equilibrium state $ (a, 0)$ for the first projectively reduced system} (7.1).
\vspace{0.75ex}

{\bf Property 9.8.}
\vspace{0.15ex}
{\it
Suppose an infinitely removed equilibrium state of the projectively non\-sin\-gu\-lar system} (D) 
{\it is lying on <<extremities>> of the straight line $x=ay.$ 
\vspace{0.15ex}
Then the type of this infinitely removed equilibrium state 
\vspace{0.15ex}
{\rm(}to within a direction of movement along trajectories  adjoining it{\rm)}
\vspace{0.15ex}
is the same as the type of the equilibrium state 
$(0, a)$ for the second projectively reduced system} (7.2).
\vspace{0.75ex}

{\bf Property 9.9.}
\vspace{0.15ex}
{\it
Suppose an infinitely removed equilibrium state of the projectively sin\-gu\-lar system} (D) 
\vspace{0.15ex}
{\it is lying on <<extremities>> of the straight line $y=ax.$ 
Then the type of this infinitely removed equilibrium state 
\vspace{0.15ex}
{\rm(}to within a direction of movement along trajectories  adjoining it{\rm)}
is\! the\! same\! as\! the\! type\! of\! the\! equilibrium\! state 
\vspace{0.75ex}
$\!(a, 0)\!\!$ of the\! first\! projectively\! reduced\! system}\! (7.3).

{\bf Property 9.10.}
\vspace{0.15ex}
{\it
Suppose an infinitely removed equilibrium state of the projectively sin\-gu\-lar system} (D) 
{\it is lying on <<extremities>> of the straight line $x=ay.$ 
\vspace{0.15ex}
Then the type of this infinitely removed equilibrium state 
\vspace{0.15ex}
{\rm(}to within a direction of movement along trajectories  adjoining it{\rm)}
\vspace{0.15ex}
is the same as the type of the equilibrium state 
$(0, a)$ for the second projectively reduced system} (7.4).
\vspace{0.75ex}

We say that points of the projective phase plane of system (D) are {\it regular points} of system (D) 
if they are not the equilibrium states of this system. Thus, each point of the projective phase plane 
is a regular point or an equilibrium state. Each trajectory is an equilibrium state or consists of regular points.
\vspace{0.25ex}

If the system (D) is projectively nonsingular, then  infinitely removed straight line of the projective phase plane 
consists of trajectories (Property 9.1) and among these trajectories can be a finite number of equilibrium states 
(since the right members $X$ and $Y$ of the differential system (D) are polynomials, we see that 
the number of equilibrium states is finite). 
Therefore the behaviour of trajectories of the projectively nonsingular system (D) 
in an neighbourhood of infinitely removed straight line is defined by the equilibrium states lying on this line.
\vspace{0.25ex}

If the system (D) is projectively singular, then infinitely removed straight line of the projective phase plane 
does not consist of trajectories (Property 9.2), but on this line can lie a finite number of equilibrium states.
For research of behaviour of trajectories for the projectively singular system (D), 
among regular infinitely removed points we will select points such that 
trajectories on the projective phase plane concern the infinitely removed straight line in these points. 
We say that these points are {\it equatorial contact points} of the projectively singular system (D). 
It can be assumed that points in which trajectories of system (D) concern the infinitely removed straight line are
preimages of points of Poincar\'{e}'s sphere  in which images of trajectories concern equator.
Let's find equatorial contact points of the projectively singular system (D) on 
the basis of contact points of the coordinate axes of the P-reduced systems.
\vspace{0.25ex}

A point in which a trajectory of system  (D) 
\vspace{0.15ex}
concerns the curve $\gamma$ (this curve is lying on the phase plane $Oxy)$ 
\vspace{0.15ex}
is said to be {\it a contact point} of the curve $\gamma.$
In particular, each point of the curve $ \gamma $ is a contact point 
\vspace{0.75ex}
if and only if  the curve $\gamma $ is a trajectory of system (D).

{\bf Property 9.11.}
{\it 
The following statements are equivalent}\;\!:
\vspace{0.15ex}

1. 
{\it
The infinitely removed point, 
\vspace{0.15ex}
which is lying on <<extremities>> of the straight line $y=ax,$ 
is an equatorial contact point of the projectively singular system} (D); 
\vspace{0.35ex}

2. 
{\it 
The point} $ (a, 0) $ {\it is a contact point of the straight line} $ \theta=0$ 
\vspace{0.15ex}
{\it  for the first projectively reduced system } (7.3);
\vspace{0.25ex}

3. 
$\xi=a$ {\it is a solution of the equation} $X_n^{}(1,\xi)=0,$ and $W_{n-1}^{}(1,a)\ne0.$
\vspace{0.5ex}

{\bf Property 9.12.}
{\it 
The following statements are equivalent}\;\!:
\vspace{0.2ex}

1. 
{\it 
The infinitely removed point, 
\vspace{0.15ex}
which is lying on <<extremities>> of the straight line $x=ay, $ 
is an equatorial contact point of the projectively singular system} (D); 
\vspace{0.35ex} 

2. 
{\it 
The point} $(0, a)$ {\it is a contact point of the straight line} $ \eta=0$ 
\vspace{0.15ex}
{\it for the first projectively reduced system } (7.4);
\vspace{0.35ex}

3. 
$\zeta=a$ {\it is a solution of the equation} $Y_n^{}(\zeta,1)=0,$ аnd $W_{n-1}^{}(a,1)\ne 0.$
\vspace{0.75ex}

A trajectory of system (D), 
\vspace{0.15ex}
which is passing through the contact point $A$ of the curve $\gamma,$ 
is said to be {\it the contact $A\!$-trajectory} of the curve $\gamma.$ 
\vspace{0.15ex}
{\it The equatorial contact $A\!$-trajectory} of the projectively singular system (D) is
\vspace{0.15ex}
a trajectory, which is lying in the projective phase plane and passing through the contact point $A$ of the infinitely removed straight line.
\vspace{0.75ex}

{\bf Property 9.13.}\!
\vspace{0.15ex}
{\it 
The\! equatorial\! contact $\!A\!$-trajectory of the projectively singular system}\! (D) 
{\it in an enough small neighbourhood of the infinitely removed point $A$ 
\vspace{0.15ex}
{\rm(}this point is lying on 
<<extremities>> of the straight line $y=ax)$ is located}\;\!:
\vspace{0.25ex}

{\it 
a{\rm)} in the coordinate halfplane $x> 0;$
\vspace{0.35ex}

b{\rm)} in the coordinate halfplane $x <0;$
\vspace{0.35ex}

c{\rm)} as in the coordinate halfplane $x> 0, $ and in the coordinate halfplane $x <0,$
\\[0.75ex]
if and only if the contact $A ^ {(1)}\!\!$-trajectory of the straight line} $\theta=0$ 
\vspace{0.35ex}
{\it for the first projectively reduced system } (7.3) {\it in an enough small neighbourhood of the point $A ^ {(1)} (a, 0)$ 
is located}\;\!:
\vspace{0.5ex}

{\it  
a{\rm)} in the coordinate halfplane $ \theta> 0; $
\vspace{0.35ex}

b{\rm)} in the coordinate halfplane $ \theta <0; $
\vspace{0.35ex}

c{\rm)} as in the coordinate halfplane $ \theta> 0,$ and in the coordinate halfplane $ \theta <0. $
}
\vspace{1ex}

{\bf Property 9.14.}\!
\vspace{0.25ex}
{\it 
The\! equatorial\! contact $\!A\!$-trajectory of the projectively singular system}\! (D) 
{\it in an enough small neighbourhood of the infinitely removed point $A$ 
\vspace{0.25ex}
{\rm(}this point is lying on <<extremities>> of the straight line $x=ay)$ is located}\;\!:
\vspace{0.35ex}

{\it 
a{\rm)} in the coordinate halfplane $y> 0;$
\vspace{0.35ex}

b{\rm)} in the coordinate halfplane $y <0;$
\vspace{0.35ex}

c{\rm)} as in the coordinate halfplane $x> 0,$ and in the coordinate halfplane $y <0,$
\\[0.75ex]
if and only if the contact $A^{(1)}\!\!$-trajectory of the straight line} $\eta=0$ 
\vspace{0.35ex}
for  {\it the second projectively reduced system } (7.4) {\it in an enough small neighbourhood of the point $A ^ {(2)} (0, a)$ is located}\;\!:
\vspace{0.5ex}

{\it  
a{\rm)} in the coordinate halfplane $ \eta> 0; $
\vspace{0.35ex}

b{\rm)} in the coordinate halfplane $ \eta <0; $
\vspace{0.35ex}

c{\rm)} as in the coordinate halfplane $ \eta> 0, $ and in the coordinate halfplane $ \eta <0. $
}
\vspace{0.75ex}

Disposition indications 
\vspace{0.75ex}
of equatorial contact trajectories are given in Properties 9.15, 9.16.

{\bf Property 9.15.}
{\it 
If $ \xi=a $ is a solution of the equation $X_n ^ {} (1, \xi) =0$ and besides}:
\\[1.75ex]
\mbox{}\hfill
{\it 
a{\rm)} 
$W_{n-1}^{}(1,a)\,\partial_y^{} X_n^{}(1,a)<0;$
\qquad
b{\rm)}  
$W_{n-1}^{}(1,a)\,\partial_y^{} X_n^{}(1,a)>0;$
\hfill\mbox{}
\\[2.25ex]
\mbox{}\hfill
c{\rm)}  
$
\partial_y^{} X_n^{}(1,a)=0,
\quad
W_{n-1}^{}(1,a)\,\partial_{yy}^{} X_n^{}(1,a)\ne0,
$
\hfill\mbox{}
\\[1.75ex]
then the equatorial contact $A\!\!$-trajectory of the projectively singular system} (D) 
\vspace{0.15ex}
{\it in an enough small neighbourhood of the infinitely removed point $A$ 
\vspace{0.15ex}
{\rm(}this point is lying on <<extremities>> of the straight line $y=ax)$ is located, respectively}\;\!:
\vspace{0.25ex}

{\it  
a{\rm)} in the coordinate halfplane $x> 0; $
\vspace{0.25ex}

b{\rm)} in the coordinate halfplane $x <0; $
\vspace{0.25ex}

c{\rm)} as in the coordinate halfplane $x> 0, $ and in the coordinate halfplane $x <0. $
}
\vspace{0.75ex}

On Fig. 9.1 the projective circle $\P\K (x, y)$ of the projectively singular system (D) is constructed. On this circle
the behaviour of trajectories in an neighbourhood of the regular equatorial contact point $A$ is represented,  
when the equatorial contact 
$A\!$-tra\-je\-c\-to\-ry is lying in the first coordinate quarter of the projective phase plane $\R\P (x, y). $ 

Taking into account Property 9.13, from Property 9.15 
\vspace{0.15ex}
we obtain disposition indications of trajectories (which are contact with the straight line $ \theta=0)$ 
\vspace{0.15ex}
for the first projectively  reduced system (7.3) in coordinate halfplanes. 
\\[3.75ex]
\mbox{}\hfill
{\unitlength=1mm
\begin{picture}(40,40)
\put(0,0){\includegraphics[width=40mm,height=40mm]{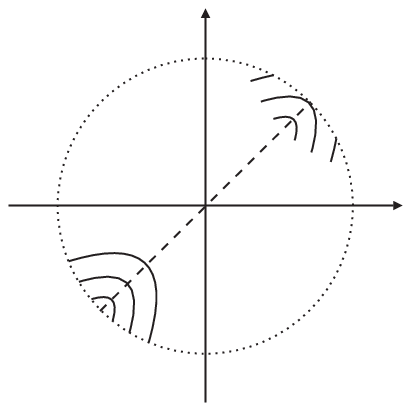}}
 
\put(22,18){\makebox(0,0)[cc]{\scriptsize $O$}}
\put(40,18){\makebox(0,0)[cc]{\scriptsize $x$}}
\put(22,40){\makebox(0,0)[cc]{\scriptsize $y$}}
\put(32.5,32.5){\makebox(0,0)[cc]{\scriptsize $A^{\star}$}}
\put(8,8){\makebox(0,0)[cc]{\scriptsize $A_{\star}^{}$}}

\put(20,-6){\makebox(0,0)[cc]{\rm Fig. 9.1}}
\end{picture}}
\hfill\mbox{}
\\[6.75ex]
\indent
{\bf Property 9.16.}
{\it 
If $ \zeta=a $ is a solution of the equation $Y_n ^ {} (\zeta,1) =0$ and besides}\;\!:
\\[2ex]
\mbox{}\hfill
{\it 
a{\rm)}  
$W_{n-1}^{}(a,1)\,\partial_x^{} Y_n^{}(a,1)>0;$
\qquad
b{\rm)}  
$W_{n-1}^{}(a,1)\,\partial_x^{} Y_n^{}(a,1)<0;$
\hfill\mbox{}
\\[2.5ex]
\mbox{}\hfill
c{\rm)}  
$
\partial_x^{} Y_n^{}(a,1)=0,\quad
W_{n-1}^{}(a,1)\,\partial_{xx}^{} Y_n^{}(a,1)\ne0,$
\hfill\mbox{}
\\[2ex]
then the equatorial contact 
\vspace{0.25ex}
$A\!\!$-trajectory of the projectively singular system} (D) 
{\it in an enough small neighbourhood of the infinitely removed point $A,$ 
\vspace{0.25ex}
which is lying on <<extremities>> of the straight line $x=ay, $ is located, respectively}\;\!:
\vspace{0.35ex}

{\it  
a{\rm)} in the coordinate halfplane $y> 0;$
\vspace{0.35ex}

b{\rm)} in the coordinate halfplane $y <0;$
\vspace{0.35ex}

c{\rm)} as in the coordinate halfplane $y> 0,$ and in the coordinate halfplane $y <0.$
}
\vspace{1ex}

Taking into account Property 9.14, from Property 9.16 
\vspace{0.15ex}
we obtain disposition indications of trajectories (which are contact with the straight line $ \eta=0)$ 
\vspace{0.15ex}
for the second projectively  reduced system (7.4) in coordinate halfplanes. 
\\[4.25ex]
\centerline{
{\bf  10. Trajectories of the first projectively reduced system
}
}
\\[1.5ex]
\indent
The behaviour of trajectories 
\vspace{0.15ex}
of the first projectively  reduced system (6.3) essentially depends 
on behaviour of trajectories of system (D) and, on the contrary, 
\vspace{0.15ex}
by trajectories of system (6.3) we can establish the course of trajectories of system (D).  
\vspace{0.25ex}
These connections are schematically represented on Fig. 5.1 
\vspace{1ex}
by means of the projective circles $\P\K (x, y), \ \P\K (\xi, \theta).$

{\bf Property 10.1.}
\vspace{0.25ex}
{\it 
The images of the equilibrium states of system {\rm (D)}, which are lying}\;\!:
\vspace{0.35ex}

1)
{\it
in the final part of the projective phase plane $ \R\P (x, y), $ but not lying on the axis $Oy; $
}
\vspace{0.35ex}

2)\!
{\it
on\! <<extremities>>\! of the straight lines $\!y\!=\!ax,\!$ 
\vspace{0.15ex}
where the parameter $\!a\!$ is any real number
\\[0.25ex]
are the equilibrium states of the first projectively  reduced system {\rm (6.3)}, 
\vspace{0.35ex}
which are lying in the final part of the projective phase plane $\R\P (\xi,\theta).$ 
\vspace{0.25ex}
Thus the type of the equilibrium states remains to within a direction of movement along trajectories adjoining them.
}
\vspace{1ex}

{\bf Property 10.2.}
\vspace{0.5ex}
{\it 
If the equilibrium state $M (x, y) $ of system {\rm (D)} lies in}\;\!:
{\it a})
{\it 
the first};
{\it b})
{\it 
the second};
{\it c})
{\it 
the third};
\vspace {0.75ex}
{\it d})
{\it 
the fourth open coordinate quarter of the phase plane $Oxy,$ 
then its image $M _ {\phantom1} ^ {(1)} \Bigl (\dfrac{y}{x}\,,\, \dfrac{1}{x}\Bigr)$ 
lies, respectively, in}\;\!:
{\it a})
{\it 
the first};
{\it b})
{\it 
the third};
\vspace {0.5ex}
{\it c})
{\it 
the fourth};
{\it d})
{\it 
the second open coordinate quarter of the phase plane $O _ {\phantom1} ^ {(1)} \xi\theta. $
\vspace {0.75ex}
If the equilibrium state $M (x, 0) $ lies in the halfplane $x> 0\  (x <0),$ 
then its image $M _ {\phantom1} ^ {(1)} \Bigl (0, \dfrac {1} {x} \Bigr) $ lies in the halfplane $ \theta> 0\ (\theta <0).$
}

\newpage

The inverse statements to Properties 10.1 and 10.2, when equilibrium states of system (D) are defined by corresponding 
equilibrium states of the first projectively  reduced system (6.3), also take place.

The infinitely removed straight line 
\vspace{0.25ex}
of  the projective phase plane for the first projectively reduced system (6.3) 
and also the straight line 
\vspace{0.35ex}
$\zeta=0$ (the coordinate axis $O_{\phantom1}^{(2)} \eta)$ of the phase plane 
of the second projectively  reduced system (6.6) are the image of the straight line $x=0$ (the coordinate axis $Oy)$ 
\vspace{0.25ex}
of the phase plane for the differential system (D) (see, for example, Fig. 5.1). 
\vspace{0.15ex}
This implies that we have the following properties about trajectories for the differential systems (D), (6.3), and (6.6).
\vspace{0.75ex}

{\bf Property 10.3.}
{\it 
The following statements are equivalent}\;\!:
\vspace{0.35ex}

1. 
{\it
The first projectively  reduced system {\rm (6.3)} is projectively nonsingular}\;\!;
\vspace{0.35ex}

2. 
{\it 
The straight line $x=0$ consists from trajectories of systems} (D);
\vspace{0.35ex}

3. 
{\it
The straight line $ \zeta=0$ consists from trajectories of the  second projectively reduced system} (6.6).
\vspace{0.75ex}

{\bf Property 10.4.}
{\it 
The following statements are equivalent}\;\!:
\vspace{0.35ex}

1. 
{\it
The first projectively  reduced system {\rm (6.3)} is projectively singular}\;\!;
\vspace{0.35ex}

2. 
{\it 
The straight line $ x=0$ 
\vspace{0.15ex}
does not consists from trajectories of the  second projectively reduced system} (D);
\vspace{0.35ex}

3. 
{\it
The straight line $ \zeta=0$ 
\vspace{0.15ex}
does not consists from trajectories of the  second projectively reduced system} (6.6).
\vspace{0.75ex}

The infinitely removed point, 
\vspace{0.25ex}
which is lying on <<extremities>> of the straight line $ \xi=a\theta$ of the projective plane $ \R\P (\xi, \theta),$ 
\vspace{0.25ex}
is the image of the point $A (0, a),$ which is lying on the straight line $x=0.$
It allows to discover infinitely removed equilibrium states of the first projectively  reduced system (6.3) 
\vspace{0.15ex}
on the equilibrium states of system (D), which are lying on the straight line $x=0.$ 
In the case when the system (6.3) is projectively singular, we may find equatorial contact points of system (6.3) on the basis 
\vspace{0.15ex}
of contact points of the straight line $x=0$ for system (D).
\vspace{0.75ex}

{\bf Property 10.5.}
\vspace{0.25ex}
{\it
The point $A (0, a)$ is the equilibrium state of system {\rm (D)} if and only if
the point, which is lying on <<extremities>> of the straight line $ \xi=a\theta,$ 
\vspace{0.25ex}
is infinitely removed equilibrium state of the  first projectively  reduced system {\rm (6.3)}. 
\vspace{0.25ex}
Thus the equilibrium states have an identical type to within a direction of movement along trajectories adjoining them.}
\vspace{0.75ex}

{\bf Property 10.6.}
\vspace{0.25ex}
{\it 
The point $A (0, a) $ is a contact point of the straight line $x=0$ of system {\rm (D)} if and only if 
\vspace{0.25ex}
the point, which is lying on <<extremities>> of the straight line $\xi=a\theta,$ 
is an equatorial contact point of the first projectively  reduced projectively singular system} (6.3).
\vspace{0.5ex}

The disposition of equatorial contact trajectories 
\vspace{0.25ex}
of the projectively singular system (6.3) rather 
the infinitely removed straight line of the projective phase plane $\R\P (\xi, \theta)$ can be established 
as on the basis of Property 9.15, applied to the projectively singular system (6.3), 
and on the basis of a disposition of contact trajectories of the straight line $x=0$ of system (D) 
concerning the straight line $x=0. $ In the second case it is not required to know the analytical representation 
of the  first projectively reduced projectively singular system (6.3), it is enough to use Property 10.6 
\vspace{0.25ex}
and map of the projective circles $ \P\K (x, y) $ and $ \P\K (\xi, \theta), $ reduced on Fig. 5.1, and also the following statement.
\vspace{0.75ex}

{\bf Property 10.7.}
\vspace{0.35ex}
{\it 
The point $A (0, a) $ is the contact point of the straight line $x=0$ of system {\rm (D)} 
if and only if 
\vspace{0.35ex}
$y=a $ is a solution of the equation $X (0, y) =0, $ and $Y (0, a) \ne 0. $
If $X (0, a) =0$ and besides}\;\!: 
\vspace{0.5ex}

{\it a}) $Y (0, a)\,\partial_y ^ {} X (0, a)> 0; $
\vspace{0.75ex}

{\it b}) $Y (0, a) \, \partial_y ^ {} X (0, a) <0; $
\vspace{0.75ex}

{\it c}) $\partial_y ^ {} X (0, a) =0,\ \, Y (0, a)\,\partial_{yy}^{} X (0, a) \ne 0,$
\\[0.5ex]
{\it then the contact 
\vspace{0.35ex}
$A\!\!$-trajectory of the straight line $x=0$ of system {\rm (D)} 
in an enough small neighbourhood of  the point $A (0, a)$ lies, respectively}\;\!:
\vspace{0.35ex}

{\it  
a{\rm)} in the halfplane $x\geq 0;$
\vspace{0.5ex}

b{\rm)} in the halfplane $x\leq 0;$
\vspace{0.5ex}

c{\rm)} as in the halfplane $x\geq 0, $ and in the halfplane $x\leq 0.$
}
\vspace{0.75ex}

Thus, using the map of the projective circles 
\vspace{0.35ex}
$\P\K (x, y)$ and $\P\K (\xi, \theta),$ reduced on Fig. 5.1, 
on the basis of phase portrait on the projective circle $\P\K (x, y)$ 
\vspace{0.35ex}
of behaviour of trajectories of system (D) on the projective phase plane $\R\P (x, y),$ 
\vspace{0.35ex}
we can construct phase portrait on the projective circle $\P\K (\xi, \theta)$ 
\vspace{0.35ex}
of behaviour of trajectories of the first projectively  reduced system (6.3) 
on the projective phase plane $ \R\P (\xi, \theta).$ 
\vspace{0.35ex}
And, on the contrary, on the basis of phase portrait on the projective circle $\P\K (\xi, \theta)\!$ of system (6.3) 
\vspace{0.35ex}
it is possible to construct phase portrait on the projective circle $\P\K (x, y)\!$ of system (D).
\\[4.25ex]
\centerline{
{\bf  11. Trajectories of the second projectively reduced system
}
}
\\[1.5ex]
\indent
Connections between trajectories of system (D)  
\vspace{0.15ex}
and trajectories of the second projectively  reduced system (6.6) 
\vspace{0.25ex}
are represented schematically on Fig. 5.1 by means of the projective circles  $\P\K(x,y)$ and $\P\K(\eta,\zeta).$
\vspace{0.75ex}

{\bf Property 11.1.}
\vspace{0.35ex}
{\it 
The images of the equilibrium states of system {\rm (D)}, which are lying}\;\!:

1)
{\it
in the final part 
\vspace{0.35ex}
of the projective phase plane $ \R\P (x, y), $ but not lying on the axis $Ox;$
}

2)\!
{\it on\! <<extremities>>\! of the straight lines $\!x\!=\!ay,\!$ where the parameter $\!a\!$ is any real number
\\[0.25ex]
are the equilibrium states of the second projectively  reduced system {\rm (6.6)}, 
\vspace{0.35ex}
which are lying in the final part of the projective phase plane $ \R\P (\xi, \theta).$ 
\vspace{0.25ex}
Thus the type of the equilibrium states remains to within a direction of movement along trajectories adjoining them.
}
\vspace{1ex}

{\bf Property 11.2.}
\vspace {0.5ex}
{\it 
If the equilibrium state $M(x, y) $ of system {\rm (D)} lies in}\;\!:
{\it a})
{\it 
the first};
{\it b})
{\it 
the second};
{\it c})
{\it 
the third};
\vspace {0.75ex}
{\it d})
{\it 
the fourth open coordinate quarter of the phase plane $Oxy,$ 
then its image $M _ {\phantom2} ^ {(1)} \Bigl (\dfrac{1}{y}\,, \, \dfrac{x}{y}\Bigr) $ lies, respectively, in}\;\!:
{\it a})
{\it 
the first};
{\it b})
{\it the fourth};
\vspace {0.35ex}
{\it c})
{\it 
the second};
{\it d})
{\it the third open coordinate quarter of the phase plane $O _ {\phantom2} ^ {(1)} \eta\zeta. $
\vspace {1ex}
If the equilibrium state $M (0, y) $ lies in the halfplane $y> 0\ (y <0), $ 
then its image 
\vspace {1.25ex}
$M_{\phantom2}^{(1)} \Bigl(\dfrac{1}{y}\,,\, 0\Bigr) $ lies in the halfplane $\eta>0\ (\eta <0).$
}

The inverse statements to Properties 11.1 and 11.2, 
\vspace{0.15ex}
when equilibrium states of system (D) are defined by corresponding 
\vspace{0.15ex}
equilibrium states of the second  projectively  reduced system (6.6), also take place.
\vspace{0.15ex}

The infinitely removed straight line 
\vspace{0.25ex}
of  the projective phase plane for the second projectively reduced system (6.6) 
and also the straight line 
\vspace{0.35ex}
$\xi=0$ (the coordinate axis $O_{\phantom1}^{(1)} \theta)$ of the phase plane 
\vspace{0.35ex}
of the first projectively  reduced system (6.3) are the image of the straight line $y=0$ (the coordinate axis $Ox)$ 
of the phase plane for system (D) (see, for example, Fig.~5.1). 
\vspace{0.15ex}

Thus, we have
\vspace{0.5ex}

{\bf Property 11.3.}
{\it 
The following statements are equivalent}\;\!:
\vspace{0.15ex}

1. 
{\it
The second projectively  reduced system {\rm (6.6)} is projectively nonsingular}\;\!;
\vspace{0.15ex}

2. 
{\it 
The straight line $y=0$ consists from trajectories of systems} (D);
\vspace{0.25ex}

3.\! 
{\it
The\! straight\! line $\!\xi\!=\!0\!$ 
\vspace{0.75ex}
consists\! from\! trajectories\! of\! the\! first\! projectively\! reduced\! system}\! (6.3).

{\bf Property 11.4.}
{\it 
The following statements are equivalent}\;\!:
\vspace{0.15ex}

1. 
{\it
The second projectively  reduced system {\rm (6.6)} is projectively singular}\;\!;
\vspace{0.15ex}

2. 
{\it 
The straight line $ y=0$ does not consists from trajectories of system} (D);
\vspace{0.25ex}

3. 
{\it
The straight line $\xi=0$ does not consists from trajectories of the  first projectively reduced system} (6.3).
\vspace{0.5ex}

The infinitely removed point, 
\vspace{0.25ex}
which is lying on <<extremities>> of the straight line $ \zeta=a\eta$ of the projective plane $ \R\P (\eta, \zeta),$ 
\vspace{0.25ex}
is the image of the point $A (a, 0),$ which is lying on the straight line $y=0.$
It allows to discover infinitely removed equilibrium states of the second projectively  reduced system (6.6) 
\vspace{0.15ex}
on the equilibrium states of system (D), which are lying on the straight line $y=0.$ 
In the case when the system (6.6) is projectively singular, we may find equatorial contact points of system (6.6) on the basis 
\vspace{0.15ex}
of contact points of the straight line $y=0$ for system (D).
\vspace{0.75ex}

{\bf Property 11.5.}
\vspace{0.15ex}
{\it
The point $A (a, 0) $ is the equilibrium state of the system {\rm (D)} if and only if 
the point, which is lying on <<extremities>> of the straight line $ \zeta=a\eta,$ 
\vspace{0.25ex}
is infinitely removed equilibrium state of the second projectively  reduced system {\rm (6.6)}. 
\vspace{0.15ex}
Thus the equilibrium states have an identical type to within a direction of movement along trajectories adjoining them.}
\vspace{0.75ex}

{\bf Property 11.6.}
\vspace{0.25ex}
{\it 
The point $A(a, 0) $ is a contact point of the straight line $y=0$ of system {\rm (D)} 
if and only if the point, which is lying on <<extremities>> of the straight line $ \zeta=a\eta,$ 
\vspace{0.25ex}
is an equatorial contact point of the second projectively reduced projectively singular system} (6.6).
\vspace{0.75ex}

The disposition of equatorial contact trajectories 
\vspace{0.25ex}
of the projectively singular system (6.6) rather 
the infinitely removed straight line of the projective phase plane $\R\P (\eta, \zeta)$ can be established 
as on the basis of Property 9.16, applied to the projectively singular system (6.6), 
and on the basis of a disposition of contact trajectories of the straight line $y=0$ of system (D) 
concerning the straight line $y=0. $ In the second case it is not required to know the analytical representation 
of the second  projectively reduced projectively singular system (6.6), it is enough to use Property 11.6 
\vspace{0.25ex}
and map of the projective circles $ \P\K (x, y) $ and $ \P\K (\eta, \zeta), $ reduced on Fig. 5.1, and also the following statement.
\vspace{0.75ex}

{\bf Property 11.7.}
\vspace {0.35ex}
{\it 
The point $A(a, 0) $ is the contact point of the straight line $y=0$ of system {\rm (D)} 
if and only if 
\vspace {0.35ex}
$x=a $ is a solution of the equation $Y (x, 0) =0,$ and $X (a, 0) \ne 0.$
If $Y (a, 0) =0$ and besides}\;\!: 
\vspace {0.5ex}

{\it a}) $X (a, 0)\,\partial_x ^ {} Y (a, 0)> 0;$
\vspace {0.75ex}

{\it b}) $X (a, 0)\,\partial_x ^ {} Y (a, 0) <0;$
\vspace {0.75ex}

{\it c}) $ \partial_x ^ {} Y (a, 0) =0,\ \ X (a, 0)\,\partial _ {xx} ^ {} Y (a, 0) \ne0,$
\\[0.75ex]
{\it then the contact $A\!\!$-trajectory 
\vspace {0.35ex}
of  the straight line  $y=0$ of system {\rm (D)} in an enough small neighbourhood 
of the point $A(a, 0)$ lies, respectively}\;\!:
\vspace {0.35ex}

{\it a{\rm)} in the halfplane $y\geq 0; $
\vspace {0.5ex}

b{\rm)} in the halfplane $y\leq 0; $
\vspace {0.5ex}

c{\rm)} as in the halfplane $y\geq 0, $ and in the halfplane $y\leq 0.$
}
\vspace{0.75ex}

Thus, using the map of the projective circles 
\vspace{0.35ex}
$\P\K (x, y)$ and $\P\K (\eta, \zeta),$ reduced on Fig. 5.1, 
on the basis of phase portrait on the projective circle $\P\K (x, y)$ 
\vspace{0.35ex}
of behaviour of trajectories of system (D) on the projective phase plane $\R\P (x, y),$ 
\vspace{0.35ex}
we can construct phase portrait on the projective circle $\P\K (\eta, \zeta)$ 
\vspace{0.35ex}
of behaviour of trajectories of the second projectively  reduced system (6.6) 
on the projective phase plane $ \R\P (\eta, \zeta).$ 
\vspace{0.35ex}
And, on the contrary, on the basis of phase portrait on the projective circle $\P\K (\eta, \zeta)\!$ of system (6.6) 
\vspace{0.35ex}
it is possible to construct phase portrait on the projective circle $\P\K (x, y)\!$ of system (D).
\\[4.25ex]
\centerline{
{\bf  12.  Linear and open limit cycles
}
}
\\[1.5ex]
\indent
A closed trajectory of system (D), 
\vspace{0.25ex}
which is lying in the final part of the projective phase plane $\R\P (x, y),$ 
is called [13, p. 22] {\it a cycle} of system (D).
\vspace{0.35ex}

If for a cycle of system (D) 
\vspace{0.15ex}
there exists an neighbourhood such that this neighbourhood  has not cycles 
which distinct from are the given cycle (isolated cycle), 
\vspace{0.25ex}
then this cycle is called [1, pp. 71 -- 91] {\it a limit cycle} of system (D).

To a cycle of system (D) there corresponds 
\vspace{0.35ex}
the closed trajectory (cycle on the projective circle) 
on the open projective circle $\R\K (x, y) \backslash\partial \R\K (x, y),$ 
\vspace{0.35ex}
and on the projective sphere $\P {\mathbb S} (x, y)$ there corresponds 
\vspace{0.25ex}
pair of the diametrically opposite closed trajectories which do not have common points 
with equator (antipodal cycles on the projective sphere). 
\vspace{0.25ex}

There are trajectories of system (D), 
\vspace{0.25ex}
which are closed on the projective phase plane $\R\P (x, y), $ but not closed on its final part $(x, y).$
\vspace{1ex}

{\bf Definition 12.1.}
\vspace{0.35ex}
{\it 
The straight line, which is a trajectory of system {\rm (D)} on the projective phase plane $ \R\P (x, y),$ 
is called \textit{\textbf{a linear cycle}} of system {\rm (D)} or {\boldmath $\ell\!$}\textit{\textbf{-cycle}} of system} (D). 
\vspace{0.75ex}

Linear cycle, being a trajectory on $\!\R\P (x, y),\!$ 
\vspace{0.25ex}
consists of the regular points, which are ly\-ing as in the final part of the projective phase plane, 
\vspace{0.35ex}
and on the infinitely removed straight line. 

If the infinitely removed straight line of the projective phase plane $\R\P (x, y)$ is 
\vspace{0.25ex}
a trajectory of system (D), then this line is called {\it $\ell_{\infty}^ {}\!$-cycle} of system (D) .
\vspace{0.35ex}

To linear cycles of system (D) there correspond circles of the big radius
\vspace{0.25ex}
on the projective sphere $\P{\mathbb S}(x, y);$
\vspace{0.35ex}
and to $\ell _ {\infty}^{}\!$-cycle there corresponds the equator. 
On the projective circle $\P\K (x, y)$ to $\ell_{\infty}^{}\!$-cycle there corresponds the boundary circle $\partial \P\K(x,y).$
\vspace{0.75ex}

{\bf Definition 12.2.}
{\it 
We say that 
\vspace{0.15ex}
\textit {\textbf {limit linear cycle}} of system {\rm (D)} is a linear cycle such that 
\vspace{0.15ex}
on the projective sphere $\P {\mathbb S} (x, y)$ to this linear cycle there corresponds 
a circle of the big radius having an neighbourhood in which 
\vspace{0.15ex}
there are no closed trajectories of system {\rm (D)}, distinct from this circle.
} 
\vspace{0.75ex}

Each trajectory, which is distinct from the equilibrium state $O (0,0),$ 
\vspace{0.25ex}
of the projectively nonsingular system (8.3) with the projective atlas of trajectories on Fig. 8.5 
\vspace{0.25ex}
is a cycle. 
Mo\-re\-o\-ver, the infinitely removed straight line of the projective phase plane $\R\P (x, y)$ is 
\vspace{0.25ex}
$\ell_ {\infty}^{}\!$-cycle of system (8.3).
\vspace{0.25ex}
Thus $\ell_{\infty}^{}\!$-cycle and the cycles, which are lying in the final part of the projective phase plane $ \R\P (x, y),$ 
are not limit cycles of this system.
\vspace{0.75ex}

The projective atlas of trajectories 
\vspace{0.15ex}
for the projectively nonsingular system (8.4) is constructed on Fig. 8.6. 
\vspace{0.15ex}
This projectively nonsingular system has the limit $\ell _ {\infty}^{}\!$-cycle. Since the system (8.4) is linear, we see that  
\vspace{0.15ex}
this system has not limit cycles in the final part of the projective phase plane $ \R\P (x, y).$ 
\vspace{0.5ex}

Jacobi's system [14, pp. 111 -- 117; 10, pp. 14 -- 25]
\\[2.5ex]
\mbox{}\hfill
$
\dfrac{dx}{dt}=x-y+x(x+y),
\qquad
\dfrac{dy}{dt}=(y+1)(x+y)
$
\hfill (12.1)
\\[2.25ex]
has the general autonomous integral 
\\[2ex]
\mbox{}\hfill
$
F\colon (x, y) \to\  \dfrac {x^2+y^2} {(y+1) ^2} \, \exp \Bigl ({}-2 \arctan \dfrac {y} {x} \Bigr)
$ 
\ for all 
$
(x, y) \in G
\hfill
 $
\\[2.25ex]
on any domain from the set $G=\{(x,y)\colon x\ne 0 \ \&\  y\ne {}-1\,\}.$
\vspace{0.75ex}

The straight line-trajectory $y = {}-1$ 
\vspace{0.35ex}
is a limit linear cycle of the projectively singular system (12.1).
\vspace{0.5ex}
In addition, the system (12.1) has not [10, p. 55] limit cycles, which are lying in the final part of the projective phase plane $\R\P (x, y).$ 
\vspace{0.75ex}

The first projectively reduced system for system (12.1) is Jacobi's system 
\\[2.25ex]
\mbox{}\hfill                   %(12.2)
$
\dfrac{d\xi}{dt}=1+\xi^2,
\qquad
\dfrac{d\theta}{dt}={}-1-\xi-\theta +\xi \theta.
$
\hfill (12.2)
\\[2.5ex]
\indent
The system (12.2) has not [10, p. 55] limit cycles 
\vspace{0.5ex}
in the final part of the projective phase plane $\R\P (\xi, \theta).$ 
The straight line-trajectory $\xi = {}-\theta $ is a limit linear cycle of system (12.2).
\vspace{0.5ex}

The transcendental function 
\\[2ex]
\mbox{}\hfill
$
F^{(1)}_{\phantom 1} \colon (\xi,\theta)\to\
\dfrac{1+\xi^2}{(\xi+\theta)^2} \, \exp ({}-2\arctan \xi)
$ 
\ for all 
$
(\xi,\theta)\in G^{(1)}_{\phantom 1},
\quad
G^{(1)}_{\phantom 1}=\{(\xi,\theta) \colon \xi\ne {}-\theta\},
\hfill
$ 
\\[1.5ex]
is the general autonomous integral of system (12.2) on any domain from the set $G^{(1)}_{\phantom 1}.$
\vspace{0.75ex}

The second projectively  reduced system for system (12.1) is Jacobi's system 
\\[2ex]
\mbox{}\hfill                  %(12.3)
$
\dfrac{d\eta}{dt}={}-(1+\eta)(1+\zeta),
\qquad
\dfrac{d\zeta}{dt}={}-(1+\zeta^2).
$
\hfill (12.3)
\\[2.5ex]
\indent
The system (12.3) has not [10, p. 55] limit cycles 
\vspace{0.5ex}
in the final part of the projective phase plane  $\R\P(\eta,\zeta).$ 
The straight line-trajectory $\eta={}-1$ 
\vspace{0.75ex}
is a limit linear cycle of system (12.3).

The function
\\[2ex]
\mbox{}\hfill
$
F^{(2)}_{\phantom 1} \colon (\eta,\zeta)\to \
\dfrac{1+\zeta^2}{(1+\eta)^2} \, \exp \Bigl({}-2\arctan \dfrac1\zeta\Bigr)
$ 
\ for all 
$
(\eta,\zeta)\in G^{(2)}_{\phantom 1},
\hfill
$
\\[2ex]
where  $G^{(2)}_{\phantom 1}=\{(\eta,\zeta) \colon \eta\ne {}-1 \ \&\  \zeta\ne 0\},$
\vspace{0.75ex}
is the general autonomous integral of system (12.3) on any domain from the set $G^{(2)}_{\phantom 1}.$
\vspace{0.5ex}

The projective atlas of trajectories for Jacobi's system (12.1) is constructed on Fig. 12.1.
\\[3.75ex]
\mbox{}\hfill
{\unitlength=1mm
\begin{picture}(42,42)
\put(0,0){\includegraphics[width=42mm,height=42mm]{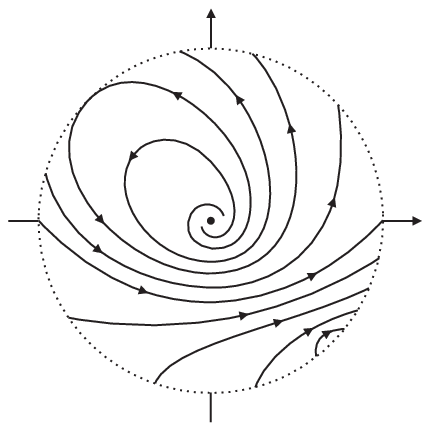}}
\put(18,41){\makebox(0,0)[cc]{ $y$}}
\put(40.2,18.2){\makebox(0,0)[cc]{ $x$}}
\end{picture}}
\quad
{\unitlength=1mm
\begin{picture}(42,42)
\put(0,0){\includegraphics[width=42mm,height=42mm]{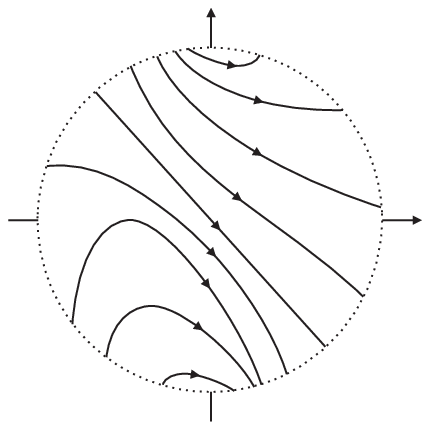}}
\put(18,41){\makebox(0,0)[cc]{ $\theta$}}
\put(40.2,17.8){\makebox(0,0)[cc]{ $\xi$}}
\put(21,-7){\makebox(0,0)[cc]{Fig. 12.1}}
\end{picture}}
\quad
{\unitlength=1mm
\begin{picture}(42,42)
\put(0,0){\includegraphics[width=42mm,height=42mm]{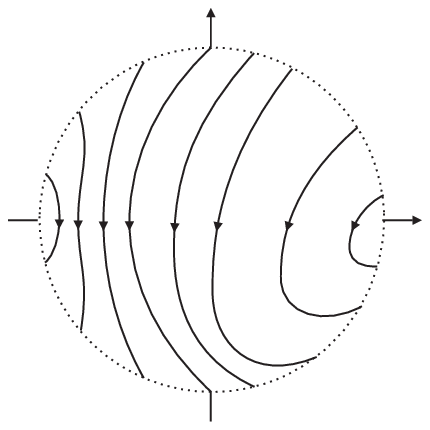}}
\put(18,41){\makebox(0,0)[cc]{ $\zeta$}}
\put(40.2,18){\makebox(0,0)[cc]{ $\eta$}}
\end{picture}}
\hfill\mbox{}
\\[7.5ex]
\indent
Closed trajectories of system (D) 
\vspace{0.15ex}
on the projective sphere $ \P {\mathbb S}(x, y),$ 
distinct from circles of the big radius, are subdivided into two classes: 
\vspace{0.15ex}

1) not having common points with the equator of the projective sphere $ \P {\mathbb S} (x, y); $
\vspace{0.5ex}

2) having 
\vspace{0.75ex}
at least one common point with the equator of the projective sphere $\P{\mathbb S}(x,y).$

{\bf Definition 12.3 (12.4).}
\vspace{0.15ex}
{\it 
An \textit{\textbf{open cycle}} {\rm (}\textit{\textbf{open limit cycle}}{\rm)} of system {\rm (D)} 
\vspace{0.25ex}
is a trajectory of system {\rm (D)} on the projective phase plane $\R\P (x, y)$ such that
to this trajectory on the pro\-j\-ec\-ti\-ve sphere $\P {\mathbb S} (x, y)$ there corresponds 
\vspace{0.25ex}
pair of the antipodal closed trajectories {\rm (}pair of the antipodal isolated closed trajectories{\rm)}, 
\vspace{0.15ex}
which are having common points with the equator of the projective sphere $\P{\mathbb S} (x, y).$
}
\vspace{0.75ex}

Then, 
\vspace{0.15ex}
to a cycle of system (D), which is lying in the final part of the projective phase plane $ \R\P (x, y),$ 
\vspace{0.15ex}
there corresponds pair of the antipodal closed trajectories of system (D) on the projective sphere $ \P {\mathbb S} (x, y),$ 
\vspace{0.15ex}
not having common points with the equator of the projective sphere $\P{\mathbb S}(x, y).$
\vspace{0.35ex}

The projectively singular Jacobi system  
\\[2ex]
\mbox{}\hfill                             %(12.4)
$
\dfrac{d\xi}{dt}={}-(1+\xi^2),
\qquad
\dfrac{d\theta}{dt}={}-\xi\;\!\theta
$
\hfill (12.4)
\\[2ex]
is the first projectively  reduced system for system (8.3).
\vspace{0.25ex}

The straight line-trajectory $ \theta=0$ 
\vspace{0.25ex}
is a linear limit cycle of system (12.4), 
and any other trajectory of system (12.4) is an open cycle (Fig. 8.5).
\vspace{0.25ex}

The system (12.4) hasn't 
\vspace{0.25ex}
linear limit cycles, open limit cycles, and limit cycles, which are lying in the final part of the projective phase plane $\R\P(\xi,\theta).$ 
\vspace{0.35ex}

The function
\\[2ex]
\mbox{}\hfill
$
F^{(1)}_{\phantom 1}\colon (\xi,\theta)\to\ 
\dfrac{1+\xi^2}{\theta^2}
$
\ for all 
$
(\xi,\theta)\in G^{(1)}_{\phantom 1}, 
\quad
G^{(1)}_{\phantom 1}=\{(\xi,\theta)\colon \theta\ne 0\},
\hfill
$
\\[1.5ex]
is the general autonomous integral of system (12.4) on any domain from the set $G^{(1)}_{\phantom 1}.$
\vspace{0.5ex}

The projectively singular Jacobi system  
\\[2ex]
\mbox{}\hfill                %(12.5)
$
\dfrac{d\eta}{dt}=\eta\;\!\zeta,
\qquad
\dfrac{d\zeta}{dt}=1+\zeta^2
$
\hfill (12.5)
\\[2ex]
is the second projectively reduced system for system (8.3).
\vspace{0.25ex}

The straight line-trajectory $\eta=0$ 
\vspace{0.25ex}
is a linear limit cycle of system (12.5), and any other trajectory of system (12.5) 
is an open cycle (Fig. 8.5).
\vspace{0.35ex}

The system (12.5) hasn't 
\vspace{0.25ex}
linear limit cycles, open limit cycles, and limit cycles, which are lying in the final part of the projective phase plane $\R\P(\eta,\zeta).$ 
\vspace{0.35ex}

The function
\\[2ex]
\mbox{}\hfill
$
F^{(2)}_{\phantom 1}\colon (\eta,\zeta)\to\ 
\dfrac{1+\zeta^2}{\eta^2}
$ 
\ for all 
$
(\eta,\zeta)\in G^{(2)}_{\phantom 1},
\quad
G^{(2)}_{\phantom 1}=\{(\eta,\zeta)\colon \eta\ne 0\},
\hfill
$
\\[1.5ex]
is the general autonomous integral of system (12.5) on any domain from the set $G^{(2)}_{\phantom 1}.$
\vspace{0.75ex}

The projectively singular Jacobi system 
\\[2ex]
\mbox{}\hfill             %(12.6)
$
\dfrac{d\xi}{dt}=1+\xi^2,
\qquad
\dfrac{d\theta}{dt}=\theta(\xi -1)
$
\hfill (12.6)
\\[2ex]
is the first projectively reduced system for system (8.4).
\vspace{0.25ex}

The straight line-trajectory $\theta=0$ is a linear limit cycle of system (12.6)  (Fig. 8.6).
\vspace{0.25ex}

The system (12.6) hasn't (including limit cycles)
\vspace{0.25ex}
open cycles and cycles, which are lying in the final part of the projective phase plane $\R\P (\xi, \theta).$  
\vspace{0.25ex}
 
The function
\\[2ex]
\mbox{}\hfill
$
F^{(1)}_{\phantom 1}\colon (\xi,\theta)\to\ 
\dfrac{1+\xi^2}{\theta^2}\, \exp ({}-2\arctan \xi)
$ 
\ for all 
$
(\xi,\theta)\in G^{(1)}_{\phantom 1},
\quad 
G^{(1)}_{\phantom 1}=\{(\xi,\theta)\colon \theta\ne 0\},
\hfill
$
\\[1.5ex]
is the general autonomous integral of system (12.6) on any domain from the set  $G^{(1)}_{\phantom 1}.$
\vspace{0.5ex}

The projectively singular Jacobi system
\\[2ex]
\mbox{}\hfill             %(12.7)
$
\dfrac{d\eta}{dt}={}-\eta(1+ \zeta),
\qquad
\dfrac{d\zeta}{dt}={}-(1+\zeta^2)
$
\hfill (12.7)
\\[2.5ex]
is the second projectively  reduced system for system (8.4).
\vspace{0.25ex}

The straight line-trajectory $\eta=0$ is a linear limit cycle of system (12.7) (Fig. 8.6).
\vspace{0.25ex}

The system (12.7) hasn't (including limit cycles)
\vspace{0.25ex}
open cycles and cycles, which are lying in the final part of the projective phase plane $\R\P(\eta, \zeta).$  
\vspace{0.25ex}

The function
\\[2ex]
\mbox{}\hfill
$
F^{(2)}_{\phantom 1}\colon (\eta,\zeta)\to\ 
\dfrac{1+\zeta^2}{\eta^2}\,\exp \Bigl({}-2\arctan \dfrac1\zeta\Bigr)
$  
\ for all 
$
(\eta,\zeta)\in G^{(2)}_{\phantom 1},
\quad 
G^{(2)}_{\phantom 1}=\{(\eta,\zeta)\colon \eta \zeta\ne 0\},
\hfill
$
\\[1.5ex]
is the general autonomous integral of system (12.7) on any domain from the set $G^{(2)}_{\phantom 1}.$
\vspace{0.75ex}

{\bf Property 12.1.}
\vspace{0.25ex}
{\it 
In the class of projectively singular systems {\rm (D)} there exist 
an open and a linear limit cycles, which are distinct from $\ell_ {\infty}^{}\!$-cycle.
\vspace{0.25ex}
In the class of projectively nonsingular systems {\rm (D)} there exists 
the limit $\ell _ {\infty}^{}\!$-cycle.
}

{\sl Proof} 
is based on that the projective type of system (D) depends on that consists or does not consist 
the infinitely removed straight line of the projective phase plane $ \R\P (x, y) $ from 
trajectories of system (D) (Properties 9.1 and 9.2). \k
\vspace{0.5ex}

{\bf Property 12.2.}
{\it 
A projectively nonsingular system {\rm (D)} with even degree $n$ 
has at least one infinitely removed equilibrium state.
}
\vspace{0.25ex}

{\sl Indeed}, 
\vspace{0.15ex}
if $n$ is an even number, then either the first projectively  reduced system (7.1) has an equilibrium state on the straight line $\theta=0$ 
\vspace{0.25ex}
or the second projectively  reduced system (7.2) has an equilibrium state in the origin of coordinates for 
the phase plane $O ^ {(2)} _ {\phantom 1}\eta\zeta.$\k
\vspace{0.75ex}

From Property 12.1, taking into account Property 12.2,  we obtain  
\vspace{0.75ex}

{\bf Property 12.3.}
\vspace{0.15ex}
{\it 
If the system {\rm (D)} has the limit $\ell_{\infty}^{}\!$-cycle, then 
this system is a projectively nonsingular with odd degree $n.$
}
\vspace{0.75ex}

The projectively singular Jacobi  system (12.1) hasn't open limit cycles.
\vspace{0.75ex}

The Darboux system [10]                          
\\[2.25ex]
\mbox{}\hfill                      %(12.8)
$
\dfrac{dx}{dt}={}-y-x(x^2+y^2-1),
\qquad
\dfrac{dy}{dt}=x-y(x^2+y^2-1),
$
\hfill (12.8)
\\[2.5ex]
in the final part of the projective phase plane $ \R\P (x, y)$ has one equilibrium state which is an unstable focus.
\vspace{0.25ex}

The circle $x^2+y^2=1$ 
\vspace{0.35ex}
is unique limit cycle of system (12.8) in the final part of the projective phase plane $ \R\P (x, y) $ 
[4, pp. 257 -- 260; 15, pp. 39 -- 46]. 
\vspace{0.35ex}

The function
\\[2ex]
\mbox{}\hfill
$
F\colon (x,y) \to\ 
\dfrac{x^2+y^2}{1-x^2-y^2}\ \exp\Bigl({}-2\arctan\dfrac{y}{x}\Bigr)
$
\ for all 
$
(x,y)\in G,
\hfill
$
\\[2ex]
where the set $G=\{(x,y)\colon x\ne 0\ \, \&\,\ x^2+y^2\ne 1\},$
\vspace{0.35ex}
is the general autonomous integral of the Darboux system (12.8) on any domain from the set $G.$
\vspace{0.5ex}

Straight lines 
\vspace{0.15ex}
are not trajectories of the projectively singular Darboux system (12.8) on the projective phase plane $\R\P (x, y).$ 
\vspace{0.15ex}
The system (12.8) hasn't equatorial contact points. 
The system (12.8) hasn't linear cycles and open cycles. 
\vspace{0.15ex}
Moreover, this system hasn't limit linear cycles and limit open cycles.
\vspace{0.5ex}

The first projectively  reduced system for system (12.8) is the projectively singular system 
\\[2ex]
\mbox{}\hfill                      %(12.9)
$
\dfrac{d\xi}{d\tau}=\theta (1+\xi^2),
\qquad
\dfrac{d\theta}{d\tau}=1+\xi^2-\theta^2+\xi\theta^2, 
\qquad 
\theta d\tau=dt.
$
\hfill (12.9)
\\[2.25ex]
\indent
This system hasn't 
\vspace{0.35ex}
linear cycles and cycles, which are lying in the final part of the projective phase plane $ \R\P (\xi, \theta) $. 
\vspace{0.5ex}

The hyperbola $\theta^2-\xi^2=1$ is a limit open cycle of system (12.9). 
\vspace{0.35ex}
The system (12.9) in the projective phase plane $\R\P(\xi,\theta) $ hasn't other limit cycles. 
\vspace{0.5ex}

The transcendental function 
\\[2ex]
\mbox{}\hfill
$
F^{(1)}_{\phantom 1}\colon (\xi,\theta)\to\ 
\dfrac{1+\xi^2}{\theta^2-\xi^2-1}\ \exp({}-2\arctan\xi)
$
\ for all 
$
(\xi,\theta)\in G^{(1)}_{\phantom 1},
\hfill
$
\\[1.5ex]
where 
\vspace{0.5ex}
$G^{(1)}_{\phantom 1}=\{(\xi,\theta)\colon \theta^2-\xi^2\ne 1\},$ 
is the general autonomous integral of system (12.9) on any domain from the set  $G^{(1)}_{\phantom 1}.$
\vspace{0.75ex}

The second projectively  reduced system for the Darboux system (12.8) is the projectively singular system
\\[2ex]
\mbox{}\hfill                      %(12.10)
$
\dfrac{d\eta}{d\nu}=1-\eta^2+\zeta^2-\eta^2\zeta, 
\qquad
\dfrac{d\zeta}{d\nu}={}-\eta(1+\zeta^2), 
\qquad 
\eta d\nu=dt.
$
\hfill (12.10)
\\[2.5ex]
\indent
The hyperbola $ \eta^2-\zeta^2=1$ is a limit open cycle of system (12.10).
\vspace{0.35ex}
The system (12.10) in the projective phase plane $ \R\P (\eta, \zeta)$ hasn't
other cycles (including limit cycles). 
\vspace{0.35ex}

The function
\\[2ex]
\mbox{}\hfill
$
F^{(2)}_{\phantom 1}\colon (\eta,\zeta) \to\ 
\dfrac{1+\zeta^2}{\eta^2-\zeta^2-1}\ \exp\Bigl({}-2\arctan\dfrac{1}{\zeta}\,\Bigr)
$
\ for all 
$
(\eta,\zeta)\in G^{(2)}_{\phantom 1},
\hfill
$ 
\\[1.75ex]
where the set
\vspace{0.5ex}
$G^{(2)}_{\phantom 1}=\{(\eta,\zeta)\colon \zeta\ne 0\,\ \&\,\ \eta^2-\zeta^2\ne 1\},$
is the general autonomous integral of system (12.10) on any domain from the set $G^{(2)}_{\phantom 1}.$
\vspace{0.35ex}

The projective atlas of trajectories for system (12.8) is constructed on Fig. 12.2. 
\\[4.25ex]
\mbox{}\hfill
{\unitlength=1mm
\begin{picture}(42,42)
\put(0,0){\includegraphics[width=42mm,height=42mm]{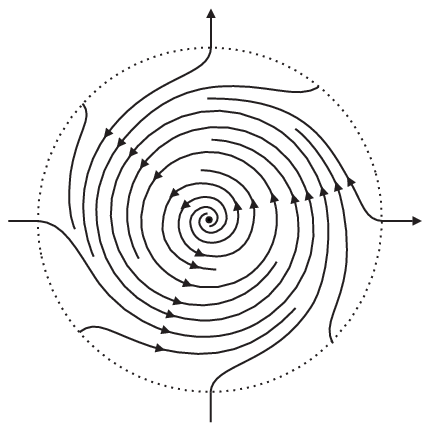}}
\put(18,41){\makebox(0,0)[cc]{ $y$}}
\put(40.2,18.2){\makebox(0,0)[cc]{ $x$}}
\end{picture}}
\quad
{\unitlength=1mm
\begin{picture}(42,42)
\put(0,0){\includegraphics[width=42mm,height=42mm]{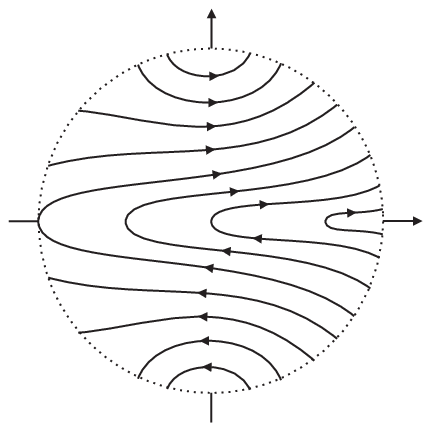}}
\put(18,41){\makebox(0,0)[cc]{ $\theta$}}
\put(40.2,17.8){\makebox(0,0)[cc]{ $\xi$}}
\put(21,-7){\makebox(0,0)[cc]{Fig. 12.2}}
\end{picture}}
\quad
{\unitlength=1mm
\begin{picture}(42,42)
\put(0,0){\includegraphics[width=42mm,height=42mm]{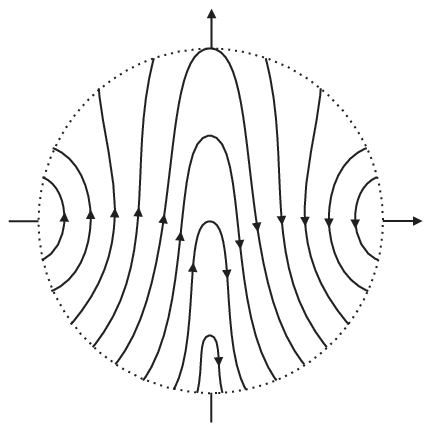}}
\put(18,41){\makebox(0,0)[cc]{ $\zeta$}}
\put(40.2,18){\makebox(0,0)[cc]{ $\eta$}}
\end{picture}}
\hfill\mbox{}
\\[8.5ex]
\indent
The projectively singular system
\\[2ex]
\mbox{}\hfill                         %(12.11)
$
\dfrac{dx}{dt}={}-2x-y+3x^2+y^2 -x(x^2+y^2),
\qquad
\dfrac{dy}{dt}={}-1+x+2xy -y(x^2+y^2)
$
\hfill (12.11)
\\[2.5ex]
is received from the Darboux system  (12.8) 
\vspace{0.25ex}
by parallel transposition of the origin of coordinates for the phase plane $Oxy $ to the point $A (1,0).$ 
\vspace{0.25ex}
Therefore, the system (12.11) in the final part of the projective phase plane $\R\P (x, y)$ 
\vspace{0.35ex}
has one equilibrium state $A (1,0),$ 
being an unstable focus, and has one limit cycle $(x-1)^2+y^2=1.$
\vspace{0.25ex}

The function 
\\[2ex]
\mbox{}\hfill
$
F\colon (x,y) \to\ 
\dfrac{(x-1)^2+y^2}{1-(x-1)^2-y^2}\ \exp\Bigl({}-2\arctan\dfrac{y}{x-1}\Bigr)
$
\ for all 
$
(x,y)\in G,
\hfill
$
\\[2ex]
where the set
\vspace{0.5ex}
$G=\{(x,y)\colon x \ne 1\ \, \&\,\ (x-1)^2+y^2\ne 1\},$
is the general autonomous integral of system (12.11) on any domain from the set $G.$
\vspace{0.75ex}

The first projectively  reduced system for system (12.12) is the projectively singular system
\\[2ex]
\mbox{}\hfill                      %(12.12)
$
\dfrac{d\xi}{d\tau}={}-\xi +\theta+2\xi\theta - \theta^2-\xi^3+\xi^2\theta,
\quad
\dfrac{d\theta}{d\tau}=1-3\theta+\xi^2+2\theta^2-\xi^2\theta+\xi\theta^2, 
\quad 
\theta d\tau=dt.
$
\hfill (12.12)
\\[2.5ex]
\indent
The parabola $\xi^2-2\theta+1=0$ is a limit open cycle of system (12.12). 
\vspace{0.35ex}
The system (12.12) hasn't other cycles (including limit cycles) in the projective phase plane $\R\P(\xi, \theta).$ 
\vspace{0.35ex}

The function 
\\[2ex]
\mbox{}\hfill
$
F^{(1)}_{\phantom 1}\colon (\xi,\theta) \to\ 
\dfrac{\xi^2+(\theta-1)^2}{2\theta-1-\xi^2}\ \exp\Bigl(2\arctan\dfrac{\xi}{\theta-1}\Bigr)
$
\ for all 
$
(\xi,\theta)\in G^{(1)}_{\phantom 1},
\hfill
$ 
\\[2ex]
where the set 
\vspace{0.5ex}
$G^{(1)}_{\phantom 1}=\{(\xi,\theta)\colon \theta\ne 1\ \,\&\,\ 1-2\theta+\xi^2\ne 0\},$ 
is the general autonomous integral of system (12.12) on any domain from the set $G^{(1)}_{\phantom 1}.$ 

\newpage

The second projectively reduced system for system (12.11) is the projectively singular system
\\[2.2ex]
\mbox{}\hfill                      %(12.13)
$
\dfrac{d\eta}{d\nu}=1-2\eta\zeta+\zeta^2+\eta^3-\eta^2\zeta\,, 
\quad
\dfrac{d\zeta}{d\nu}=1-\eta-2\eta\zeta +\zeta^2+ \eta^2\zeta-\eta\zeta^2, \ \ \ \eta d\nu=dt.
$
\hfill (12.13)
\\[2.5ex]
\indent
The hyperbola $\zeta^2-2\eta\zeta+1=0$ is a limit open cycle of system (12.13).
\vspace{0.35ex}
The system (12.13) hasn't other cycles (including limit cycles) in the projective phase plane $\R\P(\eta,\zeta).$ 
\vspace{0.35ex}
  
The function
\\[2ex]
\mbox{}\hfill
$
F^{(2)}_{\phantom 1}\colon (\eta,\zeta) \to\ 
\dfrac{1+(\eta-\zeta)^2}{2\eta\zeta-\zeta^2-1}\ 
\exp\Bigl(2\arctan\dfrac{1}{\eta-\zeta}\,\Bigr)
$
\ for all 
$
(\eta,\zeta)\in G^{(2)}_{\phantom 1},
\hfill
$ 
\\[1.5ex]
where the set 
\vspace{0.35ex}
$G^{(2)}_{\phantom 1}=\{(\eta,\zeta)\colon \eta\ne \zeta\,\ \&\,\ 2\eta\zeta-\zeta^2- 1\ne0\},$
is the general autonomous integral of system (12.13) on any domain from the set $G^{(2)}_{\phantom 1}.$
\vspace{0.5ex}

The projective atlas of trajectories for system (12.11) is constructed on Fig. 12.3.  
\\[3.75ex]
\mbox{}\hfill
{\unitlength=1mm
\begin{picture}(42,42)
\put(0,0){\includegraphics[width=42mm,height=42mm]{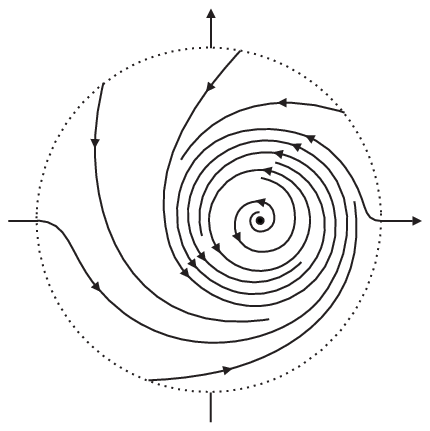}}
\put(18,41){\makebox(0,0)[cc]{ $y$}}
\put(40.2,18.2){\makebox(0,0)[cc]{ $x$}}
\end{picture}}
\quad
{\unitlength=1mm
\begin{picture}(42,42)
\put(0,0){\includegraphics[width=42mm,height=42mm]{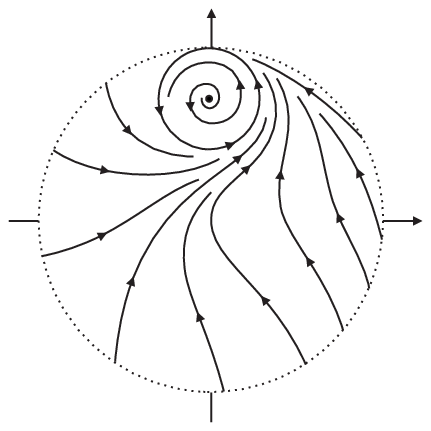}}
\put(18,41){\makebox(0,0)[cc]{ $\theta$}}
\put(40.2,17.8){\makebox(0,0)[cc]{ $\xi$}}
\put(21,-7){\makebox(0,0)[cc]{Fig. 12.3}}
\end{picture}}
\quad
{\unitlength=1mm
\begin{picture}(42,42)
\put(0,0){\includegraphics[width=42mm,height=42mm]{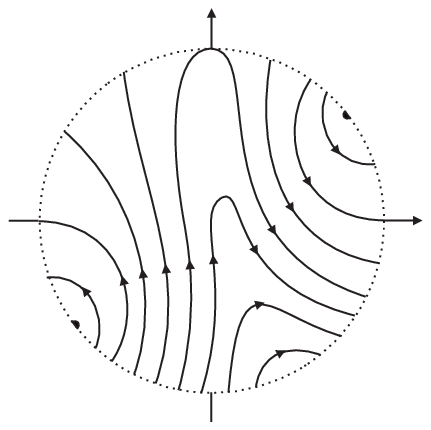}}
\put(18,41){\makebox(0,0)[cc]{ $\zeta$}}
\put(40.2,18){\makebox(0,0)[cc]{ $\eta$}}
\end{picture}}
\hfill\mbox{}
\\[9ex]
\indent
Since the projectively singular system
\\[2.2ex]
\mbox{}\hfill                                       %(12.14)
$
\displaystyle
\dfrac{dx}{dt}=1-(15+9a-8b)x+y-(5+a)x^2-cxy-x^2(x-y),
\hfill
$
\\[0.1ex]
\mbox{}\hfill (12.14)
\\[0.1ex]
\mbox{}\hfill
$
\displaystyle
\dfrac{dy}{dt}=10y-(5+a)xy-cy^2-xy(x-y),
\hfill
$
\\[2.75ex]
where
\vspace{0.15ex}
$
a={}-10^{-13},\ 
b={}-10^{-52},\ 
c={}-10^{-200}, 
$
has one equilibrium state, which is saddle, 
we see that 
\vspace{0.15ex}
this system hasn't limit cycles, which are lying in the final part of the projective phase plane $ \R\P (x, y).$ 
\vspace{0.35ex}

The second projectively  reduced system for system (12.14) is the system
\\[2.2ex]
\mbox{}\hfill                           %(12.15)
$
\displaystyle
\dfrac{d\eta}{d\nu}=c\eta-\zeta -10\eta^2+(5+a)\eta \zeta +\zeta^2,
\hfill
$
\\[0.1ex]
\mbox{}\hfill (12.15)
\\[0.1ex]
\mbox{}\hfill
$
\displaystyle
\dfrac{d\zeta}{d\nu}=\eta(1+\eta-(25+9a-8b)\zeta),
\qquad 
\eta\, d\nu =dt,
\hfill
$
\\[2.75ex]
where the numbers
$
a={}-10^{-13},\ 
b={}-10^{-52},\ 
c={}-10^{-200}.
$
\vspace{0.15ex}

The system (12.15) in the final part of the projective phase plane $\R\P (\eta, \zeta)$ 
\vspace{0.25ex}
has two equilibrium states $O ^ {(2)}(0,0) $ and $A^{(2)} (0,1),$ which are focuses.
\vspace{0.15ex}

The system (12.15) in the final part of the projective phase plane $ \R\P (\eta, \zeta) $ 
\vspace{0.25ex}
has not less than four limit cycles [16]. 
\vspace{0.25ex}
The system (12.15) is projectively nonsingular and has infinitely removed equilibrium state (saddle [17, p. 182]). 
\vspace{0.25ex}
The system (12.15) does not have linear and open cycles (including limit cycles).

Using 
\vspace{0.25ex}
the phase portrait of behaviour of trajectories of system (12.15) on the projective circle $ \P\K (\eta, \zeta)$ 
\vspace{0.25ex}
(constructed on Fig. 7 in [17, p. 182]  on the case of four limit cycles) and 
the rule of map of the projective circles $\P\K (x, y), \, \P\K (\xi, \theta), \, \P\K (\eta, \zeta)$ (see Fig. 5.1), 
\vspace{0.25ex}
we have the projective atlas of trajectories for system (12.15) is constructed on Fig. 12.4.
\vspace{0.25ex}

Consequently, 
\vspace{0.15ex}
the system (12.14) hasn't linear limit cycles and limit cycles, which are lying in the final part of the projective phase plane 
\vspace{0.15ex}
$\R\P (x, y),$ and in addition this system has not less than four open limit cycles.
\vspace{0.5ex}

The first projectively  reduced system for system (12.14) is the system
\\[2ex]
\mbox{}\hfill                      %(12.16)
$
\dfrac{d\xi}{d\tau}=\xi\;\!\theta(25+9a -8b-\xi -\theta),
\hfill
$
\\[0.25ex]
\mbox{}\hfill(12.16)
\\[0.25ex]
\mbox{}\hfill
$
\dfrac{d\theta}{d\tau}=1-\xi +(5+a)\theta+c\xi \theta+(15+9a-8b) \theta^2-\xi\theta^2-\theta^3, 
\quad 
\theta d\tau=dt,
\hfill
$
\\[2.5ex]
where the numbers
$
a={}-10^{-13},\ 
b={}-10^{-52},\ 
c={}-10^{-200}. 
$
\vspace{0.25ex}

The system (12.16) hasn't linear cycles. 
\vspace{0.15ex}
At the same time the system (12.16) has open limit cycles and limit cycles, 
\vspace{0.15ex}
which are lying in the final part of the projective phase plane $\R\P (\xi, \theta)$ 
(note that the sum of these cycles not less than four). 
\\[3.75ex]
\mbox{}\hfill
{\unitlength=1mm
\begin{picture}(42,42)
\put(0,0){\includegraphics[width=42mm,height=42mm]{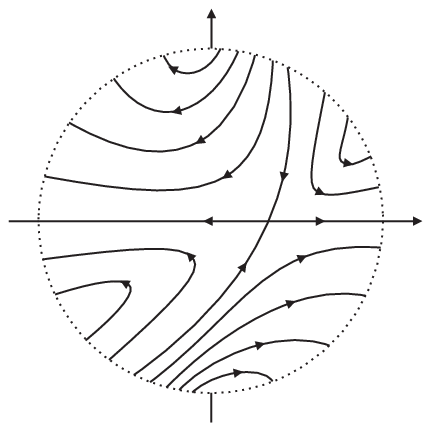}}
\put(18,41){\makebox(0,0)[cc]{ $y$}}
\put(40.2,18.5){\makebox(0,0)[cc]{ $x$}}
%\put(21,-3){\makebox(0,0)[cc]{ $Oxy$}}
%\put(22.5,-6){\makebox(0,0)[cc]{Рис. 1}}
\end{picture}}
\quad
{\unitlength=1mm
\begin{picture}(42,42)
\put(0,0){\includegraphics[width=42mm,height=42mm]{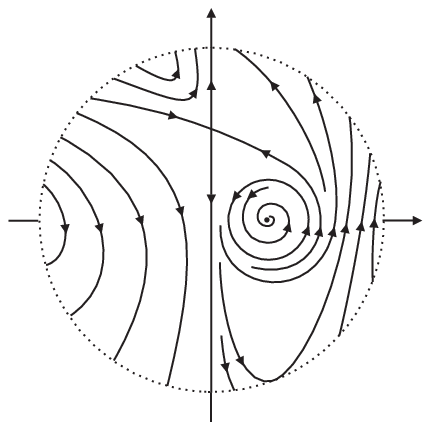}}
\put(18,41){\makebox(0,0)[cc]{ $\theta$}}
\put(40.2,18){\makebox(0,0)[cc]{ $\xi$}}
%\put(21,-3){\makebox(0,0)[cc]{ $O^{{}^{(1)}}uz$}}
\put(21,-7){\makebox(0,0)[cc]{Fig. 12.4}}
\end{picture}}
\quad
{\unitlength=1mm
\begin{picture}(42,42)
\put(0,0){\includegraphics[width=42mm,height=42mm]{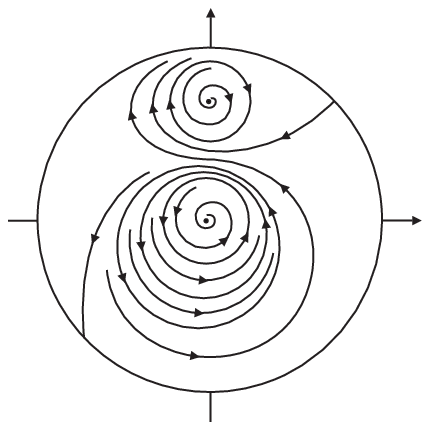}}
\put(18,41){\makebox(0,0)[cc]{ $\zeta$}}
\put(40.2,18){\makebox(0,0)[cc]{ $\eta$}}
%\put(21,-3){\makebox(0,0)[cc]{ $O^{{}^{(2)}}zv$}}
%\put(22.5,-6){\makebox(0,0)[cc]{Рис. 1}}
\end{picture}}
\hfill\mbox{}
\\[8.5ex]
\centerline{
{\bf  13. 
Symmetry of the phase field of directions
}
}
\\[1.5ex]
\indent
Through the regular point $M (x, y) $ 
\vspace{0.25ex}
of the phase plane $Oxy $ we will spend a segment of the straight line with the directing vector 
$\vec{a}(x,y)=(X(x,y), Y(x,y)).$
\vspace{0.25ex}

It can be assume that the segment is not directed and has unit length with the middle in the point $M.$ 
A trajectory of system  (D) concerns this segment in its middle-point (i.e. $M$ is contact point of the segment).
The set of such segments constructed in each regular point of system (D) is called {\it the phase directional field} for system (D). 
The phase directional field in an equilibrium state of system (D) is not defined.
\vspace{0.35ex}

From geometrical reasons, we get the following criteria of symmetry of the phase directional field 
for the differential system (D).
\vspace{0.5ex}

{\bf Property 13.1}  
\vspace{0.5ex}
[18, pp. 4 -- 5].
%\marginpar{[2150]}
{\it 
The phase directional field of system {\rm (D)} is  symmetric with respect to the straight line 
$Ax + By + C = 0\ (|A | + |B | \ne 0) $ if and only if
\\[2.25ex]
\mbox{}\hfill
$
R\Bigl(x - \dfrac{2A(Ax + By + C)}{A^{2} + B^{2}}\,,\ y
- \dfrac{2B(Ax + By + C)}{A^{2} + B^{2}}\, \Bigr) \, =\,
\dfrac{2AB - \bigl(A^2 - B^2\,\bigr)R(x,y)}{A^{2} - B^{2} +2AB\,R(x,y)}
\hfill
$
\\[2.75ex]
\mbox{}\hfill
for all 
$
(x,y) \in \Omega,
\qquad
\Omega=\{(x,y)\colon X(x,y)\ne 0\},
\hfill
$
\\[2.5ex]
where the rational function 
$R\colon (x,y)\to\, \dfrac{Y(x,y)}{X(x,y)}$ for all $(x,y) \in \Omega.$
}
\vspace{0.75ex}

Suppose the system (D) has the phase directional field, symmetric relative to some straight line $l$.  
Then, for each trajectory of this system there exists a trajectory, symmetric to it relative to the straight line $l.$ 
Besides, the possibility of presence of the trajectories, which are symmetric relative to the straight line $l,$ is not excluded.
\vspace{0.15ex}

For example, the symmetry of the phase directional field for system (D) relative to bisectrixes of coordinate angles 
\vspace{0.75ex}
of the phase plane $Oxy $ is established by following properties. 

{\bf Property 13.2.} 
{\it 
The phase directional field of system {\rm (D)} 
\vspace{0.15ex}
is  symmetric with respect to the straight line $y = x $ if and only if
\\[1.5ex]
\mbox{}\hfill
$
X(x,y)\;\!X(y,x)\;\!  -\;\!  Y(x,y)\;\!Y(y,x)  =  0
$ 
\ for all 
$
(x,y)\in \R^2.
\hfill
$
}
\\[1.75ex]
\indent
{\bf Property 13.3.} 
{\it 
The phase directional field of system {\rm (D)} 
\vspace{0.15ex}
is  symmetric with respect to the straight line $y ={} -x $ if and only if
\\[1.5ex]
\mbox{}\hfill
$
X({}-x,{}-y)\;\!X(y,x)\;\!  -\;\! Y({}-x,{}-y)\;\!Y(y,x)  =  0
$
\ for all 
$
(x,y)\in \R^2.
\hfill
$
}
\\[1.5ex]
\indent
The next statements 
\vspace{0.15ex}
are analytical criteria of symmetry of  the phase directional field for system (D) 
\vspace{0.75ex}
relative to coordinate axes and the origin of coordinates of the phase plane $Oxy.$

{\bf Property 13.4.} [12, pp. 50 -- 51].
%\marginpar{[296]}
{\it 
The phase directional field of system {\rm (D)} 
\vspace{0.25ex}
is  symmetric with respect to the axis $Ox$ if and only if
\\[1.5ex]
\mbox{}\hfill
$
X(x,y)\;\!Y(x,{}-y)\;\!  +\;\!  X(x,{}-y)\;\!Y(x,y)  =  0
$
\ for all 
$
(x,y)\in \R^2.
\hfill
$
}
\\[1.75ex]
\indent
{\bf  Property 13.5.} 
{\it 
The phase directional field of system {\rm (D)} 
\vspace{0.15ex}
is  symmetric with respect to the axis $Oy$ if and only if
\\[1.5ex]
\mbox{}\hfill
$
X(x,y)\;\!Y({}-x,y)\;\! +\;\! X({}-x,y)\;\!Y(x,y)  =  0
$
\ for all 
$
(x,y)\in\R^2.
\hfill
$
}
\\[1.75ex]
\indent
{\bf  Property 13.6.} 
{\it 
The phase directional field of system {\rm (D)} 
\vspace{0.15ex}
is  symmetric with respect to the origin of coordinate $O(0,0)$ 
if and only if
\\[1.5ex]
\mbox{}\hfill
$
X(x,y)\;\!Y({}-x,{}-y)\;\!  -\;\! X({}-x,{}-y)\;\!Y(x,y)  = 0
$
\ for all 
$
(x,y)\in \R^2.
\hfill
$
}
\\[1.75ex]
\indent
For the first and the second projectively  reduced systems,
the symmetry of the phase directional fields  relative to coordinate axes and the origin of coordinates 
is established by means of Properties 13.4 -- 13.6, and also the symmetry can be established on the basis 
of the connections containing in Properties 13.7 -- 13.9.
\vspace{0.75ex}

{\bf Property 13.7.} 
{\it 
The following statements are equivalent\;\!{\rm:}
\vspace{0.35ex}

{\rm 1.} 
The phase directional field of system {\rm (D)} 
\vspace{0.25ex}
is  symmetric with respect to the origin of coordinate $O(0,0)$ 
of the phase plane $Oxy;$ 
\vspace{0.75ex}

{\rm 2.} 
\vspace{0.75ex}
The identity 
$
X(x,y)\;\!Y({}-x,{}-y)  - X({}-x,{}-y)\;\!Y(x,y)  = 0
$
for all 
$(x,y)\!\in\! \R^2$
is true{\rm;}

{\rm 3.} 
The phase directional field of the first projectively reduced system  {\rm (6.3)} 
\vspace{0.35ex}
is  symmetric with respect 
to the coordinate axes $O_{\phantom1}^{(1)}\xi;$
\vspace{0.75ex}
 
{\rm 4.} 
\vspace{0.75ex}
The identity
$
\Xi(\xi,\theta)\;\!\Theta(\xi,{}-\theta) + \Xi(\xi,{}-\theta)\;\!\Theta(\xi,\theta)  = 0
$
for all 
$(\xi,\theta)\in \R^2$
is true{\rm;}

{\rm 5.} 
The phase directional field of the second projectively reduced system {\rm (6.6)} 
\vspace{0.35ex}
is  symmetric with respect to the coordinate axes $O_{\phantom1}^{(2)}\zeta;$
\vspace{0.75ex}
 
{\rm 6.} 
The identity
$
H(\eta,\zeta)\;\! Z({}-\eta,\zeta) + H({}-\eta,\zeta)\;\! Z(\eta,\zeta)  = 0
$
for all 
$(\eta,\zeta)\in \R^2$
is true.
}
\vspace{1.25ex}

{\bf Property 13.8.} 
{\it 
The following statements are equivalent\;\!{\rm:}
\vspace{0.35ex}

{\rm 1.} 
\vspace{0.25ex}
The phase directional field of system {\rm (D)} is  symmetric with respect to the coordinate axes $Ox;$ 
\vspace{0.75ex}

{\rm 2.} 
\vspace{0.75ex}
The identity
$
X(x,y)\;\!Y(x,{}-y)  + X(x,{}-y)\;\!Y(x,y)  = 0
$
for all 
$(x,y)\in \R^2$
is true{\rm;}

{\rm 3.} 
The phase directional field of the first projectively reduced system {\rm (6.3)} 
\vspace{0.35ex}
is  symmetric with respect to the coordinate axes  $O_{\phantom1}^{(1)}\theta;$
\vspace{0.75ex}
 
{\rm 4.} 
\vspace{0.75ex}
The identity
$
\Xi(\xi,\theta)\;\!\Theta({}-\xi,\theta) + \Xi({}-\xi,\theta)\;\!\Theta(\xi,\theta)  = 0
$
for all 
$(\xi,\theta)\in \R^2$
is true{\rm;}

{\rm 5.} 
The phase directional field of the second projectively reduced system  {\rm (6.6)} 
\vspace{0.25ex}
is  symmetric with respect to the origin of coordinate $O(0,0)$ of the phase plane $O_{\phantom1}^{(2)}\eta\zeta\,;$
\vspace{0.75ex}

{\rm 6.} 
The identity
\vspace{1.25ex}
$
H(\eta,\zeta)Z({}-\eta,{}-\zeta) - H({}-\eta,{}-\zeta)Z(\eta,\zeta)  = 0
$
for all 
$(\eta,\zeta)\in \R^2$
is true.
}

{\bf Property 13.9.} 
{\it 
The following statements are equivalent\;\!{\rm:}
\vspace{0.35ex}

{\rm 1.} 
The phase directional field of system {\rm (D)} is  symmetric with respect to the coordinate axes $Oy;$ 
\vspace{0.5ex}

{\rm 2.} 
\vspace{0.75ex}
The identity
$
X(x,y)\;\!Y({}-x,y)  + X({}-x,y)\;\!Y(x,y)  = 0
$
for all 
$(x,y)\in \R^2$
is true{\rm;}

{\rm 3.} 
The phase directional field of the first projectively reduced system {\rm (6.3)} 
\vspace{0.25ex}
is  symmetric with respect to the origin of coordinate $O(0,0)$ of the phase plane $O_{\phantom1}^{(1)}\xi\theta\, ;$
\vspace{0.75ex}
 
{\rm 4.} 
\vspace{0.75ex}
The identity
$
\Xi(\xi,\theta)\;\!\Theta({}-\xi,{}-\theta) - \Xi({}-\xi,{}-\theta)\;\!\Theta(\xi,\theta)  = 0
$
for all 
$
(\xi,\theta)\in \R^2
$
is true{\rm;}

{\rm 5.} 
The phase directional field of the second projectively reduced system {\rm (6.6)} 
\vspace{0.25ex}
is  symmetric with respect to the coordinate axes $O_{\phantom1}^{(2)}\eta\, ;$
\vspace{0.5ex}
 
{\rm 6.} 
The identity
$
H(\eta,\zeta)\;\! Z(\eta,{}-\zeta) + H(\eta,{}-\zeta)\;\! Z(\eta,\zeta)  = 0
$
for all 
$
(\eta,\zeta)\in \R^2
$
is true.
}
\\[4.25ex]
\centerline{
{\bf  14. Sets of projectively nonsingular and projectively singular systems}
}
\\[1.5ex]
\indent
Suppose $l$ is some straight line of the plane $Oxy.$ 
Then, there is a linear nondegenerate transformation of the plane $Oxy$ 
such that in an new coordinate system the straight line $l$ will be the axis of ordinate (abscissas).
Using this fact, on the basis of Properties 10.3 and 11.3, we get the next statement.
\vspace{0.5ex}

{\bf Property 14.1.}
{\it 
Suppose some straight line of the phase plane $Oxy$ consists from trajectories of system {\rm (D)}.  
Then, using a linear nondegenerate transformation and the first {\rm (}the second{\rm)} Poincar\'{e}'s transformation of  system {\rm (D)},
we obtain the system {\rm(D)} is reduced to a projectively nonsingular system.
}
\vspace{0.5ex}

All possible straight lines of the phase plane of system (D) cannot consist simultaneously of trajectories.
Therefore always it is possible to specify linear nondegenerate transformation of the phase plane $Oxy$ 
such that the system (D) is reduced to a system that in an new coordinate system the axis of ordinates (abscissas) will not consist of its trajectories. 
By Properties 10.3, 10.4, and 11.3, 11.4, we have the statement.
\vspace{0.5ex}

{\bf Property 14.2.}
{\it 
For the system {\rm (D)} always there exists a linear nondegenerate tran\-s\-for\-ma\-ti\-on 
that converts it into a system such that 
the first {\rm (}the second{\rm)} projectively  reduced system for this new system is projectively singular}. 
\vspace{0.5ex}

In other words, by means of a linear nondegenerate transformation and by means 
of the first (the second) Poincar\'{e}'s transformation, any system (D) is reduced to a projectively singular system.

Let $D$ be the set of all ordinary autonomous differential systems of the second order with polynomial right members of any degree;

$A\subset D$ is  the set of all systems, which are projectively singular or which by means 
of a linear nondegenerate transformation and by the first or the second Poincar\'{e}'s transformation are reduced to projectively singular systems;

$B\subset D $ is the set of all systems, which are projectively nonsingular or which by means of a linear nondegenerate transformation 
and by the first or the second Poincar\'{e}'s transformation are reduced to projectively nonsingular systems;

$C\subset D $ is the set of all systems, which are projectively singular or by means of a linear nondegenerate transformation 
and by Poincar\'{e}'s transformations are not reduced to projectively nonsingular systems.
\vspace{0.35ex}

Then, $B\cap C=\O,\ B\cup C=D.$ 
\vspace{0.35ex}
By Property 14.2, we get $A=D.$ So, $B\subset A,\ C\subset A.$ Thus the complemen of the set $\!B\!$ 
to the set $\!D\!=\!A\!$ is the set $\!C,\!$ disjunctive with the set $\!B.$

Let's prove existence of systems belonging to the set $C.$ 
\vspace{0.5ex}
We will use the next statement.

{\bf Property 14.3.}
{\it 
A system is belonging to the set $C$ 
if and only if 
the projective phase plane of this system hasn't
straight lines consisting of trajectories of this system.
} 
\vspace{0.5ex}

{\sl Indeed}, 
the infinitely removed straight line is not consists from trajectories of system (D) if and only if 
the system (D) is projectively singular (Property 9.2). 
Therefore each system of the set $C$ such that the infinitely removed straight line does not consist of its trajectories.
Otherwise, i.e. when the infinitely removed straight line consists of trajectories, 
system is projectively nonsingular (Property 9.1) and, hence, this system does not belong to the set $C.$

According to Property 10.4 (to Property 11.4) the first (the second) reduced system is projectively singular 
if and only if 
the axis of ordinates (abscissas) does not consist of its trajectories.
If on the phase plane $Oxy $ 
\vspace{0.15ex}
there exists a straight line consisting of trajectories of system (D), 
then the system (D) is reduced to a projectively nonsingular system (Property~14.1). 
Therefore each system of the set $C$ such that on the phase plane is not present straight lines, which are consist of its trajectories. 
Existence of straight lines consisting of trajectories of system (D), means that the system (D) is 
\vspace{0.15ex}
reduced to a projectively nonsingular system. So, the system (D) $\!\not\in C.$\k
\vspace{0.5ex}

Everyone projectively  singular differential system (12.8) -- (12.13) hasn't 
\vspace{0.15ex}
straight lines on the projective phase plane consisting of its trajectories. 
\vspace{0.15ex}
By Property 14.3, the differential systems (12.8) -- (12.13) belong to the set $C,$ and therefore, $C\ne\O.$
\vspace{0.75ex}

{\bf Theorem 14.1.}
{\it 
A system} (D) {\it at $n=0,1,2$ does not belong to the set} $C.$
\vspace{0.35ex}

{\sl Proof}.
\vspace{0.25ex}
The system (D) at $n=0$ has the form (7.8) and this system is projectively nonsingular (Example~7.1). 
So, the system (D) at $n=0$ is not belong to the set $C. $
\vspace{0.25ex}

At $n=1$ the system (D) is the linear autonomous system (7.11) 
\vspace{0.25ex}
which will be projectively singular system at $a_2 ^ {} = b_1 ^ {} = 0, \, b_2 ^ {} = a_1 ^ {} \ne 0$ (Example 7.2). 
\vspace {0.25ex}
However, in this case any straight line of the set 
$C_1 ^ {} (a_1 ^ {} y+b _ {_ 0} ^ {}) +C_2 ^ {} (a_1 ^ {} x+a _ {_ 0} ^ {}) =0$ 
\vspace{0.25ex}
consists of trajectories of system (7.11). 
\linebreak
By Property 14.1, the system (7.11) is not belong to the set $C.$
\vspace{0.25ex}

If $n=2, $ thet the system (D) has the form (7.16) 
\vspace{0.25ex}
and this system is projectively singular at 
$a_5^{}=b_3^{}=0, \ b_4^{}=a_3^{},\ b_5^{}=a_4^{}, \ |a_3^{}|+|a_4^{}|\ne 0$ (Example 7.3). 
\vspace{0.35ex}
The equation of trajectories of the projectively singular system (7.16) 
\vspace{0.25ex}
is the Jacobi equation which always has a straight line consisting of its trajectories [10, pp. 15 -- 16]. 
By Property 14.1, 
\vspace{0.25ex}
the system (7.16) is not belong to the set $C.$ \k
\\[3.75ex]
\centerline{
{\bf  15. Topological equivalence of differential systems }
}
\\[0.25ex]
\centerline{
{\bf on projective circle and on projective sphere 
}
}
\\[1.5ex]
\indent
We will distinguish behaviour of trajectories 
on projective circle or on projective sphere to within topological equivalence.
\vspace{0.75ex}

{\bf Definition\! 15.1}\! [11, p. 34].\!\!
\vspace{0.15ex}
{\it 
Two autonomous polynomial differential systems of the se\-cond order are 
\textit{\textbf{topologically equivalent on projective circle}} 
\vspace{0.15ex}
if exists a homeomorphism of
the projective circles, translating the trajectories of one system 
\vspace{0.5ex}
to the trajectories of other system}.

Since a linear nondegenerate transformation of  the plane is a homeomorphism, we see that next statement is takes place.
\vspace{0.5ex}

{\bf Property 15.1.}\!
{\it 
At a linear nondegenerate transformation of the phase plane $\!(x, y)\!$ re\-ma\-in{\rm:}
{\rm a)} topological equivalence of system {\rm (D);} 
{\rm b)} degree $n $ of system {\rm (D);} 
{\rm c)} type of system {\rm (D)} on the projective phase plane.
}
\vspace{0.75ex}

For example, the system (7.8) and the system
\\[2.25ex]
\mbox{}\hfill
$
\dfrac{dx}{dt}=\widetilde{a}_{_0},
\qquad 
\dfrac{dy}{dt}=\widetilde{b}_{_0},
\qquad
|\widetilde{a}_{_0}|+|\widetilde{b}_{_0}|\ne0,
\hfill
$
\\[2.25ex]
are topologically  equivalent on the projective circle.

Really, by the parallel transposition $v=x+\widetilde{a}_{_0}-a_{_0},\ w=y+\widetilde{b}_{_0}-b_{_0},$ we obtain
\\[2.25ex]
\mbox{}\hfill
$
\dfrac{dv}{dt}=\dfrac{dx}{dt}+\widetilde{a}_{_0}-a_{_0}=\widetilde{a}_{_0},
\qquad 
\dfrac{dw}{dt}=\dfrac{dy}{dt}+\widetilde{b}_{_0}-b_{_0}=\widetilde{b}_{_0}.
\hfill
$
\\[2.25ex]
\indent
Therefore the phase portrait of trajectories of system (7.8) 
\vspace{0.15ex}
to within topological equ\-i\-va\-len\-ce on the projective circle is constructed both on the circle 
\vspace{0.15ex}
$\P\K (x, y) $ from Fig.~8.1, and on the circle $\P\K (x, y) $ from Fig. 8.2.
\vspace{0.5ex}

{\bf Definition 15.2.}
\vspace{0.15ex}
{\it 
Two autonomous polynomial differential systems of the second order are 
\textit{\textbf{topologically equivalent on projective sphere}} 
\vspace{0.15ex}
if exists the homeomorphism of projective spheres translating the trajectories of one system to the trajectories of other system}.
\vspace{0.5ex}

Topological equivalence of behaviour of trajectories of systems on the projective sphere 
remains at superposition of  a linear nondegenerate transformation and the transformations of Poincar\'{e}. 
Therefore the following two properties are true.
\vspace{0.5ex}

{\bf Property 15.2.}
{\it 
The systems {\rm (D), (6.3)}, and {\rm (6.6)}  are  topologically equivalent on the projective sphere.
}
\vspace{0.5ex}

{\bf Property 15.3.}
{\it 
The sum of linear, open, and usual limit cycles is identical for to\-po\-lo\-gi\-cal\-ly equivalent on the projective sphere systems.
}
\vspace{0.5ex}

For example, according to Property 15.2,
the system (8.1) with the projective atlas on Fig. 8.3 and 
the system (8.2) with the projective atlas on Fig. 8.4 are topologically equivalent on the projective sphere. 
In addition, the system (8.2) is the first projectively reduced system of system (8.1), 
\vspace{0.25ex}
and the system (8.1) is the second projectively reduced system of system (8.2).

If ${\rm deg}\,({\rm D})=n,$ then the system (D) is denoted by $({\rm D}^n).$
\vspace{0.75ex}

{\bf Property 15.4.}
\vspace{0.25ex}
{\it 
A projectively nonsingular system $\!({\rm D}^n)\!$ is topologically equivalent on the projective sphere 
to a projectively singular system $ ({\rm D}^{n+1}),$ 
\vspace{0.15ex}
and the system $({\rm D}^{n+1})$ has a straight line consisting of its trajectories.
}
\vspace{0.35ex}

{\sl Proof}\, is the consequence of Properties 7.3, 10.4, 11.3, 15.1, and Theorem 4.1. \k
\vspace{0.5ex}

This\! property\! allows\! instead\! of\! qualitative\! research\! as\! a\! whole of projectively singular system 
in the presence of a straight line consisting of its trajectories, 
to fulfil qualitative research as a whole of system corresponding to projectively nonsingular system. 
For the class of systems $A$ such possibility of passage to research of projectively nonsingular systems is excluded.

For example, 
\vspace{0.25ex}
the Jacobi system (7.21) has the straight line consisting of its trajectories [10, pp. 15 -- 16].
Then, according to Property 15.4, we have the next statement.
\vspace{0.75ex}

{\bf Property 15.5.}
\vspace{0.15ex}
{\it 
The Jacobi system {\rm (7.21)} is topologically equivalent on the projective sphere to 
the linear stationary system {\rm (7.11)}.
}
\vspace{0.5ex}

Thereupon, we say that Jacobi's system is {\it projectively linear differential system}.
\vspace{0.5ex}

A system $({\rm D}^{2})$ will be projectively singular if and only if 
\vspace{0.25ex}
this system is Jacobi's system (Property 7.4).
Therefore on the basis of property 15.5 we can state the following.
\vspace{0.75ex}

{\bf Property 15.6.}
\vspace{0.15ex}
{\it 
The projectively singular system $({\rm D}^{2})$ is topologically equivalent on the projective sphere 
to the system $({\rm D}^{1}).$
}
\vspace{0.75ex}

{\bf Property 15.7.}
{\it 
The projectively singular system $({\rm D}^{2})$ is projectively linear differential system.
}
\vspace{0.75ex}

The linear-fractional transformations 
\\[2.25ex]
\mbox{}\hfill               %(15.1)
$
\xi=\dfrac{{}-\beta x+\alpha y}{\alpha x+\beta y+\gamma}\,,
\qquad 
\theta=\dfrac{\sqrt{\alpha^2+\beta^2}}{\alpha x+\beta y+\gamma}
$
\hfill (15.1)
\\[1.5ex]
and 
\\[1.5ex]
\mbox{}\hfill               %(15.2)
$
\eta=\dfrac{\sqrt{\alpha^2+\beta^2}}{\alpha x+\beta y+\gamma}\,,
\qquad 
\zeta= \dfrac{\beta x-\alpha y}{\alpha x+\beta y+\gamma}\,,
$
\hfill (15.2)
\\[2.5ex]
can be used 
\vspace{0.25ex}
at passage to topologically equivalent on the projective sphere systems, 
where $\alpha^2+\beta^2\ne 0.$ By these linear-fractional transformations, 
\vspace{0.25ex}
the straight line $ \alpha x +\beta y +\gamma=0$ is reduced to the infinitely removed straight lines 
\vspace{0.25ex}
of the projective phase planes $\R\P(\xi,\theta)$ and $\R\P(\eta,\zeta),$ respectively.
\vspace{0.35ex}

In particular, 
if the straight line $ \alpha x +\beta y +\gamma=0$ consists from trajectories of system (D), 
then the system (D) by the transformations (15.1) and (15.2) is reduced to projectively nonsingular systems.
\\[4.5ex]
\centerline{
{\bf  16. Examples of global qualitative research of trajectories
}
}
\\[0.5ex]
\centerline{
{\bf
for differential systems on projective phase plane}
}
\\[1.75ex]
\indent
{\bf 16.1.
Trajectories of Darboux's differential system} [14, pp. 109 -- 111; 10]
\\[2ex]
\mbox{}\hfill        % (16.1)
$
\dfrac{dx}{dt}={}-y+x^3\equiv X(x,y),
\qquad 
\dfrac{dy}{dt}=x(1+xy)\equiv Y(x,y).
$
\hfill(16.1)
\\[2ex]
\indent
{\it Integral basis of system} (16.1). 
The functionally independent first integrals
\\[2ex]
\mbox{}\hfill
$
F_1^{}\colon (t,x,y)\to\ 
{}-t+\arctan \dfrac{y}{x}
$
\, \ for all 
$
(t,x,y) \in \{(t,x,y)\colon x\ne 0\}
\hfill
$
\\[1ex]
and
\\[1ex] 
\mbox{}\hfill
$
F_2^{}\colon (t,x,y)\to\ 
t+ \dfrac{1+xy}{x^2+y^2}
$
\, \ for all 
$
(t,x,y) \in \{(t,x,y)\colon |x|+|y|\ne 0\}
\hfill
$
\\[2ex]
are integral basis [19] of system (16.1) on any domain from the set $\{(t, x, y) \colon x\ne 0\}.$ 
\vspace{0.75ex}

{\it The general autonomous integral of system} (16.1).  
The transcendental function
\\[2ex]
\mbox{}\hfill
$
F\colon (x,y)\to\
\dfrac{1+xy}{x^2+y^2}+\arctan \dfrac{y}{x}
$
\ \, 
for all 
$
(x,y)\in\{(x,y)\colon x\ne 0\}
\hfill
$
\\[2ex]
is the general autonomous integral [20, pp. 112 -- 114] of the differential system (16.1) 
\vspace{0.5ex}
on any domain from the set $\{(x, y)\colon x\ne 0\}.$
\vspace{0.5ex}

{\it Equilibrium states of system} 
\vspace{0.35ex}
(16.1) {\it in the final part of the projective phase plane} $\R\P(x,y).$ 
The cubic parabola $y=x^3$ hasn't common points with the hyperbola $xy+1=0, $ 
\vspace{0.35ex}
and this parabola intersects the straight line $x=0$ in the point $O(0,0).$
\vspace{0.35ex}
The system (16.1) has one equilibrium state $O (0,0) $ in the final part of the projective phase plane $\R\P (x, y) $.
\vspace{0.35ex}

Using the characteristic equation 
\vspace{0.25ex}
$ \lambda^2+1=0,$ we get the equilibrium state $O (0,0) $ is a centre or a focus
[3, pp. 139 -- 145].
\vspace{0.25ex}

The set of trajectories of system (16.1) in the polar coordinate system $O\rho\varphi $ is defined by the equation 
\\[1.25ex]
\mbox{}\hfill
$
\rho^{{}-2}=C-\varphi-\dfrac{1}{2}\, \sin 2\varphi\;\!.
\hfill
$
\\[1.5ex]
\indent
Therefore, $O (0,0) $ is an unstable focus.
\vspace{0.5ex}

{\it Movement of radius-vector of organising point along trajectories of system } (16.1).
The function 
\\[1.25ex]
\mbox{}\hfill
$
W\colon (x,y)\to\ 
x\;\!Y(x,y)-y\;\!X(x,y)=x^2+y^2>0
$
\ for all 
$
(x,y)\in \R^2\backslash\{(0,0)\}.
\hfill
$
\\[1.5ex]
\indent
Then, the angle between the radius-vector of  the organising point and the positive direction of the axis $Ox$ 
at movement along trajectories of system (16.1) is increasing.
\vspace{0.5ex}

{\it Symmetry of the phase directional field of system}  (16.1).
The phase directional field of system (16.1) is 
symmetric with respect to the origin of coordinates of the phase plane $Oxy $ (Property 13.6).
For each trajectory of system (16.1) there exists a symmetric trajectory with respect to the origin of coordinates of the phase plane $Oxy. $ 
\vspace{0.5ex}

{\it Zero and orthogonal isoclines of system}  (16.1).
\vspace{0.35ex}
From the equation $Y (x, y) =0$ at $X (x, y) \ne 0$ it follows that 
\vspace{0.35ex}
zero isoclines of system (16.1) are the hyperbola $xy+1=0$ and the straight line $x=0$ without the point $O (0,0).$ 
\vspace{0.35ex}

The tangent to a trajectory of system (16.1) 
\vspace{0.35ex}
at each point of the hyperbola $xy+1=0$ and at each point of the straight line $x=0,$ 
\vspace{0.35ex}
which is distinct from the origin of coordinates of phase plane  $Oxy,$ 
is parallel to the axis $Ox.$
\vspace{0.5ex}

From the equation $X (x, y) =0$ at $Y (x, y) \ne 0$ it follows that 
\vspace{0.5ex}
an orthogonal isocline of system (16.1) is the cubic parabola $y=x^3$ without the point $O (0,0).$
\vspace{0.35ex}

The tangent to a trajectory of system (16.1) at each point of  the parabola $y=x^3,$ 
\vspace{0.35ex}
which is distinct from the origin of coordinates of phase plane  $Oxy,$ is parallel to the axis $Oy.$
\vspace{0.75ex}

{\it Definite domains of the phase directional field of system}  (16.1).
From the inequality 
\\[1.5ex]
\mbox{}\hfill
$
X (x, y)\;\!Y (x, y)>0
\hfill
$ 
\\[1.25ex]
it follows that domains of positivity for the phase directional field of system (16.1) are  
\\[2ex]
\mbox{}\hfill
$
\Omega_{1}^{+}=\Bigl\{(x,y)\colon x>0\ \&\ {}-\dfrac1{x}<y<x^3\Bigr\}
$ 
\ \ and \ \  
$
\Omega_{2}^{+}=\Bigl\{(x,y)\colon x<0\ \&\ x^3<y<{}-\dfrac1{x}\Bigr\}. 
\hfill
$
\\[2ex]
\indent
The tangent to a trajectory of system (16.1) at each point of the set 
\vspace{0.35ex}
$
\Omega_{}^{+}=\Omega_{1}^{+}\sqcup\Omega_{2}^{+}
$
organises an acute angle with the positive direction of the axis $Ox.$
\vspace{0.35ex}

From the inequality 
\\[1ex]
\mbox{}\hfill
$
X (x, y)\;\!Y (x, y) <0
\hfill
$ 
\\[1.5ex]
it follows that domains of negativity for the phase directional field of system (16.1) are 
\\[2ex]
\mbox{}\hfill
$
\Omega_{1}^{-}=\{(x,y)\colon x>0\ \&\ y>x^3\},
\qquad
\Omega_{2}^{-}=\Bigl\{(x,y)\colon x<0\ \&\ y>{}-\dfrac1{x}\Bigr\},
\hfill
$
\\[2.5ex]
\mbox{}\hfill
$ 
\Omega_{3}^{-}=\{(x,y)\colon x<0\ \&\ y<x^3\}, 
\qquad 
\Omega_{4}^{-}=\Bigl\{(x,y)\colon x>0\ \&\ y<{}-\dfrac1{x}\Bigr\}.
\hfill
$
\\[2ex]
\indent
The tangent to a trajectory of system (16.1) at each point of the set 
\\[1.5ex]
\mbox{}\hfill
$
\Omega_{}^{-}=\Omega_{1}^{-}\sqcup\Omega_{2}^{-}\sqcup\Omega_{3}^{-}\sqcup\Omega_{4}^{-}
\hfill
$
\\[1.25ex]
organises an obtuse angle with the positive direction of the axis $Ox.$
\vspace{0.75ex}

{\it  Contact points of system} (16.1) {\it for the coordinate axes of the phase plane} $Oxy.$ 
\vspace{0.35ex}
The zero isoclines 
\vspace{0.5ex}
$\{(x,y)\colon x=0\ \&\ y\ne 0\}$ and $\{ (x,y) \colon  xy+1=0\}$
haven't common points with the axis $Ox.$
The system (16.1) hasn't contact points on the axis $Ox.$ 
\vspace{0.5ex}

The orthogonal isocline 
\vspace{0.5ex}
$\{ (x,y) \colon  y=x^3\ \&\ x\ne 0\}$ hasn't common points with the axis $Oy. $ 
The system (16.1) hasn't contact points on the axis $Oy.$
\vspace{0.75ex}

{\it  Projective type of system} (16.1).  
The function 
\\[1.75ex]
\mbox{}\hfill
$
W_3^{}\colon  (x,y)\to\ 
x\;\!Y_3^{}(x,y)-y\;\!X_3^{}(x,y)=0
$
\ for all 
$
(x,y)\in \R^2.
\hfill
$
\\[1.75ex]
\indent
The system (16.1) is projectively singular.\!
\vspace{0.35ex}
The infinitely removed straight line of the pro\-j\-ec\-ti\-ve phase plane $\!\R\P (x, y)\!$ 
\vspace{0.75ex}
doesn't consist of trajectories of system (16.1) (Property 9.2). 

{\it The first projectively reduced system for system } (16.1). 
\vspace{0.5ex}
The projectively singular system (16.1) by 
the first transformation of Poincar\'{e} $x =\dfrac{1}{\theta}\,, \ y =\dfrac{\xi}{\theta}$ 
\vspace{0.75ex}
is reduced to the first projectively reduced system (Property 9.2)
\\[2.5ex]
\mbox{}\hfill        % (16.2)
$
\dfrac{d\xi}{d\tau}= \theta+\xi^2\;\!\theta\equiv \Xi(\xi,\theta),
\qquad 
\dfrac{d\theta}{d\tau}={}-1+\xi\;\!\theta^2\equiv \Theta(\xi,\theta),
\qquad
\theta\,d\tau=dt.
$
\hfill (16.2)
\\[-2ex]

\newpage

{\it The second projectively reduced system for system} (16.1). 
\vspace{0.5ex}
The projectively singular system (16.1) by the second transformation of Poincar\'{e} 
$x =\dfrac{1}{\theta}\,, \ y =\dfrac{\xi}{\theta}$ 
\vspace{0.5ex}
is reduced to the second projectively reduced system (Property 9.2)
\\[2.5ex]
\mbox{}\hfill        % (16.3)
$
\dfrac{d\eta}{d\nu}= {}-\zeta^2-\eta^2\zeta\equiv H(\eta,\zeta),
\qquad 
\dfrac{d\zeta}{d\nu}={}-\eta-\eta\zeta^2 \equiv Z(\eta,\zeta),
\qquad
\eta\,d\nu=dt.
$
\hfill (16.3)
\\[2.5ex]
\indent
{\it Equilibrium states of system} (16.2) {\it on the coordinate axis $O ^ {(1)} _ {\phantom1} \xi $ 
\vspace{0.35ex}
in the final part of the projective phase plane} $\R\P(\xi,\theta).$ 
Since
\\[2ex]
\mbox{}\hfill
$
\Theta(\xi,0)={}-1\ne 0  
$ 
\ for all  
$
\xi \in \R,
\hfill
$
\\[1.75ex]
we see that the system (16.2) hasn't equilibrium states 
\vspace{0.35ex}
on the coordinate axis $O ^ {(1)} _ {\phantom1} \xi $ in the final part of the projective phase plane $ \R\P (\xi, \theta).$ 
\vspace{0.5ex}

{\it Equilibrium state of system} (16.3) 
\vspace{0.25ex}
{\it in the origin of coordinates of the phase plane} $O^{(2)}_{\phantom1} \eta\zeta.$
The origin of coordinates of the phase plane $O ^ {(2)} _ {\phantom1} \eta\zeta $ 
\vspace{0.5ex}
is a complicated equilibrium state of system (16.3) with the characteristic equation $\lambda^2=0.$ 
\vspace{0.5ex}
By Theorem 6.2.1 from the monograph [21, pp. 128 -- 129], 
\vspace{0.25ex}
the equilibrium state $O ^ {(2)} _ {\phantom1} (0,0) $ of system (16.3) is the two-separatrix saddle 
which separatrixes adjoin in the direction of coordinate axis $O^{(2)}_{\phantom1} \zeta.$
\vspace{0.75ex}

{\it Equilibrium\! states\! of\! system}\! (16.1)\! {\it in\! the\! projective\! phase\! plane} $\!\R\P(x,y).\!\!$ 
\vspace{0.35ex}
The\! system\! (16.1) in the projective phase plane $ \R\P (x, y) $ has two equilibrium states:  
\vspace{0.35ex}
the unstable focus $O (0,0)$ and the two-separatrix saddle, 
\vspace{0.75ex}
which is lying on <<extremities>> of the straight line $x=0.$

{\it Equatorial contact points of system} (16.1).
\vspace{0.35ex}
The equation $X_3 ^ {} (1, \xi) =0$ hasn't roots. 
\linebreak
The system (16.1) hasn't 
\vspace{0.35ex}
equatorial contact points on <<extremities>> of the straight lines $y=ax,\ a\in\R.$
\vspace{0.35ex}
The two-separatrix saddle is lying on <<extremities>> of the axis $Oy$ and its separatrixes are orthogonal to the axis $Oy. $
\vspace{0.15ex}
The projectively singular system (16.1) hasn't equatorial contact points.
\vspace{0.75ex}

{\it  Limit cycles of system} (16.1) {\it in the projective phase plane $\R\P(x,y).$}
The divergence 
\\[1.75ex]
\mbox{}\hfill
$
{\rm div}\ \vec{a}(x,y)=4x^2
$
\ for all 
$
(x,y)\in\R^2
\hfill
$
\\[1ex]
of the vector field 
\\[1.25ex]
\mbox{}\hfill
$
\vec{a}\colon (x,y)\to\ 
({}-y+x^3,x+x^2y)
$
\ for all 
$
(x,y)\in\R^2
\hfill
$
\\[1.75ex]
is positive on the simply connected domain $\R^2\backslash \{(0,0) \}.$
\vspace{0.35ex}

By the Bendixon criterion [4, p. 120], 
\vspace{0.25ex}
the system (16.1) hasn't limit cycles in the final part of the projective phase plane $\R\P(x,y).$
\vspace{0.25ex}

The set of all trajectories of the projectively singular system (16.1) 
\vspace{0.25ex}
on the projective phase plane $ \R\P (x, y) $ hasn't straight lines.
The system (16.1) hasn't limit linear cycles. 
\vspace{0.25ex}

The projectively singular system (16.1) has one infinitely removed equilibrium state, 
\vspace{0.25ex}
which is the two-separatrix saddle, and this system hasn't equatorial contact points.
\vspace{0.25ex}
The limit cycle cannot surround the two-separatrix saddle.
The system (16.1) hasn't open limit cycles.
\vspace{0.75ex}

{\it The general autonomous integral of system} (16.2). 
The function
\\[2ex]
\mbox{}\hfill
$
F^{(1)}_{\phantom1} \colon (\xi,\theta)\to\
\dfrac{\xi +\theta^2}{1+\xi^2}+\arctan \xi 
$
\ for all 
$
(\xi,\theta)\in \R^2
\hfill
$
\\[2ex]
is the general autonomous integral of system (16.2) on $\R^2.$
\vspace{0.75ex}

{\it Equilibrium states of system} (16.2) {\it in the projective phase plane} $\R\P(\xi,\theta).$ 
\vspace{0.35ex}
The system (16.2) in the projective phase plane $ \R\P (\xi, \theta) $ has 
two equilibrium state: the unstable focus on <<extremities>> of axis $O ^ {(1)} _ {\phantom1} \theta $ 
and the two-separatrix saddle on <<extremities>> of axis $O ^ {(1)} _ {\phantom1} \xi,$ 
which separatrixes adjoin in the direction of axis  $O^{(1)}_{\phantom1}\xi.$ 
\vspace{0.5ex}

{\it Movement of radius-vector of organising point along trajectories of system} (16.2).
The function 
\\[1.25ex]
\mbox{}\hfill
$
W^{(1)}_{\phantom1}\colon (\xi,\theta)\to\ 
\xi \,\Theta(\xi,\theta)-\theta \,\Xi(\xi,\theta)={}-\xi +\theta^2
$
\ for all 
$
(\xi,\theta)\in \R^2
\hfill
$
\\[2ex]
is negative if $ \xi> \theta^2$ 
and this function is positive if $\xi<\theta^2.$
\vspace{0.5ex}

The angle between radius-vector 
\vspace{0.35ex}
of the organising point and the positive direction of the axis 
\vspace{0.5ex}
$O^{(1)}_{\phantom1}\xi$ at movement along the parts of trajectories of system (16.2) from the domains: 
\linebreak
a) $\{(\xi,\theta)\colon \xi>\theta^2\};$ b) $\{(\xi,\theta)\colon \xi<\theta^2\}$ 
is a) decreasing; b) increasing, respectively.
\vspace{0.75ex}

Through each point, which is lying on the parabola $\xi =\theta^2,$ 
\vspace{0.35ex}
the trajectory of system (16.2) passes in the direction of the radius-vector of this point.
\vspace{0.75ex}

{\it Symmetry of the phase directional field of system}  (16.2).
\vspace{0.25ex}
The phase directional field of system (16.2) is symmetric with respect to the coordinate axis 
$O ^ {(1)} _ {\phantom1} \xi $ (Property 13.6).
\vspace{0.35ex}
For each trajectory of system (16.2) there exists a symmetric trajectory 
\vspace{0.25ex}
with respect the coordinate axis $O ^ {(1)} _ {\phantom1} \xi. $
\vspace{0.35ex}
Each trajectory of system (16.2), which is intersecting the coordinate axis $O ^ {(1)} _ {\phantom1} \xi,$ 
is symmetric with respect to this coordinate axis.
\vspace{0.75ex}

{\it Zero and orthogonal isoclines of system}  (16.2).
\vspace{0.5ex}
From the equation $ \Theta (\xi, \theta) =0$ 
it follows that a zero isocline of system (16.2) is the square hyperbola $\xi\! =\!\dfrac{1}{\theta^2}\,.\!$ 
\vspace{0.75ex}
The tangent to a trajec\-to\-ry of system (16.2) at each point of the hyperbola $\xi\! =\!\dfrac{1}{\theta^2}$ 
\vspace{0.75ex}
is parallel to the axis $O_{\phantom1}^{(1)}\xi.$

From the equation $ \Xi (\xi, \theta) =0$ 
\vspace{0.35ex}
it follows that an orthogonal isocline of system (16.2) is $ \theta=0.$
Trajectories of system (16.2) intersect the coordinate axis $O _ {\phantom1} ^ {(1)} \xi $ at right angle.
\vspace{0.75ex}

{\it Definite domains of the phase directional field of system}  (16.2).
From the inequality 
\\[1.5ex]
\mbox{}\hfill
$
\Xi (\xi, \theta)\;\!\Theta (\xi, \theta)> 0
\hfill
$ 
\\[1.5ex]
it follows that domains of positivity of the phase directional field of system (16.2) are  
\\[1.75ex]
\mbox{}\hfill
$
\Omega_{1}^{{}+(1)}=\Bigl\{(\xi,\theta)\colon \xi>\dfrac1{\theta^2}\ \ \& \ \theta>0\Bigr\}
$ 
\ \ 
and 
\ \ 
$
\Omega_{2}^{{}+(1)}=\Bigl\{(\xi,\theta)\colon \xi<\dfrac1{\theta^2}\ \ \& \ \theta<0\Bigr\}.
\hfill
$
\\[2ex]
\indent
The tangent to a trajectory of system (16.2) at each point of the set 
\vspace{0.25ex}
$
\Omega_{\phantom1}^{{}+(1)}=\Omega_{1}^{{}+(1)}\sqcup\Omega_{2}^{{}+(1)}
$
organises an acute angle with the positive direction of the axis  $O_{\phantom1}^{(1)}\xi.$
\vspace{0.5ex}

From the inequality 
\\[1ex]
\mbox{}\hfill
$
X (\xi, \theta)\;\!Y (\xi, \theta) <0 
\hfill
$ 
\\[1.5ex]
it follows that domains of negativity of the phase directional field of system (16.2) are 
\\[1.75ex]
\mbox{}\hfill
$
\Omega_{1}^{{}-(1)}=\Bigl\{(\xi,\theta)\colon \xi<\dfrac1{\theta^2}\ \ \& \ \theta>0\Bigr\}
$
\ \ and \ \ 
$ 
\Omega_{2}^{{}-(1)}=\Bigl\{(\xi,\theta)\colon \xi>\dfrac1{\theta^2}\ \ \& \ \theta<0\Bigr\}.
\hfill
$
\\[2ex]
\indent
The tangent to a trajectory of system (16.2) at each point of the set  
\vspace{0.25ex}
$
\Omega_{\phantom1}^{{}-(1)}=\Omega_{1}^{{}-(1)}\sqcup\Omega_{2}^{{}-(1)}$
organises an obtuse angle with the positive direction of the axis $O_{\phantom1}^{(1)}\xi.$
\vspace{0.5ex}

{\it Contact points of system} (16.2) 
\vspace{0.35ex}
{\it for the coordinate axes of the phase plane} $O^{(1)}_{\phantom1} \xi\theta.$
The zero isocline $ \xi =\dfrac {1} {\theta^2} $ 
\vspace{0.75ex}
does not intersect the axis $O _ {\phantom1} ^ {(1)} \xi. $
The system (16.2) hasn't contact points on the axis $O _ {\phantom1} ^ {(1)} \xi.$ 
\vspace{0.35ex}

The orthogonal isocline $ \theta=0$ 
\vspace{0.35ex}
intersects the axis $O _ {\phantom1} ^ {(1)} \theta $ in one point $O _ {\phantom1} ^ {(1)} (0,0).$ 
The origin of coordinates $O _ {\phantom1} ^ {(1)} (0,0) $ 
is a unique contact point of the axis $O _ {\phantom1} ^ {(1)} \theta $ of system (16.2).
Since $\Theta(0,0)\,\partial_{\theta}^{}\,\Xi(0,0)={}-1<0,$ we see that
\vspace{0.35ex}
the contact $O _ {\phantom1} ^ {(1)} $\!-trajectory of system (16.2)
in enough small neighbourhood of the point $O _ {\phantom1} ^ {(1)} (0,0)$ 
is lying in the half-plane $ \xi\leq 0.$ 
\vspace{0.75ex}

{\it Projective type of system} (16.2). 
\vspace{0.35ex}
The straight line $x=0$ does not consist of trajectories of system (16.1). 
The system (16.2) is projectively singular (Property 10.4). 
\vspace{0.5ex}

The infinitely removed straight line of the projective phase plane $ \R\P (\xi, \theta)$ 
\vspace{0.35ex}
does not consist of trajectories of system (16.2) (Property 9.2).
\vspace{0.75ex}

{\it Equatorial contact points of system} (16.2). 
\vspace{0.35ex}
The system (16.1) hasn't contact points on the axis $Oy.$
\vspace{0.35ex}
The systems (16.1) has the equilibrium state (two-separatrix saddle), which is lying
on <<extremities>> of axes $Oy.$ 
The system (16.2) hasn't equatorial contact points. 
\vspace{0.75ex}

{\it Limit cycles of system} (16.2) {\it in the projective phase plane $\R\P(\xi,\theta).$} 
\vspace{0.25ex}
The system (16.1) hasn't linear limit cycles, open limit cycles, and 
\vspace{0.25ex}
limit cycles, which are lying in the final part of the projective phase plane $\R\P (x, y).$  
\vspace{0.5ex}

The system (16.2) hasn't linear limit cycles, open limit cycles, and 
\vspace{0.35ex}
limit cycles, which are lying in the final part of the projective phase plane $\R\P (\xi, \theta).$
\vspace{1ex}

{\it The general autonomous integral of system} (16.3). 
The function
\\[1.5ex]
\mbox{}\hfill
$
F^{(2)}_{\phantom1} \colon (\eta,\zeta)\to\
\dfrac{\zeta+\eta^2}{1+\zeta^2} + \arctan\dfrac{1}{\zeta} 
$
\ for all
$ 
(\eta,\zeta)\in G^{(2)}_{\phantom1},
\quad   
G^{(2)}_{\phantom1}=\{(\eta,\zeta)\colon \zeta\ne0\},
\hfill
$
\\[1.25ex]
is the general autonomous integral of system (16.3) on any domain from the set $G^{(2)}_{\phantom1}.$
\vspace{0.75ex}

{\it Equilibrium\! states of system}\! (16.3)\! {\it in\! the\! projective\! phase\! plane} $\!\R\P(\eta,\zeta).\!$ 
\vspace{0.35ex}
The system (16.3) in the projective phase plane $\R\P(\eta,\zeta)$  
has two equilibrium states: 
\vspace{0.35ex}
the two-separatrix saddle  $O^{(2)}_{\phantom1} (0,0),$ which separatrixes adjoin in the direction of 
axis $O^{(2)}_{\phantom1}\zeta,$ and
\vspace{0.35ex}
the unstable focus, which is lying on <<extremities>> of axis $O^{(2)}_{\phantom1}\eta.$ 
\vspace{0.75ex}

{\it Movement of radius-vector of organising point along trajectories of system} (16.3).
The function 
\\[1.25ex]
\mbox{}\hfill
$
W^{(2)}_{\phantom1}\colon (\eta,\zeta)\to\ 
\eta \,Z(\eta,\zeta)-\zeta \,H(\eta,\zeta)={}-\eta^2 +\zeta^3
$ 
\ for all 
$
(\eta,\zeta)\in \R^2
\hfill
$
\\[2.25ex]
is negative if $\zeta<\sqrt[\scriptstyle3]{\eta^2}$ and 
this function is positive if $\zeta>\sqrt[\scriptstyle3]{\eta^2}\,.$
\vspace{0.75ex}

The angle between 
\vspace{0.35ex}
radius-vector of the organising point and the positive direction of the axis $O^{(2)}_{\phantom1}\eta$
\vspace{0.5ex}
at movement along the parts of trajectories of system (16.3) from the domains: 
a) $\{(\eta,\zeta)\colon \zeta<\sqrt[\scriptstyle3]{\eta^2}\,\};$  
b) $\{(\eta,\zeta)\colon \zeta>\sqrt[\scriptstyle3]{\eta^2}\,\}$
is a) decreasing; b) increasing, respectively.
\vspace{1ex}

Through each point, which is lying on the curve $\zeta=\sqrt[\scriptstyle3]{\eta^2}\,,$ 
\vspace{0.5ex}
the trajectory of system (16.3) passes in the direction of radius-vector of this point.
\vspace{0.75ex}

{\it Symmetry of the phase directional field of system}  (16.3).
\vspace{0.25ex}
The phase directional field of system (16.3) is symmetric with respect to the coordinate axis $O^{(2)}_{\phantom1} \zeta$ (Property 13.6).
\vspace{0.25ex}
For each trajectory of system (16.3) there exists 
\vspace{0.25ex}
a symmetric trajectory with respect to the coordinate axis $O^{(2)}_{\phantom1} \zeta.$
\vspace{0.35ex}
Each trajectory of system (16.3), which is intersecting the coordinate axis $O^{(2)}_{\phantom1} \zeta$  
(the point of intersection is distinct 
\vspace{0.35ex}
from the origin of coordinates of the phase plane $O^{(2)}_{\phantom1} \eta\zeta),$
is symmetric with respect to this coordinate axis.
\vspace{0.75ex}

{\it Zero and orthogonal isoclines of system}  (16.3).
\vspace{0.35ex}
From the equation $Z(\eta,\zeta)=0$ it follows that
a zero isocline of system (16.3) is
\vspace{0.35ex}
the straight line $ \eta=0$ without the point $O _ {\phantom1} ^ {(2)} (0,0).$ 

Trajectories of system (16.3) intersect the coordinate axis 
\vspace{0.5ex}
$O_{\phantom1}^{(2)}\zeta$ at right angle in each point $(0,\zeta),$ where $\zeta\ne0.$
\vspace{0.5ex}

From the equation  $H(\eta,\zeta)=0$ it follows that 
\vspace{0.75ex}
orthogonal isoclines of system (16.3) are
$\zeta=0$ at $\eta\ne0$ and $\zeta={}-\eta^2$ at $\eta\ne0$.
\vspace{0.5ex}

Trajectories of system (16.3) 
\vspace{0.35ex}
intersect the coordinate axis  $O_{\phantom1}^{(2)}\eta$
at right angle in each point $(\eta,0),$ where $\eta\ne0.$
\vspace{0.35ex}

The tangent to a trajectory of system (16.3)  at each point of the parabola 
\vspace{0.35ex}
$\zeta={}-\eta^2,$ which is distinct from the origin of coordinates of the phase plane 
\vspace{0.35ex}
$O_{\phantom1}^{(2)}\eta\zeta,$  is parallel to the coordinate  axis $O_{\phantom1}^{(2)}\zeta\;\!.$
\vspace{0.75ex}

{\it Definite domains of the phase directional field of system}  (16.3).
From the inequality 
\\[1.5ex]
\mbox{}\hfill
$
H(\eta,\zeta)\,Z(\eta,\zeta)>0
\hfill
$
\\[1.5ex]
it follows that domains of positivity of the phase directional field of system (16.3) are
\\[1.75ex]
\mbox{}\hfill
$
\Omega_{1}^{{}+(2)}=\{(\eta,\zeta)\colon \eta>0\ \ \& \ \zeta>0\},
\qquad  
\Omega_{2}^{{}+(2)}=\{(\eta,\zeta)\colon \eta<0\ \ \& \ {}-\eta^2<\zeta<0\},
\hfill
$
\\[2.5ex]
\mbox{}\hfill
$  
\Omega_{3}^{{}+(2)}=\{(\eta,\zeta)\colon \eta>0\ \ \& \ \zeta<{}-\eta^2\}
\hfill
$
\\[2.5ex]
\indent
The tangent to a trajectory of system (16.3) in each point of the set
\\[1.5ex]
\mbox{}\hfill
$
\Omega_{\phantom1}^{{}+(2)}=\Omega_{1}^{{}+(2)}\sqcup\Omega_{2}^{{}+(2)}\sqcup\Omega_{3}^{{}+(2)}
\hfill
$
\\[1.25ex]
organises an acute angle with the positive direction of the coordinate axis $O_{\phantom1}^{(2)}\eta.$
\vspace{0.35ex}

From the inequality  
\\[0.75ex]
\mbox{}\hfill
$
H(\eta,\zeta)\,Z(\eta,\zeta)<0
\hfill
$ 
\\[1.5ex]
it follows that domains of negativity of the phase directional field of system (16.3) are
\\[1.75ex]
\mbox{}\hfill
$
\Omega_{1}^{{}-(2)}=\{(\eta,\zeta)\colon \eta<0\ \ \& \ \zeta>0\},
\qquad 
\Omega_{2}^{{}-(2)}=\{(\eta,\zeta)\colon \eta<0\ \ \& \ \zeta<{}-\eta^2\},
\hfill
$
\\[2.5ex]
\mbox{}\hfill
$  
\Omega_{3}^{{}-(2)}=\{(\eta,\zeta)\colon \eta>0\ \ \& \ {}-\eta^2<\zeta<0\}.
\hfill
$
\\[2.5ex]
\indent
The tangent to a trajectory of system (16.3) in each point of the set 
\\[1.5ex]
\mbox{}\hfill
$
\Omega_{\phantom1}^{{}-(2)}=\Omega_{1}^{{}-(2)}\sqcup\Omega_{2}^{{}-(2)}\sqcup\Omega_{3}^{{}-(2)}
\hfill
$
\\[1.25ex]
organises an obtuse angle with the positive direction of the coordinate axis $O_{\phantom1}^{(2)}\eta.$
\vspace{0.75ex}

{\it Contact points of system} (16.3) {\it for the coordinate axes of the phase plane} $O^{(2)}_{\phantom1} \eta\zeta.$
\vspace{0.25ex}
The zero isocline $\{(\eta,\zeta)\colon \eta=0\ \,\&\,\ \zeta\ne 0\}$ doesn't intersect the axis $O_{\phantom1}^{(2)}\eta.$
\vspace{0.5ex}
The system (16.3) hasn't contact points on the axis $O_{\phantom1}^{(2)}\eta.$ 
\vspace{0.35ex}

The orthogonal isoclines 
\vspace{0.5ex}
$\{(\eta,\zeta)\colon \zeta=0\ \,\&\,\ \eta\ne 0\}$
and
$\{(\eta,\zeta)\colon \zeta={}-\eta^2\ \,\&\,\ \eta\ne 0\}$
don't intersect the axis  $O_{\phantom1}^{(2)}\zeta.$
The system (16.3) hasn't contact points on the axis $O _ {\phantom1} ^ {(2)} \zeta.$ 
\vspace{0.75ex}

{\it Projective type of system} (16.3). 
\vspace{0.25ex}
The straight line $y=0$ doesn't consist of trajectories of system (16.1). 
The system (16.3) is projectively singular (Property 10.4). 
\vspace{0.25ex}

The infinitely removed straight line of the projective phase plane $\R\P(\eta,\zeta)$ 
\vspace{0.25ex}
doesn't consist of trajectories of system (16.3) (Property 9.2).
\vspace{0.75ex}

{\it Equatorial contact points of system} (16.3). 
\vspace{0.35ex}
The system (16.1) hasn't contact points on the coordinate axis $Ox.$
\vspace{0.25ex}
The system (16.3) hasn't equatorial contact points on <<extremities>> of the straight lines $\zeta=a\eta,$ 
where  $a$ is any real coefficient.
\vspace{0.5ex}
  
The point $O _ {\phantom1} ^ {(1)}$ 
\vspace{0.35ex}
is a contact point of straight line $ \xi=0$ for the differential system (16.2) and 
the contact $O _ {\phantom1} ^ {(1)}\!$-trajectory in enough small neighbourhood of the point $\!O _ {\phantom1} ^ {(1)}\!\!$ 
\vspace{0.35ex}
is lying in the half-plane $\!\xi\leq 0.$ 
\vspace{0.35ex}

The equatorial contact point $O _ {\phantom1} ^ {(1)} $ of system (16.3) 
\vspace{0.35ex}
is lying on <<extremities>> of the coordinate axis $O _ {\phantom1} ^ {(2)} \zeta $  
\vspace{0.35ex}
and the equatorial contact $O _ {\phantom1}^{(1)}\!$-trajectory of system (16.3) 
\vspace{0.35ex}
in enough small neighbourhood 
\vspace{0.35ex}
of the infinitely removed straight line of the projective phase plane $ \R\P (\eta, \zeta)$ 
is lying in the half-plane $\zeta<0.$
\vspace{0.75ex}

{\it Limit cycles of system} (16.3) {\it in the projective phase plane $\R\P(\eta,\zeta).$} 
\vspace{0.25ex}
The system (16.1) hasn't linear limit cycles, open limit cycles, and 
\vspace{0.25ex}
limit cycles, which are lying in the final part of the projective phase plane $\R\P(x,y).$
\vspace{0.25ex}
  
The system (16.3) hasn't linear limit cycles, open limit cycles, and 
\vspace{0.25ex}
limit cycles, which are lying in the final part of the projective phase plane $\R\P(\eta,\zeta).$
\vspace{0.75ex}

{\it Intersection of trajectories of system} (16.1) 
\vspace{0.25ex}
{\it by the infinitely removed straight line of the projective phase plane  $\R\P(x,y).$}
\vspace{0.25ex}
The two-separatrix saddle (its separatrixes are orthogonal to the axis $Oy)$ 
is lying on <<extremities>> of the axis $Oy.$
\vspace{0.25ex}
The straight line $ \theta=0$ is an orthogonal isocline of system (16.2). 
\vspace{0.25ex}
Trajectory of system (16.1) in each point, which 
is not lying on the axis $Oy$
\vspace{0.25ex}
and belong to the boundary circle of the projective circle $ \P\K (x, y),$ 
orthogonally intersects the boundary circle.
\vspace{0.5ex}

{\it Intersection of trajectories of system} (16.2) 
\vspace{0.25ex}
{\it by the infinitely removed straight line of the projective phase plane $\R\P(\xi,\theta).$}
\vspace{0.25ex}
Equilibrium states is lying on <<extremities>> of the coordinate axes $O _ {\phantom1} ^ {(1)} \xi $ and $O _ {\phantom1} ^ {(1)} \theta.$ 
\vspace{0.25ex}
The straight line $x=0$ without the point $O (0,0)$ is a zero isocline of system (16.1). 
\vspace{0.5ex}
Trajectory of system (16.2) in each point, which is not lying 
on the coordinate axes $O _ {\phantom1} ^ {(1)} \xi, \ O _ {\phantom1} ^ {(1)} \theta$ 
\vspace{0.35ex}
and belong to the boundary circle of the projective circle $ \P\K (\xi, \theta),$ 
orthogonally intersects of the boundary circle.
\vspace{0.75ex}

{\it The projective atlas of trajectories for system} (16.1) is constructed on Fig. 16.1.
\\[3.75ex]
\mbox{}\hfill
{\unitlength=1mm
\begin{picture}(42,42)
\put(0,0){\includegraphics[width=42mm,height=42mm]{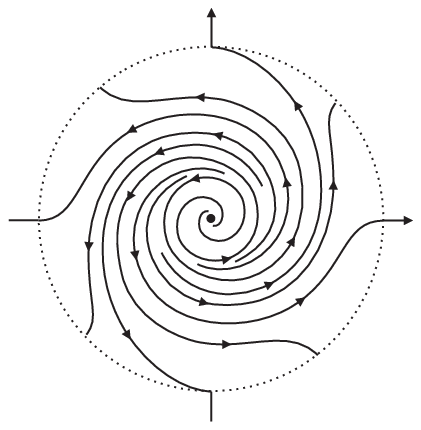}}
\put(18,41){\makebox(0,0)[cc]{ $y$}}
\put(40.2,18.2){\makebox(0,0)[cc]{ $x$}}
\end{picture}}
\qquad
{\unitlength=1mm
\begin{picture}(42,42)
\put(0,0){\includegraphics[width=42mm,height=42mm]{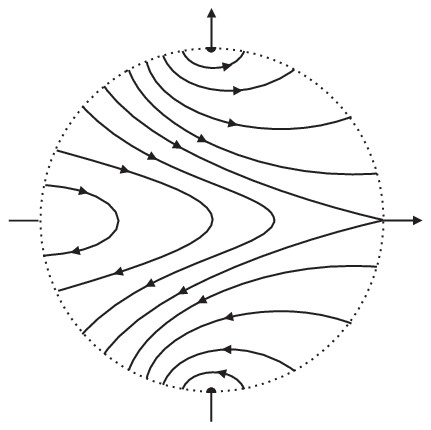}}
\put(18,41){\makebox(0,0)[cc]{ $\theta$}}
\put(40.2,17.8){\makebox(0,0)[cc]{ $\xi$}}
\put(21,-7){\makebox(0,0)[cc]{Fig. 16.1}}
\end{picture}}
\qquad
{\unitlength=1mm
\begin{picture}(42,42)
\put(0,0){\includegraphics[width=42mm,height=42mm]{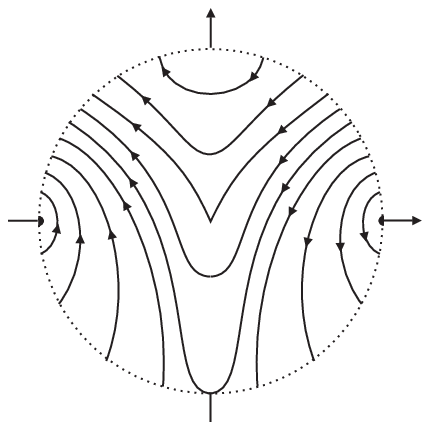}}
\put(18,41){\makebox(0,0)[cc]{ $\zeta$}}
\put(40.2,18){\makebox(0,0)[cc]{ $\eta$}}
\end{picture}}
\hfill\mbox{}
\\[9ex]
\indent
{\bf 16.2.
Trajectories of differential system} [1, pp. 84 -- 85]
\\[2ex]
\mbox{}\hfill        % (16.4)
$
\dfrac{dx}{dt}=1-x^2-y^2\equiv X(x,y),\ \ \dfrac{dy}{dt}=xy-1\equiv Y(x,y).
$
\hfill (16.4)
\\[2.5ex]
\indent
Since the circle $x^2+y^2=1$ and the hyperbola $xy=1$ don't intersect,
\vspace{0.25ex}
we see that the differential system (16.4) 
\vspace{0.15ex}
hasn't equilibrium states in the final part of the projective phase plane $\R\P (x, y).$ 
 
The phase directional field of system (16.4) is symmetric with respect to the origin of coordinates of the phase plane $Oxy.$ 
For each trajectory of system (16.4) there exists a symmetric trajectory with respect to the origin of  coordinates 
of the phase plane $Oxy.$ 
Trajectory of system (16.4), which is passing through the origin of coordinates of the phase plane $Oxy,$ 
is symmetric with respect to the origin of coordinates of this phase plane.
\vspace{0.25ex}

A zero isocline of system (16.4) is the hyperbola $xy=1.$ 
\vspace{0.25ex}
The tangent to a trajectory of system (16.4) at each point of the hyperbola $xy=1$ is parallel to the axis $Ox.$
\vspace{0.25ex}

An orthogonal isocline of system (16.4) is the circle $x^2+y^2=1.$ 
\vspace{0.25ex}
The tangent to a trajectory of system (16.4) at each point of the circle $x^2+y^2=1$ is parallel to the axis $Oy.$

The tangent to a trajectory of system (16.4) at each point of the domain of positivity 
\\[1.75ex]
\mbox{}\hfill
$
\Omega_{\phantom1}^{+}=\{(x,y)\colon  x^2+y^2>1\ \&\ xy<1\}
\hfill
$ 
\\[1.5ex]
of phase directional field   
\vspace{0.35ex}
organises an acute angle with the positive direction of the axis $Ox.$

Domains of negativity for the phase directional field of system (16.4) are the domains
\\[1.5ex]
\mbox{}\hfill
$
\Omega_{1}^{-}=\{(x,y)\colon x^2+y^2<1\},
\  
\Omega_{2}^{-}=\{(x,y)\colon x<0\ \&\ xy>1\},
\ 
\Omega_{3}^{-}=\{(x,y)\colon x>0\ \&\ xy>1\}.
\hfill
$
\\[1.75ex]
\indent
The tangent to a trajectory of system (16.4) at each point of the set 
\vspace{0.25ex}
$\Omega_{\phantom1}^{-}=\Omega_1^-\sqcup \Omega_2^-\sqcup \Omega_3^-$
organises an obtuse angle with the positive direction of the axis $Ox.$
\vspace{0.35ex}

The zero isocline $xy=1$ hasn't common points with the axis $Ox.$ 
\vspace{0.25ex}
The system (16.4) hasn't contact points on the axis $Ox.$
\vspace{0.35ex}

The orthogonal isocline $x^2+y^2=1$ intersects the axis $Oy $ at the points $(0, {} \pm1).$ 
\vspace{0.25ex}
The system (16.4) has two contact points $A_1 ^ {} (0, {}-1) $ and $A_2 ^ {} (0,1)$ on the axis $Oy.$
\vspace{0.75ex}
 
If $y = {}-1,$ then the product 
\vspace{0.5ex}
$Y (0, y)\;\!\partial_y^{} X(0, y) = 2\;\!y$ 
is negative, 
and if $y=1,$ then this product is positive.
\vspace{0.35ex}
The contact $A_1 ^ {}\!$-trajectory in enough small neighbourhood of the point $A_1 ^ {}$ 
is lying  in the half-plane $x\leq 0,$ and 
\vspace{0.35ex}
the contact $A_2 ^ {}\!$-trajectory in enough small neighbourhood of the point $A_2 ^ {}$ 
is lying in the half-plane $x\geq 0.$ 
\vspace{0.25ex}

The function 
\\[1ex]
\mbox{}\hfill
$
W_2^{}\colon (x,y)\to\ 
y(2x^2+y^2)
$
\ for all 
$
(x,y)\in \R^2
\hfill
$
\\[1.5ex]
is not identically zero on the plane $\R^2.$
\vspace{0.15ex}
The differential system (16.4) is projectively nonsingular. 
\vspace{0.15ex}
The infinitely removed straight line of the projective phase plane $\R\P(x, y)$ 
consists of trajectories of system (16.4). 
\vspace{0.35ex}

The first and the second projectively reduced systems for system (16.4) are the systems  
\\[2ex]
\mbox{}\hfill        % (16.5)
$
\dfrac{d\xi}{d\tau}= 2\xi-\theta^2+\xi\bigl(\xi^2-\theta^2\bigr)\equiv \Xi(\xi,\theta),
\quad 
\dfrac{d\theta}{d\tau}=\theta\bigl(1+\xi^2-\theta^2\bigr)\equiv \Theta(\xi,\theta),
\quad
\theta\, d\tau=dt,
$
\hfill (16.5)
\\[2ex]
and
\\[1.75ex]
\mbox{}\hfill        % (16.6)
$
\dfrac{d\eta}{d\nu}= {}-\eta\bigl(\zeta-\eta^2\bigr)\equiv H(\eta,\zeta),
\quad 
\dfrac{d\zeta}{d\nu}={}-1+\eta^2-2\zeta^2 + \eta^2\zeta\equiv Z(\eta,\zeta),
\quad
\eta\,d\nu=dt.
$
\hfill (16.6)
\\[3ex]
\indent
The system of equations 
\\[1ex]
\mbox{}\hfill
$\Xi(\xi,0)=0,
\quad 
\Theta(\xi,0)=0
\hfill
$ 
\\[1.5ex]
has one solution $\xi=0.$
\vspace{0.25ex}
The system (16.5) has one equilibrium state $O ^ {(1)} _ {\phantom1} (0,0)$ 
on the coordinate axis $O ^ {(1)} _ {\phantom1} \xi $ in the final part of the projective phase plane $ \R\P (\xi, \theta).$
\vspace{0.5ex}
This equilibrium state is an unstable node with the characteristic equation $\lambda^2-3\lambda+2=0.$ 
\vspace{0.75ex}

Since $Z(0,0)={}-1\ne 0,$
\vspace{0.35ex}
we see that the origin of coordinates of the phase plane $O ^ {(2)} _ {\phantom1} \eta\zeta $ 
is not an equilibrium state of system (16.6).  
\vspace{0.35ex}

The differential system (16.4) has one equilibrium state $O _ {\phantom1} ^ {(1)}$ 
\vspace{0.35ex}
in the projective phase plane $ \R\P (x, y).$
\vspace{0.25ex}
The equilibrium state $O _ {\phantom1} ^ {(1)}$ is lying on <<extremities>> of the straight line $y=0$ and 
the equilibrium state $O _ {\phantom1} ^ {(1)}$ is an unstable node.
\vspace{0.35ex}

The system (16.4) hasn't equilibrium states 
\vspace{0.25ex}
in the final part of the projective phase plane $ \R\P (x, y).$
Hence the system (16.4) hasn't limit cycles.
\vspace{0.35ex}

The system (16.4) is projectively nonsingular 
\vspace{0.25ex}
and this systen has infinitely removed equilibrium state. 
\vspace{0.5ex}
Thus the system (16.4) hasn't  linear limit cycles and open limit cycles.

The system (16.5) has one equilibrium state $O _ {\phantom1} ^ {(1)} (0,0)$ 
\vspace{0.35ex}
in the projective phase plane $ \R\P (\xi, \theta).$
This equilibrium state is an unstable node.
\vspace{0.35ex}

\newpage

The system (16.5) hasn't linear cycles, open cycles, and cycles, which are lying 
\vspace{0.35ex}
in the final part of the projective phase plane $ \R\P (\xi, \theta),$ (including limit cycles). 
\vspace{0.5ex}

The phase directional field of system (16.5) 
\vspace{0.35ex}
is symmetric with respect to the coordinate axis $O ^ {(1)} _ {\phantom1} \xi.$ 
\vspace{0.35ex}
For each trajectory of system (16.5) there exists a symmetric trajectory with respect 
to the coordinate axis $O^{(1)}_{\phantom1} \xi.$
\vspace{0.5ex}

From the equation $ \Theta (\xi, \theta) =0,$ we get two equations $ \theta=0$ and $ \theta^2-\xi^2=1.$
\vspace{0.5ex}
Zero isoclines of system (16.5) are the hyperbola $ \theta^2-\xi^2=1$ and 
\vspace{0.35ex}
the straight line $ \theta=0$ without the point $O _ {\phantom1} ^ {(1)} (0,0).$
\vspace{0.35ex}
The straight line $ \theta=0$ consists of trajectories of system (16.5).
The tangent to a trajectory of system (16.5) at each point of the hyperbola $\theta^2-\xi^2=1$ 
\vspace{0.35ex}
is parallel to the coordinate axis $O_{\phantom1}^{(1)}\xi.$ 
\vspace{0.5ex}

The equation $\Xi (0, \theta) =0$ has one root $ \theta=0,$ but $ \Theta (0,0) =0. $
\vspace{0.35ex}
The system (16.5) hasn't contact points on the coordinate axis $O^{(1)} _ {\phantom1} \theta.$ 
\vspace{0.35ex}

The straight line $x=0$ does not consist of trajectories of system (16.4). 
\vspace{0.25ex}
The system (16.5) is projectively singular. 
\vspace{0.35ex}
The infinitely removed straight line of the projective phase plane $ \R\P (\xi, \theta) $ does not consist 
of trajectories of system (16.5).
\vspace{0.5ex}

The equatorial contact points of system (16.5), which are lying on <<extremities>> 
\vspace{0.35ex}
of the straight lines $ \theta = {}-\xi $ and $ \theta =\xi,$ corresponding to 
\vspace{0.35ex}
the contact points $A_1 ^ {} (0, {}-1) $ and $A_2 ^ {} (0,1)$ 
of the axis $Oy $ for system (16.4). 
\vspace{0.35ex}
The equatorial contact $A_1 ^ {}\!$- and $A_2 ^ {}\!$-trajectories of system (16.5) 
\vspace{0.35ex}
in enough small neighbourhood of the infinitely removed straight line 
of the projective phase plane $\R\P(\xi, \theta)$ are lying in the half-plane $\xi>0.$
\vspace{0.5ex}

The system (16.6) has one equilibrium state $O _ {\phantom1} ^ {(1)}$ 
\vspace{0.25ex}
in the projective phase plane $ \R\P (\eta, \zeta).$
The equilibrium state $O _ {\phantom1} ^ {(1)}$ is lying on 
\vspace{0.25ex}
 <<extremities>> of the straight line $ \eta=0$ 
and this equilibrium state is an unstable node.
\vspace{0.5ex}

The system (16.6) hasn't linear cycles, open cycles, and cycles, which are lying 
\vspace{0.35ex}
in the final part of the projective phase plane $\R\P (\eta, \zeta),$ (including limit cycles). 
\vspace{0.5ex}

The phase directional field of system (16.6) 
\vspace{0.35ex}
is symmetric with respect to the axis $O ^ {(2)} _ {\phantom1} \zeta.$ 
For each trajectory of system (16.6) 
\vspace{0.35ex}
there exists a symmetric trajectory with respect to the axis $O ^ {(2)} _ {\phantom1} \zeta. $
\vspace{0.35ex}
The straight line-trajectory $ \eta=0$ of system (16.6) is symmetric with respect to the origin of coordinates 
of the phase plane $O ^ {(2)} _ {\phantom1} \eta\zeta.$
\vspace{0.75ex}

From the equation $H (\eta, \zeta) =0,$ we get two equations $\eta=0$ and $\zeta =\eta^2.$
\vspace{0.5ex}
Orthogonal isoclines of system (16.6) are the straight line $\eta=0$ and the parabola $\zeta =\eta^2.$
\vspace{0.35ex}
The orthogonal isocline $\eta=0$ of system (16.6) is a trajectory of system (16.6). 
\vspace{0.35ex}
The tangent to a trajectory of system (16.6) at each point of the parabola $ \zeta =\eta^2$ 
is parallel to the axis $O^{(2)}_{\phantom1} \zeta.$
\vspace{0.75ex}

The equation $Z (\eta, 0) =0$ has two roots 
\vspace{0.75ex}
$\eta_1 ^ {} = {}-1$ and $ \eta_2 ^ {} = 1,$ and $H ({} \pm1,0) = \pm1\ne 0. $
The system (16.6) has two contact points $A_1 ^ {(2)} ({}-1,0) $ and $A_2 ^ {(2)} (1,0) $
\vspace{0.75ex}
on the axis $O ^ {(2)} _ {\phantom1} \eta.$ 
The product $H (\eta, 0) \; \! \partial_\eta ^ {} Z (\eta, 0) =2\eta^4$  
\vspace{0.35ex}
is positive as at $\eta = {}-1, $ and at $ \eta=1.$
The trajectories of system (16.6), which are passing through the contact points $A_1 ^ {(2)} $ and $A_2 ^ {(2)},$ 
\vspace{0.75ex}
in enough small neighbourhood of each of these points are lying in the half-plane $\zeta\geq 0.$ 
\vspace{0.75ex}

The straight line $y=0$ does not consist of trajectories of system (16.4). 
\vspace{0.35ex}
The system (16.6) is projectively singular. 
\vspace{0.35ex}
The infinitely removed straight line of the projective phase plane $\R\P (\eta, \zeta)$ 
does not consist of trajectories of system (16.6).
\vspace{0.75ex}

The differential system (16.4) hasn't contact points on the axis $Ox.$ 
\vspace{0.35ex}
The system (16.4) has an equilibrium state on <<extremities>> of the straight line $y=0.$
\vspace{0.35ex}
Therefore the differential system (16.6) hasn't  equatorial contact points.

The direction of movement along trajectories of systems (16.4) -- (16.6) is defined by instability of the node  
$O^{(1)}_{\phantom1}.$
\vspace{0.35ex}

The projective atlas of trajectories for system (16.4) is constructed on Fig. 16.2.
\\[3.75ex]
\mbox{}\hfill
{\unitlength=1mm
\begin{picture}(42,42)
\put(0,0){\includegraphics[width=42mm,height=42mm]{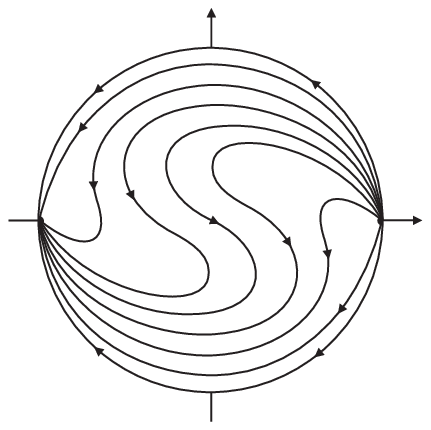}}
\put(18,41){\makebox(0,0)[cc]{ $y$}}
\put(40.2,18.5){\makebox(0,0)[cc]{ $x$}}
%\put(21,-3){\makebox(0,0)[cc]{ $Oxy$}}
%\put(22.5,-6){\makebox(0,0)[cc]{Рис. 1}}
\end{picture}}
\qquad
{\unitlength=1mm
\begin{picture}(42,42)
\put(0,0){\includegraphics[width=42mm,height=42mm]{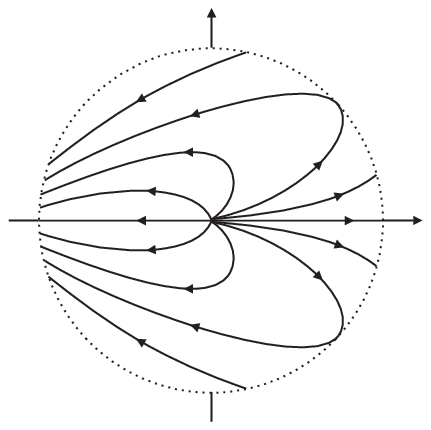}}
\put(18,41){\makebox(0,0)[cc]{ $\theta$}}
\put(40.2,17.8){\makebox(0,0)[cc]{ $\xi$}}
%\put(21,-3){\makebox(0,0)[cc]{ $O^{{}^{(1)}}uz$}}
\put(21,-7){\makebox(0,0)[cc]{Fig. 16.2}}
\end{picture}}
\qquad
{\unitlength=1mm
\begin{picture}(42,42)
\put(0,0){\includegraphics[width=42mm,height=42mm]{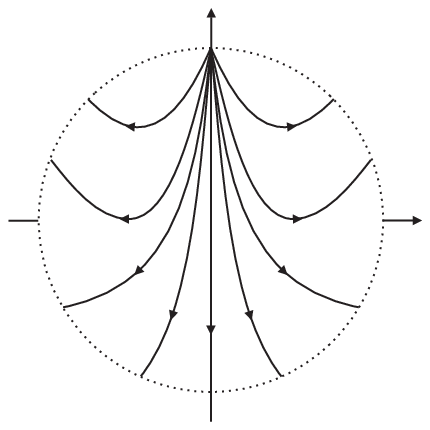}}
\put(18,41){\makebox(0,0)[cc]{ $\zeta$}}
\put(40.2,18){\makebox(0,0)[cc]{ $\eta$}}
%\put(21,-3){\makebox(0,0)[cc]{ $O^{{}^{(2)}}zv$}}
%\put(22.5,-6){\makebox(0,0)[cc]{Рис. 1}}
\end{picture}}
\hfill\mbox{}
\\[9.5ex]
\indent
{\bf 16.3.
Trajectories of differential system} [1, pp. 85 -- 87; 2, pp. 209 -- 212] 
\\[2.25ex]
\mbox{}\hfill        % (16.7)
$
\dfrac{dx}{dt}={}-1+x^2+y^2\equiv X(x,y),
\qquad 
\dfrac{dy}{dt}={}-5+5xy\equiv Y(x,y).
$
\hfill (16.7)
\\[2.75ex]
\indent
The differential system (16.7) hasn't equilibrium states in the final part of the projective phase plane $\R\P (x, y).$
\vspace{0.35ex}

The phase directional field of system (16.7) is symmetric with respect to the origin of coordinates of the phase plane $Oxy.$  
For each trajectory of system (16.7) there exists a symmetric trajectory with respect to 
the origin of  coordinates of the phase plane $Oxy.$ 
The trajectory of system (16.7), which is passing through the origin of coordinates of the phase plane $Oxy,$ 
is symmetric with respect to the origin of coordinates of this phase plane.
\vspace{0.35ex}

A zero isocline of system (16.7) is the hyperbola $xy=1.$ 
\vspace{0.35ex}
The tangent to a trajectory of system (16.7) at each point of the hyperbola $xy=1$ is parallel to the axis $Ox.$
\vspace{0.5ex}

An orthogonal isocline of system (16.7) is the circle $x^2+y^2=1.$ 
\vspace{0.35ex}
The tangent to a trajectory of system (16.7) at each point of the circle $x^2+y^2=1$ is parallel to the axis $Oy.$
\vspace{0.5ex}

Domains of positivity for the phase directional field of system (16.7) are the domains 
\\[1.5ex]
\mbox{}\hfill
$
\Omega_{1}^{+}=\{(x,y)\colon  x^2+y^2<1\},
\ 
\Omega_{2}^{+}=\{(x,y)\colon  x<0\ \&\ xy>1\},
\ 
\Omega_{3}^{+}=\{(x,y)\colon  x>0\ \&\ xy>1\}.
\hfill
$
\\[1.75ex]
\indent
The tangent to a trajectory of system (16.7) at each point of the set 
\vspace{0.35ex}
$
\Omega_{}^{+}=
\Omega_{1}^{+}\sqcup
\Omega_{2}^{+}\sqcup
\Omega_{3}^{+}
$
organises an acute angle with the positive direction of the axis $Ox.$
\vspace{0.35ex}

The tangent to a trajectory of system (16.7) at each point of the domain of negativity 
\\[1.5ex]
\mbox{}\hfill
$
\Omega_{}^{-}=\{(x,y)\colon x^2+y^2>1\ \&\ xy<1\}
\hfill
$
\\[1.5ex]
for the phase directional field of system (16.7)
organises an obtuse angle with the positive direction of the axis $Ox.$
\vspace{0.5ex}

The zero isocline $xy=1$ hasn't common points with the axis $Ox.$ 
\vspace{0.35ex}
The system  (16.7) hasn't contact points on the axis $Ox.$ 
\vspace{0.5ex}

The orthogonal isocline $x^2+y^2=1$ intersects the axis $Oy $ at the points $ (0, {} \pm1).$ 
\vspace{0.35ex}
The system (16.7) has two contact points $A_1 ^ {} (0, {}-1) $ and $A_2 ^ {} (0,1)$ on the axis $Oy.$
\vspace{0.5ex}
If $y={}-1,$ then the product $Y(0,y)\;\!\partial_y X(0,y)= {}-10y$  
\vspace{0.5ex}
is negative, and if $y=1,$ then this product is positive.
\vspace{0.5ex}
The contact $A_1^{}\!$-trajectory in enough small neighbourhood of the point $A_1 ^ {} $ 
is lying in the half-plane $x\geq 0,$ and 
\vspace{0.5ex}
the contact $A_2^{}\!$-trajectory in enough small neighbourhood of the point $A_2 ^ {}$ 
is lying in the half-plane $x\leq 0.$ 
\vspace{0.5ex}

The function 
\\[1ex]
\mbox{}\hfill
$
W_2^{}\colon (x,y)\to\ 
y(4x^2-y^2)
$
\ for all 
$
(x,y)\in \R^2
\hfill
$
\\[1.75ex]
is not identically zero on the plane $\R^2.$  
\vspace{0.15ex}
The differential system (16.7) is projectively nonsingular. 
\vspace{0.15ex}
The infinitely removed straight line of the projective phase plane $ \R\P (x, y)$ consists of trajectories of system (16.7). 
\vspace{0.35ex}

The first and the second projectively reduced systems for system (16.4) are the systems    
\\[2.5ex]
\mbox{}\hfill        % (16.8)
$
\dfrac{d\xi}{d\tau}= 4\xi-5\theta^2-\xi\bigl(\xi^2-\theta^2\bigr)\equiv \Xi(\xi,\theta),
\quad 
\dfrac{d\theta}{d\tau}={}-\theta\bigl(1+\xi^2-\theta^2\bigr)\equiv \Theta(\xi,\theta),
\quad
\theta\, d\tau=dt,
$
\hfill (16.8)
\\[2.25ex]
and
\\[2.25ex]
\mbox{}\hfill        % (16.9)
$
\dfrac{d\eta}{d\nu}= {}-5\eta\bigl(\zeta-\eta^2\bigr)\equiv H(\eta,\zeta),
\quad 
\dfrac{d\zeta}{d\nu}=1-\eta^2-4\zeta^2 + 5\eta^2\zeta\equiv Z(\eta,\zeta),
\quad
\eta\,d\nu=dt.
$
\hfill (16.9)
\\[3ex]
\indent
The differential system (16.8) has three equilibrium states  
\vspace{0.25ex}
on the coordinate axis $O ^ {(1)} _ {\phantom1} \xi$ in the final part of the projective phase plane 
$\R\P (\xi, \theta)\colon\!$ 
\vspace{0.35ex}
the saddle  $O^{(1)}_{\phantom1}(0,0);$
the stable no\-des $A^{(1)}_{\phantom1}({}-2,0)$ and $B^{(1)}_{\phantom1}(2,0).$
The characteristic equation of these equilibrium points is
\\[1.75ex]
\mbox{}\hfill
$
\lambda^2+(4\xi^2-3)\lambda+(\xi^2+1)(3\xi^2-4)=0.
\hfill
$
\\[1.75ex]
\indent
Since $Z(0,0)=1\ne 0,$ we see that 
\vspace{0.25ex}
the origin of coordinates of the phase plane $O^{(2)}_{\phantom1} \eta\zeta$
is not an equilibrium state of system (16.9).  
\vspace{0.5ex}

The system (16.7) 
\vspace{0.35ex}
has three equilibrium states in the projective phase plane $\R\P(x,y)\colon$ 
the saddle $O_{\phantom1}^{(1)},$ which is lying on <<extremities>> of the straight line $y=0;$ 
\vspace{0.35ex}
the stable node $A^{(1)}_{\phantom1},$
which is lying on <<extremities>> of the straight line $y={}-2x;$  
\vspace{0.35ex}
the stable node $B^{(1)}_{\phantom1},$
which is lying on <<extremities>> of the straight line $y=2x.$  
\vspace{0.5ex}

The system (16.4) 
\vspace{0.25ex}
in the final part of the projective phase plane $ \R\P (x, y)$ hasn't equilibrium states. 
Therefore this system hasn't limit cycles.
\vspace{0.35ex}

The system (16.7) is projectively nonsingular 
\vspace{0.25ex}
and this system has infinitely removed equilibrium state. 
The system (16.7) hasn't linear limit cycles and open limit cycles.
\vspace{0.5ex}

The system (16.8) 
\vspace{0.35ex}
has three equilibrium states in the projective phase plane $\R\P(\xi,\theta)\colon$ 
$O_{\phantom1}^{(1)}(0,0)$ is a saddle;
$A^{(1)}_{\phantom1}({}-2, 0)$  and $B^{(1)}_{\phantom1}(2, 0)$ are stable nodes.
\vspace{0.75ex}

The system (16.8) hasn't linear limit cycles, open limit cycles, and 
\vspace{0.35ex}
limit cycles, which are lying in the final part of the projective phase plane $\R\P(\xi,\theta).$
\vspace{0.5ex}

The phase directional field of system (16.8) 
\vspace{0.35ex}
is symmetric with respect to the axis $O^{(1)}_{\phantom1} \xi.$ 
For each trajectory of the differential system (16.8) 
\vspace{0.25ex}
there exists a symmetric trajectory with respect to the axis $O^{(1)}_{\phantom1} \xi.$
\vspace{0.75ex}

Using the equation $ \Theta (\xi, \theta) =0,$ we obtain 
\vspace{0.35ex}
the system (16.8) has one zero isocline $ \theta=0.$ This zero isocline consists of trajectories of system (16.8).
\vspace{0.75ex}

The equation $\Xi(0,\theta)=0$ has one root $\theta=0,$ but $\Theta(0,0)=0.$ 
\vspace{0.35ex}
The system (16.8) hasn't contact points on the coordinate axis $O^{(1)}_{\phantom1} \theta.$
\vspace{0.5ex}
 
The straight line $x=0$ does not consist of trajectories of system (16.7). 
\vspace{0.35ex}
The differential system (16.8) is projectively singular.  
\vspace{0.35ex}
The infinitely removed straight line of the projective phase plane $\R\P(\xi,\theta)$ 
does not consist of trajectories of system (16.8).
\vspace{0.5ex}

The equatorial contact points of system (16.8), which are lying on <<extremities>> 
\vspace{0.35ex}
of the straight lines $\theta={}-\xi$ and $\theta=\xi,$ 
\vspace{0.35ex}
corresponding to the contact points $A_1^{}(0,{}-1)$ and $A_2^{}(0,1)$ of the axis $Oy$ for system (16.7). 

\newpage

The equatorial contact trajectories of system (16.8) in an 
\vspace{0.5ex}
neighbourhood of the infinitely removed straight line of the projective phase plane $\R\P(\xi,\theta)$ 
\vspace{0.75ex}
are lying in the half-plane $\xi<0.$ 

The system (16.9) has 
\vspace{0.5ex}
three equilibrium states in the projective phase plane $\R\P(\eta,\zeta)\colon$
the saddle $O_{\phantom1}^{(1)}$ on  <<extremities>> of the straight line $\eta=0,$
\vspace{0.25ex}
the unstable node $A_{\phantom1}^{(2)}\Bigl(0, {}-\dfrac{1}{2}\Bigr),$ and
the stable node $B_{\phantom1}^{(2)}\Bigl(0, \dfrac{1}{2}\Bigr).$
\vspace{0.75ex}

The system (16.9) hasn't linear limit cycles, 
\vspace{0.25ex}
open limit cycles, and 
limit cycles, which are lying in the final part of the projective phase plane $\R\P(\eta,\zeta).$ 
\vspace{0.5ex}

The phase directional field of system (16.9) 
\vspace{0.25ex}
is symmetric with respect to the axis $O^{(2)}_{\phantom1} \zeta.$ 
For each trajectory of the differential system (16.9) 
\vspace{0.25ex}
there exists a symmetric trajectory with respect to 
the axis $O^{(2)}_{\phantom1} \zeta.$
\vspace{0.5ex}

From the equation 
\vspace{0.35ex}
$H(\eta,\zeta)=0,$ we get two equations $\eta=0$ and $\zeta=\eta^2.$
Orthogonal isoclines of system (16.9) are the straight line $\eta=0$ and the parabola $\zeta=\eta^2.$
\vspace{0.35ex}
The orthogonal isocline $\eta=0$ of system (16.9) consists of trajectories of system (16.9). 
\vspace{0.25ex}
The tangent to a trajectory of system (16.9) at each point of the parabola 
\vspace{0.75ex}
$\zeta=\eta^2$ is parallel to the axis $O^{(2)}_{\phantom1} \zeta.$

The equation $Z(\eta,0)=0$
\vspace{0.5ex}
has two roots $\eta_1^{}={}-1$ and $\eta_2^{}=1,$ and 
$H({}\pm1,0)=\pm1\ne 0.$
The system (16.9) has two contact points $A_1^{(2)}({}-1,0)$ и $A_2^{(2)}(1,0)$
\vspace{0.5ex}
on the coordinate axis $O^{(2)}_{\phantom1} \eta.$ 
The product 
$H(\eta,0)\;\!\partial_\eta^{} Z(\eta,0)={}-10\eta^4$
\vspace{0.5ex}
is negative as at $\eta={}-1,$ and at $\eta=1.$
The trajectories of system (16.9), 
\vspace{0.5ex}
which are passing through the contact points $A_1^{(2)}$ and $A_2^{(2)},$ 
in enough small neighbourhood of each of these points 
\vspace{0.75ex}
are lying in the half-plane $\zeta\leq 0.$ 

The straight line $y=0$ does not consist of trajectories of system (16.7). 
\vspace{0.25ex}
The differential system (16.9) is projectively singular. 
\vspace{0.25ex}
The infinitely removed straight line of the projective phase plane $\R\P(\eta,\zeta)$ 
does not consist of trajectories of system (16.9).
\vspace{0.5ex}

The system (16.7) hasn't contact points on the axis $Ox.$
\vspace{0.25ex}
The system (16.7) has the equilibrium state, which is lying
on <<extremities>> of the straight line $y=0.$
\vspace{0.25ex}
The differential system (16.9) hasn't equatorial contact points.
\\[3.75ex]
\mbox{}\hfill
{\unitlength=1mm
\begin{picture}(42,42)
\put(0,0){\includegraphics[width=42mm,height=42mm]{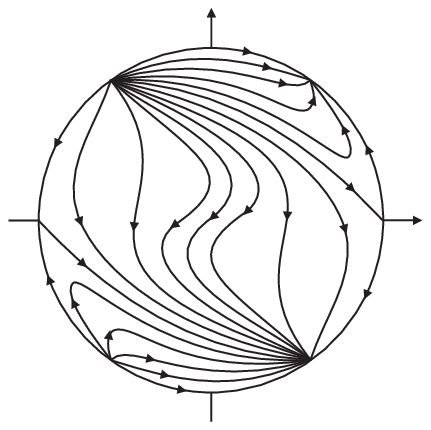}}
\put(18,41){\makebox(0,0)[cc]{ $y$}}
\put(40.2,18.5){\makebox(0,0)[cc]{ $x$}}
%\put(21,-3){\makebox(0,0)[cc]{ $Oxy$}}
%\put(22.5,-6){\makebox(0,0)[cc]{Рис. 1}}
\end{picture}}
\qquad
{\unitlength=1mm
\begin{picture}(42,42)
\put(0,0){\includegraphics[width=42mm,height=42mm]{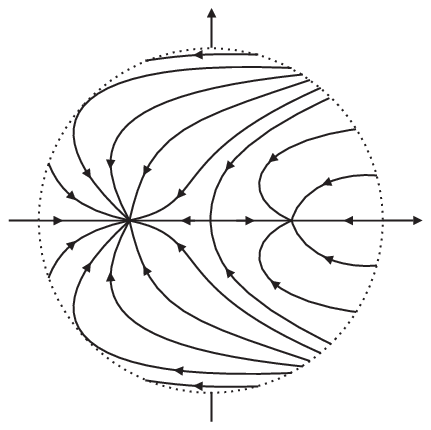}}
\put(18,41){\makebox(0,0)[cc]{ $\theta$}}
\put(40.2,17.8){\makebox(0,0)[cc]{ $\xi$}}
%\put(21,-3){\makebox(0,0)[cc]{ $O^{{}^{(1)}}uz$}}
\put(21,-7){\makebox(0,0)[cc]{Fig. 16.3}}
\end{picture}}
\qquad
{\unitlength=1mm
\begin{picture}(42,42)
\put(0,0){\includegraphics[width=42mm,height=42mm]{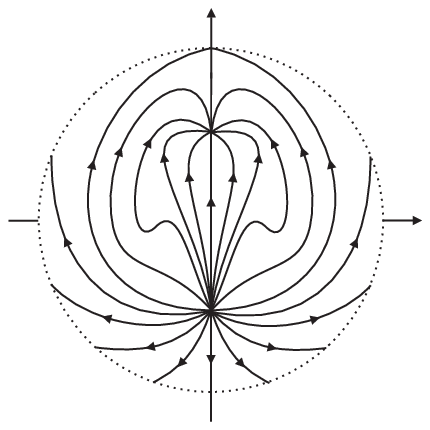}}
\put(18,41){\makebox(0,0)[cc]{ $\zeta$}}
\put(40.2,18){\makebox(0,0)[cc]{ $\eta$}}
%\put(21,-3){\makebox(0,0)[cc]{ $O^{{}^{(2)}}zv$}}
%\put(22.5,-6){\makebox(0,0)[cc]{Рис. 1}}
\end{picture}}
\hfill\mbox{}
\\[8ex]
\indent
The projective atlas of trajectories of system (16.7) is constructed on Fig. 16.3.
\vspace{0.35ex}
First we established the behaviour of trajectories 
\vspace{0.35ex}
for the projectively singular system (16.8) on the projective circle $ \P\K (\xi, \theta),$ 
\vspace{0.35ex}
and secondly we used the maps of the projective circles 
$\P\K(x, y),\ \P\K(\xi, \theta),\ \P\K(\eta, \zeta).$
\vspace{0.5ex}

The behaviour of trajectories on the projective circle $ \P\K (\xi, \theta) $ 
\vspace{0.35ex}
for the differential system (16.8) 
is defined unambiguously. Indeed, two separatrixes of the saddle $O ^ {(1)} _ {\phantom1}$ 
\vspace{0.35ex}
are lying on the axis $O ^ {(1)} _ {\phantom1} \xi,$ 
\vspace{0.35ex}
and other two separatrixes are symmetric with respect to the axis $O ^ {(1)} _ {\phantom1} \xi $ 
and lying in the half-plane $\xi>0.$ 
\vspace{0.35ex}
The last statement follow from that 
the system (16.8) hasn't contact points on the axis $O ^ {(1)} _ {\phantom1} \theta,$ 
\vspace{0.5ex}
and equatorial contact trajectories of the differential system (16.8) are lying in the half-plane $ \xi <0.$ 
\vspace{0.75ex}

The direction of movement along trajectories of systems (16.7) -- (16.9) 
\vspace{0.35ex}
is defined by stability of the nodes  
$A^{(1)}_{\phantom1}$ and $B^{(1)}_{\phantom1}.$
\\[3.25ex]
\indent
{\bf 16.4. Trajectories of differential system} [1, p. 88]
\\[2.25ex]
\mbox{}\hfill        % (16.10)
$
\begin{array}{l}
\dfrac{dx}{dt}=x(x^2+y^2-1)-y(x^2+y^2+1)
\equiv X(x,y),
\\[4.25ex] 
\dfrac{dy}{dt}=y(x^2+y^2-1)+x(x^2+y^2+1)
\equiv Y(x,y).
\end{array}
$
\hfill (16.10)
\\[3ex]
\indent
Integral basis of system (16.10) on any domain from the set $D=\{(t,x,y)\colon x\ne 0\}$
is the first integrals
\\[2ex]
\mbox{}\hfill
$
F_1^{}\colon (t,x,y)\to\ 
(x^2+y^2)\exp \Bigl(4t-2\arctan \dfrac{y}{x}\Bigr)
$
\ for all 
$
(t,x,y) \in D
\hfill
$
\\[1.75ex]
and 
\\[1.75ex]
\mbox{}\hfill
$
F_2^{}\colon (t,x,y)\to\ 
(x^2+y^2-1)\exp \Bigl(2t-2\arctan \dfrac{y}{x}\Bigr)
$
\ for all 
$
(t,x,y) \in D.
\hfill
$
\\[2.25ex]
\indent
The general autonomous integral of system (16.10) on any domain from the set $G$ is
\\[2.25ex]
\mbox{}\hfill
$
F\colon (x,y)\to\ 
\dfrac{x^2+y^2}{(x^2+y^2-1)^2}\,\exp \Bigl(2\arctan \dfrac{y}{x}\Bigr)
$
\ for all 
$
(x,y) \in G,
\hfill
$
\\[2.25ex]
where 
$
G=\{(x,y)\colon x^2+y^2\ne 1\ \, \& \ \, x\ne 0\}.
$
\vspace{1ex}

The system of equations 
\\[1ex]
\mbox{}\hfill
$
X(x,y)=0,
\quad 
Y(x,y)=0
\hfill
$ 
\\[1.75ex]
has one solution $x=0, \ y=0.$ 
\vspace{0.35ex}
The differential system (16.10) has one equilibrium state $O(0,0)$ 
in the final part of the projective phase plane $\R\P(x, y).$
\vspace{0.35ex}
This equilibrium state is a stable focus with the characteristic equation $\lambda^2+2\lambda+2=0.$
\vspace{0.35ex}

The function 
\\[1.5ex]
\mbox{}\hfill
$
W\colon (x,y)\to\ 
x Y(x,y)-yX(x,y)=(x^2+y^2)(x^2+y^2+1)>0 
$
\ for all 
$
(x,y)\in\R^2\backslash\{(0,0)\}.
\hfill
$
\\[1.75ex]
\indent
The angle between the radius-vector 
\vspace{0.35ex}
of the organising point and the positive direction of the axis $Ox$ 
at movement along trajectories of system (16.10) is increasing.
\vspace{0.5ex}
 
The phase directional field of system (16.10) 
\vspace{0.25ex}
is symmetric with respect to the origin of coordinates of the phase plane $Oxy.$
\vspace{0.25ex}
For each trajectory of system (16.10) there exists a symmetric trajectory with respect to 
the origin of coordinates of the phase plane $Oxy. $ 
\vspace{0.5ex}

The equation $Y (x, 0) =0$ has one root $x=0, $ and $X(0,0)=0.$
\vspace{0.5ex}

The equation $X(0,y)=0$ has one root $y=0,$ and $Y(0,0)=0.$ 
\vspace{0.5ex}

Thus the system (16.10) hasn't contact points on the coordinate axes $Ox$ and $Oy.$ 
\vspace{0.35ex}

The function 
\\[1ex]
\mbox{}\hfill
$
W_3^{}\colon (x,y)\to\ 
(x^2+y^2)^2
$
\ for all 
$
(x,y)\in \R^2
\hfill
$ 
\\[2ex]
is not identically zero on the plane $\R^2.$ The system (16.10) is projectively nonsingular.
\vspace{0.5ex}

The projectively nonsingular differential system (16.10) 
\vspace{0.5ex}
by the first transformation of Poincar\'{e} 
$x =\dfrac{1}{\theta}\,, \ y =\dfrac{\xi}{\theta}$ 
is reduced to the first projectively reduced system 
\\[2ex]
\mbox{}\hfill        % (16.11)
$
\begin{array}{c}
\dfrac{d\xi}{d\tau}= (1+\xi^2)(1+\xi^2+\theta^2)\equiv \Xi(\xi,\theta),
\\[4ex] 
\dfrac{d\theta}{d\tau}={}-\theta(1+\xi^2-\theta^2-\xi(1+\xi^2+\theta^2)) \equiv \Theta(\xi,\theta),
\quad \
\theta^2\, d\tau=dt.
\end{array}
$
\hfill (16.11)
\\[2.5ex]
\indent
The projectively nonsingular differential system (16.10) 
\vspace{0.5ex}
by the second transformation of Poincar\'{e}
$x=\dfrac{\zeta}{\eta}\,,\ y=\dfrac{1}{\eta}$
is reduced to the second projectively reduced system
\\[2.25ex]
\mbox{}\hfill        % (16.12)
$
\begin{array}{c}
\dfrac{d\eta}{d\nu}= {}-\eta(1-\eta^2+\zeta^2+\zeta(1+\eta^2+\zeta^2))\equiv H(\eta,\zeta),
\\[4ex]  
\dfrac{d\zeta}{d\nu}={}-(1+\zeta^2)(1+\eta^2+\zeta^2) \equiv Z(\eta,\zeta),
\quad \
\eta^2\,d\nu=dt.
\end{array}
$
\hfill (16.12)
\\[2.5ex]
\indent
Since $\Xi(\xi,\theta)>0$ for all  $(\xi,\theta)\in \R^2,$ we see that
\vspace{0.5ex}
the system (16.11) hasn't equilibrium states  
in the final part of the projective phase plane $\R\P(\xi,\theta).$ 
\vspace{0.5ex}

Since $Z(\eta,\zeta)<0$ for all $(\eta,\zeta)\in \R^2,$ we see that
\vspace{0.5ex}
the system (16.12) hasn't  equilibrium states
in the final part of the projective phase plane $\R\P(\eta,\zeta).$ 
\vspace{0.5ex}

The differential system (16.10) has one equilibrium state $O(0,0)$
\vspace{0.35ex}
in the projective phase plane $\R\P(x,y).$ 
This equilibrium state is a stable focus.
\vspace{0.5ex}

In [1, p. 88], 
\vspace{0.35ex}
H. Poincar\'{e} proved that the circle $x^2+y^2=1$ is a limit cycle of the differential system (16.10).
\vspace{0.35ex}
From the equation $F (x, y) =C $ of sets of trajectories for system (16.10) it follows that 
\vspace{0.35ex}
the circle $x^2+y^2=1$ is unique limit cycle of the differential system (16.10) 
in the final part of the projective phase plane $\R\P(x,y).$
\vspace{0.5ex}

The system (16.10) 
\vspace{0.35ex}
is projectively nonsingular and hasn't equilibrium states on 
the infinitely removed straight line of the projective phase plane $\R\P(x,y).$  
\vspace{0.5ex}
The system (16.10) hasn't open limit cycles. 
Unique linear limit cycle of system (16.10) is the limit $\ell_\infty ^ {}\!$-cycle.
\vspace{0.5ex}

The differential system (16.11) 
\vspace{0.5ex}
in the projective phase plane $\R\P(\xi,\theta)$ has one equlibrium state $O,$ 
which is lying on <<extremities>> of the straight line $\xi=0.$
\vspace{0.25ex}
This equlibrium state $O$ is a stable focus.
\vspace{0.35ex}

The differential system (16.11) in the projective phase plane $\R\P(\xi,\theta)$ 
\vspace{0.5ex}
has two limit cycles:
the linear limit cycle $\theta=0$
and the open limit cycle $\theta^2-\xi^2=1.$
\vspace{0.5ex}

The phase directional field of system (16.11) 
\vspace{0.35ex}
is symmetric with respect to the axis $O^{(1)}_{\phantom1}\xi.$ 
For each trajectory of system (16.11) 
\vspace{0.35ex}
there exists a symmetric trajectory with respect to 
the coordinate axis $O^{(1)}_{\phantom1} \xi.$
\vspace{0.5ex}
The straight line-trajectory $\theta=0$ passes
through the origin of coordinates of the phase plane $O ^ {(1)} _ {\phantom1} \xi\theta$ 
and coincides with the coordinate axis $O^{(1)}_{\phantom1} \xi.$
\vspace{0.75ex}

The straight line $\theta=0$  is a trajectory of system (16.11). 
\vspace{0.5ex}
The equation $\Xi(0,\theta)=0$ hasn't solutions. 
The system (16.11) hasn't contact points on the coordinate axis $O^{(1)}_{\phantom1} \theta.$  
\vspace{0.5ex}

The straight line $x=0$ 
\vspace{0.25ex}
does not consist of trajectories of system (16.10). 
The differential system  (16.11) is projectively singular. 
\vspace{0.35ex}
The infinitely removed straight line of the projective phase plane $\R\P(\xi,\theta)$ 
does not consist of trajectories of system (16.11).
\vspace{0.5ex}

The differential system (16.10) hasn't contact points on the axis $Oy.$ 
\vspace{0.25ex}
The straight line $\theta=0$ is a trajectory of system (16.11). 
\vspace{0.25ex}
The projectively singular differential system (16.11) hasn't equatorial contact points.
\vspace{0.5ex}

The differential system (16.12) in the projective phase plane $\R\P(\eta,\zeta)$
\vspace{0.25ex}
has one equlibrium state $O,$ which is lying on <<extremities>> of  the straight line $\eta=0.$
\vspace{0.25ex}
And this equlibrium state $O$ is a stable focus.

The system (16.12) in the projective phase plane $\R\P(\eta,\zeta)$ 
\vspace{0.35ex}
has two limit cycles: 
the linear limit cycle $\eta=0$
and the open limit cycle $\eta^2-\zeta^2=1.$
\vspace{0.35ex}

The phase directional field of system (16.12) is symmetric with respect to 
\vspace{0.35ex}
the coordinate axis $O^{(2)}_{\phantom1} \zeta.$ 
\vspace{0.35ex}
For each trajectory of system (16.12) there exists a symmetric trajectory with respect to 
the coordinate axis $O^{(2)}_{\phantom1} \zeta.$
\vspace{0.35ex}
The straight line-trajectory $\eta=0$ passes 
through the origin of coordinates of the phase plane $O^{(2)}_{\phantom1} \eta\zeta$
\vspace{0.5ex}
and coincides with the coordinate axis $O^{(2)}_{\phantom1} \zeta.$

The equation $Z(\eta,0)=0$ hasn't solutions. 
\vspace{0.35ex}
The system (16.12) hasn't contact points on the coordinate axis $O^{(2)}_{\phantom1} \eta.$
The straight line $\eta=0$  is a trajectory of system (16.12). 
\vspace{0.5ex}

The straight line $y=0$ 
\vspace{0.25ex}
does not consist of trajectories of system (16.10). 
The differential system  (16.12) is projectively singular.
\vspace{0.25ex}
The infinitely removed straight line of the projective phase plane $\R\P(\eta,\zeta)$ 
does not consist of trajectories of system (16.12).
\vspace{0.5ex}

The differential system (16.10) hasn't contact points on the axis $Ox.$
\vspace{0.25ex}
The straight line $\eta=0$ is a trajectory of system (16.12). 
\vspace{0.25ex}
The projectively singular differential system (16.12) hasn't equatorial contact points.
\vspace{0.5ex}

The projective atlas of trajectories of system (16.10) is constructed on Fig. 16.4.
\\[3.75ex]
\mbox{}\hfill
{\unitlength=1mm
\begin{picture}(42,42)
\put(0,0){\includegraphics[width=42mm,height=42mm]{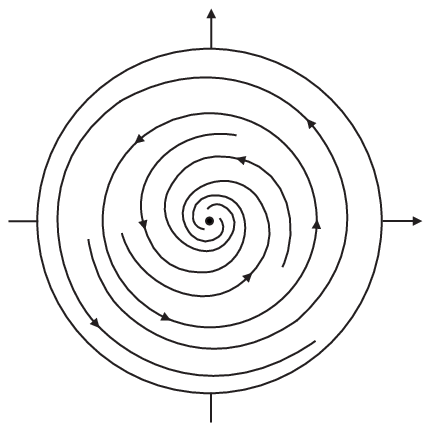}}
\put(18,41){\makebox(0,0)[cc]{ $y$}}
\put(40.2,18.2){\makebox(0,0)[cc]{ $x$}}
\end{picture}}
\qquad
{\unitlength=1mm
\begin{picture}(42,42)
\put(0,0){\includegraphics[width=42mm,height=42mm]{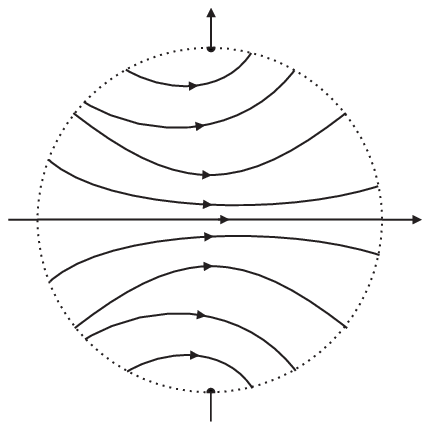}}
\put(18,41){\makebox(0,0)[cc]{ $\theta$}}
\put(40.2,17.8){\makebox(0,0)[cc]{ $\xi$}}
\put(21,-7){\makebox(0,0)[cc]{Fig. 16.4}}
\end{picture}}
\qquad
{\unitlength=1mm
\begin{picture}(42,42)
\put(0,0){\includegraphics[width=42mm,height=42mm]{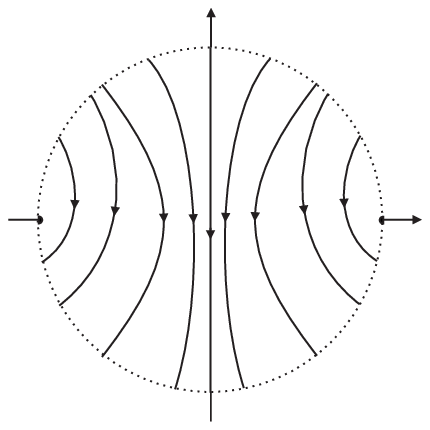}}
\put(18,41){\makebox(0,0)[cc]{ $\zeta$}}
\put(40.2,18){\makebox(0,0)[cc]{ $\eta$}}
\end{picture}}
\hfill\mbox{}
\\[9.75ex]
\indent
{\bf 16.5. Trajectories of differential system} [1, pp. 88 -- 90]
\\[2ex]
\mbox{}\hfill                                       %(16.13)
$
\begin{array}{l}
\dfrac{dx}{dt}=x(x^2+y^2-1)(x^2+y^2-9)-y(x^2+y^2-2x-8)
\equiv X(x,y),
\\[3.5ex]
\dfrac{dy}{dt}=y(x^2+y^2-1)(x^2+y^2-9)+x(x^2+y^2-2x-8)
\equiv Y(x,y).
\end{array}
$
\hfill (16.13)
\\[2.5ex]
\indent
From the system of equations 
\\[1.25ex]
\mbox{}\hfill
$
X(x,y)=0,
\quad 
Y(x,y)=0
\hfill
$ 
\\[1.75ex]
it follows that 
\vspace{0.5ex}
the system (16.13) in the final part of the projective phase plane $ \R\P (x, y)$  
has three equilibrium states $O(0,0),\ A\biggl(\dfrac12\,,\dfrac{\sqrt{35}}{2}\,\biggr),$ and
\vspace{0.75ex}
$B\biggl(\dfrac12\,,{}-\dfrac{\sqrt{35}}{2}\,\biggr).$
The equilibrium states $A$ and $B $ lie on the circle $x^2+y^2=9.$ 
\vspace{0.75ex}
This circle consists of trajectories of system (16.13).

The equilibrium state $O (0,0) $ has the characteristic equation $ \lambda^2-18\lambda+145=0.$ 
\vspace{0.35ex}
This equilibrium state is an unstable focus.
\vspace{0.5ex}

The equilibrium state  $A\biggl(\dfrac12\,,\dfrac{\sqrt{35}}{2}\,\biggr)$ has the characteristic equation 
\\[1.5ex]
\mbox{}\hfill
$
\lambda^2-\bigl( 144+\sqrt{35}\,\bigr)\lambda+25+144\;\!\sqrt{35}=0
\hfill
$
\\[1.25ex]
This equilibrium state is an unstable node.

Since the absolute term of characteristic equation for the 
\vspace{0.35ex}
equilibrium state $B\biggl(\dfrac12\,,{}-\dfrac{\sqrt{35}}{2}\,\biggr)$
is
$\Delta=25-144\;\!\sqrt{35}<0,$ we see that
the equilibrium state  $B\biggl(\dfrac12\,,{}-\dfrac{\sqrt{35}}{2}\,\biggr)$
\vspace{0.75ex}
is a saddle. 

If $(x-1)^2+y^2>9,$ then the function
\\[1.5ex]
\mbox{}\hfill
$
W\colon (x,y)\to\ 
x Y(x,y)-yX(x,y)=(x^2+y^2)(x^2+y^2-2x-8)
$
\ for all 
$
(x,y)\in \R^2
\hfill
$
\\[1.75ex]
is positive.  
If $(x-1)^2+y^2<9$ and $x^2+y^2\ne 0,$
then this function $W$ is negative.
\vspace{0.5ex}
 
The angle between 
\vspace{0.35ex}
radius-vector of the organising point and the positive direction of the axis $Ox$
\vspace{0.5ex}
at movement along the parts of trajectories of system (16.13) from the domains: 
\linebreak
a)  $\bigl\{ (x,y)\colon\, (x-1)^2+y^2<9 \ \& \ x^2+y^2\ne 0\bigr\};$
\vspace{0.75ex}
b)  $\bigl\{ (x,y)\colon\, (x-1)^2+y^2>9\bigr\}$
is a) decreasing; b) increasing, respectively.
\vspace{0.5ex}
Through each point, which is lying on the circle $(x-1) ^2+y^2=9, $ 
the trajectory of system (16.13) passes in the direction of the radius-vector of this point.
\vspace{0.5ex}

The equation $Y(x,0)=0$ has three roots $x=0,\ x={}-2,\ x=4$ such that 
\\[1.75ex]
\mbox{}\hfill
$
X(0,0)=0,\ \ 
X({}-2,0)=30\ne 0,\ \ 
X(4,0)=420\ne0.
\hfill
$
\\[1.75ex]
\indent
The differential system (16.13) has two contact points 
$C_1^{}({}-2,0)$ and $C_2^{}(4,0)$ on the axis $Ox.$
Since
\\[1.25ex]
\mbox{}\hfill
$
X({}-2,0)\;\!\partial_x^{} Y({}-2,0)=360>0,
\quad 
X(4,0)\;\! \partial_x^{} Y(4,0)=420\cdot 24>0,
\hfill
$
\\[1.5ex]
we see that 
\vspace{0.5ex}
the contact $C_1^{}\!$-  and $C_2 ^ {}\!$-trajectories of system (16.13) 
in enough small neighbourhoods of the points $C_1 ^ {} $ and $C_2 ^ {} $ lie in the half-plane $y\geq 0.$  
\vspace{0.75ex}

The equation $X(0,y)=0$ has three roots $y=0$ и $y={}\pm 2\;\!\sqrt{2}$ such that 
\\[1.75ex]
\mbox{}\hfill
$
Y(0,0)=0,
\quad
Y(0,{}-2\;\!\sqrt{2}\,)=14\;\!\sqrt{2}\,\ne 0,
\quad 
Y(0,2\;\!\sqrt{2}\,)={}-14\;\!\sqrt{2}\,\ne 0.
\hfill
$
\\[1.75ex]
\indent
The differential  system (16.13) has two contact points  
$D_1^{}(0,{}-2\;\!\sqrt{2}\,)$ and 
$D_2^{}(0,2\;\!\sqrt{2}\,)$  on the axis $Oy.$
\vspace{0.25ex}

The product 
\\[1.15ex]
\mbox{}\hfill
$
Y(0,{}-2\;\!\sqrt{2}\,)\;\! \partial_y^{} X(0,{}-2\;\!\sqrt{2}\,)=14\;\!\sqrt{2}\,\cdot ({}-16)<0.
\hfill
$
\\[1.75ex]
\indent
The contact $D_1^{}\!$-trajectory of system (16.13)
\vspace{0.5ex}
in enough small neighbourhoods of the point $D_1^{}$
lies in the half-plane  $x\leq 0.$
\vspace{0.5ex}
 
The product 
\\[1.15ex]
\mbox{}\hfill
$
Y(0,2\;\!\sqrt{2}\,)\;\! \partial_y^{} X(0,2\;\!\sqrt{2}\,)={}-14\;\!\sqrt{2}\,\cdot ({}-16)>0.
\hfill
$
\\[1.75ex]
\indent
The contact $D_2^{}\!$-trajectory of system (16.13)
\vspace{0.5ex}
in enough small neighbourhoods of the point $D_2^{}$
lies in the half-plane $x\geq 0.$
\vspace{0.35ex}

The function
\\[1ex]
\mbox{}\hfill
$
W_5^{}\colon (x,y)\to\
x\;\!Y_5^{}(x,y)-y\;\!X_5^{}(x,y)=0
$
\ for all 
$
(x,y)\in \R^2.
\hfill
$
\\[1.5ex]
\indent
The system  (16.13) is projectively singular. 
\vspace{0.25ex}
The infinitely removed straight line of the projective phase plane $ \R\P (x, y)$ 
\vspace{0.5ex}
does not consist of trajectories of system (16.13).

The first and the second projectively reduced systems for system (16.13) are the systems 
\\[2.15ex]
\mbox{}\hfill        % (16.14)
$
\dfrac{d\xi}{d\tau}= \theta\bigl(1+\xi^2\bigr)\bigl(1-2\theta +\xi^2-8\theta^2\bigr)\equiv \Xi(\xi,\theta),
\hfill
$
\\
\mbox{}\hfill (16.14)      
\\[0.25ex]
\mbox{}
$
\dfrac{d\theta}{d\tau}={}-1-2\xi^2+10\theta^2+\xi \theta^2-\xi^4+10\xi^2\;\!\theta^2-
2\xi \theta^3-9\theta^4+\xi\;\!\theta^2\bigl(\xi^2-8\theta^2\bigr)
 \equiv \Theta(\xi,\theta),
\quad
\theta^3\, d\tau=dt,
\hfill
$
\\[2ex]
and
\\[2ex]
\mbox{}\       % (16.15)
$
\dfrac{d\eta}{d\nu}= {}-1+10\eta^2-2\zeta^2-\eta^2\;\! \zeta-9\eta^4+10\eta^2\;\!\zeta^2
-\zeta^4+\eta^2\;\!\zeta\bigl( 8\eta^2+2\eta\;\!\zeta-\zeta^2\bigr)
\equiv H(\eta,\zeta),
\hfill
$
\\[0.75ex]
\mbox{}\hfill (16.15)      
\\[0.25ex]
\mbox{}\hfill
$
\dfrac{d\zeta}{d\nu}=\eta\bigl(1+\zeta^2\bigr)\bigl({}-1+8\eta^2+2\eta\zeta-\zeta^2\bigr) \equiv Z(\eta,\zeta),
\quad
\eta^3\,d\nu=dt.
\hfill 
$
\\[2.75ex]
\indent
Since
\\[1ex]
\mbox{}\hfill
$
\Theta(\xi,0)={}-1-2\xi^2-\xi^4\ne 0
$ 
\ for all  
$
\xi \in \R,
\hfill
$
\\[1.75ex]
we see that 
\vspace{0.35ex}
the system (16.14) hasn't equilibrium states 
on the coordinate axis $O^{(1)}_{\phantom1} \xi$
in the final part of the projective phase plane $\R\P(\xi,\theta).$ 
\vspace{0.35ex}

Since $H(0,0)={}-1\ne 0,$ we see that
\vspace{0.5ex}
the origin of coordinates of the phase plane $O^{(2)}_{\phantom1} \eta\zeta$
is not an equilibrium state of system (16.15).  
\vspace{0.5ex}
 
The differential system (16.13) has three equilibrium states 
\vspace{0.35ex}
in the projective phase plane $\R\P(x,y)\colon\!\!\!$ 
$O(0,0),\ A\biggl(\dfrac12\,,\dfrac{\sqrt{35}}{2}\,\biggr),$ and 
$B\biggl(\dfrac12\,,{}-\dfrac{\sqrt{35}}{2}\,\biggr).$
\vspace{0.5ex}

Since 
\\[1.5ex]
\mbox{}\hfill
$
X_5^{}(1,\xi)=(1+\xi^2)^2\ne 0 
$
\ for all 
$
\xi \in \R
$
\ \ and \ \ 
$
Y_5^{}(\zeta,1)= (1+\zeta^2)^2\ne 0 
$
\ for all 
$
\zeta \in \R,
\hfill
$
\\[1.75ex]
we see that 
the system (16.13) hasn't  equatorial contact points.
\vspace{0.5ex}

In [1, pp. 88 -- 90], H. Poincar\'{e} proved that 
\vspace{0.35ex}
the circle $x^2 + y^2=1$ is unique limit cycle of system (16.13) in the final part of the projective phase plane $\R\P(x,y).$
\vspace{0.5ex}

Straight lines are not trajectories of the projectively singular system (16.13) 
\vspace{0.35ex}
in the projective phase plane $\R\P (x, y).$  
The system (16.13) hasn't linear limit cycles.
\vspace{0.35ex}

The projectively singular system (16.13) hasn't 
\vspace{0.35ex}
infinitely removed equilibrium states and equatorial contact points.
The system (16.13) hasn't open limit cycles.
\vspace{0.35ex}

The differential system (16.14) in the projective phase plane $\R\P(\xi,\theta)$ 
\vspace{0.35ex}
has three equilibrium states:
\vspace{0.35ex}
the unstable node $A^{(1)}_{\phantom1} (\sqrt{35}\,, 2),$
the saddle  $B^{(1)}_{\phantom1} ({}-\sqrt{35}\,, 2),$
and the unstable  focus, which is lying on <<extremities>> of the straight line $\xi=0.$
\vspace{0.5ex}

From the equation $\Xi(\xi,\theta)=0,$
\vspace{0.5ex}
we get two equations $\theta=0$
and $8\theta^2+2\theta-\xi^2-1=0.$
The straight line $\theta=0$
\vspace{0.5ex}
and the hyperbola $8\theta^2+2\theta-\xi^2-1=0$
are orthogonal isoclines of system (16.14).
\vspace{0.35ex}
The tangent to a trajectory of system (16.14) at each point of these curves (straight line and hyperbola) 
is parallel to the axis  $O^{(1)}_{\phantom1} \theta.$
\vspace{0.5ex}

The straight line $\theta=0$  is an orthogonal isocline of system (16.14).
\vspace{0.35ex}
The system (16.14) hasn't contact points  
on the coordinate axis $O^{(1)}_{\phantom1} \xi.$ 
\vspace{0.5ex}
The orthogonal isoclines
$\theta=0$
and $8\theta^2+2\theta-\xi^2-1=0$
\vspace{0.5ex}
intersect the straight line $ \xi=0$ in three points: 
$O^{(1)}_{\phantom1} (0,0),\
C^{(1)}_{1} \Bigl(0,{}-\dfrac12\Bigr),$
and $C^{(1)}_{2} \Bigl(0,\dfrac14\Bigr).$
These points are contact points of the coordinate axis  $O^{(1)}_{\phantom1} \theta.$
\vspace{0.5ex}

The product 
\vspace{0.5ex}
$\Theta(0,0)\;\! \partial_\theta^{} \Xi(0,0)={}-1<0.$
The contact $O^{(1)}_{\phantom1}\!$-trajectory of system (16.14) 
in enough small neighbourhood of the point 
$O^{(1)}_{\phantom1}$ lies in the half-plane  $\xi\leq 0.$
\vspace{0.5ex}
Considering disposition of the contact $C_1 ^ {}\!$- and $C_2 ^ {}\!$-trajectories of system (16.13), 
\vspace{0.35ex}
we obtain: a) the contact $C_1 ^ {(1)}\!$-trajectory of system (16.14) 
\vspace{0.5ex}
in enough small neighbourhoods of the point $C ^ {(1)} _ {1}$ lies 
in the half-plane $ \xi\leq 0;$ 
\vspace{0.5ex}
b) the contact $C_2 ^ {(1)}\!$-trajectory of system (16.14) 
in enough small neighbourhood of the point $C ^ {(1)} _ {2}$ lies in the half-plane  $\xi\geq 0.$
\vspace{0.5ex}

The straight line $x=0$ does not consist of trajectories of system (16.13). 
\vspace{0.35ex}
The system (16.14) is projectively singular. 
\vspace{0.35ex}
The infinitely removed straight line of the projective phase plane $ \R\P (\xi, \theta)$ 
does not consist of trajectories of system (16.14).

A focus lies on <<extremities>> of the straight line $ \xi=0.$ 
\vspace{0.35ex}
The infinity removed point on <<extremities>> of the axis $O ^ {(1)} _ {\phantom1}\theta$ 
is not an equatorial contact point of system (16.14). 
\vspace{0.75ex}

The equation 
\\[1ex]
\mbox{}\hfill
$
\Xi_5^{}(1,a)=0,
\hfill
$ 
\\[1ex]
where 
\\[1ex]
\mbox{}\hfill
$
\Xi_5^{}(\xi,\theta)=\xi^2\;\!\theta(\xi^2-8\theta^2)
$
\ for all 
$
(\xi,\theta)\in\R^2,
\hfill
$
\\[1.75ex]
has three roots: $a=0$ and $a={}\pm \dfrac{\sqrt{2}}{4}\,,$
\vspace{0.5ex}
such that $W^{(1)}_4(1,a)=a\, \Xi_4^{}(1,a)- \Theta_4^{}(1,a)\ne 0.$

The system (16.14) has three equatorial contact points: 
\vspace{0.5ex}
$O^{(2)}_{\phantom1},\ D^{(1)}_{1},$ and $D^{(1)}_{2},$
which are lying on <<extremities>> of the straight lines 
$\theta=0,\ \theta={}-\dfrac{\sqrt{2}}{4}\ \xi,$ and 
\vspace{0.5ex}
$\theta=\dfrac{\sqrt{2}}{4}\ \xi,$ respectively.

Considering disposition of contact trajectories of the axis $Oy$ for system (16.13), 
\vspace{0.35ex}
we get equatorial contact trajectories of system (16.14) in enough small neighbourhood 
\vspace{0.35ex}
of the in\-fi\-ni\-te\-ly removed straight line of the projective phase plane $\!\R\P (\xi, \theta)\!$ 
\vspace{0.5ex}
lie in the half-plane $\!\xi\!>\! 0.$

The differential system (16.14) hasn't linear limit cycles 
\vspace{0.35ex}
and limit cycles, which are lying 
in the final part of the projective phase plane $ \R\P (\xi, \theta).$
\vspace{0.35ex}
The hyperbola $ \theta^2-\xi^2=1$ is an open limit cycle of the differential system (16.14).
\vspace{0.5ex}

The differential system (16.15) in the projective phase plane $\R\P(\eta,\zeta)$ 
\vspace{0.5ex}
has three equatorial contact points: 
\vspace{0.5ex}
the unstable focus, which is lying on <<extremities>> of the straight line $\zeta=0,$
the unstable node $A^{(2)}_{\phantom1}\biggl(\dfrac{2\sqrt{35}}{35}\,,\dfrac{\sqrt{35}}{35}\,\biggl),$ and 
the saddle 
$B^{(2)}_{\phantom1}\biggl({}-\dfrac{2\sqrt{35}}{35}\,,{}-\dfrac{\sqrt{35}}{35}\,\biggl).$  
\vspace{1ex}

From the equation $Z(\eta,\zeta)=0,$
\vspace{0.5ex}
we get two equations $\eta=0$ and $8\eta^2+2\eta\zeta-\zeta^2-1=0.$
The straight line $\eta=0$
and the hyperbola $8\eta^2+2\eta\zeta-\zeta^2-1=0$
\vspace{0.5ex}
are zero isoclines of system (16.15).
The tangent to a trajectory of system (16.15) 
\vspace{0.35ex}
at each point of these curves (straight line and hyperbola) 
are parallel to the axis $O^{(2)}_{\phantom1} \eta.$
\vspace{0.5ex}

The straight line $\eta=0$  is a zero isocline of system (16.15). 
\vspace{0.35ex}
The system (16.15) hasn't contact points on the axis $O^{(2)}_{\phantom1} \zeta.$ 
\vspace{0.75ex}

The zero isoclines
\vspace{0.75ex}
$\eta=0$ and $8\eta^2+2\eta\zeta-\zeta^2-1=0$
intersect the straight line $\zeta=0$  in three points: 
\vspace{0.35ex}
$O^{(2)}_{\phantom1} (0,0),\ D^{(2)}_{1} \Bigl({}-\dfrac{\sqrt{2}}{4}\,,0\Bigr),$
and $D^{(2)}_{2} \Bigl(\dfrac{\sqrt{2}}{4}\,,0\Bigr).$
These points are contact points of the axis  $O^{(2)}_{\phantom1} \eta.$
\vspace{0.75ex}

The product 
\vspace{0.35ex}
$H(0,0)\;\! \partial_\eta^{} Z(0,0)=1>0.$
The contact $O^{(2)}_{\phantom1}\!$-trajectory of system (16.13) in enough small neighbourhood of the point  
$O^{(2)}_{\phantom1}$ lies in the half-plane $\zeta\geq 0.$
\vspace{0.5ex}

Considering disposition of the contact
\vspace{0.35ex}
$D_1^{}\!$- and $D_2^{}\!$-trajectories of system (16.13), we obtain 
the contact $D_1^{(2)}\!$- and  $D_2^{(2)}\!$-trajectory of system (16.15) 
\vspace{0.5ex}
in enough small neighbourhoods of the points 
$D_1^{(2)}$ and  $D_2^{(2)}$ lies in the half-plane $\zeta\geq 0.$
\vspace{0.5ex}

The straight line $y=0$ does not consist of trajectories of system (16.13). 
\vspace{0.35ex}
The differential system (16.15) is projectively singular.  
\vspace{0.35ex}
The infinitely removed straight line of the projective phase plane $\R\P(\eta,\zeta)$ 
does not consist of trajectories of system (16.15).
\vspace{0.5ex}

A focus lies on <<extremities>> of the straight line  $\zeta=0.$ 
\vspace{0.35ex}
The infinitely removed point on <<extremities>> of the axis $O^{(2)}_{\phantom1} \eta$
is not an equatorial contact point of system (16.15). 
\vspace{0.5ex}

The equation 
\\[1ex]
\mbox{}\hfill
$
Z_5^{}(a,1)=0,
\hfill
$ 
\\[1ex]
where 
\\[1ex]
\mbox{}\hfill
$
Z_5^{}(\eta,\zeta)=\eta \zeta^2\;\!(8\eta^2+2\eta\zeta-\zeta^2)
$ 
\ for all 
$
(\eta,\zeta)\in\R^2,
\hfill
$
\\[2.5ex]
has three roots: $a=0,\, a={}-\dfrac{1}{2}\,,$
\vspace{0.75ex}
and $a=\dfrac{1}{4}\,,$
such that $W^{(2)}_4(a,1)\!=\!H_4^{}(a,1)-aZ_4^{}(a,1)\!\ne\! 0.$

The system (16.15) has three equatorial contact points: 
\vspace{0.5ex}
$O^{(1)}_{\phantom1},\ C^{(2)}_{1},$ and $C^{(2)}_{2},$
which are lying on <<extremities>> of the straight lines $\eta=0,\ \zeta={}-2\eta,$ and 
$\zeta=4\eta,$ respectively.
\vspace{0.5ex}

Considering disposition of the contact trajectories of the axis $Ox$ for system (16.13),
\vspace{0.35ex}
we get the equatorial contact $O^{(1)}_{\phantom1}\!$-trajectory of system (16.15) 
\vspace{0.35ex}
in enough small neighbourhood of the infinitely removed point
$O^{(1)}_{\phantom1}$ lies in the half-plane  $\zeta< 0,$
\vspace{0.35ex}
and the equatorial contact $C^{(2)}_{1}\!$- and  $C^{(2)}_{2}\!$-trajectory of system (16.15) 
\vspace{0.35ex}
in enough small neighbourhood of the infinitely removed straight line 
of the projective phase plane $\R\P(\eta,\zeta)$ 
lie in the half-plane $\eta> 0.$
\vspace{0.5ex}

The differential system (16.15) hasn't linear limit cycles and limit cycles, 
\vspace{0.35ex}
which are lying in the final part of the projective phase plane $\R\P(\eta,\zeta).$
\vspace{0.35ex}
The hyperbola $\eta^2-\zeta^2=1$ is an open limit cycle of the differential system (16.15).
\vspace{0.5ex}

The straight line $\theta=0$ is an orthogonal isocline of system (16.14), 
\vspace{0.35ex}
and the straight line $\eta=0$ is a zero isocline of system (16.15).
\vspace{0.35ex}
Each trajectory of system (16.13), which is passing through the infinitely removed straight line 
\vspace{0.35ex}
of the projective phase plane $\R\P (x, y),$ 
is orthogonal to the boundary circle of the projective circle $\R\K(x,y).$
\vspace{0.5ex}

The projective atlas of trajectories of system (16.13) is constructed on Fig. 16.5.
\\[3.75ex]
\mbox{}\hfill
{\unitlength=1mm
\begin{picture}(42,42)
\put(0,0){\includegraphics[width=42mm,height=42mm]{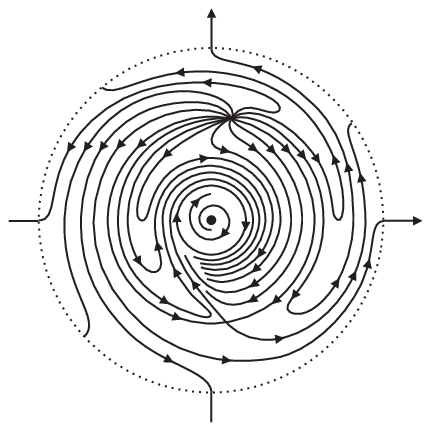}}
\put(18,41){\makebox(0,0)[cc]{ $y$}}
\put(40.2,18.2){\makebox(0,0)[cc]{ $x$}}
\end{picture}}
\qquad
{\unitlength=1mm
\begin{picture}(42,42)
\put(0,0){\includegraphics[width=42mm,height=42mm]{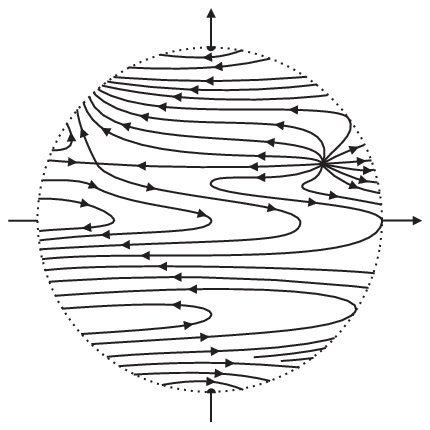}}
\put(18,41){\makebox(0,0)[cc]{ $\theta$}}
\put(40.2,17.8){\makebox(0,0)[cc]{ $\xi$}}
\put(21,-7){\makebox(0,0)[cc]{Fig. 16.5}}
\end{picture}}
\qquad
{\unitlength=1mm
\begin{picture}(42,42)
\put(0,0){\includegraphics[width=42mm,height=42mm]{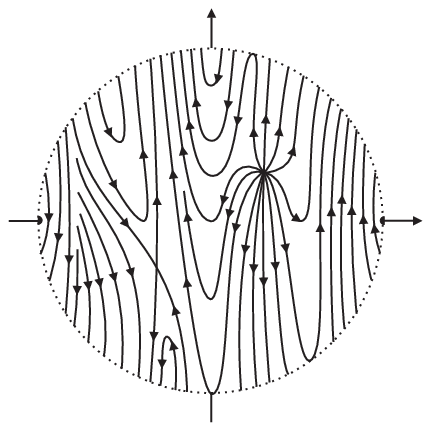}}
\put(18,41){\makebox(0,0)[cc]{ $\zeta$}}
\put(40.2,18){\makebox(0,0)[cc]{ $\eta$}}
\end{picture}}
\hfill\mbox{}
\\[10ex]
\indent
{\bf 16.6. Trajectories of differential system}  [1, pp. 90 -- 91]
\\[2ex]
\mbox{}\hfill                                       %(16.16)
$
\!\!\!
\begin{array}{l}
\dfrac{dx}{dt}=x(2x^2+2y^2+1)\Bigl((x^2+y^2)^2+x^2-y^2+\dfrac{1}{10}\Bigr)-y(2x^2+2y^2-1)
\equiv X(x,y),
\\[4.25ex]
\dfrac{dy}{dt}=y(2x^2+2y^2-1)\Bigl((x^2+y^2)^2+x^2-y^2+\dfrac{1}{10}\Bigr)+x(2x^2+2y^2+1)
\equiv Y(x,y).
\end{array}
$
\hfill (16.16)
\\[3ex]
\indent
The differential system (16.16) has three equilibrium states
\vspace{0.35ex}
in the final part of the projective phase plane $\R\P(x,y)\colon\!\!\!$ 
$O(0,0),\ A\biggl(0,\,\dfrac{\sqrt{2}}{2}\,\biggr),$ and
$B\biggl(0,{}-\dfrac{\sqrt{2}}{2}\,\biggr).$
\vspace{0.75ex}

Since the absolute term 
\vspace{0.5ex}
of the characteristic equation for the equilibrium state  $O(0,0)$ is  
$\Delta(0,0)={}-0{,}99<0,$
we see that the equilibrium state $O(0,0)$ is a saddle. 
\vspace{0.75ex}

The equilibrium states  
\vspace{0.75ex}
$A\biggl(0,\,\dfrac{\sqrt{2}}{2}\,\biggr)$ and
$B\biggl(0,{}-\dfrac{\sqrt{2}}{2}\,\biggr)$
are stable foci with the characteristic equation 
$\lambda^2+\dfrac{3}{5}\,\lambda+\dfrac{409}{100}=0.$

\newpage

The phase directional field of system (16.16) 
\vspace{0.25ex}
is symmetric with respect to the origin of coordinates 
of the phase plane $Oxy.$
\vspace{0.25ex}
For each trajectory of system (16.16) there exists a symmetric trajectory with respect to 
\vspace{0.75ex}
the origin of coordinates of the phase plane $Oxy.$ 

The equation $Y(x,0)=0$ hasn't  roots.  
The system (16.16) on the axis $Ox$ hasn't contact points. 
The equation $X(0,y)=0$ has three roots $y=0$ and $y={}\pm \dfrac{\sqrt{2}}{2}\,,$ but 
$Y(0,0)=0$ and $Y\biggl(0, {}\pm \dfrac{\sqrt{2}}{2}\,\biggr)=0.$
The system (16.16) on the axis $Oy$ hasn't contact points. 
\vspace{1ex}

The function 
\\[1.25ex]
\mbox{}\hfill
$
W_7^{}\colon (x,y)\to\ 
x\;\!Y_7^{}(x,y)-y\;\!X_7^{}(x,y)=0
$
\ for all 
$
(x,y)\in \R^2.
\hfill
$
\\[1.75ex]
\indent
The system  (16.16) is projectively singular. 
\vspace{0.35ex}
The infinitely removed straight line of the projective phase plane $\R\P(x,y)$
does not consist of trajectories of system (16.16).
\vspace{0.5ex}

The first and the second projectively reduced systems for system (16.16) are the systems  
\\[3ex]
\mbox{}\       % (16.17)
$
\dfrac{d\xi}{d\tau}= \theta\Bigl({}-2\xi +2\theta^2 -4\xi^3 -2\xi\theta^2+ 4\xi^2\theta^2+\theta^4 -2\xi^5+2\xi^3\theta^2-
\dfrac{1}{5}\, \xi\theta^4+\xi^2\theta^2\bigl(2\xi^2-\theta^2\bigr)\Bigr)\equiv \Xi(\xi,\theta),
\hfill
$
\\[1.25ex]
\mbox{}\hfill (16.17)      
\\[0.25ex]
\mbox{}\hfill
$
\dfrac{d\theta}{d\tau}={}-2-6\xi^2-3\theta^2-6\xi^4-2\xi^2\;\!\theta^2-
\dfrac{6}{5}\, \theta^4+2\xi\theta^4\ - 
\hfill
$
\\[2.5ex]
\mbox{}\hfill
$
-\ 2\xi^6+\xi^4\theta^2+
\dfrac{4}{5}\,\xi^2\theta^4-\dfrac{1}{10}\, \theta^6+\xi\theta^4\bigl(2\xi^2-\theta^2\bigr)
 \equiv \Theta(\xi,\theta),
\qquad 
\theta^5\, d\tau=dt,
\hfill
$
\\[1.5ex]
and
\\[1.5ex]
\mbox{}\hfill       % (16.18)
$
\dfrac{d\eta}{d\nu}= 
{}-2+3\eta^2-6\zeta^2-\dfrac{6}{5}\,\eta^4+2\eta^2\;\!\zeta^2
-6\zeta^4-2\eta^4\;\!\zeta\ +
\hfill
$
\\[2.5ex]
\mbox{}\hfill
$
+\ \dfrac{1}{10}\,\eta^6+ 
\dfrac{4}{5}\,\eta^4\;\!\zeta^2-
\eta^2\;\!\zeta^4-2\zeta^6-\eta^4\;\!\zeta\bigl(\eta^2+2\zeta^2\bigr)
\equiv H(\eta,\zeta),
\hfill
$
\\[0.5ex]
\mbox{}\hfill (16.18)      
\\[0.5ex]
\mbox{}\hfill
$
\dfrac{d\zeta}{d\nu}=\eta\;\!\Bigl(2\zeta-2\eta^2-2\eta^2\zeta+4\zeta^3+\eta^4
-4\eta^2\;\!\zeta^2\ +
\hfill
$
\\[2.5ex]
\mbox{}\hfill
$
+\ \dfrac{1}{5}\,\eta^4\;\!\zeta+ 
2\eta^2\;\!\zeta^3+2\zeta^5-\eta^2\;\!\zeta^2\bigl(\eta^2+2\zeta^2\bigr)\Bigr)
\equiv Z(\eta,\zeta),
\qquad
\eta^5\,d\nu=dt.
\hfill 
$
\\[2.5ex]
\indent
Since 
\\[0.75ex]
\mbox{}\hfill
$
\Theta(\xi,0)={}-2-6\xi^2-6\xi^4-2\xi^6\ne 0  
$
\ for all  
$
\xi \in \R,
\hfill
$
\\[1.75ex]
we see that
\vspace{0.35ex}
the differential system (16.17) hasn't equilibrium states
on the coordinate axis $O^{(1)}_{\phantom1} \xi$
in the final part of the projective phase plane $\R\P(\xi,\theta).$ 
\vspace{0.75ex}

From $H(0,0)={}-2\ne 0$ it follows that 
\vspace{0.35ex}
the origin of coordinates of the phase plane $O^{(2)}_{\phantom1} \eta\zeta$
is not an equilibrium state of system (16.18).
\vspace{0.75ex}

The differential system (16.16) has three equilibrium states
\vspace{0.35ex}
in the projective phase plane $\R\P(x,y)\colon\!\!\!$ 
$O(0,0),\ A\biggl(0,\,\dfrac{\sqrt{2}}{2}\,\biggr),$ and $B\biggl(0,{}-\dfrac{\sqrt{2}}{2}\,\biggr).$
\vspace{0.5ex}

Since
\\[1.25ex]
\mbox{}\hfill
$
X_7^{}(1,\xi)=2(1+\xi^2)^3\ne 0 
$
\ for all 
$
\xi \in \R
$
\ \ and \ \ 
$
Y_7^{}(\zeta,1)=2(1+\zeta^2)^3\ne 0 
$
\ for all 
$
\zeta \in \R,
\hfill
$
\\[1.5ex]
we see that the projective singular system (16.16) hasn't  equatorial contact points.
\vspace{0.5ex}

In [1, pp. 90 -- 91], H. Poincar\'{e} proved that 
\vspace{0.35ex}
limit cycles of the differential system (16.16) in 
the finite part of the projective phase plane $\R\P(x,y)$
are two not intersected closed curves given by the equation 
\\[1ex]
\mbox{}\hfill
$
(x^2+y^2)^2+x^2-y^2+\dfrac{1}{10}=0.
\hfill
$ 
\\[1.5ex]
\indent
Straight lines are not trajectories of the projectively singular system  (16.16) 
\vspace{0.35ex}
in the projective phase plane $\R\P(x,y).$ 
The system (16.16) hasn't linear limit cycles.
\vspace{0.5ex}

The projectively singular system  (16.16) 
\vspace{0.35ex}
hasn't infinitely removed equilibrium states and equatorial contact points.
The system (16.16) hasn't open limit cycles.
\vspace{0.5ex}

The differential system (16.17) 
\vspace{0.5ex}
has three equilibrium states in the projective phase plane $\R\P(\xi,\theta)\colon$
the stable foci $A^{(1)}_{\phantom1}$ and $B^{(1)}_{\phantom1},$ 
\vspace{0.5ex}
which are lying on <<extremities>> of the straight lines $\xi=\dfrac{\sqrt2}{2}\,\theta$ and 
$\xi={}-\dfrac{\sqrt2}{2}\,\theta,$ respectively; 
\vspace{1ex}
the saddle, which is lying on <<extremities>> of the straight line $\xi=0.$
\vspace{0.5ex}

The differential system (16.17) 
\vspace{0.35ex}
hasn't linear limit cycles and limit cycles, 
which are lying in the final part of the projective phase plane $ \R\P (\xi, \theta).$
\vspace{0.5ex}

The curves given by the equation 
\\[1.5ex]
\mbox{}\hfill
$
10(\xi^2+1)^2-10\theta^2(\xi^2-1)^2+\theta^4=0
\hfill
$ 
\\[1.25ex]
are two open limit cycles of the differential system (16.17).
\vspace{0.35ex}

The phase directional field of system (16.17) is symmetric with respect to the axis $O^{(1)}_{\phantom1} \xi.$
For each trajectory of system (16.17) 
\vspace{0.35ex}
there exists a symmetric trajectory with respect to the axis $O ^ {(1)} _ {\phantom1} \xi. $
\vspace{0.35ex}
Each trajectory of system (16.17), which is intersecting the axis $O ^ {(1)} _ {\phantom1} \xi,$ 
is symmetric with respect to this coordinate axis. 
\vspace{0.35ex}

The straight line $\theta=0$  is an orthogonal isocline of system (16.17).
\vspace{0.35ex}
The system (16.7) hasn't contact points on the axis $O^{(1)}_{\phantom1} \xi.$ 
\vspace{0.5ex}

The equation $\Xi(0,\theta)=0$ has one root $\theta=0,$ and $\Theta(0,0)={}-2\ne0.$
\vspace{0.5ex}

The differential system (16.17) on the axis $O^{(1)}_{\phantom1} \theta$  
\vspace{0.35ex}
has one contact point $O^{(1)}_{\phantom1} (0,0).$
The contact $O^{(1)}_{\phantom1}\!$-trajectory of system (16.17)
\vspace{0.5ex}
in enough small neighbourhood of the point $O^{(1)}_{\phantom1}(0,0)$ 
lies in the half-plane $\theta\geq 0.$
\vspace{0.75ex}

The straight line $x=0$ does not consist of trajectories of system (16.16). 
\vspace{0.35ex}
The differential system  (16.17) is projectively singular. 
\vspace{0.35ex}
The infinitely removed straight line of the projective phase plane $\R\P(\xi,\theta)$ 
does not consist of trajectories of system (16.17).
\vspace{0.5ex}

The differential system (16.16) on the axis  $Oy$ hasn't contact points. 
\vspace{0.5ex}

The differential system (16.17) 
\vspace{0.35ex}
hasn't equatorial contact points on <<extremities>> of  the straight lines $\xi=a\;\!\theta$ 
at any real coefficient $a.$ 
\vspace{0.35ex}
Since $\Xi_7^{}(1,0)=0,\
W_6^{(1)}(1,0)={}-2\ne 0,$ we see that: 
a) equatorial contact point of systems (16.17) 
\vspace{0.35ex}
lies on <<extremities>> of the axis $O ^ {(1)} _ {\phantom1} \xi;$  
b) equatorial contact trajectory 
\vspace{0.5ex}
in enough small neighbourhood of the infinitely removed straight line 
of the projective phase plane $ \R\P (\xi, \theta)$ lies in the half-plane $\xi<0.$
\vspace{0.5ex}

The differential system (16.18) 
\vspace{0.35ex}
in the projective phase plane $\R\P(\eta,\zeta)$ has three equilibrium states:
the stable foci $A^{(2)}_{\phantom1}(\sqrt2\,,0)$ and 
$B^{(2)}_{\phantom1}({}-\sqrt2\,,0);$ 
\vspace{0.35ex}
the saddle, which is lying on <<extremities>> of the straight line $\zeta=0.$
\vspace{0.5ex}

The system (16.18) hasn't linear limit cycles and limit cycles, 
\vspace{0.35ex}
which are lying in the final part of the projective phase plane $\R\P(\eta,\zeta).$
\vspace{0.35ex}

The curves given by the equation 
\\[1.25ex]
\mbox{}\hfill
$
10(\zeta^2+1)^2+10\eta^2(\zeta^2-1)^2+\eta^4=0
\hfill
$ 
\\[1.25ex]
are two open limit cycles of the differential system (16.18).
\vspace{0.35ex}

The phase directional field of system (16.18) 
\vspace{0.35ex}
is symmetric with respect to the axis $O^{(2)}_{\phantom1} \zeta.$
For each trajectory of system (16.18) there exists 
\vspace{0.35ex}
a symmetric trajectory with respect to the axis $O^{(2)}_{\phantom1} \zeta.$
\vspace{0.35ex}
Each trajectory of system (16.18), which is intersecting the axis $O^{(2)}_{\phantom1} \zeta,$ 
is symmetric with respect to this coordinate axis. 
\vspace{0.35ex}

The straight line $\eta=0$  is a zero isocline of system (16.18). 
\vspace{0.25ex}
The system (16.18) hasn't contact points on the axis $O^{(2)}_{\phantom1} \zeta.$
\vspace{0.5ex}

The equation $Z(\eta,0)=0$ has three roots $\eta=0,\ \eta={}\pm \sqrt2\,,$ and 
\\[1.25ex]
\mbox{}\hfill
$
H(0,0)={}-2,
\quad
H({}\pm \sqrt2\,,0)=0.
\hfill
$ 
\\[1.5ex]
\indent
The differential system (16.18) has one contact point $O^{(2)}_{\phantom1} (0,0)$
\vspace{0.35ex}
on the axis $O^{(2)}_{\phantom1} \eta.$
The contact $O^{(2)}_{\phantom1}\!$-trajectory of system (16.18)
\vspace{0.35ex}
in enough small neighbourhood of the point $O^{(2)}_{\phantom1}(0,0)$ 
lies in the half-plane $\zeta\leq 0.$
\vspace{0.5ex}

The straight line $y=0$ does not consist of trajectories of system (16.16). 
\vspace{0.35ex}
The differential system  (16.18) is projectively singular. 
\vspace{0.35ex}
The infinitely removed straight line of the projective phase plane $\R\P(\eta,\zeta)$ 
does not consist of trajectories of system (16.18).
\vspace{0.5ex}

The differential system (16.16) on the axis $Ox$ hasn't contact points.
\vspace{0.35ex}

The differential system (16.18) hasn't equatorial contact points on <<extremities>> 
\vspace{0.35ex}
of the straight lines $\eta=a\zeta$ at any real coefficient  $a.$
\vspace{0.75ex}

The point $O_{\phantom1}^{(1)}$ 
\vspace{0.35ex}
is a contact point of the straight line $\xi=0$ of system (16.17) 
and 
the contact $O _ {\phantom1} ^ {(1)}\!$-trajectory
\vspace{0.25ex}
in enough small neighbourhood of the point $O _ {\phantom1} ^ {(1)}$ 
lies in the half-plane $\xi\geq 0.$ 
\vspace{0.35ex}
The equatorial contact point $O_{\phantom1}^{(1)}$ of system (16.8) lies 
\vspace{0.35ex}
on <<extremities>> of the axis $O_{\phantom1}^{(2)}\zeta$
and the equatorial contact $O _ {\phantom1} ^ {(1)}\!$-trajectory of system (16.18) 
\vspace{0.5ex}
in enough small neighbourhood of infinitely removed straight line of the projective phase plane 
\vspace{0.35ex}
$\R\P (\eta, \zeta) $ lies in the half-plane $\zeta>0.$
\vspace{0.5ex}

The straight line $\theta=0$ is an orthogonal isocline of system (16.17) and 
\vspace{0.35ex}
the straight line $\eta=0$ is a zero isocline of system (16.18).
\vspace{0.35ex}
Each trajectory of system (16.16), which is passing through the infinitely removed straight line 
\vspace{0.35ex}
of the projective phase plane $ \R\P (x, y),$ 
is orthogonal to the boundary circle of the projective circle $\R\K(x,y).$
\vspace{0.5ex}

The projective atlas of trajectories of system (16.16) is constructed on Fig. 16.6.
\\[3.75ex]
\mbox{}\hfill
{\unitlength=1mm
\begin{picture}(42,42)
\put(0,0){\includegraphics[width=42mm,height=42mm]{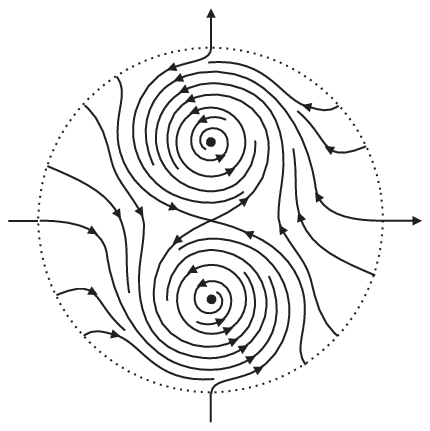}}
\put(18,41){\makebox(0,0)[cc]{ $y$}}
\put(40.2,18.2){\makebox(0,0)[cc]{ $x$}}
\end{picture}}
\qquad
{\unitlength=1mm
\begin{picture}(42,42)
\put(0,0){\includegraphics[width=42mm,height=42mm]{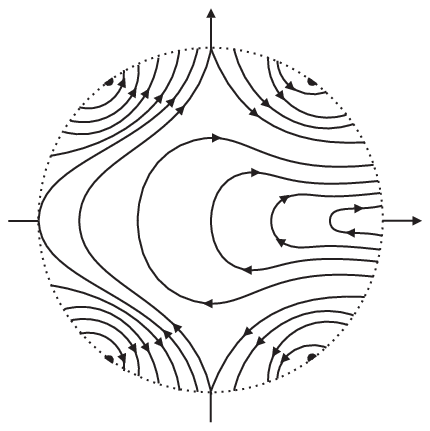}}
\put(18,41){\makebox(0,0)[cc]{ $\theta$}}
\put(40.2,17.8){\makebox(0,0)[cc]{ $\xi$}}
\put(21,-7){\makebox(0,0)[cc]{Fig. 16.6}}
\end{picture}}
\qquad
{\unitlength=1mm
\begin{picture}(42,42)
\put(0,0){\includegraphics[width=42mm,height=42mm]{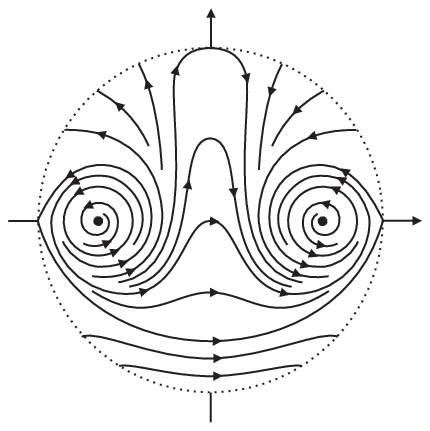}}
\put(18,41){\makebox(0,0)[cc]{ $\zeta$}}
\put(40.2,18){\makebox(0,0)[cc]{ $\eta$}}
\end{picture}}
\hfill\mbox{}
\\[11.75ex]
\indent
{\sl Remark}.
Note that the main results of this paper were originally published by the author in 
the articles [5; 22 -- 24]. 

\newpage

\mbox{}
\\[-3.5ex]

}

\newpage

{\normalsize

\Russian
\pagenumbering{arabic}

\sloppy

\lhead
    [\scriptsize В.Н. Горбузов]
    {\scriptsize В.Н. Горбузов}
\rhead
    [\it \scriptsize Проективный атлас траекторий дифференциальных систем]
    {\it \scriptsize Проективный атлас траекторий дифференциальных систем}

\thispagestyle{empty}

\mbox{}
\\[-3.75ex]
\centerline{
{\large
\bf
ПРОЕКТИВНЫЙ\;\! АТЛАС\;\! ТРАЕКТОРИЙ 
}
}
\\[0.75ex]
\centerline{
{\large
\bf
ДИФФЕРЕНЦИАЛЬНЫХ СИСТЕМ 
}
}
\\[2.25ex]
\centerline{
\bf 
В.Н. Горбузов
}
\\[2ex]
\centerline{
\it 
Гродненский государственный университет имени Янки Купалы
}
\\[0.5ex]
\centerline{
\it 
{\rm(}Ожешко 22, Гродно, Беларусь, 230023\;\!{\rm)}
}
\\[1.5ex]
\centerline{
E-mail: gorbuzov@grsu.by
}
\\[4.25ex]
\centerline{{\large\bf Резюме}}
\\[1ex]
\indent
Изложены топологические основы поведения траекторий автономных дифференциальных систем второго 
порядка на проективной фазовой плоскости. С помощью кругов Пуанкаре построен проективный атлас 
траекторий. Установлены дифференциальные связи между траекториями проективно сопряженных 
дифференциальных систем. Исследовано поведение траекторий в окрестности бесконечно удаленной 
прямой проективной фазовой плоскости и свойства замкнутых траекторий на проективной фазовой плоскости. 
Приведены примеры полного качественного исследования траекторий дифференциальных систем 
на проективной фазовой плоскости.
\\[1.5ex]
\indent
{\it Ключевые слова}:
дифференциальная система, сфера Пуанкаре, круг Пуанкаре, предельный цикл, проективная плоскость, 
атлас карт многообразия.
\\[1.25ex]
\indent
{\it 2000 Mathematics Subject Classification}: 34A26, 34C05.
\\[4.25ex]
\centerline{{\large\bf Содержание}}
\\[1.25ex]
{\bf  Введение}                   \dotfill\ 2
\\[1ex]
{\bf \S 1. 
Сфера Пуанкаре}
                                                 \dotfill \ 2
\\[0.75ex]
\mbox{}\hspace{1.35em}
1. Отображение плоскости на сферу Пуанкаре
                                                 \dotfill \ 2
\\[0.5ex]
\mbox{}\hspace{1.35em}
2. Атлас сферы Пуанкаре
                                                 \dotfill \ 4
\\[0.5ex]
\mbox{}\hspace{1.35em}
3. Круг Пуанкаре
                                                 \dotfill \ 5
\\[0.5ex]
\mbox{}\hspace{1.35em}
4. Отображения Пуанкаре
                                                 \dotfill \ 7
\\[0.5ex]
\mbox{}\hspace{1.35em}
5. Атлас проективных кругов проективной плоскости
                                                 \dotfill \ 9
\\[1ex]
\noindent
{\bf \S 2. 
Преобразования Пуанкаре дифференциальных систем
}
                                                 \dotfill \ 10
\\[0.75ex]
\mbox{}\hspace{1.35em}
6. Проективно приведенные системы
                                                 \dotfill \ 10
\\[0.5ex]
\mbox{}\hspace{1.35em}
7. Проективный тип дифференциальной системы 
                                                 \dotfill \ 11
\\[0.5ex]
\mbox{}\hspace{1.35em}
8. Проективный атлас траекторий дифференциальных систем
                                                 \dotfill \ 16
\\[1ex]
\noindent
{\bf \S 3. 
Траектории на сфере Пуанкаре
}
                                                 \dotfill \ 19
\\[0.75ex]
\mbox{}\hspace{1.35em}
9.\ \, Траектории на проективной фазовой плоскости
                                                 \dotfill \ 19
\\[0.5ex]
\mbox{}\hspace{1.35em}
10. Траектории первой проективно приведенной системы
                                                 \dotfill \ 23
\\[0.5ex]
\mbox{}\hspace{1.35em}
11. Траектории второй проективно приведенной системы
                                                 \dotfill \ 25
\\[0.5ex]
\mbox{}\hspace{1.35em}
12.  Линейные и разомкнутые  предельные циклы
                                                 \dotfill \ 27
\\[0.5ex]
\mbox{}\hspace{1.35em}
13. Симметpичность фазового поля направлений
                                                 \dotfill \ 33
\\[0.5ex]
\mbox{}\hspace{1.35em}
14. Множества проективно неособых и проективно особых систем
                                                 \dotfill \ 35
\\[0.5ex]
\mbox{}\hspace{1.35em}
15. Топологическая эквивалентность дифференциальных систем 
\\[0.25ex]
\mbox{}\hspace{3.1em}
на проективном круге и на проективной сфере 
                                                 \dotfill \ 37
\\[0.5ex]
\mbox{}\hspace{1.35em}
16. Примеры глобального качественного исследования траекторий
\\[0.25ex]
\mbox{}\hspace{3.1em}
дифференциальных систем на проективной фазовой плоскости
                                                 \dotfill \ 38
\\[1ex]
{\bf Список литературы}
                                              \dotfill \ 60

\newpage

\mbox{}
\\[-1.75ex]
\centerline{\large\bf  Введение}
\\[1.5ex]
\indent
Объектом исследования является обыкновенная автономная полиномиальная
дифференциальная система второго порядка
\\[2ex]
\mbox{}\hfill                                   %(D)
$
\displaystyle 
\dfrac{dx}{dt} =
\sum \limits_{k=0}^{n}\, X_k^{}(x,y)\equiv
 X(x,y), 
\qquad 
\dfrac{dy}{dt} =
\sum \limits_{k=0}^{n}\,Y_k^{}(x,y)\equiv
Y(x,y),
$
\hfill (D)
\\[2.25ex]
где 
\vspace{0.5ex}
$X_k^{}$ и $Y_k^{}$ --- однородные полиномы по переменным $x$ и $y$ степени $k,\ k=0,1,\ldots,n,$ такие, что 
\vspace{0.35ex}
$|X_n^{}(x,y)|+|Y_n^{}(x,y)|\not\equiv0$ на $\R^2,$ а полиномы 
$X$ и $Y$ --- взаимно простые, т.е. не имеют общих делителей, отличных от вещественных чисел.

С целью изучения поведения траекторий системы (D) А. Пуанкаре наряду с состояниями равновесия, 
лежащими в конечной части фазовой плоскости, исследовал и бесконечно удаленные состояния равновесия [1, с. 23 --- 31]. 
В этой связи фазовая плоскость $(x,y)$ пополнялась бесконечно удаленными точками, а траектории системы (D) 
проецировались на сферу с последующим изображением их на круге [1, с. 84 --- 91]. 
Тем самым было положено начало глобальной качественной теории автономных обыкновенных дифференциальных 
систем второго порядка [2 --- 4].

В данной статье с учетом исследований [5] 
рассматривается  поведение траекторий системы (D) на проективной фазовой плоскости $\R\P(x,y).$
Исходим из того, что сфера Пуанкаре (двумерная сфера единичного радиуса с отождествленными антиподальными точками [6, c. 749]) 
диффеоморфна проективной плоскости, а проективная плоскость является двумерным многообразием, 
которое аналитически описывается с помощью трех локальных прямоугольных декартовых систем координат 
[7, с. 96; 8, с. 421 --- 423].
\\[5ex]
\centerline{
{\bf\large \S\;\!1. Сфера Пуанкаре}}
\\[2ex]
\centerline{
{\bf  1. 
Отображение плоскости на сферу Пуанкаре
}
}
\\[1.5ex]
\indent
Введем 
\vspace{0.15ex}
трехмерную прямоугольную декартову систему координат
$O^{\star}x^{\star}y^{\star}z^{\star},$ 
совмещенную с правой прямоугольной декартовой системой координат $Oxy,$
\vspace{0.15ex}
соблюдая условия:
прямая $OO^{\star}$ ортогональна плоскости $Oxy,$
\vspace{0.15ex}
длина отрезка $OO^{\star}$ равна одной единице масштаба системы
координат $Oxy;$
ось $O^{\star}x^{\star}$ сонаправлена с осью $Ox,$ ось
$O^{\star}y^{\star}$ сонаправлена с осью $Oy,$ а
\vspace{0.15ex}
ось $O^{\star}z^{\star}$ направлена так, что система
\vspace{0.15ex}
координат $O^{\star}x^{\star}y^{\star}z^{\star}$ будет правой;
\vspace{0.15ex}
масштаб в системе координат $O^{\star}x^{\star}y^{\star}z^{\star}$
такой же, как и в системе координат $Oxy.$
Построим сферу с центром $O^{\star}$
единичного радиуса: 
\\[2.25ex]
\mbox{}\hfill                                   %(1.1)
$
S^2=\bigl\{
(x^{\star},y^{\star},z^{\star})
\colon
x^{\star}{}^{\,^{\scriptstyle 2}}+
y^{\star}{}^{\,^{\scriptstyle 2}}+
z^{\star}{}^{\,^{\scriptstyle 2}}=1\bigr\}.
$
\hfill (1.1)
\\[2.25ex]
\noindent
Точки 
\vspace{0.25ex}
$N(0,0,1)$ и $S(0,0,{}-1)$ ---  соответственно северный и южный полюсы этой сферы. 
При этом южный полюс $S(0,0,{}-1)$ совпадает с началом $O(0,0)$ системы координат $Oxy.$
Уравнение $z^{\star}={}-1 $
\vspace{0.25ex}
является уравнением  в системе координат $O^{\star}x^{\star}y^{\star}z^{\star}$
плоскости $Oxy.$
Плоскость $Oxy$ касается
сферы (1.1) в южном полюсе $S(0,0,{}-1).$
\vspace{0.35ex}

На плоскости $Oxy$ произвольным образом выберем точку $M$ и
проведем луч с началом $M$ через центр $O^{\star}$ сферы (1.1). 
Луч $MO^{\star}$ пересекает сферу в двух антиподальных (диаметрально противоположных)
точках  $M^{\star}$ и  $M_{\star}^{}$ (рис. 1.1).
Тем самым, каждой точке плоскости $Oxy$ соответствуют две антиподальные точки сферы (1.1).
И, наоборот, всякой паре антиподальных точек сферы (1.1), за исключением лишь точек, 
лежащих
на экваторе
сферы (окружности сферы, лежащей в координатной плоскости $O^{\star}x^{\star}y^{\star}),$ 
сопоставляется одна точка плоскости $Oxy.$

Чтобы распространить соответствие на всю сферу (1.1), поступим следующим образом.
Каждую прямую $l_{_0},$ проходящую через точку $O$ 
и лежащую в координатной плоскости $Oxy,$ 
пополним бесконечно удаленной точкой $L,$
лежащей на ее <<концах>>. 
При этом разным прямым, проходящим через точку $O$ 
и лежащим в плоскости $Oxy,$  соответствуют разные бесконечно удаленные точки,
а пучку параллельных прямых соответствует одна бесконечно удаленная точка.
Прямую $l,$ непрерывно продолженную бесконечно удаленной точкой $L,$ обозначим
$l_{_L}$ или $\overline{l}.$
Плоскость $(x,y),$ дополненная бесконечно удаленными точками 
всех лежащих на ней прямых, 
\vspace{0.15ex}
является [9, с. 247] проективной плоскостью,
которую обозначим $\R\P^2$ или $\R \P (x,y).$
Множество бесконечно удаленных точек проективной плоскости есть [9, с. 245]
бесконечно удаленная прямая этой проективной плоскости. Проективная плоскость 
$\R \P (x,y)$ есть дизъюнктивное объединение плоскости $(x,y)$ 
и бесконечно удаленной прямой.
\\[4.25ex]
\mbox{}\hfill
{\unitlength=1mm
\begin{picture}(77,50)
\put(0,0){\includegraphics[width=96.7mm,height=49.84mm]{r01-01.eps}}
 
\put(44.9,13.7){\makebox(0,0)[cc]{\scriptsize $O$}}
\put(41.5,17.7){\makebox(0,0)[cc]{\scriptsize $S$}}
\put(75,17.5){\makebox(0,0)[cc]{\scriptsize $x$}}
\put(83,17.5){\makebox(0,0)[cc]{\scriptsize $x$}}
\put(30,5){\makebox(0,0)[cc]{\scriptsize $y$}}
\put(49,19){\makebox(0,0)[cc]{\scriptsize $y$}}
%\put(41.9,3.6){\makebox(0,0)[cc]{\scriptsize $(x,y)$}}
\put(66,4){\makebox(0,0)[cc]{\scriptsize $M$}}

\put(40.9,30.7){\makebox(0,0)[cc]{\scriptsize $O^{\star}$}}
\put(52,22.8){\makebox(0,0)[cc]{\scriptsize $M^{\star}$}}
\put(34.2,36){\makebox(0,0)[cc]{\scriptsize $M_{\star}^{}$}}
%\put(32.5,27){\makebox(0,0)[cc]{\scriptsize $L^{\star}_{\scriptscriptstyle-}$}}
%\put(52.8,31.3){\makebox(0,0)[cc]{\scriptsize $L^{\star}_{\scriptscriptstyle+}$}}
%\put(9,2.7){\makebox(0,0)[cc]{\scriptsize $L^{}_{\scriptscriptstyle-}$}}
%\put(65,19.7){\makebox(0,0)[cc]{\scriptsize $L^{}_{\scriptscriptstyle+}$}}
%\put(17,10){\makebox(0,0)[cc]{\scriptsize $l^{}_{\scriptscriptstyle 0}$}}
\put(61.6,27.3){\makebox(0,0)[cc]{\scriptsize $x^{\star}$}}
\put(58,30.5){\makebox(0,0)[cc]{\scriptsize $1$}}
\put(27,30.6){\makebox(0,0)[cc]{\scriptsize ${\scriptstyle-}1$}}
\put(57.3,40){\makebox(0,0)[cc]{\scriptsize $y^{\star}$}}
\put(49.6,37.5){\makebox(0,0)[cc]{\scriptsize $1$}}
\put(35.3,25.3){\makebox(0,0)[cc]{\scriptsize ${\scriptstyle-}1$}}

\put(41.3,44.3){\makebox(0,0)[cc]{\scriptsize $N$}}
\put(44.5,44.3){\makebox(0,0)[cc]{\scriptsize $1$}}
\put(45.5,49.5){\makebox(0,0)[cc]{\scriptsize $z^{\star}$}}

\put(48,-5){\makebox(0,0)[cc]{\rm Рис. 1.1}}
\end{picture}}
\hfill\mbox{}
\\[11ex]
\mbox{}\hfill
{\unitlength=1mm
\begin{picture}(78.9,50)
\put(3,0){\includegraphics[width=75.9mm,height=49.84mm]{r01-02.eps}}
 
\put(44.9,13.7){\makebox(0,0)[cc]{\scriptsize $O$}}
\put(41.7,17.7){\makebox(0,0)[cc]{\scriptsize $S$}}
\put(60,14){\makebox(0,0)[cc]{\scriptsize $x$}}
\put(49,19){\makebox(0,0)[cc]{\scriptsize $y$}}
%\put(41.9,3.6){\makebox(0,0)[cc]{\scriptsize $(x,y)$}}
%\put(66,4){\makebox(0,0)[cc]{\scriptsize $M$}}

\put(40.9,30.7){\makebox(0,0)[cc]{\scriptsize $O^{\star}$}}
%\put(52,22.8){\makebox(0,0)[cc]{\scriptsize $M^{\star}$}}
%\put(37.6,38.9){\makebox(0,0)[cc]{\scriptsize $M_{\star}^{}$}}
\put(32.5,27){\makebox(0,0)[cc]{\scriptsize $L_{\star}$}}
\put(52.8,31.3){\makebox(0,0)[cc]{\scriptsize $L^{\star}$}}
\put(12,4){\makebox(0,0)[cc]{\scriptsize $L$}}
\put(65,20){\makebox(0,0)[cc]{\scriptsize $L$}}
\put(23,12){\makebox(0,0)[cc]{\scriptsize $l_{{}_0}$}}
\put(61.6,27.3){\makebox(0,0)[cc]{\scriptsize $x^{\star}$}}
\put(58,30.5){\makebox(0,0)[cc]{\scriptsize $1$}}
\put(27,30.6){\makebox(0,0)[cc]{\scriptsize ${\scriptstyle-}1$}}
\put(57.3,40){\makebox(0,0)[cc]{\scriptsize $y^{\star}$}}
\put(49.6,37.5){\makebox(0,0)[cc]{\scriptsize $1$}}
\put(35.3,25.3){\makebox(0,0)[cc]{\scriptsize ${\scriptstyle-}1$}}

\put(41.3,44.3){\makebox(0,0)[cc]{\scriptsize $N$}}
\put(44.5,44.3){\makebox(0,0)[cc]{\scriptsize $1$}}
\put(45.5,49.5){\makebox(0,0)[cc]{\scriptsize $z^{\star}$}}

\put(39,-5){\makebox(0,0)[cc]{\rm Рис. 1.2}}
\end{picture}}
\hfill\mbox{}
\\[6.5ex]
\indent
На координатной плоскости $Oxy$ произвольным образом выберем прямую $l_{_0}\,,$
проходящую через начало координат $O.$ Построим плоскость $\Pi_{l_{_0}}^{},$
проходящую через ось $O^{\star}z^{\star}$ и прямую $l_{_0}.$
Плоскость  $\Pi_{l_{_0}}^{}$ пересекает экватор сферы (1.1) в двух антиподальных точках 
$L^{\star}$ и  $L_{\star}^{}$ (рис. 1.2).
Бесконечно удаленной точке $L$ прямой $\overline{l}$ поставим в соответствие две антиподальные точки 
$L^{\star}$ и  $L_{\star}^{}$ экватора сферы (1.1),
лежащие на плоскости $\Pi_{l_{_0}}^{}.$
Тогда продолженной прямой $\overline{l}$ проективной плоскости $\R \P (x,y)$ на сфере (1.1) будет соответствовать
окружность большого радиуса, проходящая через точки $L^{\star}$ и  $L_{\star}^{}.$

Итак, установлено бинарное двузначное соответствие между проективной плоскостью $\R \P (x,y)$ и сферой (1.1), при котором образом каждой точки плоскости $\R \P (x,y)$ является множество, состоящее из двух антиподальных точек сферы (1.1).
Таким образом, такую сферу (1.1) будем называть {\it проективной сферой} плоскости $Oxy$ и обозначать 
\vspace{0.5ex}
$\P{\mathbb S}(x,y).$

Проективную сферу  
\vspace{0.15ex}
$\P{\mathbb S}(x,y)$ с отождествленными антиподальными точками назовем [6, c. 749]
{\it сферой Пуанкаре}  плоскости $Oxy$ и будем говорить о сфере Пуанкаре  $\P{\mathbb S}(x,y).$

Тогда введенное бинарное двузначное соответствие между проективной плоскостью $\R \P (x,y)$ и 
проективной сферой $\P{\mathbb S}(x,y)$ устанавливает 
биективное отображение проективной плоскости $\R \P (x,y)$ на сферу Пуанкаре $\P{\mathbb S}(x,y).$

Следовательно, сфера Пуанкаре является двумерным многообразием [7, c. 92 -- 93],
гомеоморфным проективной плоскости.
\\[2.5ex]
\centerline{
{\bf  2. Атлас карт  сферы Пуанкаре}
}
\\[1.5ex]
\indent
Биекция между 
сферой Пуанкаре $\P{\mathbb S}(x,y)$
и 
проективной плоскостью  $\R\P (x,y)$ позволяет построить атлас карт сферы Пуанкаре на основании атласа карт проективной плоскости. Для этого, следуя [7, c. 96], наряду с системой координат $Oxy$
введем еще две плоские прямоугольные декартовы системы координат.
\vspace{0.15ex}

На плоскости,
\vspace{0.35ex}
касающейся сферы (1.1) в точке с 
координатами
$x^{\star}=1,\, y^{\star}=0,\, z^{\star}=0,$
введем правую прямоугольную декартову систему координат $O^{(1)}_{\phantom{1}} \xi \theta$ так, что 
\vspace{0.35ex}
ее начало $O^{(1)}_{\phantom{1}}$ совпадает с точкой касания сферы (1.1), 
ось $O^{(1)}_{\phantom{1}} \xi$ сонаправлена с осью 
\vspace{0.35ex}
$O^{\star} y^{\star}$ (а значит, и с осью $O y$\!), 
ось $O^{(1)}_{\phantom{1}} \theta$ противоположно направлена с осью   $O^{\star} z^{\star}$ (рис. 2.1).
\\[3.5ex]
\mbox{}\hfill
{\unitlength=1mm
\begin{picture}(76.272,44)
\put(0,0){\includegraphics[width=81.272mm,height=44mm]{r02-01.eps}}

\put(34.5,9.7){\makebox(0,0)[cc]{\scriptsize $O$}}
%\put(36.7,13.7){\makebox(0,0)[cc]{\scriptsize $S$}}
%\put(71,13.5){\makebox(0,0)[cc]{\scriptsize $x$}}
\put(45,10.5){\makebox(0,0)[cc]{\scriptsize $x$}}
%\put(26.6,2.9){\makebox(0,0)[cc]{\scriptsize $y$}}
\put(42.7,15){\makebox(0,0)[cc]{\scriptsize $y$}}
%\put(36.9,3.6){\makebox(0,0)[cc]{\scriptsize $(x,y)$}}
\put(59.8,4.7){\makebox(0,0)[cc]{\scriptsize $M$}}
\put(46,19.5){\makebox(0,0)[cc]{\scriptsize $M^{\star}$}}
\put(50.2,15.9){\makebox(0,0)[cc]{\scriptsize $M^{(1)}_{\phantom1}$}}
\put(31,34.3){\makebox(0,0)[cc]{\scriptsize $M_{\star}$}}

%\put(19.5,27.3){\makebox(0,0)[cc]{\scriptsize $O^{(1)}_{\star}$}}
\put(52.7,30.3){\makebox(0,0)[cc]{\scriptsize $O^{(1)}_{\phantom{1}}$}}
\put(55.5,23.2){\makebox(0,0)[cc]{\scriptsize $x^{\star}$}}

\put(42,33){\makebox(0,0)[cc]{\scriptsize $O^{(2)}_{\phantom{1}}$}}
%\put(31.7,16.6){\makebox(0,0)[cc]{\scriptsize $O^{(2)}_{\star}$}}
\put(33.9,26.7){\makebox(0,0)[cc]{\scriptsize $O^{\star}$}}
%\put(14.5,22.7){\makebox(0,0)[cc]{\scriptsize $z^{\star(1)}$}}
%\put(53,26.5){\makebox(0,0)[cc]{\scriptsize $1$}}
%\put(22,26.6){\makebox(0,0)[cc]{\scriptsize $-1$}}
\put(47.4,38){\makebox(0,0)[cc]{\scriptsize $y^{\star}$}}
\put(39,42){\makebox(0,0)[cc]{\scriptsize $z^{\star}$}}

\put(48.8,42.5){\makebox(0,0)[cc]{\scriptsize $\vec{n}$}}
\put(51.8,18.8){\makebox(0,0)[cc]{\scriptsize $\theta$}}
\put(56.5,34){\makebox(0,0)[cc]{\scriptsize $\xi$}}
%\put(44.6,33.5){\makebox(0,0)[cc]{\scriptsize $1$}}
%\put(30.3,21.3){\makebox(0,0)[cc]{\scriptsize $-1$}}

%\put(36.3,40.3){\makebox(0,0)[cc]{\scriptsize $N$}}
%\put(39.5,40.3){\makebox(0,0)[cc]{\scriptsize $1$}}
%\put(34,4.5){\makebox(0,0)[cc]{\scriptsize $\theta^{\star}$}}

\put(40.6,-5){\makebox(0,0)[cc]{\rm Рис. 2.1}}
\end{picture}}
\hfill\mbox{}
\\[10ex]
\mbox{}\hfill
{\unitlength=1mm
\begin{picture}(87.153,44)
\put(0,0){\includegraphics[width=87.153mm,height=44mm]{r02-02.eps}}

\put(40.6,9.7){\makebox(0,0)[cc]{\scriptsize $O$}}
%\put(45.7,13.7){\makebox(0,0)[cc]{\scriptsize $S$}}
%\put(80,13.5){\makebox(0,0)[cc]{\scriptsize $x$}}
\put(51,10.5){\makebox(0,0)[cc]{\scriptsize $x$}}
%\put(35.6,2.9){\makebox(0,0)[cc]{\scriptsize $y$}}
\put(48.3,15){\makebox(0,0)[cc]{\scriptsize $y$}}
%\put(45.9,3.6){\makebox(0,0)[cc]{\scriptsize $(x,y)$}}
\put(65.8,4.7){\makebox(0,0)[cc]{\scriptsize $M$}}
\put(51.6,19.8){\makebox(0,0)[cc]{\scriptsize $M^{\star}$}}
\put(31.5,40.5){\makebox(0,0)[cc]{\scriptsize $M^{(2)}_{\phantom1}$}}
\put(36.9,34.6){\makebox(0,0)[cc]{\scriptsize $M_{\star}$}}

\put(39.9,26.7){\makebox(0,0)[cc]{\scriptsize $O^{\star}$}}
%\put(25.5,27.3){\makebox(0,0)[cc]{\scriptsize $O^{(1)}_{\star}$}}
\put(59.7,26.7){\makebox(0,0)[cc]{\scriptsize $O^{(1)}_{\phantom{1}}$}}
\put(47.3,33){\makebox(0,0)[cc]{\scriptsize $O^{(2)}_{\phantom{1}}$}}
\put(37.7,17){\makebox(0,0)[cc]{\scriptsize $O^{(2)}_{\star}$}}
\put(60.6,23.1){\makebox(0,0)[cc]{\scriptsize $x^{\star}$}}
%\put(53,26.5){\makebox(0,0)[cc]{\scriptsize $1$}}
%\put(22,26.6){\makebox(0,0)[cc]{\scriptsize $-1$}}
%\put(31,10){\makebox(0,0)[cc]{\scriptsize $z^{\star(2)}$}}
\put(49.4,23){\makebox(0,0)[cc]{\scriptsize $\eta$}}
\put(65,28){\makebox(0,0)[cc]{\scriptsize $\zeta$}}
\put(57.8,37){\makebox(0,0)[cc]{\scriptsize $y^{\star}$}}
\put(44.6,42){\makebox(0,0)[cc]{\scriptsize $z^{\star}$}}

%\put(44.6,33.5){\makebox(0,0)[cc]{\scriptsize $1$}}
%\put(30.3,21.3){\makebox(0,0)[cc]{\scriptsize $-1$}}

%\put(36.3,40.3){\makebox(0,0)[cc]{\scriptsize $N$}}
%\put(39.5,40.3){\makebox(0,0)[cc]{\scriptsize $1$}}
%\put(45,5){\makebox(0,0)[cc]{\scriptsize $\eta^{\star}$}}
\put(76,37){\makebox(0,0)[cc]{\scriptsize $\vec{n}$}}

\put(43.5,-5){\makebox(0,0)[cc]{\rm Рис. 2.2}}
\end{picture}}
\hfill\mbox{}
\\[6ex]
\indent
На плоскости, 
\vspace{0.35ex}
касающейся сферы (1.1) в точке с координатами $x^{\star}=0,\, y^{\star}=1,\, z^{\star}=0,$
введем правую прямоугольную декартову систему координат $O^{(2)}_{\phantom{1}} \eta \zeta$ так, что 
\vspace{0.35ex}
ее начало $O^{(2)}_{\phantom{1}}$ совпадает с точкой касания сферы (1.1), 
\vspace{0.35ex}
ось $O^{(2)}_{\phantom{1}} \eta$ противоположно направлена с осью  $O^{\star} z^{\star},$
ось $O^{(2)}_{\phantom{1}} \zeta$ сонаправлена с осью $O^{\star} x^{\star}$ (рис. 2.2).
\vspace{0.35ex}
Масштаб в системах координат $Oxy,\, O^{\star}x^{\star}y^{\star}z^{\star},\, O^{(1)}_{\phantom{1}}\xi\theta$ и 
$O^{(2)}_{\phantom{1}} \eta \zeta$
одинаковый.
\vspace{0.5ex}

Покроем сферу Пуанкаре (с учетом отождествленности антиподальных точек) тремя полусферами без края
\\[2.5ex]
\mbox{}\hfill
$
\displaystyle
U_1^{}=\Bigl\{
(x^{\star}, y^{\star}, z^{\star})\colon\,  z^{\star}={}-\sqrt{
1-{x^{\star}}^{\,^{\scriptstyle 2}}-{y^{\star}}^{\,^{\scriptstyle 2}} }\,,\ \
{x^{\star}}^{\,^{\scriptstyle 2}}+{y^{\star}}^{\,^{\scriptstyle 2}}<1\Bigr\},
\hfill
$
\\[2.25ex]
\mbox{}\hfill
$
\displaystyle
U_2^{}=\Bigl\{
(x^{\star}, y^{\star}, z^{\star})\colon\,  x^{\star}=\sqrt{
1-{y^{\star}}^{\,^{\scriptstyle 2}}-{z^{\star}}^{\,^{\scriptstyle 2}} }\,,\ \
{y^{\star}}^{\,^{\scriptstyle 2}}+{z^{\star}}^{\,^{\scriptstyle 2}}<1\Bigr\},
\hfill
$
\\[2.25ex]
\mbox{}\hfill
$
\displaystyle
U_3^{}=\Bigl\{
(x^{\star}, y^{\star}, z^{\star})\colon\, y^{\star}=\sqrt{
1-{x^{\star}}^{\,^{\scriptstyle 2}}-{z^{\star}}^{\,^{\scriptstyle 2}} }\,,\ \
{x^{\star}}^{\,^{\scriptstyle 2}}+{z^{\star}}^{\,^{\scriptstyle 2}}<1\Bigr\}
\hfill
$
\\[2.75ex]
и введем биективные отображения [7, c. 96]
\\[2.25ex]
\mbox{}\hfill                                                              % (2.1)
$
\displaystyle
\varphi_1^{}\colon (x^{\star}, y^{\star}, z^{\star})\to \
\bigl(
x(x^{\star}, y^{\star}, z^{\star}),\, y(x^{\star}, y^{\star}, z^{\star})
\bigr), 
\hfill
$
\\
\mbox{}\hfill (2.1)
\\
\mbox{}\hfill
$
x={}-\dfrac{x^{\star}}{z^{\star}}\,, 
\quad 
y={}-\dfrac{y^{\star}}{z^{\star}}
\quad
\forall (x^{\star}, y^{\star}, z^{\star})\in U_1^{}\;\!,
\hfill
$
\\[3.5ex]
\mbox{}\hfill                                                              % (2.2)
$
\displaystyle
\varphi_2^{}\colon (x^{\star}, y^{\star}, z^{\star})\to \
\bigl(
\xi(x^{\star}, y^{\star}, z^{\star}),\, \theta(x^{\star}, y^{\star}, z^{\star})
\bigr), 
\hfill
$
\\
\mbox{}\hfill (2.2)
\\
\mbox{}\hfill
$
\xi=\dfrac{y^{\star}}{x^{\star}}\,, 
\,
\quad 
\theta={}-\dfrac{z^{\star}}{x^{\star}}
\quad
\forall (x^{\star}, y^{\star}, z^{\star})\in U_2^{}\;\!,
\hfill
$
\\[3.5ex]
\mbox{}\hfill                                                              % (2.3)
$
\displaystyle
\varphi_3^{}\colon (x^{\star}, y^{\star}, z^{\star})\to \
\bigl(
\eta(x^{\star}, y^{\star}, z^{\star}),\, \zeta(x^{\star}, y^{\star}, z^{\star})
\bigr), 
\hfill
$
\\
\mbox{}\hfill (2.3)
\\
\mbox{}\hfill
$
\eta={}-\dfrac{z^{\star}}{y^{\star}}\,, 
\,\quad 
\zeta=\dfrac{x^{\star}}{y^{\star}}
\quad
\forall (x^{\star}, y^{\star}, z^{\star})\in U_3^{}\;\!.
\hfill
$
\\[3ex]
\indent
Таким образом, построены три карты  
\vspace{0.5ex}
$(U_\tau^{}, \varphi_\tau^{}),\, \tau=1,2,3,$ сферы Пуанкаре $\P{\mathbb S}(x,y).$
Множество карт 
\vspace{0.5ex}
$(U_\tau^{}, \varphi_\tau^{}),\, \tau=1,2,3,$ образуют атлас карт сферы Пуанкаре $\P{\mathbb S}(x,y).$

Заметим, что атлас карт проективной сферы $\P{\mathbb S}(x,y)$ состоит из шести карт.
\\[3.5ex]
\centerline{
{\bf  3. Круг Пуанкаре}
}
\\[1.5ex]
\indent
Пусть точка $M$ 
\vspace{0.35ex}
расположена в конечной части $(x,y)$ проективной плоскости 
$\R\P(x,y)$ и имеет координаты $M(x,y)$ в системе координат $Oxy$ (рис. 1.1).
\vspace{0.35ex}
Тогда в пространственной системе координат 
\vspace{0.35ex}
$O^{\star}x^{\star} y^{\star} z^{\star}$ эта же точка имеет координаты $M(x,y,{}-1).$
Прямая $MO^{\star}$ в системе координат
$O^{\star}x^{\star} y^{\star} z^{\star}$ задается системой уравнений
\\[2ex]
\mbox{}\hfill                                    %(3.1)
$
\dfrac{x^{\star}}{x}=
\dfrac{y^{\star}}{y}=
\dfrac{z^{\star}}{{}-1}\,.
$
\hfill (3.1)
\\[2.5ex]
Точке $M(x,y,{}-1)$ 
\vspace{0.5ex}
соответствуют точки $M^{\star}(x^{\star},y^{\star},z^{\star})$
и $M_{\star}^{}(x_{\star}^{}\,,y_{\star}^{}\,,z_{\star}^{}),$ являющиеся точками пересечения прямой 
$MO^{\star}$ и сферы (1.1).
\vspace{0.35ex}
Поэтому координаты
$x^{\star},\,y^{\star},\,z^{\star}$ и 
$x_{\star}^{}\,,\,y_{\star}^{}\,,\,z_{\star}^{}$
точек $M^{\star}$ и
$M_{\star}^{}$
суть решения алгебраической системы уравнений
\\[2.5ex]
\mbox{}\hfill                                   %(3.2)
$
\dfrac{x^{\star}}{x}=
\dfrac{y^{\star}}{y}=
\dfrac{z^{\star}}{{}-1}\,, \quad
x^{\star}{}^{\,^{\scriptstyle 2}}+
y^{\star}{}^{\,^{\scriptstyle 2}}+
z^{\star}{}^{\,^{\scriptstyle 2}}=1.
$
\hfill (3.2)
\\[2.5ex]
\indent
Будем считать, что точка 
\vspace{0.35ex}
$M^{\star} (x^{\star},y^{\star},z^{\star})$  лежит  
в южной полусфере $S^2_{\scriptscriptstyle-}.$
Тогда ее аппликата $z^{\star}\in [{}-1;0).$
\vspace{0.35ex}
Разрешив систему уравнений (3.2) относительно 
$x^{\star},\,y^{\star},\,z^{\star}$ при ${}-1\leq z^{\star}<0,$
получим биективное отражение 
\\[2.25ex]
\mbox{}\hfill                                   %(3.3)
$
\varphi_1^{{}-1}\colon (x,y)\to\ \bigl(x^{\star}(x,y), y^{\star}(x,y), z^{\star}(x,y)\bigr),
\hfill                                  
$
\\[-0.75ex]
\mbox{}\hfill {\rm (3.3)}
\\[1.25ex]
\mbox{}\hfill 
$
x^{\star}=\dfrac{x}{\sqrt{1+x^2+y^2}}\,, \ \ \
y^{\star}=\dfrac{y}{\sqrt{1+x^2+y^2}}\,, \ \ \
z^{\star}={}-\dfrac{1}{\sqrt{1+x^2+y^2}}\quad
\forall (x,y)\in \R^2
\hfill 
$
\\[2.75ex]
конечной части $(x,y)$ 
\vspace{0.35ex}
проективной плоскости $\R\P(x,y)$ на южную полусферу 
$S^2_{\scriptscriptstyle-}$ без края $\partial S^2_{\scriptscriptstyle-}$
(экватора сферы (1.1)).

Координатные функции 
\vspace{0.15ex}
отображения (3.3) непрерывно дифференцируемы.
Координатная функция $z^{\star}$
выражается через координатные функции $x^{\star},\, y^{\star}$
по формуле
\\[2.25ex]
\mbox{}\hfill                                  
$
z^{\star}(x,y)={}-\sqrt{1- x^{\star}{}^{\,^{\scriptstyle 2}}(x,y)- y^{\star}{}^{\,^{\scriptstyle 2}}(x,y) }
\quad
\forall (x,y)\in \R^2.
\hfill
$
\\[2.25ex]
Переход от координат $x,\,y$ к координатам $x^{\star},\, y^{\star}$ имеет якобиан
\\[2.25ex]
\mbox{}\hfill                                  
$
\dfrac{{\sf D} (x^{\star},y^{\star})}
{{\sf D} (x,y)}=
\dfrac{1}{(1+x^2+y^2)^2}
\ne 0\quad \forall (x,y) \in \R^2.
\hfill
$
\\[2ex]
Следовательно, биективное отображение (3.3) является диффеоморфизмом. 
\vspace{0.25ex}

Отображения (2.1) и (3.3) являются взаимно обратными.
Диффеоморфность отображения (3.3) 
\vspace{0.25ex}
означает диффеоморфность отображения (2.1)
южной полусферы без края  
$S^2_{\scriptscriptstyle-} \backslash \partial S^2_{\scriptscriptstyle-}=U_1^{}$ на $\R^2.$
\vspace{0.35ex}

С помощью биективного отображения (3.3) конечная часть 
$(x,y)$ проективной плоскости $\R\P(x,y)$ диффеоморфно отображается на южную полусферу
$S^2_{\scriptscriptstyle-}$ без края 
$\partial  S^{2}_{\scriptscriptstyle-},$
являющегося экватором 
проективной сферы $\P{\mathbb S}(x,y)$ (рис. 1.1).
Каждой бесконечно удаленной точке проективной плоскости 
$\R\P(x,y)$ соответствуют две антиподальные точки экватора
$\partial  S^{2}_{\scriptscriptstyle-}$ 
проективной сферы $\P{\mathbb S}(x,y).$
Значит, южная полусфера $S^2_{\scriptscriptstyle-}$ с отождествленными антиподальными точками 
экватора $\partial  S^{2}_{\scriptscriptstyle-}$ 
проективной сферы $\P{\mathbb S}(x,y)$
 является пространственной моделью проективной плоскости $\R\P(x,y).$
\\[3.25ex]
\mbox{}\hfill
{\unitlength=1mm
\begin{picture}(109,53)
\put(0,0){\includegraphics[width=109mm,height=52.46mm]{r03-01.eps}}

\put(51.9,13){\makebox(0,0)[cc]{\scriptsize $O$}}
%\put(36.7,13.7){\makebox(0,0)[cc]{\scriptsize $S$}}
%\put(71,13.5){\makebox(0,0)[cc]{\scriptsize $x$}}
\put(73,14){\makebox(0,0)[cc]{\scriptsize $x$}}
%\put(26.6,2.9){\makebox(0,0)[cc]{\scriptsize $y$}}
\put(59.8,25.7){\makebox(0,0)[cc]{\scriptsize $y$}}
%\put(36.9,3.6){\makebox(0,0)[cc]{\scriptsize $(x,y)$}}
\put(76,2.7){\makebox(0,0)[cc]{\scriptsize $M$}}

\put(48,34.4){\makebox(0,0)[cc]{\scriptsize $O^{\star}$}}
\put(62.8,23.3){\makebox(0,0)[cc]{\scriptsize $M^{\star}$}}
\put(55.2,11.2){\makebox(0,0)[cc]{\scriptsize $\varkappa(M)$}}
%\put(34,8.2){\makebox(0,0)[cc]{\scriptsize $D^{}_{\scriptscriptstyle-}$}}
%\put(71,18){\makebox(0,0)[cc]{\scriptsize $D^{}_{\scriptscriptstyle+}$}}
\put(36,29.5){\makebox(0,0)[cc]{\scriptsize $L^{}_{\star}$}}
\put(68.5,38.5){\makebox(0,0)[cc]{\scriptsize $L^{\star}$}}
\put(10,2){\makebox(0,0)[cc]{\scriptsize $L$}}
\put(88,22.7){\makebox(0,0)[cc]{\scriptsize $L$}}
\put(22,10){\makebox(0,0)[cc]{\scriptsize $l^{}_{\scriptscriptstyle 0}$}}
\put(73,31){\makebox(0,0)[cc]{\scriptsize $x^{\star}$}}
%\put(53,26.5){\makebox(0,0)[cc]{\scriptsize $1$}}
%\put(22,26.6){\makebox(0,0)[cc]{\scriptsize $-1$}}
\put(62,46.5){\makebox(0,0)[cc]{\scriptsize $y^{\star}$}}
%\put(54,19){\makebox(0,0)[cc]{\scriptsize $z$}}
%\put(59.3,34){\makebox(0,0)[cc]{\scriptsize $u$}}
%\put(44.6,33.5){\makebox(0,0)[cc]{\scriptsize $1$}}
%\put(30.3,21.3){\makebox(0,0)[cc]{\scriptsize $-1$}}

%\put(36.3,40.3){\makebox(0,0)[cc]{\scriptsize $N$}}
%\put(39.5,40.3){\makebox(0,0)[cc]{\scriptsize $1$}}
\put(53,51.5){\makebox(0,0)[cc]{\scriptsize $z^{\star}$}}

\put(54.5,-5){\makebox(0,0)[cc]{\rm Рис. 3.1}}
\end{picture}}
\hfill\mbox{}
\\[7ex]
\indent
Естественной проекцией южной полусферы 
\vspace{0.75ex}
$S^2_{\scriptscriptstyle-}$ 
на плоскость $(x,y)$ является круг единичного радиуса $K(x,y)=\{(x,y)\colon x^2+y^2\leq1\}$ (рис 3.1).
\vspace{0.5ex}
С помощью этой проекции устанавливается диффеоморфное отображение $p$ 
\vspace{0.5ex}
южной полусферы Пуанкаре 
$S^2_{\scriptscriptstyle-}$ на круг $K(x,y).$ Суперпозиция  отображений 
\vspace{1.5ex}
$\varkappa=p\circ \varphi_1^{{}-1}$ есть диффеоморфное отображение
%\\[2.25ex]
%\mbox{}\hfill                                 
$
\varkappa \colon (x,y)\to
\biggl( \dfrac{x}{\sqrt{1+x^2+y^2}}\,,\ \dfrac{y}{\sqrt{1+x^2+y^2}}\biggr)\ \
\forall (x,y) \in \R^2
%\hfill  
$
%\\[2.5ex]
\vspace{1.5ex}
конечной части $(x,y)$ 
проективной плоскости $\R\P(x,y)$ на открытый круг $K(x,y)\backslash\partial K(x,y).$ 
\vspace{0.5ex}
Каждой бесконечно удаленной точке проективной плоскости $\R\P(x,y)$ 
\vspace{0.5ex}
соответствуют две антиподальные точки граничной окружности 
\vspace{0.5ex}
$\partial K(x,y)=\{(x,y)\colon x^2+y^2=1\}$ круга $K(x,y).$ 
Такой круг $K(x,y)$  будем называть 
\vspace{0.5ex}
{\it проективным кругом} плоскости $Oxy$
и обозначать $\P\K(x,y).$

Проективный круг $\P\K(x,y)$ с отождествленными 
\vspace{0.35ex}
антиподальными точками граничной окружности $\partial K(x,y)$ назовем {\it кругом Пуанкаре} 
\vspace{0.35ex}
плоскости $Oxy$ 
и будем говорить о круге Пуанкаре $\P\K(x,y).$
\vspace{0.35ex}
Круг Пуанкаре  $\P\K(x,y)$ диффеоморфен проективной плоскости $\R\P(x,y),$ а значит, является плоской компактной моделью этой проективной плоскости.
\\[7.5ex]
\centerline{
{\bf  4. Отображения Пуанкаре}
}
\\[1.5ex]
\indent
Установим связи между 
\vspace{0.35ex}
локальными системами координат
$Oxy,\ O^{(1)}_{\phantom{1}}\xi\theta,\ O^{(2)}_{\phantom{1}}\eta\zeta$
атласа карт сферы Пуанкаре
$\P{\mathbb S}(x,y).$
\vspace{0.5ex}
В системе координат 
$O^{\star}x^{\star}y^{\star}z^{\star}$ плоскость $O^{(1)}_{\phantom{1}} \xi \theta$
\vspace{0.5ex}
задается уравнением $x^{\star}=1.$ 
Произвольным образом выберем в конечной части $(x,y)$
проективной плоскости $\R \P (x,y)$ точку $M(x,y)$ 
\vspace{0.5ex}
так, чтобы она не лежала на оси $Oy.$
Тогда  (рис. 2.1) прямая $MO^{\star}$ 
\vspace{0.5ex}
пересекает плоскость $O^{(1)}_{\phantom{1}} \xi \theta$
в точке $M^{(1)}.$ Прямая $MO^{\star}$ задается системой уравнений (3.1),
\vspace{0.5ex}
а плоскость  $O^{(1)}_{\phantom{1}} \xi \theta$ --- уравнением $x^{\star}=1.$ 
Поэтому в системе координат $O^{\star}x^{\star}y^{\star}z^{\star}$
у точки $M^{(1)}$ абсцисса $x^{\star}=1,$ 
\vspace{0.5ex}
а ординату
$y^{\star}=\dfrac{y}{x}$ и аппликату $z^{\star}={}-\dfrac{1}{x}$ 
\vspace{0.5ex}
находим из системы уравнений (3.1) при $x^{\star}=1.$
В системе координат
$O^{(1)}_{\phantom{1}}\xi\theta$
у точки $M_{\phantom{1}}^{(1)}$ абсцисса $\xi=y^{\star}=\dfrac{y}{x}\,,$ а 
ордината $\theta={}-z^{\star}=\dfrac{1}{x}\,.$
\vspace{1ex}
Таким образом, лежащей в конечной части $(x,y)$ 
\vspace{0.5ex}
проективной плоскости $\R \P (x,y)$ точке $M(x,y)$
с абсциссой $x\ne 0$ на  координатной плоскости 
$O^{(1)}_{\phantom{1}}\xi\theta$
\vspace{0.5ex}
соответствует точка $M_{\phantom{1}}^{(1)}(\xi,\theta),$
координаты которой через координаты точки $M(x,y)$ выражаются по формулам
\\[2.25ex]
\mbox{}\hfill       % (4.1)
$
\xi=\dfrac{y}{x}\,,
\quad \
\theta=\dfrac{1}{x}\,.
$
\hfill (4.1)
\\[2ex]
\indent
Тем самым установлено биективное отображение
\\[2ex]
\mbox{}\hfill       % (4.2)
$
P^{(1)}_{}\colon (x,y)\to \ 
\Bigl(\, \dfrac{y}{x}\,,\, \dfrac{1}{x}\;\!\Bigr)
\quad
\forall (x,y)\in\R^2\backslash \{(x,y)\colon x=0\}
$
\hfill (4.2)
\\[2ex]
плоскости $Oxy,$
\vspace{0.5ex}
из которой удалена  ось $Oy,$
на плоскость 
$O^{(1)}_{\phantom{1}}\xi\theta,$
из которой удалена ось
$O^{(1)}_{\phantom{1}}\xi.$
\vspace{0.75ex}
Отображение (4.2) назовем {\it первым отображением Пуанкаре}
плоскости $Oxy.$

Отображение (4.2) является диффеоморфизмом с якобианом
\\[2.25ex]
\mbox{}\hfill
$
\dfrac{{\sf D}(\xi,\theta)}{{\sf D}(x,y)}=\dfrac{1}{x^3}\ne 0
\quad
\forall (x,y)\in\R^2\backslash \{(x,y)\colon x=0\}.
\hfill
$
\\[2.25ex]
\indent
Разрешив равенства (4.1) относительно $x$ и $y,$ получим формулы
\\[2ex]
\mbox{}\hfill       % (4.3)
$
x=\dfrac{1}{\theta}\,,
\quad \
y=\dfrac{\xi}{\theta}\,,
$
\hfill (4.3)
\\[2.25ex]
по которым координаты точки 
\vspace{0.35ex}
$M(x,y)$ при $x\ne 0$ выражаются через координаты точки 
$M^{(1)}$ при $\theta\ne 0.$
\vspace{0.35ex}
Формулы (4.3) назовем
{\it первым преобразованием Пуанкаре}
плоскости $Oxy$ [1, с. 31].
\vspace{0.35ex}

Функции (4.1) 
\vspace{0.35ex}
есть функции перехода от координат $(\xi,\theta)$ к координатам $(x,y),$
а функции (4.3) есть функции перехода от координат $(x,y)$ к координатам $(\xi,\theta).$
\vspace{0.5ex}

Произвольным образом возьмем прямую $y=ax$ 
\vspace{0.35ex}
с параметром $a\in\R$ и преобразуем ее по формулам (4.3).
В результате получим прямую $\xi=a$ на плоскости 
$O^{(1)}_{\phantom{1}}\xi\theta.$
\vspace{0.35ex}
При этом бесконечно удаленной точке $L,$ 
\vspace{0.5ex}
лежащей на <<концах>> прямой $y=ax,$ соответствует точка 
$L_{\phantom{1}}^{(1)}(a,0),$ лежащая на оси 
$O^{(1)}_{\phantom{1}}\xi$
\vspace{0.5ex}
плоскости $O^{(1)}_{\phantom{1}}\xi\theta.$
И, наоборот, каждой точке 
$L_{\phantom{1}}^{(1)}(a,0),$ лежащей на координатной оси 
\vspace{0.5ex}
$O^{(1)}_{\phantom{1}}\xi$
плоскости $O^{(1)}_{\phantom{1}}\xi\theta,$
соответствует бесконечно удаленная точка 
$L,$ лежащая на <<концах>> прямой 
$y=ax.$ 
\vspace{0.5ex}

Стало быть, имеет место
\vspace{0.5ex}

{\bf Свойство 4.1.} 
{\it
Первое преобразование Пуанкаре {\rm (4.3)} 
\vspace{0.5ex}
устанавливает диффеоморфное отображение проективной плоскости $\!\R \P (x,y),\!$ 
\vspace{0.5ex}
у которой удалена прямая $x=0,$
на сов\-ме\-щен\-ную координатную плоскость 
$\!O^{(1)}_{\phantom{1}}\xi\theta.\!\!$
\vspace{0.5ex}
При этом образом бесконечно удален\-ной прямой
проективной плоскости $\R \P (x,y),$ 
\vspace{0.5ex}
из которой удалена точка, лежащая на 
<<концах>> прямой $\!x\!=\!0,\!$ есть координатная ось 
$\!O^{(1)}_{\phantom{1}}\xi\!$ {\rm (}прямая $\theta\!=\!0)\!$
\vspace{1ex}
плоскости $\!O^{(1)}_{\phantom{1}}\xi\theta.$}

В системе координат
\vspace{0.5ex}
$O^{\star}x^{\star}y^{\star}z^{\star}$ плоскость 
$O^{(2)}_{\phantom{1}}\eta\zeta$ задается уравнением 
$y^{\star}=1.$
Произвольным образом выберем в конечной части $(x,y)$ 
\vspace{0.5ex}
проективной плоскости 
$\R \P (x,y)$ точку $M(x,y)$ так, чтобы она не лежала на оси $Ox.$ 
\vspace{0.5ex}
Тогда (рис. 2.2) 
прямая $MO^{\star}$ пересекает плоскость 
$O^{(2)}_{\phantom{1}}\eta\zeta$ в точке $M^{(2)}.$
\vspace{0.75ex}
Прямая $MO^{\star}$ задается системой уравнений (2.1), а плоскость  
\vspace{0.5ex}
$O^{(2)}_{\phantom{1}}\eta\zeta$ --- уравнением $y^{\star}=1.$ 
\vspace{0.75ex}
Поэтому в системе координат 
$O^{\star}x^{\star}y^{\star}z^{\star}$ у точки 
$M_{\phantom{1}}^{(2)}$ ордината $y^{\star}=1,$ 
а абсциссу
\vspace{0.75ex}
$x^{\star}=\dfrac{x}{y}$ и аппликату $z^{\star}={}-\dfrac{1}{y}$
находим из системы уравнений (2.1), полагая 
$y^{\star}=1.$
\vspace{1ex}
В системе координат 
$O^{(2)}_{\phantom{1}}\eta\zeta$
у точки $M_{\phantom{1}}^{(2)}$ абсцисса 
$\eta={}-z^{\star}=\dfrac{1}{y}\,,$ а ордината $\zeta=x^{\star}=\dfrac{x}{y}\,.$
\vspace{1ex}
Таким образом, лежащей в конечной части $(x,y)$ проективной плоскости $\R \P (x,y)$ точке $M(x,y)$
\vspace{0.75ex}
с ординатой $y\ne 0$ на 
координатной плоскости 
$O^{(2)}_{\phantom{1}}\eta\zeta$
соответствует точка
\vspace{0.75ex}
$M_{\phantom{1}}^{(2)}(\eta,\zeta),$
координаты которой че\-рез координаты точки $M(x,y)$ выражаются по формулам
\\[2.25ex]
\mbox{}\hfill       % (4.4)
$
\eta=\dfrac{1}{y}\,,
\quad \ \ 
\zeta=\dfrac{x}{y}\,.
$
\hfill (4.4)
\\[2ex]
\indent
Тем самым установлено биективное отображение
\\[2ex]
\mbox{}\hfill       % (4.5)
$
P^{(2)}_{}\colon
(x,y)\to \ 
\Bigl(\,\dfrac{1}{y}\,,\,\dfrac{x}{y}\;\!\Bigr)
\quad
\forall (x,y)\in\R^2\backslash\{(x,y)\colon y=0\}
$
\hfill (4.5)
\\[2ex]
плоскости $Oxy,$
\vspace{0.5ex}
из которой удалена  ось $\!Ox,\!$
на плоскость 
$\!O^{(2)}_{\phantom{1}}\eta\zeta,\!$
из которой удалена  ось
$O^{(2)}_{\phantom{1}}\zeta.$
Отображение (4.5) назовем {\it вторым отображением Пуанкаре}
\vspace{0.75ex}
плоскости $Oxy.$

Отображение  (4.5) является диффеоморфизмом с якобианом
\\[2.25ex]
\mbox{}\hfill
$
\dfrac{{\sf D}(\eta,\zeta)}{{\sf D}(x,y)}=\dfrac{1}{y^3}\ne 0
\quad
\forall (x,y)\in\R^2\backslash \{(x,y)\colon y=0\}.
\hfill
$
\\[2.25ex]
\indent
Разрешив равенства (4.4) относительно $x$ и $y,$ получим формулы
\\[2ex]
\mbox{}\hfill       % (4.6)
$
x=\dfrac{\zeta}{\eta}\,,
\quad \ \ 
y=\dfrac{1}{\eta}\,,
$
\hfill (4.6)
\\[2.25ex]
по которым координаты точки $M(x,y)$ 
\vspace{0.35ex}
с ординатой $y\ne 0$ выражаются через координаты точки 
$M_{\phantom{1}}^{(2)}(\eta,\zeta)$ при $\eta\ne 0.$
\vspace{0.75ex}
Формулы (4.6) назовем
{\it вторым преобразованием Пуанкаре}
плоскости $Oxy$  [1, с. 31].
\vspace{0.5ex}

Функции (4.4) есть 
\vspace{0.5ex}
функции перехода от координат $(\eta,\zeta)$ к координатам $(x,y),$
а функции (4.6) есть функции перехода от координат $(x,y)$ к координатам $(\eta,\zeta).$
\vspace{0.75ex}

Произвольным 
\vspace{0.5ex}
образом возьмем прямую $x=b\;\!y$ с параметром $b\in\R$ и преобразуем ее по формулам (4.6).
В результате получим прямую $\zeta=b$ на плоскости 
\vspace{0.5ex}
$O^{(2)}_{\phantom{1}}\eta\zeta.$
При этом бесконечно удаленной точке $L,$ лежащей на <<концах>> прямой 
\vspace{0.5ex}
$x=b\;\!y,$ 
соответствует точка 
$L_{\phantom{1}}^{(2)}(0,b),$ лежащая на оси 
$O^{(2)}_{\phantom{1}}\zeta$
плоскости $O^{(2)}_{\phantom{1}}\eta\zeta.$
\vspace{0.5ex}
И, наоборот, каждой точке 
$L_{\phantom{1}}^{(2)}(0,b),$ лежащей на координатной оси 
$O^{(2)}_{\phantom{1}}\zeta$
\vspace{0.5ex}
плоскости $O^{(2)}_{\phantom{1}}\eta\zeta,$
соответствует бесконечно удаленная точка 
$L,$ лежащая на <<концах>> прямой 
$x=b\;\!y.$ 
\vspace{0.5ex}

Стало быть справедливо

{\bf Свойство 4.2.} 
\vspace{0.5ex}
{\it
Второе преобразование Пуанкаре {\rm(4.6)} устанавливает 
диффеоморфное отображение проективной плоскости $\R \P (x,y),$
\vspace{0.5ex}
у которой удалена прямая $y=0,$ на совмещенную 
координатную плоскость
$O^{(2)}_{\phantom{1}}\eta\zeta.$
\vspace{0.5ex}
При этом 
образом бесконечно удаленной прямой
проективной плоскости $\R \P (x,y),$
\vspace{0.5ex}
у которой удалена точка, 
лежащая на <<концах>> 
\vspace{0.5ex}
прямой $y=0,$ является координатная ось
$O^{(2)}_{\phantom{1}}\zeta$
$(\!$прямая $\eta=0)$ плоскости 
$O^{(2)}_{\phantom{1}}\eta\zeta.$
}
\vspace{1ex}

{\bf Теорема 4.1.} 
\vspace{0.5ex}
{\it
Тождественное отображение $I\colon (x,y)\to (x,y)\;\; \forall (x,y)\in\R^2,$
первое {\rm(4.2)} и вто\-рое {\rm(4.5)}
отображения Пуанкаре образуют группу третьего порядка}:
\\[2.25ex]
\mbox{}\hfill                                   %(4.7)
$
I\circ I=I,
\quad 
P^{(1)}_{}\circ I=I\circ P^{(1)}_{}=P^{(1)}_{},
\quad
P^{(2)}_{}\circ I=I\circ P^{(2)}_{}=P^{(2)}_{},
\hfill        
$
\\[-0.75ex]
\mbox{}\hfill        
\hfill (4.7)       
\\[1ex]
\mbox{}\hfill        
$
P^{(1)}_{}\circ P^{(2)}_{}=P^{(2)}_{}\circ P^{(1)}_{}=I,
\quad
P^{(1)}_{}\circ P^{(1)}_{}=P^{(2)}_{},
\quad
P^{(2)}_{}\circ P^{(2)}_{}=P^{(1)}_{}.
\hfill
$
\\[2ex]
\indent
{\sl Действительно,} первое и второе преобразования Пуанкаре
являются взаимно обратными:
\\[2ex]
\mbox{}\hfill
$
P^{(1)}_{}\circ P^{(2)}_{}=(x,y)\stackrel{\stackrel{\scriptstyle P^{(2)}_{}}
{\mbox{}}}{\rightarrow}\Bigl(\dfrac{1}{y}\,,\dfrac{x}{y}\Bigr)
\stackrel{\stackrel{\scriptstyle P^{(1)}_{}}{\mbox{}}}
{\rightarrow}\Bigl(\dfrac{x}{y}\,:\dfrac{1}{y}\,,1:
\dfrac{1}{y}\Bigr)=(x,y)=I;
\hfill
$
\\[2.75ex]
\mbox{}\hfill
$
P^{(2)}_{}\circ P^{(1)}_{}=(x,y)\stackrel{\stackrel{\scriptstyle P^{(1)}_{}}
{\mbox{}}}{\rightarrow}\Bigl(\dfrac{y}{x}\,,\dfrac{1}{x}\Bigr)
\stackrel{\stackrel{\scriptstyle P^{(2)}_{}}{\mbox{}}}
{\rightarrow}\Bigl(1:\dfrac{1}{x}\,,\dfrac{y}{x}\,:
\dfrac{1}{x}\Bigr)=(x,y)=I.
\hfill
$
\\[2.25ex]
\indent
Кроме этого,
\\[2ex]
\mbox{}\hfill        
$
P^{(1)}_{}\circ P^{(1)}_{}=(x,y)\stackrel{\stackrel{\scriptstyle P^{(1)}_{}}{\mbox{}}}
{\rightarrow}\Bigl(\dfrac{y}{x}\,,\dfrac{1}{x}\Bigr)
\stackrel{\stackrel{\scriptstyle P^{(1)}_{}}{\mbox{}}}
{\rightarrow}\Bigl(\dfrac{1}{x}:\dfrac{y}{x}\,,1:
\dfrac{y}{x}\Bigr)=\Bigl(\dfrac{1}{y}\,,\dfrac{x}{y}\Bigr)=P^{(2)}_{},
\hfill 
$
\\[2.75ex]
\mbox{}\hfill      
$
P^{(2)}_{}\circ P^{(2)}_{}=(x,y)\stackrel{\stackrel{\scriptstyle P^{(2)}_{}}{\mbox{}}}
{\rightarrow}\Bigl(\dfrac{1}{y}\,,\dfrac{x}{y}\Bigr)
\stackrel{\stackrel{\scriptstyle P^{(2)}_{}}{\mbox{}}}{\rightarrow}
\Bigl(1:\dfrac{x}{y}\,,\dfrac{1}{y}:\dfrac{x}{y}\Bigr)=
\Bigl(\dfrac{y}{x}\,,\dfrac{1}{x}\Bigr)=P^{(1)}_{}.
\hfill 
$
\\[2ex]
\indent
Вполне очевидно, что
\\[2ex]
\mbox{}\hfill
$
P^{(1)}_{}\circ I=I\circ P^{(1)}_{}=P^{(1)}_{},
\quad
P^{(2)}_{}\circ I=I\circ P^{(2)}_{}=P^{(2)}_{},
\quad
I\circ I=I.\
\k
\hfill
$
\\[4.25ex]
\centerline{
{\bf  5. Атлас проективных кругов проективной плоскости}
}
\\[1.5ex]
\indent
Возьмем атлас карт $(U_\tau^{}, \varphi_\tau^{}),\ \tau=1,2,3,$ сферы Пуанкаре 
\vspace{0.75ex}
$\P{\mathbb S}(x,y).$
На локальных прямоугольных декартовых системах координат 
\vspace{0.75ex}
$Oxy,\ O^{(1)}\xi\theta,\ O^{(2)}\eta\zeta$ построим (рис.~5.1) проективные круги
$\P\K(x,y),\, \P\K(\xi,\theta),\, \P\K(\eta,\zeta).$
\\[4.75ex]
\mbox{}\hfill
{\unitlength=1mm
\begin{picture}(42,42)
\put(0,0){\includegraphics[width=42mm,height=42mm]{r05-01.eps}}

\put(18,41){\makebox(0,0)[cc]{ $y$}}
\put(40.2,18.2){\makebox(0,0)[cc]{ $x$}}

\put(28.5,22.3){\makebox(0,0)[cc]{\scriptsize $1$}}
\put(31.3,31.5){\makebox(0,0)[cc]{\scriptsize $2$}}
\put(22.3,29){\makebox(0,0)[cc]{\scriptsize $3$}}

\put(19,29){\makebox(0,0)[cc]{\scriptsize $4$}}
\put(9.7,31.5){\makebox(0,0)[cc]{\scriptsize $5$}}
\put(12,22.3){\makebox(0,0)[cc]{\scriptsize $6$}}

\put(12,18.6){\makebox(0,0)[cc]{\scriptsize $7$}}
\put(9.3,9.5){\makebox(0,0)[cc]{\scriptsize $8$}}
\put(19,12.3){\makebox(0,0)[cc]{\scriptsize $9$}}

\put(22.3,12.3){\makebox(0,0)[cc]{\scriptsize $10$}}
\put(31.5,9.5){\makebox(0,0)[cc]{\scriptsize $11$}}
\put(28.4,18.8){\makebox(0,0)[cc]{\scriptsize $12$}}

\put(28,28){\makebox(0,0)[cc]{\scriptsize $M_1^{}$}}
\put(14,28.8){\makebox(0,0)[cc]{\scriptsize $M_2^{}$}}
\put(12.5,13.5){\makebox(0,0)[cc]{\scriptsize $M_3^{}$}}
\put(28.8,13.5){\makebox(0,0)[cc]{\scriptsize $M_4^{}$}}

%\put(22.5,-6){\makebox(0,0)[cc]{\rm $K(x,y)$}}
\end{picture}
}
\quad
{\unitlength=1mm
\begin{picture}(42,42)
\put(0,0){\includegraphics[width=42mm,height=42mm]{r05-01.eps}}

\put(18,41){\makebox(0,0)[cc]{ $\theta$}}
\put(40.2,17.8){\makebox(0,0)[cc]{ $\xi$}}

\put(28.5,22.3){\makebox(0,0)[cc]{\scriptsize $2$}}
\put(31.3,31.5){\makebox(0,0)[cc]{\scriptsize $3$}}
\put(22.3,29){\makebox(0,0)[cc]{\scriptsize $1$}}

\put(18.8,29){\makebox(0,0)[cc]{\scriptsize $12$}}
\put(9.5,31.5){\makebox(0,0)[cc]{\scriptsize $10$}}
\put(12,22.1){\makebox(0,0)[cc]{\scriptsize $11$}}

\put(12,18.6){\makebox(0,0)[cc]{\scriptsize $5$}}
\put(9.3,9.5){\makebox(0,0)[cc]{\scriptsize $4$}}
\put(19,12.3){\makebox(0,0)[cc]{\scriptsize $6$}}

\put(22.3,12.3){\makebox(0,0)[cc]{\scriptsize $7$}}
\put(31.5,9.5){\makebox(0,0)[cc]{\scriptsize $9$}}
\put(28.4,18.8){\makebox(0,0)[cc]{\scriptsize $8$}}

\put(28,28){\makebox(0,0)[cc]{\scriptsize $M_1^{}$}}
\put(14,28.8){\makebox(0,0)[cc]{\scriptsize $M_4^{}$}}
\put(12.5,13.5){\makebox(0,0)[cc]{\scriptsize $M_2^{}$}}
\put(28.8,13.5){\makebox(0,0)[cc]{\scriptsize $M_3^{}$}}

%\put(22.5,-6){\makebox(0,0)[cc]{\rm $K(\xi,\theta)$}}
\put(21,-6){\makebox(0,0)[cc]{\rm Рис. 5.1}}
\end{picture}
}
\quad
{\unitlength=1mm
\begin{picture}(42,42)
\put(0,0){\includegraphics[width=42mm,height=42mm]{r05-01.eps}}

\put(18,41){\makebox(0,0)[cc]{ $\zeta$}}
\put(40.2,18){\makebox(0,0)[cc]{ $\eta$}}

\put(28.5,22.3){\makebox(0,0)[cc]{\scriptsize $3$}}
\put(31.3,31.5){\makebox(0,0)[cc]{\scriptsize $1$}}
\put(22.3,29){\makebox(0,0)[cc]{\scriptsize $2$}}

\put(19,29){\makebox(0,0)[cc]{\scriptsize $8$}}
\put(9.7,31.5){\makebox(0,0)[cc]{\scriptsize $7$}}
\put(12,22.3){\makebox(0,0)[cc]{\scriptsize $9$}}

\put(12,18.7){\makebox(0,0)[cc]{\scriptsize $10$}}
\put(9.3,9.7){\makebox(0,0)[cc]{\scriptsize $12$}}
\put(18.8,12.3){\makebox(0,0)[cc]{\scriptsize $11$}}

\put(22.3,12.3){\makebox(0,0)[cc]{\scriptsize $5$}}
\put(31.5,9.5){\makebox(0,0)[cc]{\scriptsize $6$}}
\put(28.4,18.8){\makebox(0,0)[cc]{\scriptsize $4$}}

\put(28,28){\makebox(0,0)[cc]{\scriptsize $M_1^{}$}}
\put(14,28.8){\makebox(0,0)[cc]{\scriptsize $M_3^{}$}}
\put(12.5,13.5){\makebox(0,0)[cc]{\scriptsize $M_4^{}$}}
\put(28.8,13.5){\makebox(0,0)[cc]{\scriptsize $M_2^{}$}}

%\put(22.5,-6){\makebox(0,0)[cc]{\rm $K(\eta,\zeta)$}}
\end{picture}
}
\hfill\mbox{}
\\[-7ex]
%\indent

\newpage 

Групповое свойство (теорема 4.1) 
\vspace{0.35ex}
отображений Пуанкаре (4.2) и (4.5) позволяет установить соответствия между 
проективными кругами 
\vspace{0.35ex}
$\P\K(x,y),\, \P\K(\xi,\theta),\, \P\K(\eta,\zeta),$ которые показаны на рис. 5.1. С помощью точек
$M_1^{},\,M_2^{},\,M_3^{},\,M_4^{}$
\vspace{0.15ex}
показано в какой откры\-той координатной четверти проективных кругов они лежат при переходе от одного 
проек\-тив\-но\-го круга к другому. 
Числами $1,\ldots,12$ отражены соответствия между полу\-ок\-рес\-т\-нос\-тя\-ми точек, 
лежащих на границах координатных четвертей 
\vspace{0.5ex}
этих проективных кругов.

Упорядоченную тройку $(\P\K(x,y),\, \P\K(\xi,\theta),\, \P\K(\eta,\zeta))$ назовем 
\vspace{0.5ex}
{\it атласом проективных кругов} проективной плоскости $\R\P(x,y).$ 
\vspace{0.5ex}
Тогда упорядоченная тройка $(\P\K(\xi,\theta),\, \P\K(\eta,\zeta),\, \P\K(x,y))$ --- атлас 
\vspace{0.5ex}
проективных кругов проективной плоскости $\R \P(\xi,\theta),$ 
\vspace{0.5ex}
а упорядоченная тройка $(\P\K(\eta,\zeta),\, \P\K(x,y),\, \P\K(\xi,\theta))$ --- 
атлас проективных кругов проективной плоскости $\R \P(\eta,\zeta).$
\\[3.75ex]
\centerline{
{\bf\large 
\S\;\!2. Преобразования Пуанкаре дифференциальных систем}
}
\\[2.25ex]
\centerline{
{\bf  6. Проективно приведенные системы
}
}
\\[1.5ex]
\indent
Первым преобразованием Пуанкаре (4.3) дифференциальную систему (D) приводим
к системе
\\[2ex]
\mbox{}\hfill                             % (6.1)
$
\begin{array}{l}
\dfrac{d\;\!\xi}{dt}\, ={}-\xi\;\! \theta\;\! 
X\!\Bigl(\;\!\dfrac{1}{\theta}\,,\;\!\dfrac{\xi}{\theta}\;\!\Bigr) +\, 
\theta\;\! Y\!\Bigl(\;\!\dfrac{1}{\theta}\,,\;\!\dfrac{\xi}{\theta}\;\!\Bigr) \equiv\;\!
\Xi_{\phantom1}^{(1)}(\xi,\theta),
\\[4.25ex]
\dfrac{d\;\!\theta}{dt}={}-\theta^{\;\!2}
X\Bigl(\;\!\dfrac{1}{\theta}\,,\;\!\dfrac{\xi}{\theta}\;\!\Bigr)\equiv\;\! 
\Theta_{\phantom1}^{(1)}(\xi,\theta).
\end{array}
$
\hfill (6.1)
\\[2.5ex]
\indent
Поскольку $X$ и $Y$ --- полиномы, то систему (6.1) можно записать в виде
\\[2ex]
\mbox{}\hfill                               %(6.2)
$
\dfrac{d\;\!\xi}{dt}\;\!=\;\!\dfrac{\Xi(\xi,\theta)}{\theta^{m}}\,,
\qquad
\dfrac{d\;\!\theta}{dt}\;\!=\;\!\dfrac{\Theta(\xi,\theta)}{\theta^{m}}\,,
$
\hfill (6.2)
\\[2.25ex]
где $\Xi$ и $\Theta$ --- полиномы, не делящиеся
\vspace{0.25ex}
одновременно на $\theta,$ а число $m$ --- целое неотрицательное.
На основании системы (6.2) составим систему
\\[2ex]
\mbox{}\hfill                                     %(6.3)
$
\dfrac{d\;\!\xi}{d\tau}=\Xi(\xi,\theta),
\qquad
\dfrac{d\;\!\theta}{d\tau}=\Theta(\xi,\theta),
$
\hfill (6.3)
\\[2.25ex]
где $\theta^{m}\;\!d\tau=dt,$ у которой правые части $\Xi$ и $\Theta$ суть взаимно простые полиномы.
\vspace{0.5ex}

Вторым преобразованием Пуанкаре (4.6)
дифференциальную систему (D) приводим
к системе
\\[1.75ex]
\mbox{}\hfill                     %(6.4)
$
\begin{array}{l}
\dfrac{d\;\!\eta}{dt}={}-\eta^2\,Y\Bigl(\;\!\dfrac{\zeta}{\eta}\,,\;\!\dfrac{1}{\eta}\;\!\Bigr)
\equiv H_{\phantom1}^{(2)}(\eta,\zeta),
\\[4.25ex]
\dfrac{d\;\!\zeta}{dt}=\eta\,X\Bigl(\;\!\dfrac{\zeta}{\eta}\,,\;\!\dfrac{1}{\eta}\;\!\Bigr)
-\eta\;\!\zeta\,Y\Bigl(\;\!\dfrac{\zeta}{\eta}\,,\;\!\dfrac{1}{\eta}\;\!\Bigr) \equiv
Z_{\phantom1}^{(2)}(\eta,\zeta).
\end{array}
$
\hfill (6.4)
\\[2.5ex]
\indent
Поскольку $X$ и $Y$ --- полиномы, то систему (6.4) можно записать
в виде
\\[2ex]
\mbox{}\hfill                    %(6.5)
$
\dfrac{d\eta}{dt}\;\!=\;\!\dfrac{H(\eta,\zeta)}{\eta^{m}}\,,
\qquad
\dfrac{d\zeta}{dt}\;\!=\;\!\dfrac{Z(\eta,\zeta)}{\eta^{m}}\,,
$
\hfill (6.5)
\\[2ex]
где $H$ и $Z$ --- полиномы, не делящиеся
\vspace{0.25ex}
одновременно на $\eta,$ а число $m$ --- целое неотрицательное.
На основании системы (6.5)
составим систему
\\[1.75ex]
\mbox{}\hfill                    %(6.6)
$
\dfrac{d\eta}{d\nu}=H(\eta,\zeta),
\qquad
\dfrac{d\zeta}{d\nu} =Z(\eta,\zeta),
$
\hfill (6.6)
\\[2ex]
где $\eta^{m}\;\!d\nu=dt,$ у которой правые части $H$ и $Z$ суть взаимно простые полиномы.

Полученные на основании  систем  (6.1) и (6.4) автономные полиномиальные дифференциальные системы (6.3) и (6.6) назовем 
{\it проективно приведенными системами} или P-{\it приведенными системами} системы (D). 

Систему (6.3) назовем {\it первой проективно приведенной системой} или (P-1)-{\it при\-веденной системой} системы (D), 
а систему (6.6) --- {\it второй проективно приведенной сис\-темой} или
(P-2)-{\it приведенной системой} системы (D). 
\vspace{0.5ex}

Плоскости $(\xi,\theta)$ и $(\eta,\zeta)$ 
\vspace{0.25ex}
будут фазовыми плоскостями P-приведенных систем (6.3) и (6.6) соответственно.
\\[3.25ex]
\mbox{}\hfill
{\unitlength=1mm
\begin{picture}(35,26)
\put(0,0){\includegraphics[width=34.78mm,height=25.02mm]{r06-01.eps}}

\put(4.6,4.3){\makebox(0,0)[cc]{\scriptsize (6.3)}}
\put(30.5,4.3){\makebox(0,0)[cc]{\scriptsize (6.6)}}
\put(17.5,20.5){\makebox(0,0)[cc]{\scriptsize (D)}}

\put(8,14.2){\makebox(0,0)[cc]{\scriptsize $P^{(1)}_{}$}}
\put(28.4,14.2){\makebox(0,0)[cc]{\scriptsize $P^{(2)}_{}$}}
\put(17.5,6.5){\makebox(0,0)[cc]{\scriptsize $P^{(1)}_{}$}}

\put(18,-6){\makebox(0,0)[cc]{\rm Рис. 6.1}}
\end{picture}}
\qquad
{\unitlength=1mm
\begin{picture}(35,26)
\put(0,0){\includegraphics[width=34.78mm,height=25.02mm]{r06-02.eps}}

\put(4.6,4.3){\makebox(0,0)[cc]{\scriptsize (6.3)}}
\put(30.5,4.3){\makebox(0,0)[cc]{\scriptsize (6.6)}}
\put(17.5,20.5){\makebox(0,0)[cc]{\scriptsize (D)}}

\put(8,14.2){\makebox(0,0)[cc]{\scriptsize $P^{(1)}_{}$}}
\put(28.4,14.2){\makebox(0,0)[cc]{\scriptsize $P^{(2)}_{}$}}
\put(17.5,6.5){\makebox(0,0)[cc]{\scriptsize $P^{(2)}_{}$}}

\put(18,-6){\makebox(0,0)[cc]{\rm Рис. 6.2}}
\end{picture}}
\qquad
{\unitlength=1mm
\begin{picture}(61.3,8.5)
\put(0,9){\includegraphics[width=61.3mm,height=8.5mm]{r06-03.eps}}

\put(4.6,13.3){\makebox(0,0)[cc]{\scriptsize (6.3)}}
\put(57,13.3){\makebox(0,0)[cc]{\scriptsize (6.6)}}
\put(30.6,13.3){\makebox(0,0)[cc]{\scriptsize (D)}}

\put(17.5,20){\makebox(0,0)[cc]{\scriptsize $P^{(1)}_{}$}}
\put(44,5.7){\makebox(0,0)[cc]{\scriptsize $P^{(1)}_{}$}}
\put(44,20){\makebox(0,0)[cc]{\scriptsize $P^{(2)}_{}$}}
\put(17.5,5.7){\makebox(0,0)[cc]{\scriptsize $P^{(2)}_{}$}}

\put(30.5,-6){\makebox(0,0)[cc]{\rm Рис. 6.3}}
\end{picture}}
\hfill\mbox{}
\\[7.75ex]
\indent
Групповым свойством (теорема 4.1) 
\vspace{0.35ex}
отображений Пуанкаре устанавливаем связи между системами (D), (6.3), (6.6). 
\vspace{0.75ex}

Схема последовательности  $P^{(1)}_{}\circ P^{(1)}_{}=P^{(2)}_{}$
\vspace{0.75ex}
отображений Пуанкаре 
систем (D), (6.3), (6.6) изображена на рис. 6.1,
\vspace{0.75ex}
схема последовательности 
$P^{(2)}_{}\circ P^{(2)}_{}=P^{(1)}_{}$
 --- на рис. 6.2, 
а схемы последовательностей 
$P^{(1)}_{}\circ P^{(2)}_{}=I$ и 
$P^{(2)}_{}\circ P^{(1)}_{}=I$  
--- на рис. 6.3. 
\vspace{1ex}

{\bf Свойство 6.1.} 
\vspace{0.75ex}
{\it 
Первая проективно приведенная система {\rm (6.3)} с помощью первого преобразования Пуанкаре 
$
\xi=\dfrac{1}{\zeta}\,,\ 
\theta=\dfrac{\eta}{\zeta}
$
приводится к второй проективно приведенной системе {\rm (6.6)}. 
}
\vspace{1.25ex}

{\bf Свойство 6.2.} 
\vspace{0.75ex}
{\it 
Первая проективно приведенная система {\rm (6.3)} с помощью второго преобразования Пуанкаре 
$
\xi=\dfrac{y}{x}\,,\ 
\theta=\dfrac{1}{x}
$
приводится к системе {\rm (D)}. 
}
\vspace{1.25ex}

{\bf Свойство 6.3.} 
\vspace{0.75ex}
{\it 
Вторая проективно приведенная система {\rm (6.6)} с помощью первого преобразования Пуанкаре 
$
\eta=\dfrac{1}{y}\,,\ 
\zeta=\dfrac{x}{y}
$
приводится к системе {\rm (D)}. 
}
\vspace{1.25ex}

{\bf Свойство 6.4.} 
\vspace{0.75ex}
{\it 
Вторая проективно приведенная система {\rm (6.6)} с помощью второго преобразования Пуанкаре 
$
\eta=\dfrac{\theta}{\xi}\,,\ 
\zeta=\dfrac{1}{\xi}
$
приводится к первой проективно приведенной системе {\rm (6.3)}. 
}
\\[4.25ex]
\centerline{
{\bf  7.
Проективный тип дифференциальной системы 
}
}
\\[1.5ex]
\indent
Для дифференциальной системы (D) 
\vspace{0.25ex}
вид первой и второй проективно приведенных систем (6.3) и (6.6) зависит от того, является ли полином
\\[2ex]
\mbox{}\hfill                           
$
W_n^{}\colon (x,y)\to\  x\;\!Y_n^{}(x, y)\;\!-\;\!y\;\!X_n^{}(x, y)
\quad \forall (x,y)\in \R^2
\hfill
$
\\[2ex]
тождественным нулем на плоскости $\R^2$ или нет. 
\vspace{0.5ex}

Если $W_n^{}(x,y) \not\equiv 0$ на $\R^2,$ то первая проективно приведенная система (6.3) будет иметь вид
\\[2ex]
\mbox{}\hfill                            % (7.1)
$
\begin{array}{l}
\displaystyle
\dfrac{d\;\!\xi}{d\tau}\,=\ \sum\limits_{i=0}^{n}\,\theta^{n-i}\,Y_{i}^{}(1,\xi)-
\xi\;\!\sum\limits_{i=0}^{n}\,\theta^{n - i}X_{i}^{}(1, \xi)
\,\equiv\,
\sum\limits_{i=0}^{n}\,\theta^{n - i}\,W_{i}^{}(1, \xi)
\equiv\;\!
\widetilde{\Xi}(\xi,\theta),
\\[4.75ex]
\displaystyle
\dfrac{d\;\!\theta}{d\tau}\, =
{}-\theta\;\! \sum\limits_{i=0}^{n}\,\theta^{n - i}X_{i}^{}(1, \xi)
\equiv\;\!
\widetilde{\Theta}(\xi,\theta),
\end{array}
$
\hfill (7.1)
\\[2.5ex]
где 
\vspace{0.75ex}
$\theta^{n-1}\;\!d\tau=dt,\ \, W_i^{}(x,y)=x\;\!Y_i^{}(x,y)-yX_i^{}(x,y)\;\; \forall (x,y)\in\R^2,\ i=0,1,\ldots,n,$ причем
$W_n^{}(1,\xi)\not\equiv0$ на $\R,$ а вторая проективно приведенная система (6.6) будет иметь вид
\\[2.5ex]
\mbox{}\hfill                             % (7.2)
$
\begin{array}{l}
\displaystyle
\dfrac{d\;\!\eta}{d\nu}\,=
{}-\eta\;\!\sum\limits_{i=0}^{n}\, \eta^{n - i}\, Y_{i}^{}(\zeta,1)
\equiv\;\!
\widetilde{H}(\eta,\zeta),
\\[4.75ex]
\displaystyle
\dfrac{d\;\!\zeta}{d\nu}\, =\
\sum\limits_{i = 0}^{n}\, \eta^{n-i}X_{i}^{}(\zeta,1)\;\! -\;\!
\zeta\;\!\sum\limits_{i = 0}^{n}\,\eta^{n - i}\, Y_{i}^{}(\zeta,1)
\;\!\equiv
{}-\sum\limits_{i = 0}^{n}\, \eta^{n - i}\,W_{i}^{}(\zeta,1)
\equiv\;\!
\widetilde{Z}(\eta,\zeta),
\end{array}
$
\hfill (7.2)
\\[2.75ex]
где $\eta^{n-1}\,d\nu=dt,$ причем $W_n^{}(\zeta,1)\not\equiv 0$ на поле $\R.$
\vspace{1.25ex}

Если $W_n^{}(x,y) \equiv 0\ \forall (x,y)\in\R^2,$ 
\vspace{0.5ex}
то первая проективно приведенная система (6.3) будет иметь вид
\\[2.25ex]
\mbox{}\hfill                           % (7.3)
$
\begin{array}{l}
\displaystyle
\dfrac{d\;\!\xi}{d\tau}\,=\
\sum\limits_{j=0}^{n-1}\,\theta^{n-j-1}\,Y_{j}^{}(1, \xi)\;\!-\;\!
\xi\;\!\sum\limits_{j = 0}^{n-1}\,\theta^{n-j-1}X_{j}^{}(1, \xi)
\;\!\equiv\;\!
\sum\limits_{j = 0}^{n-1}\,\theta^{n-j-1}\,W_{j}^{}(1, \xi)
\equiv\;\!
\widehat{\Xi}(\xi,\theta),
\\[5ex]
\displaystyle
\dfrac{d\;\!\theta}{d\tau}\,=
{}-\sum\limits_{i = 0}^{n}\,\theta^{n - i}X_{i}^{}(1, \xi)
\equiv\;\!
\widehat{\Theta}(\xi,\theta),
\end{array}
$
\hfill (7.3)
\\[2ex]
где  $\theta^{n-2}\,d\tau=dt,$ причем $X_n^{}(1,\xi) \not\equiv0$ на поле $\R,$ 
\vspace{0.5ex}
а вторая проективно приведенная система (6.6) будет иметь вид
\\[2ex]
\mbox{}\hfill                             % (7.4)
$
\begin{array}{l}
\displaystyle
\dfrac{d\;\!\eta}{d\nu}\,=
{}-\sum\limits_{i=0}^{n}\,\eta^{n - i}\,Y_{i}^{}(\zeta,1)
\equiv \;\!
\widehat{H}(\eta,\zeta),
\\[5ex]
\displaystyle
\dfrac{d\zeta}{d\nu}\,=\
\sum\limits_{j = 0}^{n - 1}\,\eta^{n-j-1}X_{j}^{}(\zeta,1)\;\!-\;\!
\zeta\;\!\sum\limits_{j = 0}^{n - 1}\, \eta^{n-j-1}\,Y_{j}^{}(\zeta,1)\;\!
\equiv
{}-\sum\limits_{j = 0}^{n - 1}\, \eta^{n-j-1}\,W_{j}^{}(\zeta,1)
\equiv\;\!
\widehat{Z}(\eta,\zeta),
\end{array}
\!\!\!\!\!\!\!\!\!\!
$
\hfill (7.4)
\\[2.5ex]
где $\eta^{n - 2}\,d\nu=dt,$ причем $Y_n^{}(\zeta,1)\not\equiv0$ на поле $\R.$
\vspace{0.75ex}

Если $W_n^{}(x,y) \not\equiv 0$ на $\R^2,$ 
\vspace{0.5ex}
то систему (D) назовем {\it проективно неособой} или  P-{\itнеособой}. 
В противном случае, т.е. когда 
\vspace{0.5ex}
$W_n^{}(x,y)=0\;\;\forall(x,y)\in\R^2,$ систему (D) будем называть {\it проективно особой}
или  P-{\itособой}. Дифференциальные системы (D) делятся на два класса:  P-особые и P-неособые. 
Вид проективно приведенных систем (6.3) и (6.6) зависит от принадлежности системы (D) тому или иному классу.
\vspace{0.5ex}

{\bf Свойство 7.1}.
{\it У проективно неособой системы {\rm(D)} первая проективно приведенная система {\rm(6.3)} имеет вид {\rm(7.1),} а вторая проективно приведенная система {\rm(6.6)} имеет вид} (7.2).
\vspace{0.35ex}

{\bf Свойство 7.2}.
{\it У проективно особой системы {\rm(D)} первая проективно приведенная система {\rm(6.3)} имеет вид {\rm(7.3),} а вторая проективно приведенная система {\rm(6.6)} имеет вид} (7.4).
\vspace{0.35ex}

Число $n=\max\{{\rm deg}\,X,\, {\rm deg}\,Y\}$ назовем
\vspace{0.35ex}
{\it степенью системы} (D) и обозначим ${\rm deg}$(D).

Степени P-приведенных систем (7.1), (7.2), (7.3) и (7.4) 
\vspace{0.25ex}
зависят от степени $n$ системы (D), 
а также от того состоят или нет прямые $x=0$ и $y=0$ из траекторий системы (D). 

Введем три числа $\delta,\ \delta^{(1)}$ и $\delta^{(2)}.$ 
\vspace{0.25ex}
Если система (D)  проективно особая, то $\delta=0;$ а если система (D) проективно неособая, то $\delta=1.$
\vspace{0.25ex}
Если прямая $x=0$ не состоит из траекторий системы (D), то $\delta^{(1)}=0;$ 
\vspace{0.25ex}
а если прямая $x=0$ состоит из траекторий системы (D), то $\delta^{(1)}=1.$
\vspace{0.25ex}
Если прямая $y=0$ не состоит из траекторий системы (D), то $\delta^{(2)}=0;$ 
\vspace{0.25ex}
а если прямая $y=0$ состоит из траекторий системы (D), то $\delta^{(2)}=1.$ 

Тогда имеет место
\vspace{0.35ex}

{\bf Свойство 7.3}.
{\it 
Степень первой проективно приведенной системы} (6.3) {\it находится по формуле}
\\[1.5ex]
\mbox{}\hfill                              %(7.5)
$
{\rm deg}\,(6.3) =n+\delta-\delta^{(1)},
$
\hfill (7.5)
\\[2ex]
{\it а степень второй проективно приведенной системы} (6.6) {\it находится по формуле}
\\[2ex]
\mbox{}\hfill                              %(7.6)
$
{\rm deg}\,(6.6) =n+\delta-\delta^{(2)}.
$
\hfill (7.6)
\\[2ex]
\indent
Степени P-приведенных систем (7.1) --- (7.4) в соответствии с формулами (7.5) и (7.6) находятся следующим образом:
\\[2ex]
\mbox{}\hfill                              %(7.7)
$
{\rm deg}\,(7.1) =n-\delta^{(1)}+1,
\quad \ \
{\rm deg}\,(7.2) =n-\delta^{(2)}+1,
\hfill
$
\\
\mbox{}\hfill (7.7)
\\
\mbox{}\hfill
$
{\rm deg}\,(7.3) =n-\delta^{(1)},
\quad \ \
{\rm deg}\,(7.4) =n-\delta^{(2)}.
\hfill
$
\\[2.5ex]
\indent
{\bf Примеры}
\vspace{0.5ex}

{\bf 7.1.}
Рассмотрим автономную систему
\\[2ex]
\mbox{}\hfill                           %(7.8)
$
\dfrac{dx}{dt}=a_{_0}\equiv X(x,y),
\quad \ \ 
\dfrac{dy}{dt}=b_{_{0}}\equiv Y(x,y),
\qquad
|a_{_0}|+|b_{_0}|\ne 0.
$
\hfill (7.8)
\\[2.5ex]
\indent
Условие 
\vspace{0.75ex}
$|a_{_0}|+|b_{_0}|\ne 0$ обосновано лишь тем, что 
$X_{_0}(x,y)=a_{_0},\ Y_{_0}(x,y)=b_{_0},$
а система (D) при $n=0$ должна быть такой, что 
$|X_{_0}(x,y)|+|Y_{_0}(x,y)|\not\equiv0$ на $\R^2.$
\vspace{1ex}

Полином 
\vspace{0.5ex}
$W_{_0}(x,y)=b_{_0}x-a_{_0}y\not\equiv0$ на $\R^2$ при 
$|a_{_0}|+|b_{_0}|\ne 0.$ 
Система (7.8) проективно неособая.
\vspace{0.35ex}

У системы (7.8) первой проективно приведенной системой является система
\\[2ex]
\mbox{}\hfill                           % (7.9)
$
\dfrac{d\xi}{d\tau}=b_{_0}-a_{_0}\xi
\equiv\widetilde{\Xi}(\xi,\theta),
\quad
\dfrac{d\theta}{d\tau}={}-a_{_0}\theta
\equiv\widetilde{\Theta}(\xi,\theta),
$
\ где \  $d\tau=\theta dt,\ |a_{_0}|+|b_{_0}|\ne 0;$ 
\hfill\mbox{} (7.9)
\\[2.35ex]
а второй проективно приведенной системой является система
\\[2.1ex]
\mbox{}\hfill                             % (7.10)
$
\dfrac{d\eta}{d\nu}={}-b_{_0}\eta
\equiv \widetilde{H}(\eta,\zeta),
\quad
\dfrac{d\zeta}{d\nu}=a_{_0}-b_{_0}\zeta
\equiv \widetilde{Z}(\eta,\zeta),
$
\ где \ $d\nu=\eta dt,\ |a_{_0}|+|b_{_0}|\ne 0.$ 
\hfill\mbox{} (7.10)
\\[3.25ex]
\indent
{\bf 7.2.}
Рассмотрим линейную стационарную систему
\\[2ex]
\mbox{}\hfill                                           %(7.11)
$
\dfrac{dx}{dt}  =  a_{_0}  +  a^{}_{1}x  +  a^{}_{2}y\;\!
\equiv X(x,y),
\qquad
\dfrac{dy}{dt} =  b_{_0}  +  b_{1}^{}x  +  b_{2}^{}y\;\!
\equiv Y(x,y),
$
\hfill (7.11)
\\[2.5ex]
с коэффициентами такими, что $|a_{1}^{}| + |a_{2}^{}| +|b_1^{}| + |b_2^{}|\ne 0.$
\vspace{1.25ex}

Полином 
\vspace{1ex}
$W_1^{}(x,y)=b_1^{}x^2+(b_2^{}-a_1^{})xy-a_2^{}y^2\;\; \forall (x,y)\in\R^2.$
Если $|a_2^{}|+|b_1^{}|+|a_1^{}-b_2^{}|\ne 0,$ то система (7.11) P-неособая;
а если $a_2^{}=b_1^{}=0,\ a_1^{}=b_2^{}\ne 0,$ 
\vspace{1ex}
то система (7.11) P-особая.

Пусть
$|a_{2}^{}|+|b_{1}^{}| +|a_{1}^{}-b_{2}^{}|\ne 0,$
\vspace{0.75ex}
т.е. система (7.11) P-неособая.

С помощью первого преобразования Пуанкаре 
\vspace{0.75ex}
$x=\dfrac{1}{\theta}\,,\ y=\dfrac{\xi}{\theta}$
P-не\-о\-со\-бую систему (7.11) приводим к
(P-1)-приведенной системе
\\[2.5ex]
\mbox{}\hfill                                 %(7.12)
$
\begin{array}{l}
\dfrac{d\xi}{dt} = 
b_{1}^{} - (a_{1}^{} - b_{2}^{})\xi + b_{_0}\theta -
a_{2}^{}\xi^{2} - a_{_0}\xi\theta
\,\equiv \widetilde{\Xi} (\xi,\theta),
\\[3.75ex]
\dfrac{d\theta}{dt} = {} - a_{1}^{}\theta - a_{2}^{}\xi\theta - a_{_0}\theta^{2}
\,\equiv  \widetilde{\Theta}(\xi,\theta),
\end{array}
$
\hfill\mbox{} (7.12)
\\[2.25ex]
где $|a_{2}^{}|+|b_{1}^{}| +|a_{1}^{}-b_{2}^{}|\ne 0.$
%\vspace{0.5ex}

С помощью второго пpеобpазования 
Пуанкаре
$x=\dfrac{\zeta}{\eta}\,,\ y=\dfrac{1}{\eta}$
P-не\-о\-со\-бую систему (7.11) приводим к
(P-2)-приведенной системе
\\[2ex]
\mbox{}\hfill                                         %(7.13)
$
\begin{array}{l}
\dfrac{d\eta}{dt} = {} - b_{2}^{}\eta -b_{_0}\eta^{2} -b_{1}^{}\eta \zeta
\, \equiv \widetilde{H}(\eta,\zeta),
\\[3.75ex]
\dfrac{d\zeta}{dt} = a_{2}^{} +a_{_0}\eta+ (a_{1}^{} - b_{2}^{})\zeta
- b_{_0}\eta\zeta - b_{1}^{}\zeta^{2}
\,\equiv \widetilde{Z} (\eta,\zeta),
\end{array}
$
\hfill\mbox{} (7.13)
\\[2.5ex]
где $|a_{2}^{}|+|b_{1}^{}| +|a_{1}^{}-b_{2}^{}|\ne 0.$
\vspace{1.25ex}

Пусть 
\vspace{0.75ex}
$a_{2}^{} =  b_{1}^{} = 0,\ b_{\;\!2}^{}=a_{1}^{}\ne0,$
т.е. система (7.11) P-особая.

С помощью первого преобразования Пуанкаре 
\vspace{0.5ex}
$x=\dfrac{1}{\theta}\,,\ y=\dfrac{\xi}{\theta}$
P-особую систему (7.11) приводим к
(P-1)-приведенной системе
\\[2.25ex]
\mbox{}\hfill                                      %(7.14)
$
\dfrac{d\xi}{d\tau} = b_{_0} -  a_{_0}\xi
\;\!\equiv\;\! 
\widehat{\Xi}(\xi,\theta),
\qquad
\dfrac{d\theta}{d\tau} = {}-a_{1}^{}-a_{_0}\theta
\;\!\equiv\;\! 
\widehat{\Theta}(\xi,\theta),
$
\hfill\mbox{} (7.14)
\\[2.75ex]
где $d\tau=\theta dt,$ коэффициент $a_1^{}\ne 0.$
\vspace{0.5ex}

С помощью второго пpеобpазования Пуанкаре
$x=\dfrac{\zeta}{\eta}\,,\ y=\dfrac{1}{\eta}$
P-особую систему (7.11) приводим к
(P-2)-приведенной системе
\\[2ex]
\mbox{}\hfill                                         %(7.15)
$
\dfrac{d\eta}{d\nu} = {}-a_{1}^{} - b_{_0}\eta
\;\!\equiv\;\! \widehat{H}(\eta,\zeta),
\qquad
\dfrac{d\zeta}{d\nu} = a_{_0} - b_{_0}\zeta
\;\!\equiv\;\! \widehat{Z}(\eta,\zeta),
$
\hfill\mbox{} (7.15)
\\[2.5ex]
где $d\nu=\theta dt,$ коэффициент $a_1^{}\ne 0.$
\vspace{1ex}

{\bf 7.3.}
Рассмотрим автономную квадpатичную систему
\\[2.25ex]
\mbox{}\hfill                                         %(7.16)
$
\begin{array}{l}
\dfrac{dx}{dt} = a_{_0}+a_1^{} x+a_2^{} y+a_3^{}x^2+a_4^{} xy+a_5^{} y^2
\equiv X(x,y),
\\[4ex]
\dfrac{dy}{dt} = 
b_{_0}+b_1^{} x+b_2^{} y+b_3^{}x^2+b_4^{} xy+b_5^{} y^2
\equiv Y(x,y),
\end{array}
$
\hfill\mbox{} (7.16)
\\[2.75ex]
где
$|a_3^{}| +|a_{4}^{}|+|a_{5}^{}| +|b_{3}^{}|+|b_{4}^{}| + |b_{5}^{}|\ne  0.$
\vspace{1.25ex}

Полином 
\vspace{1ex}
$W_2^{}(x,y)=b_3^{}x^3+(b_4^{}-a_3^{})x^2y+(b_5^{}-a_4^{})xy^2 -a_5^{}y^3\;\; \forall (x,y)\in\R^2.$ 

Если 
$|a_5^{}| +|b_{3}^{}|+|a_{3}^{}-b_{4}^{}|\,+
|a_{4}^{}-b_{5}^{}|\ne  0,$ то система (7.16) 
\vspace{0.75ex}
P-неособая, а если
\linebreak
$a_5^{} =b_{3}^{}=0,\ a_{3}^{}=b_{4}^{},\ a_{4}^{}=b_{5}^{},\ |b_4^{}|+|b_5^{}|\ne  0,$ то система (7.16) 
P-особая.
\vspace{0.75ex}

Пусть
\vspace{0.5ex}
$|a_5^{}| +|b_{3}^{}|+|a_{3}^{}-b_{4}^{}|+|a_{4}^{}-b_{5}^{}|\ne  0,$ т.е. система (7.16) 
P-неособая.
У P-неособой системы (7.16)
первая проективно приведенная система имеет вид
\\[2ex]
\mbox{}\hfill                                         %(7.17)
$
\begin{array}{l}
\dfrac{d\;\!\xi}{d\tau} \ =\,  b_{3}^{} + (b_{4}^{} - a_{3}^{})\xi +
b_{1}^{}\theta + (b_{5}^{} - a_{4}^{})\xi^{2} +
(b_{2}^{} - a_{1}^{})\xi\theta\ +
\\[2.25ex]
\mbox{}\quad\  \ 
+\, b_{_0}\theta^{2} - a_{5}^{}\xi^{3} - a_{2}^{}\xi^{2}\theta - a_{_0}\xi\theta^{2}
\;\! \equiv\;\!  \widetilde{\Xi}(\xi,\theta),
\\[3ex]
\dfrac{d\;\!\theta}{d\tau}\, = {} -  a_{3}^{}\theta - a_{4}^{}\xi\theta - a_{1}^{}\theta^{2} -
a_{5}^{}\xi^{2}\theta - a_{2}^{}\xi\theta^{2} - a_{_0}\theta^{3}
\;\!\equiv\;\!  \widetilde{\Theta} (\xi,\theta),
\end{array}
$
\hfill\mbox{} (7.17)
\\[3ex]
где 
\vspace{0.5ex}
$z\,d\tau = dt,\ 
|a_5^{}| +|b_{3}^{}|+|a_{3}^{}-b_{4}^{}|+|a_{4}^{}-b_{5}^{}|\ne  0,$ 
а вторая проективно приведенная система имеет вид
\\[2.5ex]
\mbox{}\hfill                                         %(7.18)
$
\begin{array}{l}
\dfrac{d\;\!\eta}{d\nu}\  = {} -b_{5}^{}\eta - b_{2}^{}\eta^{2}- b_{4}^{}\eta\;\!\zeta  - b_{_0}\eta^{3}-
b_{1}^{}\eta^{2}\zeta-b_{3}^{}\eta\;\!\zeta^{2}
\;\!\equiv \, \widetilde{H}(\eta,\zeta),
\\[4ex]
\dfrac{d\zeta}{d\nu}\,  = \, a_{5}^{} + a_{2}^{}\eta
+ (a_{4}^{} - b_{5}^{})\zeta
+ a_{_0}\eta^{2}
+(a_{1}^{} - b_{2}^{})\eta\;\!\zeta  
+  (a_{3}^{} - b_{4}^{})\zeta^{2}\ -
\\[2.5ex]
\mbox{}\quad\ \
-\, b_{_0}\eta^{2}\;\!\zeta-b_{1}^{}\eta\;\!\zeta^{2}
- b_{3}^{}\zeta^{3}
\;\!\equiv \, \widetilde{Z} (\eta,\zeta),
\end{array}
$
\hfill\mbox{} (7.18)
\\[2.5ex]
где $z\,d\nu = dt,\ 
|a_5^{}| +|b_{3}^{}|+|a_{3}^{}-b_{4}^{}|+|a_{4}^{}-b_{5}^{}|\ne  0.$ 
\vspace{1ex}

Пусть 
\vspace{0.75ex}
$a_{5}^{} = b_{3}^{}=0,\ b_{4}^{} =  a_{3}^{},\ b_{5}^{} =  a_4^{},\ |a_{3}^{}|+|a_{4}^{}|\ne0,$ т.е. система (7.16) P-особая. 

У P-особой си\-с\-те\-мы (7.16) первая проективно приведенная система имеет вид
\\[2ex]
\mbox{}\hfill                                         %(7.19)
$
\begin{array}{l}
\dfrac{d\;\!\xi}{dt}\;\! = b_{1}^{}+(b_{2}^{}- a_{1}^{})\xi + b_{_0}\;\!\theta -
a_{2}^{}\;\!\xi^2 - a_{_0}\xi\;\!\theta
\, \equiv \, \widehat{\Xi}(\xi,\theta),
\\[4ex]
\dfrac{d\;\!\theta}{dt}\;\! ={}- a_{3}^{} - a_{4}^{}\;\!\xi -
a_{1}^{}\;\!\theta  -  a_{2}^{}\;\!\xi\;\!\theta  - a_{_0}\;\!\theta^2
\, \equiv\;\! \widehat{\Theta}(\xi,\theta),
\ \ 
\text{где} \,\  |a_3^{}|+|a_4^{}|\ne 0;
\end{array}
$
\hfill\mbox{} (7.19)
\\[2.75ex]
а вторая проективно приведенная система имеет вид
\\[2.5ex]
\mbox{}\hfill                                         %(7.20)
$
\begin{array}{l}
\dfrac{d\;\!\eta}{dt}\;\! = 
{}- a_{4}^{}  - b_{2}^{}\;\!\eta - a_{3}^{}\;\!\zeta - b_{_0}\;\!\eta^2
-b_{1}^{}\;\!\eta\;\!\zeta
\;\!\equiv\;\! \widehat{H}(\eta,\zeta),
\\[4ex]
\dfrac{d\;\!\zeta}{dt}\;\! =\;\! a_{2}^{} + a_{_0}\;\!\eta + (a_{1}^{} - b_{2}^{})\zeta
- b_{_0}\eta\;\!\zeta - b_{1}^{}\zeta^2
\;\!\equiv\;\! 
\widehat{Z} (\eta,\zeta),
\ \ 
\text{где} \,\  |a_3^{}|+|a_4^{}|\ne 0.
\end{array}
$
\hfill\mbox{} (7.20)
\\[2.75ex]
\indent
{\bf 7.4.}
Дифференциальную систему
\\[2ex]
\mbox{}\hfill                                         %(7.21)
$
\begin{array}{l}
\dfrac{dx}{dt}=a_{_0}+a_1^{}x+a_2^{}y+x(c_{_0}+c_1^{}x+c_2^{}y)\equiv X(x,y),
\\[3.5ex]
\dfrac{dy}{dt}=b_{_0}+b_1^{}x+b_2^{}y+y(c_{_0}+c_1^{}x+c_2^{}y)\equiv Y(x,y),
\end{array}
$
\hfill\mbox{} (7.21)
\\[2.5ex]
где вещественные постоянные 
\vspace{0.5ex}
$a_i^{},\,b_i^{},\,c_i^{},\,i=1,2,3,$ такие, что полиномы $X$ и $Y$ являются взаимно простыми, назовем 
[10]
%\marginpar{[33]}
{\it системой Якоби}.
\vspace{0.35ex}

Учитывая, что у дифференциальной системы (7.21) 
\\[1.5ex]
\mbox{}\hfill
$
X_2^{}(x,y)=x(c_1^{}x+c_2^{}y),
\quad 
Y_2^{}(x,y)=y(c_1^{}x+c_2^{}y)
\quad 
\forall (x,y)\in\R^2,
\hfill
$
\\[0.5ex]
получаем
\vspace{0.5ex}

{\bf Свойство 7.4.} 
{\it 
Система {\rm (7.16)} будет проективно особой в том и только в том случае, когда она является системой Якоби.
}
\vspace{0.75ex}

{\bf 7.5.}
Дифференциальную систему
\\[2ex]
\mbox{}\hfill                                         %(7.22)
$
\dfrac{dx}{dt}=A_{m}^{}(x,y)+x\;\!C_{n-1}^{}(x,y)\equiv X(x,y),
\ \ \
\dfrac{dy}{dt}=B_{m}^{}(x,y)+y\;\!C_{n-1}^{}(x,y)\equiv Y(x,y),
$
\hfill\mbox{} (7.22)
\\[2.5ex]
где $A_m^{},\,B_m^{}$ и $C_{n-1}^{}$ 
\vspace{0.5ex}
--- однородные полиномы по переменным $x,\ y$ степеней $m$ и $n-1,$ $m\leq n-1,$ такие, 
\vspace{0.5ex}
что полиномы $X$ и $Y$ являются взаимно простыми, назовем [10]
%\marginpar{[33]}
{\it системой Дарбу}.
Поскольку $m\leq n-1,$ то у системы (7.22)  
\\[1.5ex]
\mbox{}\hfill
$
X_n^{}(x,y)=x\;\!C_{n-1}^{}(x,y),
\quad 
Y_n^{}(x,y)=y\;\!C_{n-1}^{}(x,y)
\quad 
\forall (x,y)\in\R^2,
\hfill
$
\\[1ex]
а значит, имеет место
\vspace{0.5ex}

{\bf Свойство 7.5.} 
{\it 
Система Дарбу {\rm (7.22)} является проективно особой.
}
\\[3.5ex]
\centerline{
{\bf  8. Проективный атлас траекторий дифференциальных систем
}
}
\\[1.75ex]
\indent
Поведение траекторий системы (D) на проективной сфере 
\vspace{0.25ex}
$\P{\mathbb S}(x,y)$  определяется поведением траекторий систем (D), (6.3), (6.6) 
\vspace{0.5ex}
в конечных частях проективных фазовых плоскостей $\R \P(x,y),\, \R \P(\xi,\theta),\, \R \P(\eta,\zeta)$
соответственно. 
\vspace{0.5ex}

Учитывая связи проективных кругов 
\vspace{0.35ex}
$\P\K(x,y), \  \P\K(\xi,\theta),\ \P\K(\eta,\zeta),$ указанные на рис.~5.1, 
по поведению траекторий систем (D), (6.3), (6.6) в конечных частях их проективных фазовых плоскостей 
определяется ход траекторий каждой из этих систем на всей ее проективной фазовой плоскости.
\vspace{0.35ex}

Траектории систем (D), (6.3), (6.6)  
\vspace{0.5ex}
на проективных кругах атласа $(\P\K(x,y),\;\! \P\K(\xi,\theta),\;\! \P\K(\eta,\zeta))$ 
\vspace{0.5ex}
назовем {\it проективным атласом} траекторий системы~(D).

Тогда траектории систем (6.3), (6.6), (D)  
\vspace{0.35ex}
на проективных кругах атласа $(\P\K(\xi,\theta),\;\! \P\K(\eta,\zeta),\;\! \P\K(x,y))$ суть 
\vspace{0.5ex}
проективный атлас траекторий системы  (6.3); 
а траектории систем (6.6), (D), (6.3)  на  проективных кругах атласа 
\vspace{0.5ex}
$(\P\K(\eta,\zeta), \P\K(x,y), \P\K(\xi,\theta))$ --- 
проективный атлас траекторий системы  (6.6).
\vspace{0.35ex}

Таким образом, построив проективный атлас траекторий одной из систем (D), (6.3), (6.6), 
получаем  проективные атласы траекторий двух других систем.
\vspace{1ex}

{\bf Примеры}
\vspace{0.5ex}

{\bf 8.1.}
\vspace{0.5ex}
На рис. 8.1 построен проективный атлас траекторий $a_{_0}y-b_{_0}x=C$ системы (7.8) при $a_{_0}=0,\,b_{_0}\ne0,$ а 
на рис. 8.2 --- при $a_{_0}=b_{_0}\ne0.$
\vspace{0.5ex}
Направление движения вдоль траекторий определяется постоянным вектором 
$\vec{a}(x,y)=a_{_0}\vec{\imath}+b_{_0}\vec{\jmath}\,\;\; \forall (x,y)\in \R^2.$
\\[5.5ex]
\mbox{}\hfill
{\unitlength=1mm
\begin{picture}(42,42)
\put(0,0){\includegraphics[width=42mm,height=42mm]{r08-01a.eps}}
\put(18,41){\makebox(0,0)[cc]{ $y$}}
\put(40.2,18.2){\makebox(0,0)[cc]{ $x$}}
\end{picture}}
\quad
{\unitlength=1mm
\begin{picture}(42,42)
\put(0,0){\includegraphics[width=42mm,height=42mm]{r08-01b.eps}}
\put(18,41){\makebox(0,0)[cc]{ $\theta$}}
\put(40.2,17.8){\makebox(0,0)[cc]{ $\xi$}}
\put(21,-7){\makebox(0,0)[cc]{Рис. 8.1}}
\end{picture}}
\quad
{\unitlength=1mm
\begin{picture}(42,42)
\put(0,0){\includegraphics[width=42mm,height=42mm]{r08-01c.eps}}
\put(18,41){\makebox(0,0)[cc]{ $\zeta$}}
\put(40.2,18){\makebox(0,0)[cc]{ $\eta$}}
\end{picture}}
\hfill\mbox{}
\\[12.5ex]
\mbox{}\hfill
{\unitlength=1mm
\begin{picture}(42,42)
\put(0,0){\includegraphics[width=42mm,height=42mm]{r08-02a.eps}}
\put(18,41){\makebox(0,0)[cc]{ $y$}}
\put(40.2,18.2){\makebox(0,0)[cc]{ $x$}}
\end{picture}}
\quad
{\unitlength=1mm
\begin{picture}(42,42)
\put(0,0){\includegraphics[width=42mm,height=42mm]{r08-02b.eps}}
\put(18,41){\makebox(0,0)[cc]{ $\theta$}}
\put(40.2,17.8){\makebox(0,0)[cc]{ $\xi$}}
\put(21,-7){\makebox(0,0)[cc]{Рис. 8.2}}
\end{picture}}
\quad
{\unitlength=1mm
\begin{picture}(42,42)
\put(0,0){\includegraphics[width=42mm,height=42mm]{r08-02c.eps}}
\put(18,41){\makebox(0,0)[cc]{ $\zeta$}}
\put(40.2,18){\makebox(0,0)[cc]{ $\eta$}}
\end{picture}}
\hfill\mbox{}
\\[7.75ex]
\indent
{\bf 8.2.}
На рис. 8.3 --- 8.8 построены проективные атласы траекторий системы (7.11) в случаях, указанных в таблице 8.1 (соответствуют случаям, приведенным в теореме 2.15 из [11, c. 39]).
%\marginpar{[103]}

\newpage

\mbox{}
\\[-3ex]
\mbox{}\hfill
{
\renewcommand{\arraystretch}{2.5} %регулировка высоты; 
\newcommand{\PreserveBackslash}[1]{\let\temp=\\#1\let\\=\temp}
\let\PBS=\PreserveBackslash
\begin{tabular}{l}
\;\;\;Таблица 8.1
\\
\begin{tabular}                   
              {|>{\PBS\centering\hspace{0pt}}p{3.5cm}
               |>{\PBS\centering\hspace{0pt}}p{5.3cm}
               |>{\PBS\centering\hspace{0pt}}p{3.3cm}
               |>{\PBS\centering\hspace{0pt}}p{0.7cm}|}
\hline 
Система   & Семейство траекторий & Состояния равновесия   &  Рис. \\ 
\hline  
$
\begin{array}{l}
\dfrac{dx}{dt}=x,
\\[-0.3ex]
\dfrac{dy}{dt}=2y 
\end{array}
$
\hfill (8.1) 
& 
$C_1^{}y+C_2^{}x^2=0$
&
\!\!\!\!\!
\begin{tabular}{l}
$O$ --- узел,\\[-3ex]
$O^{(1)}$ --- седло,\\[-3ex]
$O^{(2)}$ --- узел
\end{tabular}
\hfill\mbox{}
& 8.3\\ \hline  
$
\begin{array}{l}
\dfrac{dx}{dt}=x,
\\[-0.3ex] 
\dfrac{dy}{dt}={}-y 
\end{array}
$
\hfill (8.2) 
& 
$xy=C$
&
\!\!\!\!\!
\begin{tabular}{l}
$O$ --- седло,\\[-2.75ex]
$O^{(1)}$ --- узел,\\[-2.75ex]
$O^{(2)}$ --- узел
\end{tabular}
\hfill\mbox{}
& 8.4\\ \hline  
$
\begin{array}{l}
\dfrac{dx}{dt}=y,
\\[-0.3ex] 
\dfrac{dy}{dt}={}-x 
\end{array}
$
\hfill (8.3) 
& 
$x^2+y^2=C$
&
$O$ --- центр
\hfill\mbox{}
& 8.5\\ \hline  
$
\begin{array}{l}
\dfrac{dx}{dt}=x-y,
\\[-0.3ex] 
\dfrac{dy}{dt}=x +y
\end{array}
$
\hfill (8.4) 
& 
$(x^2+y^2)\exp\Bigl({}-2\arctg\dfrac{y}{x}\Bigr)=C$
&
$O$ --- фокус
\hfill\mbox{}
& 8.6\\ \hline  
$
\begin{array}{l}
\dfrac{dx}{dt}=x+y,
\\[-0.3ex] 
\dfrac{dy}{dt}=y
\end{array}
$
\hfill (8.5) 
& 
$y\exp\Bigl({}-\dfrac{x}{y}\Bigr)=C$
&
\!\!\!\!\!\!
\begin{tabular}{l}
$O$ --- вырожден-\\[-4ex]
ный узел,\\[-3ex]
$O^{(1)}$ --- седло-узел
\end{tabular}
& 8.7\\ \hline  
$
\begin{array}{l}
\dfrac{dx}{dt}=x,
\\[-0.3ex] 
\dfrac{dy}{dt}=y 
\end{array}
$
\hfill (8.6) 
&$C_1^{}y+C_2^{}x=0$& 
\!\!\!\!\!
\begin{tabular}{l}
$O$ --- дикрити-
\\[-4ex]
ческий узел
\end{tabular}
\hfill\mbox{}
& 8.8\\ \hline  
\end{tabular}
\end{tabular}
}
\hfill\mbox{}
\\[6ex]
\mbox{}\hfill
{\unitlength=1mm
\begin{picture}(42,42)
\put(0,0){\includegraphics[width=42mm,height=42mm]{r08-03a.eps}}
\put(18,41){\makebox(0,0)[cc]{ $y$}}
\put(40.2,18.5){\makebox(0,0)[cc]{ $x$}}
\end{picture}}
\quad
{\unitlength=1mm
\begin{picture}(42,42)
\put(0,0){\includegraphics[width=42mm,height=42mm]{r08-03b.eps}}
\put(18,41){\makebox(0,0)[cc]{ $\theta$}}
\put(40.2,18){\makebox(0,0)[cc]{ $\xi$}}
\put(21,-7){\makebox(0,0)[cc]{Рис. 8.3}}
\end{picture}}
\quad
{\unitlength=1mm
\begin{picture}(42,42)
\put(0,0){\includegraphics[width=42mm,height=42mm]{r08-03c.eps}}
\put(18,41){\makebox(0,0)[cc]{ $\zeta$}}
\put(40.2,18){\makebox(0,0)[cc]{ $\eta$}}
\end{picture}}
\hfill\mbox{}
\\[10ex]
\mbox{}\hfill
{\unitlength=1mm
\begin{picture}(42,42)
\put(0,0){\includegraphics[width=42mm,height=42mm]{r08-04a.eps}}
\put(18,41){\makebox(0,0)[cc]{ $y$}}
\put(40.2,18.5){\makebox(0,0)[cc]{ $x$}}
\end{picture}}
\quad
{\unitlength=1mm
\begin{picture}(42,42)
\put(0,0){\includegraphics[width=42mm,height=42mm]{r08-04b.eps}}
\put(18,41){\makebox(0,0)[cc]{ $\theta$}}
\put(40.2,18){\makebox(0,0)[cc]{ $\xi$}}
\put(21,-7){\makebox(0,0)[cc]{Рис. 8.4}}
\end{picture}}
\quad
{\unitlength=1mm
\begin{picture}(42,42)
\put(0,0){\includegraphics[width=42mm,height=42mm]{r08-04c.eps}}
\put(18,41){\makebox(0,0)[cc]{ $\zeta$}}
\put(40.2,18){\makebox(0,0)[cc]{ $\eta$}}
\end{picture}}
\hfill\mbox{}
\\[10ex]
\mbox{}\hfill
{\unitlength=1mm
\begin{picture}(42,42)
\put(0,0){\includegraphics[width=42mm,height=42mm]{r08-05a.eps}}
\put(18,41){\makebox(0,0)[cc]{ $y$}}
\put(40.2,18.2){\makebox(0,0)[cc]{ $x$}}
\end{picture}}
\quad
{\unitlength=1mm
\begin{picture}(42,42)
\put(0,0){\includegraphics[width=42mm,height=42mm]{r08-05b.eps}}
\put(18,41){\makebox(0,0)[cc]{ $\theta$}}
\put(40.2,17.8){\makebox(0,0)[cc]{ $\xi$}}
\put(21,-7){\makebox(0,0)[cc]{Рис. 8.5}}
\end{picture}}
\quad
{\unitlength=1mm
\begin{picture}(42,42)
\put(0,0){\includegraphics[width=42mm,height=42mm]{r08-05c.eps}}
\put(18,41){\makebox(0,0)[cc]{ $\zeta$}}
\put(40.2,18){\makebox(0,0)[cc]{ $\eta$}}
\end{picture}}
\hfill\mbox{}
\\[10ex]
\mbox{}\hfill
{\unitlength=1mm
\begin{picture}(42,42)
\put(0,0){\includegraphics[width=42mm,height=42mm]{r08-06a.eps}}
\put(18,41){\makebox(0,0)[cc]{ $y$}}
\put(40.2,18.2){\makebox(0,0)[cc]{ $x$}}
\end{picture}}
\quad
{\unitlength=1mm
\begin{picture}(42,42)
\put(0,0){\includegraphics[width=42mm,height=42mm]{r08-06b.eps}}
\put(18,41){\makebox(0,0)[cc]{ $\theta$}}
\put(40.2,17.8){\makebox(0,0)[cc]{ $\xi$}}
\put(21,-7){\makebox(0,0)[cc]{Рис. 8.6}}
\end{picture}}
\quad
{\unitlength=1mm
\begin{picture}(42,42)
\put(0,0){\includegraphics[width=42mm,height=42mm]{r08-06c.eps}}
\put(18,41){\makebox(0,0)[cc]{ $\zeta$}}
\put(40.2,18){\makebox(0,0)[cc]{ $\eta$}}
\end{picture}}
\hfill\mbox{}
\\[10ex]
\mbox{}\hfill
{\unitlength=1mm
\begin{picture}(42,42)
\put(0,0){\includegraphics[width=42mm,height=42mm]{r08-07a.eps}}
\put(18,41){\makebox(0,0)[cc]{ $y$}}
\put(40.2,18.2){\makebox(0,0)[cc]{ $x$}}
\end{picture}}
\quad
{\unitlength=1mm
\begin{picture}(42,42)
\put(0,0){\includegraphics[width=42mm,height=42mm]{r08-07b.eps}}
\put(18,41){\makebox(0,0)[cc]{ $\theta$}}
\put(40.2,17.8){\makebox(0,0)[cc]{ $\xi$}}
\put(21,-7){\makebox(0,0)[cc]{Рис. 8.7}}
\end{picture}}
\quad
{\unitlength=1mm
\begin{picture}(42,42)
\put(0,0){\includegraphics[width=42mm,height=42mm]{r08-07c.eps}}
\put(18,41){\makebox(0,0)[cc]{ $\zeta$}}
\put(40.2,18){\makebox(0,0)[cc]{ $\eta$}}
\end{picture}}
\hfill\mbox{}
\\[10ex]
\mbox{}\hfill
{\unitlength=1mm
\begin{picture}(42,42)
\put(0,0){\includegraphics[width=42mm,height=42mm]{r08-08a.eps}}
\put(18,41){\makebox(0,0)[cc]{ $y$}}
\put(40.2,18.2){\makebox(0,0)[cc]{ $x$}}
\end{picture}}
\quad
{\unitlength=1mm
\begin{picture}(42,42)
\put(0,0){\includegraphics[width=42mm,height=42mm]{r08-08b.eps}}
\put(18,41){\makebox(0,0)[cc]{ $\theta$}}
\put(40.2,17.8){\makebox(0,0)[cc]{ $\xi$}}
\put(21,-7){\makebox(0,0)[cc]{Рис. 8.8}}
\end{picture}}
\quad
{\unitlength=1mm
\begin{picture}(42,42)
\put(0,0){\includegraphics[width=42mm,height=42mm]{r08-08c.eps}}
\put(18,41){\makebox(0,0)[cc]{ $\zeta$}}
\put(40.2,18){\makebox(0,0)[cc]{ $\eta$}}
\end{picture}}
\hfill\mbox{}
\\[7.5ex]
\indent
{\bf 8.3.}
Используя топологические картины траекторий на проективном круге  $\P\K(x,y)$ системы
\\[2.25ex]
\mbox{}\hfill        %(8.7)
$
2\,\dfrac{dx}{dt}=2y+i\;\!(x-i\;\!y)^q-i\;\!(x+i\;\!y)^q,
\quad \ \
2\,\dfrac{dy}{dt}={}-2x+(x-i\;\!y)^q+ (x+i\;\!y)^q
$
\hfill (8.7)
\\[2.5ex]
при $i=\sqrt{{}-1}\,,\, q=4$ и $q=5$  
\vspace{0.5ex}
из [12, c. 61 --- 65]), на рис. 8.9 построен проективный атлас траекторий 
системы (8.7) при $q=4,$ а на рис. 8.10 --- при $q=5.$ 
\\[5ex]
\mbox{}\hfill
{\unitlength=1mm
\begin{picture}(42,42)
\put(0,0){\includegraphics[width=42mm,height=42mm]{r08-09a.eps}}
\put(18,41){\makebox(0,0)[cc]{ $y$}}
\put(40.2,18.2){\makebox(0,0)[cc]{ $x$}}
\end{picture}}
\quad
{\unitlength=1mm
\begin{picture}(42,42)
\put(0,0){\includegraphics[width=42mm,height=42mm]{r08-09b.eps}}
\put(18,41){\makebox(0,0)[cc]{ $\theta$}}
\put(40.2,17.8){\makebox(0,0)[cc]{ $\xi$}}
\put(21,-7){\makebox(0,0)[cc]{Рис. 8.9}}
\end{picture}}
\quad
{\unitlength=1mm
\begin{picture}(42,42)
\put(0,0){\includegraphics[width=42mm,height=42mm]{r08-09c.eps}}
\put(18,41){\makebox(0,0)[cc]{ $\zeta$}}
\put(40.2,18){\makebox(0,0)[cc]{ $\eta$}}
\end{picture}}
\hfill\mbox{}
\\[10ex]
\mbox{}\hfill
{\unitlength=1mm
\begin{picture}(42,42)
\put(0,0){\includegraphics[width=42mm,height=42mm]{r08-10a.eps}}
\put(18,41){\makebox(0,0)[cc]{ $y$}}
\put(40.2,18.2){\makebox(0,0)[cc]{ $x$}}
\end{picture}}
\quad
{\unitlength=1mm
\begin{picture}(42,42)
\put(0,0){\includegraphics[width=42mm,height=42mm]{r08-10b.eps}}
\put(18,41){\makebox(0,0)[cc]{ $\theta$}}
\put(40.2,17.8){\makebox(0,0)[cc]{ $\xi$}}
\put(21,-7){\makebox(0,0)[cc]{Рис. 8.10}}
\end{picture}}
\quad
{\unitlength=1mm
\begin{picture}(42,42)
\put(0,0){\includegraphics[width=42mm,height=42mm]{r08-10c.eps}}
\put(18,41){\makebox(0,0)[cc]{ $\zeta$}}
\put(40.2,18){\makebox(0,0)[cc]{ $\eta$}}
\end{picture}}
\hfill\mbox{}
\\[10ex]
\centerline{
{\bf\large \S\;\!3. Траектории на сфере Пуанкаре}
}
\\[2.25ex]
\centerline{
{\bf  9. Траектории на проективной фазовой плоскости
}
}
\\[1.5ex]
\indent
Качественно поведение траекторий системы (D) на проективной фазовой плоскости определяется состояниями равновесия, предельными циклами, а также специфическими свойствами траекторий конкретной системы (симметричностью фазового поля направлений, наличием контактных точек, известными траекториями и т.д.). 
\vspace{0.25ex}

Прежде всего сформулируем свойства, связанные с тем, что образом  на проективной сфере 
бесконечно удаленной прямой проективной фазовой плоскости  является экватор сферы.
\vspace{0.25ex}

{\bf Свойство 9.1.}
{\it 
Равносильными являются следующие утверждения}\,:

1. 
{\it
Бесконечно удаленная прямая проективной фазовой плоскости $\R\P(x,y)$ состоит из траекторий системы} (D);

2. 
{\it 
Полином $W^{}_n $ не является тождественным нулем на $\R^2;$} 
\vspace{0.25ex}

3. 
{\it
Система} (D) --- {\it проективно  неособая}\;\!;

4. 
{\it 
У системы} (D) {\it первой проективно приведенной является система} (7.1);

5. 
{\it
Прямая $\theta=0$ состоит из траекторий первой проективно приведенной системы дифференциальной системы} (D);

6. 
{\it
У системы} (D)  {\it второй проективно приведенной является система} (7.2);

7. 
{\it
Прямая $\eta=0$ состоит из траекторий второй проективно приведенной системы дифференциальной системы} (D).
\vspace{0.5ex}

{\bf Свойство 9.2.}
{\it 
Равносильными являются следующие утверждения}\,:

1. 
{\it
Бесконечно удаленная прямая проективной фазовой плоскости $\R\P(x,y)$ не состоит из траекторий системы} (D);

2. 
{\it
Полином} $W_n^{} (x,y)=0\;\; \forall (x,y)\in \R^2;$
\vspace{0.25ex}

3. 
{\it
Система} (D) --- {\it проективно  особая}\;\!;

4. 
{\it 
У системы} (D) {\it первой проективно приведенной является система} (7.3);

5. 
{\it
Прямая $\theta=0$ не состоит из траекторий первой проективно приведенной системы дифференциальной системы} (D);

6. 
{\it
У системы} (D)  {\it второй проективно приведенной является система } (7.4);

7. 
{\it
Прямая $\eta=0$ не состоит из траекторий второй  приведенной системы дифференциальной системы} (D).
\vspace{0.5ex}

Наличие бесконечно удаленных состояний равновесия на проективной фазовой плоскости у системы (D) устанавливается на основании свойств 9.3 --- 9.6.
\vspace{0.5ex}

{\bf Свойство 9.3.}
{\it 
Равносильными являются следующие утверждения}\,:

1. 
{\it Проективно неособая система}  (D) {\it имеет бесконечно удаленное состояние равновесия, лежащее на 
<<концах>> прямой} $y=ax;$
\vspace{0.25ex}

2. 
{\it 
Точка} $(a,0)$ {\it является состоянием равновесия первой проективно приведенной системы} (7.1);

3. 
$\xi=a$ {\it является корнем уравнения} $W_n^{}(1,\xi)=0.$
\vspace{0.75ex}

{\bf Свойство 9.4.}
{\it 
Равносильными являются следующие утверждения}\,:

1. 
{\it Проективно неособая система}  (D) {\it имеет бесконечно удаленное состояние равновесия, лежащее на 
<<концах>> прямой} $x=ay;$
\vspace{0.25ex}

2. 
{\it 
Точка} $(0,a)$ {\it является состоянием равновесия второй проективно приведенной системы} (7.2);

3. 
$\zeta=a$ {\it является корнем уравнения} $W_n^{}(\zeta,1)=0.$
\vspace{0.75ex}

{\bf Свойство 9.5.}
{\it 
Равносильными являются следующие утверждения}\,:

1. 
{\it Проективно особая система}  (D) {\it имеет бесконечно удаленное состояние равновесия, лежащее на 
<<концах>> прямой} $y=ax;$
\vspace{0.25ex}

2. 
{\it 
Точка} $(a,0)$ {\it является состоянием равновесия первой проективно приведенной системы} (7.3);

3. 
$\xi=a$ {\it является решением системы уравнений} $X_n^{}(1,\xi)=0,\ W_{n-1}^{}(1,\xi)=0.$
\vspace{0.75ex}

{\bf Свойство 9.6.}
{\it 
Равносильными являются следующие утверждения}\,:
\vspace{0.2ex}

1. 
{\it Проективно особая система}  (D) {\it имеет бесконечно удаленное состояние равновесия, лежащее на 
<<концах>> прямой} $x=ay;$
\vspace{0.25ex}

2. 
{\it 
Точка} $(0,a)$ {\it является состоянием равновесия второй проективно приведенной системы} 
\vspace{0.2ex}
(7.4);

3. 
$\zeta=a$ {\it является решением системы уравнений} $Y_n^{}(\zeta,1)=0,\ W_{n-1}^{}(\zeta,1)=0.$
\vspace{0.75ex}

Вид бесконечно удаленных состоояний равновесия на проективной фазовой плоскости 
системы (D) устанавливается на основании свойств 9.7 --- 9.12.
\vspace{0.5ex}

{\bf Свойство 9.7.}
{\it
Вид бесконечно удаленного состояния равновесия проективно  неособой системы} (D), {\it
лежащего на <<концах>> прямой $y=ax,$ с точностью до направления движения вдоль примыкающих к нему траекторий такой же как и у состояния равновесия} $(a,0)$ {\it первой проективно приведенной системы} (7.1).
\vspace{0.5ex}

{\bf Свойство 9.8.}
{\it
Вид бесконечно удаленного состояния равновесия проективно  неособой системы} (D), {\it
лежащего на <<концах>> прямой $x=ay,$ с точностью до направления движения вдоль примыкающих к нему траекторий такой же как и у состояния равновесия} $(0,a)$  {\it второй проективно приведенной системы} (7.2).
\vspace{0.5ex}

{\bf Свойство 9.9.}
{\it
Вид бесконечно удаленного состояния равновесия проективно  особой системы} (D), {\it
лежащего на <<концах>> прямой $y=ax,$ с точностью до направления движения вдоль примыкающих к нему траекторий такой же как и у состояния равновесия} $(a,0)$ {\it первой проективно приведенной 
\vspace{0.5ex}
системы} (7.3).

{\bf Свойство 9.10.}
{\it
Вид бесконечно удаленного состояния равновесия проективно особой системы} (D), {\it
лежащего на <<концах>> прямой $x=ay,$ с точностью до направления движения вдоль примыкающих к нему траекторий такой же как и у состояния равновесия} $(0,a)$ {\it второй проективно приведенной 
\vspace{0.65ex}
системы} (7.4).

Точки проективной фазовой плоскости системы (D) назовем {\it регулярными точками}
системы (D), если они не являются состояниями равновесия этой системы. 

Таким образом, каждая точка проективной фазовой плоскости регулярная или состояние равновесия. Каждая траектория является состоянием равновесия или состоит из  регулярных точек.
\vspace{0.2ex}

Если система (D) проективно неособая, то бесконечно удаленная прямая проективной фазовой плоскости состоит из траекторий (свойство 9.1), среди которых может быть конечное число состояний равновесия (число состояний равновесия конечно ввиду того, что правые части $X$ и $Y$ системы (D) суть полиномы). 
Поэтому поведение траекторий проективно неособой системы (D) в окрестности бесконечно удаленной прямой 
определяется лежащими на ней состояниями равновесия.
\vspace{0.2ex}

Если система (D) проективно особая, то бесконечно удаленная прямая проективной фазовой плоскости не состоит из траекторий (свойство 9.2), но на ней может лежать конечное число состояний равновесия.
Для исследования поведения траекторий проективно  особой системы (D)
среди регулярных бесконечно удаленных точек будем выделять те, в которых траектории на проективной фазовой плоскости касаются бесконечно удаленной прямой. Такие точки назовем
{\it экваториально контактными точками} проективно  особой системы (D). При этом будем
исходить из того, что точки касания траекторий системы (D) бесконечно удаленной прямой являются прообразами точек сферы Пуанкаре, в которых образы траекторий касаются экватора.
Находить экваториально контактные точки проективно  особой системы (D)
будем на основании контактных точек координатных осей P-приведенных систем.
\vspace{0.35ex}

Точку, в которой траектория системы (D) касается кривой $\gamma,$ лежащей на фазовой плоскости $Oxy,$
назовем 
{\it контактной точкой} кривой $\gamma.$
В частности, каждая точка кривой $\gamma$ будет контактной точкой в том и только в том случае, когда кривая \vspace{0.3ex}
$\gamma$ является траекторией системы (D).
\vspace{0.5ex}

{\bf Свойство 9.11.}
{\it 
Равносильными являются следующие утверждения}\,:
\vspace{0.25ex}

1. 
{\it 
Бесконечно удаленная точка, 
\vspace{0.15ex}
лежащая на <<концах>> прямой $y=ax,$ 
является экваториально контактной точкой проективно  особой системы} (D); 
\vspace{0.35ex}

2. 
{\it 
Точка} $(a,0)$ {\it является контактной точкой прямой} $\theta=0$ 
\vspace{0.25ex}
{\it первой проективно приведенной системы} (7.3);
\vspace{0.35ex}

3. 
$\xi=a$ {\it является корнем уравнения} $X_n^{}(1,\xi)=0,$ а $W_{n-1}^{}(1,a)\ne0.$

\newpage

{\bf Свойство 9.12.}
{\it 
Равносильными являются следующие утверждения}\,:
\vspace{0.15ex}

1. 
{\it 
Бесконечно удаленная точка, лежащая на <<концах>> прямой $x=ay,$ 
\vspace{0.15ex}
является экваториально контактной точкой проективно  особой системы} (D); 
\vspace{0.15ex}

2. 
{\it 
Точка} $(0,a)$ {\it является контактной точкой прямой} $\eta=0$ 
\vspace{0.15ex}
{\it второй проективно приведенной системы} (7.4);
\vspace{0.15ex}

3. 
$\zeta=a$ {\it является корнем уравнения} $Y_n^{}(\zeta,1)=0,$ а $W_{n-1}^{}(a,1)\ne0.$
\vspace{0.75ex}

Траекторию системы (D), 
\vspace{0.15ex}
проходящую через контактную точку $A$ кривой $\gamma,$ назовем 
{\it
контактной $A\!$-траекторией
}
кривой $\gamma.$ 
\vspace{0.15ex}

{\it
Экваториально контактной $A\!$-траекторией
}
проективно особой системы (D) назо\-вем траекторию, лежащую в проективной фазовой плоскости и
проходящую через
кон\-так\-т\-ную точку $A$ бесконечно удаленной прямой.
\vspace{0.75ex}

{\bf Свойство 9.13.}
{\it 
Экваториально контактная
$A\!$-траектория проективно особой системы} (D)
{\it
в достаточно малой окрестности бесконечно удаленной точки $A,$ 
лежащей на <<концах>> прямой $y=ax,$  расположена}:
\vspace{0.25ex}

{\it 
а{\rm)} в координатной полуплоскости $x>0;$
\vspace{0.35ex}

б\;\!{\rm)} в координатной полуплоскости $x<0;$
\vspace{0.35ex}

в{\rm)} как в координатной полуплоскости $x>0,$ 
\vspace{0.35ex}
так и в координатной полуплоскости $x<0,$
если и только если контактная $A^{(1)}$\!-траектория прямой} 
\vspace{0.25ex}
$\theta=0$ {\it первой проективно приведенной системы} (7.3)
{\it в достаточно малой окрестности  точки $A^{(1)}(a,0)$ соответственно лежит}:
\vspace{0.25ex}

{\it  
а{\rm)} в координатной полуплоскости $\theta>0;$
\vspace{0.35ex}

б\;\!{\rm)} в координатной полуплоскости $\theta<0;$
\vspace{0.35ex}

в{\rm)}\! как\! в\! координатной\! полуплоскости 
\vspace{0.75ex}
$\!\theta\!>\!0,\!$ так\! и\! в\! координатной\! полуплоскости $\!\theta\!<\!0.\!$}

{\bf Свойство 9.14.}
{\it 
Экваториально контактная
$A$\!-траектория проективно  особой системы} (D)
{\it
в достаточно малой окрестности бесконечно удаленной точки $A,$ 
лежащей на <<концах>> прямой $x=ay,$  расположена}:
\vspace{0.25ex}

{\it 
а{\rm)} в координатной полуплоскости $y>0;$
\vspace{0.35ex}

б\;\!{\rm)} в координатной полуплоскости $y<0;$
\vspace{0.35ex}

в{\rm)} как в координатной полуплоскости $y>0,$ 
\vspace{0.5ex}
так и в координатной полуплоскости $y<0,$
если и только если контактная $A^{(2)}$\!-траектория прямой} 
\vspace{0.35ex}
$\eta=0$ {\it второй проективно приведенной системы} (7.4)
{\it в достаточно малой окрестности  точки $A^{(2)}(0,a)$ соответственно лежит}:
\vspace{0.25ex}

{\it  
а{\rm)} в координатной полуплоскости $\eta>0;$
\vspace{0.35ex}

б\;\!{\rm)} в координатной полуплоскости $\eta<0;$
\vspace{0.35ex}

в{\rm)}\! как\! в\! координатной\! полуплоскости 
\vspace{1ex}
$\!\eta\!>\!0,\!$ так\! и\! в\! координатной\! полуплоскости\! 
$\!\eta\!<\!0.\!$}

Признаки расположения экваториально контактных траекторий даны 
\vspace{0.2ex}
в свойствах 9.15 и 9.16.
\\[4.25ex]
\mbox{}\hfill
{\unitlength=1mm
\begin{picture}(40,40)
\put(0,0){\includegraphics[width=40mm,height=40mm]{r09-01.eps}}
 
\put(22,18){\makebox(0,0)[cc]{\scriptsize $O$}}
\put(40,18){\makebox(0,0)[cc]{\scriptsize $x$}}
\put(22,40){\makebox(0,0)[cc]{\scriptsize $y$}}
\put(32.5,32.5){\makebox(0,0)[cc]{\scriptsize $A^{\star}$}}
\put(8,8){\makebox(0,0)[cc]{\scriptsize $A_{\star}^{}$}}

\put(20,-5){\makebox(0,0)[cc]{\rm Рис. 9.1}}
\end{picture}}
\hfill\mbox{}
\\

\newpage

{\bf Свойство 9.15.}\!
{\it 
Если
$\xi\!=\!a\!$ является корнем уравнения $X_n^{}(1,\xi)\!=\!0$ и кроме этого}:
\\[2.25ex]
\mbox{}\hfill
{\it 
а{\rm)} 
$W_{n-1}^{}(1,a)\,\partial_y^{} X_n^{}(1,a)<0;$
\qquad
б\;\!{\rm)}  
$W_{n-1}^{}(1,a)\,\partial_y^{} X_n^{}(1,a)>0;$
\hfill\mbox{}
\\[2.75ex]
\mbox{}\hfill
в{\rm)}  
$
\partial_y^{} X_n^{}(1,a)=0,\
W_{n-1}^{}(1,a)\,\partial_{yy}^{} X_n^{}(1,a)\ne0,
$
\hfill\mbox{}
\\[2.25ex]
то
экваториально контактная 
\vspace{0.25ex}
$A$\!-траектория проективно  особой системы} (D) 
{\it в достаточно малой окрестности 
бесконечно удаленной точки $A,$ 
\vspace{0.15ex}
лежащей на <<концах>> прямой $y=ax,$ соответственно расположена}:
\vspace{0.35ex}

{\it  
а{\rm)}  в координатной полуплоскости $x>0;$
\vspace{0.35ex}

б\;\!{\rm)}  в координатной полуплоскости $x<0;$
\vspace{0.35ex}

в{\rm)}\!  как\! в\! координатной\! полуплоскости $\!x\!>\!0,\!$ 
\vspace{1.25ex}
так\! и\! в\! координатной\! полуплоскости $\!x\!<\!0.\!$}

На рис. 9.1 построен проективный круг  
\vspace{0.25ex}
$\P\K(x,y)$ проективно особой системы (D), 
на котором изо\-бра\-же\-но поведение траекторий в окрестности регулярной экваториально контактной точки $A,$ 
когда экваториально контактная 
\vspace{0.35ex}
$A$\!-траектория лежит в первой координатной четверти проективной фазовой плоскости $\R\P(x,y).$ 
\vspace{0.5ex}

Признаки расположения 
в координатных полуплоскостях контактных с прямой $\theta=0$ траекто\-рий 
первой проективно приведенной системы (7.3) получаем из свойства 9.15 
с учетом свойства 9.13. 
\vspace{0.75ex}

{\bf Свойство 9.16.}\!
{\it 
Если $\zeta\!=\!a\!$ является корнем уравнения $Y_n^{}(\zeta,1)\!=\!0$ и кроме этого}:
\\[2.25ex]
\mbox{}\hfill
{\it 
а{\rm)}  
$W_{n-1}^{}(a,1)\,\partial_x^{} Y_n^{}(a,1)>0;$
\qquad
б\;\!{\rm)}  
$W_{n-1}^{}(a,1)\,\partial_x^{} Y_n^{}(a,1)<0;$
\hfill\mbox{}
\\[2.75ex]
\mbox{}\hfill
в{\rm)}  
$
\partial_x^{} Y_n^{}(a,1)=0,\
W_{n-1}^{}(a,1)\,\partial_{xx}^{} Y_n^{}(a,1)\ne0,
$
\hfill\mbox{}
\\[2.25ex]
то экваториально контактная 
\vspace{0.25ex}
$A$\!-траектория проективно  особой системы} (D) 
{\it 
в достаточно малой окрестности 
бесконечно удаленной точки $A,$ 
\vspace{0.15ex}
лежащей на <<концах>> прямой $x=ay,$ 
соответственно расположена}:
\vspace{0.35ex}

{\it  
а{\rm)}  в координатной полуплоскости $y>0;$
\vspace{0.35ex}

б\;\!{\rm)}  в координатной полуплоскости $y<0;$
\vspace{0.35ex}

в{\rm)}\!  как\! в\! координатной\! полуплоскости $\!y\!>\!0,\!$ так\! и\! в\! координатной\! 
\vspace{1ex}
полуплоскости $\!y\!<\!0.\!$}

Признаки расположения 
\vspace{0.25ex}
в координатных полуплоскостях контактных с прямой $\eta=0$ траекто\-рий 
второй проективно приведенной системы (7.4) 
\vspace{0.15ex}
получаем из свойства 9.16 с учетом свойства 9.14. 
\\[4.75ex]
\centerline{
{\bf  10. Траектории первой проективно приведенной системы
}
}
\\[1.5ex]
\indent
Поведение траекторий 
\vspace{0.15ex}
первой проективно приведенной системы (6.3) существенным образом зависит 
от поведения траекторий системы (D) и, наоборот, 
\vspace{0.15ex}
по траекториям системы (6.3) устанавливается ход траекторий системы (D). 
\vspace{0.35ex}
Эти связи схематически изображены на рис. 5.1 с помощью проективных кругов  $\P\K(x,y)$ и $\P\K(\xi,\theta).$
\vspace{1ex}

{\bf Свойство 10.1.}
\vspace{0.15ex}
{\it 
Состояниями равновесия первой проективно приведенной системы {\rm (6.3)}, 
лежащими в конечной части проективной фазовой плоскости $\R\P(\xi,\theta),$ 
\vspace{0.25ex}
являются образы состояний равновесия системы {\rm (D)}, лежащих}:
\vspace{0.35ex}

1)
{\it
в конечной части проективной фазовой плоскости $\R\P(x,y),$ но не лежащих на оси $Oy;$
}
\vspace{0.35ex}

2)
{\it
на <<концах>> прямых $y=ax,$ где параметр $a$ --- любое вещественное число.
\\[0.35ex]
При этом вид состояний равновесия 
\vspace{0.15ex}
сохраняется с точностью до направления движения вдоль примыкающих к ним траекторий.
}

\newpage

{\bf Свойство 10.2.}\!
\vspace{0.5ex}
{\it 
Если состояние равновесия $\!M(x,y)\!$ системы {\rm (D)} лежит в}\,:
{\it а})
{\it 
пер\-вой};
{\it б}\;\!)
{\it 
второй};
{\it в}\;\!)
{\it 
третьей};
\vspace{0.75ex}
{\it г}\;\!)
{\it 
четвертой открытой координатной четверти фазовой плоскости $\!Oxy,\!$ то его образ
\vspace{0.75ex}
$\!M_{\phantom1}^{(1)}\Bigl(\dfrac{y}{x}\,,\dfrac{1}{x}\Bigr)\!$ соответственно лежит в}:
{\it а})
{\it 
первой};
{\it б})
{\it 
тре\-тьей};
\vspace{0.75ex}
{\it в}\;\!)
{\it 
четвертой};
{\it г}\;\!)
{\it 
второй открытой координатной четверти фазовой плоскости $O_{\phantom1}^{(1)}\xi\theta.$
\vspace{0.75ex}
Если состояние равновесия $M(x,0)$ лежит в полуплоскости $x>0\ (x<0),$ то его образ
$M_{\phantom1}^{(1)}\Bigl(0,\dfrac{1}{x}\Bigr)$ лежит в полуплоскости $\theta>0$ $(\theta<0).$
}
\vspace{1.25ex}

Имеют место и утверждения, обратные к сформулированным в свойствах 10.1 и 10.2, 
когда состояния равновесия системы (D) определяются по соответствующим состояниям равновесия 
первой проективно приведенной системы (6.3).
\vspace{0.25ex}

Образом прямой $x=0$ (координатной оси $Oy)$ 
\vspace{0.25ex}
фазовой плоскости системы (D) является беско\-нечно удаленная прямая 
\vspace{0.25ex}
проективной фазовой плоскости первой проективно приведенной системы (6.3), 
\vspace{0.25ex}
а также прямая $\zeta=0$ (координатная ось $O_{\phantom1}^{(2)}\eta)$ 
фазовой плоскости второй проективно приведенной системы (6.6) (см., например, рис. 5.1). 
\vspace{0.35ex}
Это позволяет указать ряд свойств траекторий систем (D), (6.3) и (6.6).
\vspace{0.75ex}

{\bf Свойство 10.3.}
{\it 
Равносильными являются следующие утверждения}:
\vspace{0.25ex}

1. 
{\it
Первая проективно приведенная система {\rm (6.3)} проективно неособая};
\vspace{0.25ex}

2. 
{\it 
Прямая $x=0$ состоит из тректорий системы} (D);
\vspace{0.25ex}

3. 
{\it
Прямая $\zeta=0$ 
\vspace{0.15ex}
состоит из траекторий второй проективно приведенной дифференциальной системы} (6.6).
\vspace{0.75ex}

{\bf Свойство 10.4.}
{\it 
Равносильными являются следующие утверждения}:
\vspace{0.25ex}

1. 
{\it
Первая проективно приведенная система {\rm (6.3)} проективно особая};
\vspace{0.25ex}

2. 
{\it 
Прямая $x=0$ не состоит из тректорий системы} (D);
\vspace{0.25ex}

3. 
{\it
Прямая $\zeta=0$ 
\vspace{0.15ex}
не состоит из траекторий второй проективно приведенной дифференциальной системы} (6.6).
\vspace{0.75ex}

Образом лежащей на прямой $x=0$ точки $A(0,a)$ 
\vspace{0.15ex}
на проективной плоскости 
$\R\P(\xi,\theta)$ является бесконечно удаленная точка, лежащая на <<концах>> прямой $\xi=a\theta.$
Это позволяет находить бесконечно удаленные состояния равновесия первой проективно приведенной системы
(6.3) по лежащим на прямой $x=0$ состояниям равновесия системы (D), а в случае, когда система (6.3) проективно особая, кроме этого находить ее экваториально контактные точки по контактным точкам прямой 
$x=0$ системы (D).
\vspace{0.5ex}

{\bf Свойство 10.5.}
{\it 
Точка $A(0,a)$ является состоянием равновесия системы {\rm (D)} 
в том и только в том случае, когда лежащая на <<концах>> прямой $\xi=a\theta$ 
точка является бесконечно удаленным состоянием равновесия первой проективно приведенной системы {\rm (6.3)}. 
При этом состояния равновесия одинакового вида с точностью до направления движения  вдоль 
примыкающих к ним траекторий}.
\vspace{0.5ex}

{\bf Свойство 10.6.}
{\it 
Точка $A(0,a)$ является контактной точкой прямой $x=0$ системы {\rm (D)} 
в том и только в том случае, когда лежащая на <<концах>> прямой $\xi=a\theta$ 
точка является экваториально контактной точкой первой проективно приведенной проективно особой системы} (6.3).
\vspace{0.5ex}

Расположение экваториально контактных траекторий проективно особой системы (6.3) 
\vspace{0.25ex}
относительно бесконечно удаленной прямой проективной фазовой плоскости $\R\P(\xi,\theta)$ 
\vspace{0.25ex}
может быть установлено как на основании свойства 9.15, 
примененного к проективно особой системе (6.3), 
\vspace{0.25ex}
так и на основании расположения контактных траекторий прямой $x=0$ системы (D) 
относительно прямой $x=0.$ Во втором случае не требуется знать аналитическое задание 
первой проективно приведенной проективно особой системы (6.3), 
\vspace{0.25ex}
достаточно использовать свойство 10.6 и отображение проективных кругов  $\P\K(x,y)$ и $\P\K(\xi,\theta),$ 
приведенное на рис. 5.1, а также
\vspace{0.5ex}

{\bf Свойство 10.7.}
\vspace{0.35ex}
{\it 
Точка $A(0,a)$ является контактной точкой прямой $x=0$ системы {\rm (D)} тогда и только тогда, когда 
\vspace{0.5ex}
$y=a$ является корнем уравнения $X(0,y)=0,$ а $Y(0,a)\ne 0.$
Если $X(0,a)=0$ и кроме этого}: 
\vspace{0.75ex}

{\it a}) $Y(0,a)\,\partial_y^{} X(0,a)>0;$
\vspace{1ex}

{\it б}\;\!) $Y(0,a)\,\partial_y^{} X(0,a)<0;$
\vspace{1ex}

{\it в}) $\partial_y^{} X(0,a)=0,$ $Y(0,a)\,\partial_{yy}^{} X(0,a)\ne0,$
\\[0.75ex]
{\it то контактная $A$\!-траектория прямой 
\vspace{0.35ex}
$x=0$ системы {\rm (D)} в достаточно малой окрестности точки 
$A(0,a)$ соответственно лежит}:
\vspace{0.35ex}
{\it  
а{\rm)} в полуплоскости $x\geq 0;$
б{\rm)}~в полуплоскости $x\leq 0;$
в{\rm)} как в полуплоскости $x\geq 0,$ так и в полуплоскости $x\leq 0.$}
\vspace{1.25ex}

Таким образом, 
\vspace{0.35ex}
используя отображение проективных кругов  $\P\K(x,y)$ и $\P\K(\xi,\theta),$ 
приведенное на рис. 5.1, по фазовому портрету на проективном  круге $\P\K(x,y)$ 
\vspace{0.35ex}
поведения траекторий системы (D) на проективной фазовой плоскости $\R\P(x,y)$ 
\vspace{0.35ex}
всякий раз можно построить фазовый портрет на проективном круге $\P\K(\xi,\theta)$ 
\vspace{0.35ex}
поведения траекторий первой проективно приведенной системы (6.3) 
\vspace{0.35ex}
на проективной фазовой плоскости $\R\P(\xi,\theta)$ и, наоборот, 
\vspace{0.35ex}
по фазовому портрету на  проективном круге  $\P\K(\xi,\theta)$ системы (6.3)
построить фазовый портрет на  проективном круге  $\P\K(x,y)$ системы (D).
\\[4.25ex]
\centerline{
{\bf  11. Траектории второй проективно приведенной системы
}
}
\\[1.5ex]
\indent
Связи между 
\vspace{0.15ex}
траекториями дифференциальной системы (D) и второй проективно приведенной системы (6.6) 
схематически изображены на рис. 5.1 
\vspace{0.35ex}
с помощью  проективных  кругов  $\P\K(x,y)$ и $\P\K(\eta,\zeta).$
\vspace{0.75ex}

{\bf Свойство 11.1.}
\vspace{0.25ex}
{\it 
Состояниями равновесия второй проективно приведенной системы {\rm (6.6)}, 
лежащими в конечной части проективной фазовой плоскости $\R\P(\eta,\zeta),$ 
являются образы состояний равновесия системы {\rm (D)}, лежащих}:
\vspace{0.5ex}

1)
{\it
в конечной части проективной фазовой плоскости $\R\P(x,y),$ но не лежащих на оси $Ox;$
}
\vspace{0.35ex}

2)
{\it
на <<концах>> прямых $x=ay,$ где параметр $a$ --- любое вещественное число.
\\[0.35ex]
При этом вид состояний равновесия 
\vspace{0.15ex}
сохраняется с точностью до направления движения вдоль примыкающих к ним траекторий.
}
\vspace{1ex}

{\bf Свойство 11.2.}\!
\vspace{0.5ex}
{\it 
Если состояние равновесия $\!M(x,y)\!$ системы {\rm (D)} лежит в}\,:
{\it а})
{\it 
пер\-вой};
{\it б}\;\!)
{\it 
второй};
{\it в}\;\!)
{\it 
третьей};
\vspace{0.75ex}
{\it г}\;\!)
{\it 
четвертой открытой координатной четверти фазовой плоскости $\!Oxy,\!$ то его образ
\vspace{0.75ex}
$\!M_{\phantom1}^{(2)}\Bigl(\dfrac{1}{y}\,,\dfrac{x}{y}\Bigr)\!$ соответственно лежит в}:
{\it а})
{\it 
первой};
{\it б})
{\it 
чет\-вер\-той};
\vspace{0.75ex}
{\it в}\;\!)
{\it 
второй};
{\it г}\;\!)
{\it 
третьей открытой координатной четверти фазовой плоскости $O_{\phantom1}^{(2)}\eta\zeta.$
\vspace{0.75ex}
Если состояние равновесия $M(0,y)$ лежит в полуплоскости $y>0\ (y<0),$ то его образ
$M_{\phantom1}^{(2)}\Bigl(\dfrac{1}{y}\,,0\Bigr)$ лежит в полуплоскости $\eta>0\ (\eta<0).$
}
\vspace{1ex}

Имеют место и утверждения, обратные к сформулированным в свойствах 11.1 и 11.2, 
когда состояния равновесия системы (D) определяются по соответствующим состояниям 
равновесия второй проективно приведенной системы (6.6).
\vspace{0.35ex}

Образом прямой $y=0$ 
\vspace{0.25ex}
(координатной оси $Ox)$ фазовой плоскости системы (D) 
является бесконечно удаленная прямая 
\vspace{0.35ex}
проективной фазовой плоскости второй проективно приведенной системы (6.6), 
а также прямая $\xi=0$ (координатная ось $O_{\phantom1}^{(1)}\theta)$ 
\vspace{0.35ex}
фазовой плоскости первой проективно приведенной системы (6.3) (см., например, рис. 5.1). 
\vspace{0.35ex}

Тогда имеют место

\newpage

{\bf Свойство 11.3.}
{\it 
Равносильными являются следующие утверждения}:
\vspace{0.35ex}

1. 
{\it
Вторая  проективно приведенная система {\rm (6.6)} проективно неособая};
\vspace{0.35ex}

2. 
{\it 
Прямая $y=0$ состоит из тректорий системы} (D);
\vspace{0.35ex}

3. 
{\it
Прямая $\xi=0$ состоит из траекторий первой проективно приведенной 
дифференциальной системы} (6.3).
\vspace{0.5ex}

{\bf Свойство 11.4.}
{\it 
Равносильными являются следующие утверждения}:
\vspace{0.35ex}

1. 
{\it
Вторая  проективно приведенная система {\rm (6.6)} проективно особая};
\vspace{0.35ex}

2. 
{\it 
Прямая $y=0$ не состоит из тректорий системы} (D);
\vspace{0.35ex}

3. 
{\it
Прямая $\xi=0$ не состоит из траекторий первой проективно приведенной дифференциальной системы} (6.3).
\vspace{0.5ex}

Образом лежащей на прямой $y=0$ точки $A(a,0)$ на проективной плоскости 
$\R\P(\eta,\zeta)$ является бес\-ко\-нечно удаленная точка, лежащая на <<концах>> прямой $\zeta=a\eta.$
Это позволяет находить бесконечно удаленные состояния равновесия второй проективно приведенной системы
(6.6) по лежащим на прямой $y=0$ состояниям равновесия системы (D), а в случае, 
когда система (6.6) проективно особая, кроме этого находить ее экваториально контактные точки по контактным точкам прямой 
$y=0$ системы (D).
\vspace{0.5ex}

{\bf Свойство 11.5.}
{\it 
Точка $A(a,0)$ является состоянием равновесия системы {\rm (D)} в том и только в том случае, 
когда лежащая на <<концах>> прямой $\zeta=a\eta$ точка является 
бесконечно удаленным состоянием равновесия второй проективно приведенной системы {\rm (6.6)}. 
При этом состояния равновесия одинакового вида с точностью до направления движения  вдоль 
примыкающих к ним траекторий}.
\vspace{0.5ex}

{\bf Свойство 11.6.}
{\it 
Точка $A(a,0)$ является контактной точкой прямой $y=0$ системы {\rm (D)} в том и только в том случае, 
когда лежащая на <<концах>> прямой $\zeta=a\eta$ точка является экваториально контактной точкой 
второй проективно приведенной проективно особой системы} (6.6).
\vspace{0.5ex}

Расположение экваториально контактных траекторий проективно особой системы (6.6) 
относительно бесконечно удаленной прямой проективной фазовой плоскости $\R\P(\eta,\zeta)$ 
может быть установлено как на основании свойства 9.16, 
примененного к  проективно особой системе (6.6), так и на основании расположения 
контактных траекторий прямой $y=0$ системы (D) относительно прямой $y=0.$ 
Во-втором случае не требуется знать аналитическое задание 
второй проективно приведенной проективно особой системы (6.6), 
\vspace{0.25ex}
достаточно использовать свойство 11.6 и отображение  проективных
кругов  $\P\K(x,y)$ и $\P\K(\eta,\zeta),$ приведенное на рис. 5.1, а также
\vspace{0.75ex}

{\bf Свойство 11.7.}
\vspace{0.35ex}
{\it 
Точка $A(a,0)$ является контактной точкой прямой $y=0$ системы {\rm (D)} 
\vspace{0.5ex}
тогда и только тогда, когда 
$x=a$ является корнем уравнения $Y(x,0)=0,$ а $X(a,0)\ne 0.$
Если $Y(a,0)=0$ и кроме этого}: 
\vspace{0.75ex}

{\it a}) $X(a,0)\,\partial_x^{} Y(a,0)>0;$
\vspace{0.75ex}

{\it б}\;\!) $X(a,0)\,\partial_x^{} Y(a,0)<0;$
\vspace{0.75ex}

{\it в}) $\partial_x^{} Y(a,0)=0,$ $X(a,0)\,\partial_{xx}^{} Y(a,0)\ne0,$
\\[0.5ex]
{\it то контактная $A$\!-траектория прямой $y=0$ системы {\rm (D)} 
\vspace{0.25ex}
в достаточно малой окрестности точки 
$A(a,0)$ соответственно лежит}:
\vspace{0.35ex}
{\it  
а{\rm)} в полуплоскости $y\geq 0;$
б{\rm)}~в полуплоскости $y\leq 0;$
в{\rm)} как в полуплоскости $y\geq 0,$ так и в полуплоскости $y\leq 0.$}
\vspace{0.5ex}

Итак, 
\vspace{0.35ex}
используя отображение проективных кругов  $\P\K(x,y)$ и $\P\K(\eta,\zeta),$ 
приведенное на рис. 5.1, по фазовому портрету на проективном круге  $\P\K(x,y)$ 
\vspace{0.35ex}
поведения траекторий системы (D) на проективной фазовой плоскости $\R\P(x,y)$ 
\vspace{0.35ex}
всякий раз можно построить фазовый портрет на  проективном круге $\P\K(\eta,\zeta)$ 
\vspace{0.35ex}
поведения траекторий второй проективно приведенной системы (6.6) 
\vspace{0.35ex}
на проективной фазовой плоскости $\R\P(\eta,\zeta)$ и, наоборот, 
\vspace{0.35ex}
по фазовому портрету на  проективном круге  
$\P\K(\eta,\zeta)$ 
системы (6.6)
построить фазовый портрет на проективном круге 
$\P\K(x,y)$ 
системы (D).
\\[3.5ex]
\centerline{
{\bf  12.  Линейные и разомкнутые  предельные циклы
}
}
\\[1.5ex]
\indent
Замкнутую траекторию системы\! (D), 
\vspace{0.25ex}
лежащую в конечной части проективной фазо\-вой плоскости $\R\P(x,y),$ 
назовем 
[13, c. 22]
%\marginpar{[2156]}
{\it циклом} системы (D).
\vspace{0.25ex}
Если существует окрестность цикла системы (D), 
\vspace{0.25ex}
в которой нет циклов, отличных от данного цикла, (изолированный цикл), то такой цикл назовем 
[1, c. 71 --- 91]
%\marginpar{[1506]}
{\it предельным циклом} 
\vspace{0.5ex}
системы (D).

Циклу системы (D) на открытом проективном круге 
\vspace{0.25ex}
$\R\K(x,y)\backslash\partial \R\K(x,y)$ соответствует замкнутая траектория (цикл на проективном круге), 
\vspace{0.35ex}
а на проективной сфере 
$\P{\mathbb S}(x,y)$
 --- пара 
\vspace{0.25ex}
диаметрально противоположных замкнутых траекторий, не имеющих общих точек с экватором 
(антиподальные циклы на проективной сфере). 
\vspace{0.25ex}

Существуют траектории системы (D), 
\vspace{0.25ex}
являющиеся замкнутыми на проективной фазовой плоскости $\R\P(x,y),$ но незамкнутые на ее конечной части $(x,y).$
\vspace{0.75ex}

{\bf Определение 12.1.}
\vspace{0.25ex}
{\it 
Прямую, являющуюся траекторией системы {\rm (D)} на проективной фазовой плоскости $\R\P(x,y),$ 
\vspace{0.25ex}
назовем \textit{\textbf{линейным циклом}} 
системы {\rm (D)} или {\boldmath$\ell\!$}\textit{\textbf{-циклом}} 
системы} (D). 
\vspace{0.5ex}

Линейный цикл, будучи траекторией на $\R\P(x,y),$  
\vspace{0.25ex}
состоит из регулярных точек, лежащих как в конечной части 
\vspace{0.15ex}
проективной фазовой плоскости, так и на бесконечно удаленной прямой. 
\vspace{0.15ex}

Если бесконечно удаленная  
\vspace{0.15ex}
прямая  проективной фазовой плоскости $\R\P(x,y)$  является траекторией системы (D), 
\vspace{0.25ex}
то она будет линейным циклом, который назовем 
\linebreak
$\ell_{\infty}^{}\!$-цик\-лом системы (D).
\vspace{0.5ex}

На проективной сфере $\P{\mathbb S}(x,y)$ 
\vspace{0.5ex}
линейным циклам системы (D) соответствуют окружности большого радиуса; 
\vspace{0.5ex}
$\ell_{\infty}^{}\!$-циклу соответствует экватор. 
На проективном круге  $\P\K(x,y)\ \,\ell_{\infty}^{}\!$-циклу соответствует граничная окружность 
$\partial \P\K(x,y).$
\vspace{0.75ex}

{\bf Определение 12.2.}
\vspace{0.25ex}
{\it 
\textit{\textbf{Предельным линейным циклом}} системы {\rm (D)} назовем такой линейный цикл, которому на 
проективной сфере 
$\P{\mathbb S}(x,y)$
\vspace{0.25ex}
соответствует окружность большого радиуса, 
\vspace{0.25ex}
имеющая окрестность, в которой нет замкнутых траекторий 
системы {\rm (D)}, отличных от этой окружности.} 
\vspace{0.75ex}

У проективно неособой системы (8.3) 
\vspace{0.25ex}
с проективным атласом траекторий на рис. 8.5 каждая  траектория,
\vspace{0.25ex}
отличная от состояния равновесия $O(0,0),$ является циклом. 
В том числе бесконечно удаленная прямая проективной фазовой плоскости 
\vspace{0.35ex}
$\R\P(x,y)$ есть 
\linebreak
 $\ell_{\infty}^{}\!$-цикл системы (8.3).
\vspace{0.5ex}
При этом  $\ell_{\infty}^{}\!$-цикл и  циклы, лежащие в конечной части проективной фазовой плоскости $\R\P(x,y),$ 
\vspace{0.5ex}
не являются предельными циклами этой  системы.

Проективно неособая система (8.4),
\vspace{0.35ex}
проективный атлас траекторий которой построен на рис. 8.6, имеет 
предельный  $\ell_{\infty}^{}$-цикл.
\vspace{0.5ex}
В конечной части проективной фазовой плоскости $\R\P(x,y)$  у системы (8.4)  предельных циклов нет, так как 
\vspace{1ex}
система (8.4) линейная.

Система Якоби [14, c. 111 --- 117; 10, c. 14 --- 25]
\\[2.5ex]
\mbox{}\hfill
$
\dfrac{dx}{dt}=x-y+x(x+y),
\qquad
\dfrac{dy}{dt}=(y+1)(x+y)
$
\hfill (12.1)
\\[2.5ex]
имеет общий автономный интеграл 
\\[2ex]
\mbox{}\hfill
$
F\colon (x,y)\to\ \dfrac{x^2+y^2}{(y+1)^2}\, \exp \Bigl( {}-2 \arctg \dfrac{y}{x}\Bigr)
\quad 
\forall (x,y)\in G
\hfill
$
\\[2ex]
на любой области из множества $G=\{(x,y)\colon x\ne 0 \ \&\  y\ne {}-1\,\}.$
\vspace{0.5ex}

Прямая-траектория $y={}-1$
\vspace{0.35ex}
является предельным линейным циклом проективно особой системы (12.1).
\vspace{0.35ex}
Предельных циклов,  лежащих в конечной части  проективной фазовой плоскости $\R\P(x,y),$
у системы (12.1) нет [10, c. 55].
\\[2ex]
\mbox{}\hfill
{\unitlength=1mm
\begin{picture}(42,42)
\put(0,0){\includegraphics[width=42mm,height=42mm]{r12-1a.eps}}
\put(18,41){\makebox(0,0)[cc]{ $y$}}
\put(40.2,18.2){\makebox(0,0)[cc]{ $x$}}
\end{picture}}
\quad
{\unitlength=1mm
\begin{picture}(42,42)
\put(0,0){\includegraphics[width=42mm,height=42mm]{r12-1b.eps}}
\put(18,41){\makebox(0,0)[cc]{ $\theta$}}
\put(40.2,17.8){\makebox(0,0)[cc]{ $\xi$}}
\put(21,-7){\makebox(0,0)[cc]{Рис. 12.1}}
\end{picture}}
\quad
{\unitlength=1mm
\begin{picture}(42,42)
\put(0,0){\includegraphics[width=42mm,height=42mm]{r12-1c.eps}}
\put(18,41){\makebox(0,0)[cc]{ $\zeta$}}
\put(40.2,18){\makebox(0,0)[cc]{ $\eta$}}
\end{picture}}
\hfill\mbox{}
\\[9ex]
\indent
Первой проективно приведенной системой системы (12.1) является система Якоби
\\[2ex]
\mbox{}\hfill                   %(12.2)
$
\dfrac{d\xi}{dt}=1+\xi^2,
\qquad
\dfrac{d\theta}{dt}={}-1-\xi-\theta +\xi \theta,
$
\hfill (12.2)
\\[2.25ex]
которая не имеет 
\vspace{0.5ex}
предельных циклов в конечной части  проективной фазовой плоскости $\R\P(\xi,\theta)$ [10, c. 55].
Прямая-траектория 
\vspace{0.5ex}
$\xi={}-\theta$ --- предельный линейный цикл  системы Якоби (12.2).
Трансцендентная функция 
\\[2ex]
\mbox{}\hfill
$
F^{(1)}_{\phantom 1} \colon (\xi,\theta)\to\
\dfrac{1+\xi^2}{(\xi+\theta)^2} \, \exp ({}-2\arctg \xi)
\quad 
\forall (\xi,\theta)\in G^{(1)}_{\phantom 1},
\qquad 
G^{(1)}_{\phantom 1}=\{(\xi,\theta) \colon \xi\ne {}-\theta\},
\hfill
$
\\[2ex]
есть 
\vspace{0.5ex}
общий автономный интеграл системы (12.2) на любой области из множества $G^{(1)}_{\phantom 1}.$

Второй проективно приведенной системой системы (12.1) является система Якоби
\\[2ex]
\mbox{}\hfill                  %(12.3)
$
\dfrac{d\eta}{dt}={}-(1+\eta)(1+\zeta),
\qquad
\dfrac{d\zeta}{dt}={}-(1+\zeta^2),
$
\hfill (12.3)
\\[2.25ex]
которая не имеет 
\vspace{0.5ex}
предельных циклов в конечной части  проективной фазовой плоскости $\R\P(\eta,\zeta)$ [10, c. 55].
Прямая-траектория 
\vspace{0.5ex}
$\eta={}-1$ --- предельный линейный цикл системы Якоби (12.3).
Трансцендентная функция 
\\[2.25ex]
\mbox{}\hfill
$
F^{(2)}_{\phantom 1} \colon (\eta,\zeta)\to \
\dfrac{1+\zeta^2}{(1+\eta)^2} \, \exp \Bigl({}-2\arctg \dfrac1\zeta\Bigr)
\quad 
\forall (\eta,\zeta)\in G^{(2)}_{\phantom 1},
\quad 
G^{(2)}_{\phantom 1}\!=\!\{(\eta,\zeta) \colon \eta\ne\! {}-1 \, \&\,  \zeta\ne 0\},
\hfill
$
\\[1.75ex]
есть общий автономный интеграл системы (12.3) на любой области из множества $G^{(2)}_{\phantom 1}.$
\vspace{0.35ex}

На рис. 12.1 построен проективный атлас траекторий системы Якоби (12.1).
\vspace{0.5ex}

Замкнутые на 
проективной сфере 
$\P{\mathbb S}(x,y)$
траектории системы (D), отличные от окружностей большого радиуса, подразделяются на два вида: 
\vspace{0.35ex}

1) не имеющие общих точек с экватором 
проективной сферы 
$\P{\mathbb S}(x,y);$
\vspace{0.35ex}

2) имеющие хотя бы одну общую точку с экватором 
проективной сферы 
$\P{\mathbb S}(x,y).$
\vspace{0.75ex}

{\bf Определение 12.3 (12.4).}
\vspace{0.35ex}
{\it
\textit{\textbf{Разомкнутым циклом}} 
{\rm(}\textit{\textbf{разомкнутым предельным  циклом}}{\rm)} 
\vspace{0.35ex}
системы {\rm (D)} назовем ее траекторию  на проективной фазовой плоскости $\R\P(x,y),$ которой на 
проективной сфере $\P{\mathbb S}(x,y)$
\vspace{0.35ex}
соответствует пара антиподальных замкнутых траекторий 
\vspace{0.35ex}
{\rm(}пара антиподальных изолированных замкнутых траекторий{\rm)}, имеющих общие точки с экватором 
проективной сферы $\P{\mathbb S}(x,y).$
}
\vspace{0.75ex}

Тогда лежащему  в конечной части  проективной фазовой плоскости $\R\P(x,y)$ 
\vspace{0.25ex}
циклу дифференциальной системы (D) на проективной сфере $\P{\mathbb S}(x,y)$
\vspace{0.25ex}
соответствует пара антиподальных замкнутых траекторий дифференциальной системы (D), 
\vspace{0.25ex}
не имеющих общих точек с экватором проективной сферы $\P{\mathbb S}(x,y).$

\newpage

Проективно особая система Якоби 
\\[2ex]
\mbox{}\hfill                             %(12.4)
$
\dfrac{d\xi}{dt}={}-(1+\xi^2),
\qquad
\dfrac{d\theta}{dt}={}-\xi\theta
$
\hfill (12.4)
\\[2ex]
является первой проективно приведенной системой системы (8.3).
\vspace{0.25ex}

Прямая-траектория 
\vspace{0.15ex}
$\theta=0$ ---  линейный цикл  системы (12.4), 
а любая другая траектория системы (12.4) ---  разомкнутый цикл (рис. 8.5).
\vspace{0.15ex}

Линейных, разомкнутых 
\vspace{0.25ex}
и лежащих в конечной части  проективной фазовой плоскости $\R\P(\xi,\theta)$
предельных циклов у системы (12.4) нет.
\vspace{0.35ex}

Функция
\\[1.75ex]
\mbox{}\hfill
$
F^{(1)}_{\phantom 1}\colon (\xi,\theta)\to\ 
\dfrac{1+\xi^2}{\theta^2}
\quad 
\forall (\xi,\theta)\in G^{(1)}_{\phantom 1},
\quad 
G^{(1)}_{\phantom 1}=\{(\xi,\theta)\colon \theta\ne 0\},
\hfill
$
\\[1.5ex]
есть
\vspace{0.75ex}
общий автономный интеграл системы (12.4) на любой области из множества $G^{(1)}_{\phantom 1}.$

Проективно особая система Якоби 
\\[2ex]
\mbox{}\hfill                %(12.5)
$
\dfrac{d\eta}{dt}=\eta \zeta,
\qquad
\dfrac{d\zeta}{dt}=1+\zeta^2
$
\hfill (12.5)
\\[2ex]
является второй  проективно приведенной системой системы (8.3).
\vspace{0.25ex}

Прямая-траектория 
\vspace{0.25ex}
$\eta=0$ ---  линейный цикл  системы (12.5), 
а любая другая траектория системы (12.5) ---  разомкнутый цикл (рис. 8.5).
\vspace{0.15ex}

Линейных, 
\vspace{0.25ex}
разомкнутых и лежащих в конечной части  проективной фазовой плоскости $\R\P(\eta,\zeta)$
предельных циклов у системы (12.5) нет.
\vspace{0.35ex}

Функция
\\[1.75ex]
\mbox{}\hfill
$
F^{(2)}_{\phantom 1}\colon (\eta,\zeta)\to\ 
\dfrac{1+\zeta^2}{\eta^2}
\quad 
\forall (\eta,\zeta)\in G^{(2)}_{\phantom 1},
\quad
G^{(2)}_{\phantom 1}=\{(\eta,\zeta)\colon \eta\ne 0\},
\hfill
$
\\[1.5ex]
есть 
\vspace{0.75ex}
общий автономный интеграл системы (12.5) на любой области из множества $G^{(2)}_{\phantom 1}.$

Проективно особая система Якоби 
\\[2ex]
\mbox{}\hfill             %(12.6)
$
\dfrac{d\xi}{dt}=1+\xi^2,
\qquad
\dfrac{d\theta}{dt}=\theta(\xi -1)
$
\hfill (12.6)
\\[2ex]
является первой проективно приведенной системой системы (8.4).
\vspace{0.25ex}

Прямая-траектория 
\vspace{0.15ex}
$\theta=0$ является предельным линейным циклом  системы (12.6)  (рис. 8.6).
\vspace{0.35ex}
Разомкнутых и лежащих в конечной части  проективной фазовой плоскости $\R\P(\xi,\theta)$
циклов (в том числе и предельных) у системы (12.6) нет.
\vspace{0.35ex}

Функция
\\[1.75ex]
\mbox{}\hfill
$
F^{(1)}_{\phantom 1}\colon (\xi,\theta)\to\ 
\dfrac{1+\xi^2}{\theta^2}\, \exp ({}-2\arctg \xi)
\quad 
\forall (\xi,\theta)\in G^{(1)}_{\phantom 1},
\quad 
G^{(1)}_{\phantom 1}=\{(\xi,\theta)\colon \theta\ne 0\},
\hfill
$
\\[1.5ex]
есть 
\vspace{0.75ex}
общий автономный интеграл системы (12.6) на любой области из множества $G^{(1)}_{\phantom 1}.$

Проективно особая система Якоби 
\\[2ex]
\mbox{}\hfill             %(12.7)
$
\dfrac{d\eta}{dt}={}-\eta(1+ \zeta),
\quad
\dfrac{d\zeta}{dt}={}-(1+\zeta^2)
$
\hfill (12.7)
\\[2ex]
является второй  проективно приведенной системой системы (8.4).
\vspace{0.25ex}

Прямая-траектория $\eta=0$
\vspace{0.25ex}
является предельным  линейным циклом  системы (12.7) (рис. 8.6).
\vspace{0.5ex}
Разомкнутых и лежащих в конечной части  проективной фазовой плоскости $\R\P(\eta,\zeta)$
циклов (в том числе и предельных) у системы (12.7) нет.
\vspace{0.25ex}

\newpage

Функция
\\[1.75ex]
\mbox{}\hfill
$
F^{(2)}_{\phantom 1}\colon (\eta,\zeta)\to\ 
\dfrac{1+\zeta^2}{\eta^2}\,\exp \Bigl({}-2\arctg \dfrac1\zeta\Bigr)
\quad  
\forall (\eta,\zeta)\in G^{(2)}_{\phantom 1},
\quad 
G^{(2)}_{\phantom 1}=\{(\eta,\zeta)\colon \eta \zeta\ne 0\},
\hfill
$
\\[1.5ex]
есть 
\vspace{1ex}
общий автономный интеграл системы (12.7) на любой области из множества $G^{(2)}_{\phantom 1}.$

{\bf Свойство 12.1.}
{\it 
Разомкнутые и линейные предельные циклы, отличные от 
\vspace{0.25ex}
$\ell_{\infty}^{}\!$-цик\-ла, существуют у проективно особых систем {\rm (D);} 
\vspace{0.25ex}
а предельный $\ell_{\infty}^{}\!$-цикл существует у проективно неособых систем {\rm (D)}.
}
\vspace{0.25ex}

{\sl Доказательство} 
основано на том, что проективный тип системы (D) зависит от того состоит 
или не состоит из траекторий системы (D) 
\vspace{0.25ex}
бесконечно удаленная прямая 
проективной фазовой плоскости $\R\P(x,y)$  (свойства 9.1 и 9.2).\k
\vspace{0.75ex}

{\bf Свойство 12.2.}
{\it 
Проективно неособая система {\rm (D)} с четной степенью $n$ 
имеет хотя бы одно бесконечно удаленное состояние равновесия.
}
\vspace{0.35ex}

{\sl Действительно}, 
\vspace{0.25ex}
если $n$ --- четное число, то или  первая проективно приведенная дифференциальная система (7.1) 
имеет состояние равновесия на прямой $\theta=0$ или  
вторая проективно приведенная  дифференциальная система (7.2) 
\vspace{0.35ex}
имеет состояние равновесия 
в начале координат фазовой плоскости $O^{(2)}_{\phantom 1} \eta \zeta.$
\vspace{0.75ex}

Из свойства 12.1 с учетом свойства 12.2 получаем 
\vspace{0.75ex}

{\bf Свойство 12.3.}
\vspace{0.25ex}
{\it 
Для того, чтобы система {\rm (D)} имела предельный $\ell_{\infty}^{}\!$-цикл, 
необходимо,  чтобы она была  проективно неособой   с  нечетной степенью $n.$
}
\vspace{0.5ex}

Проективно особая система  Якоби (12.1) 
\vspace{0.75ex}
разомкнутых предельных циклов не имеет.

Система Дарбу [10]                          
\\[2.25ex]
\mbox{}\hfill                      %(12.8)
$
\dfrac{dx}{dt}={}-y-x(x^2+y^2-1),
\qquad
\dfrac{dy}{dt}=x-y(x^2+y^2-1),
$
\hfill (12.8)
\\[2.5ex]
в конечной части  проективной фазовой плоскости $\R\P(x,y)$
\vspace{0.35ex}
имеет одно состояние равновесия, которое является неустойчивым фокусом.
\vspace{0.35ex}
Окружность $x^2+y^2=1$  --- единственный предельный цикл системы (12.8) в 
\vspace{0.5ex}
конечной части  проективной фазовой плоскости $\R\P(x,y)$ 
[4, c. 257 -- 260; 15, c. 39 -- 46]. 
Функция 
\\[2.5ex]
\mbox{}\hfill
$
F\colon (x,y) \to\ \dfrac{x^2+y^2}{1-x^2-y^2}\ \exp\Bigl({}-2\arctg\dfrac{y}{x}\Bigr)
\quad
\forall (x,y)\in G
\hfill
$
\\[2ex]
является 
\vspace{0.5ex}
общим автономным интегралом системы (12.8) на любой области из множества 
$G=\{(x,y)\colon x\ne 0\ \, \&\,\ x^2+y^2\ne 1\}.$
\vspace{0.75ex}

Среди траекторий проективно особой системы Дарбу (12.8) 
\vspace{0.25ex}
на проективной фазовой плоскости $\R\P(x,y)$ нет прямых. 
\vspace{0.35ex}
У системы (12.8) нет экваториально контактных точек. 
Система (12.8) не имеет линейных и разомкнутых циклов, в том числе, 
у нее нет предельных линейных и предельных разомкнутых циклов.
\vspace{0.5ex}

Первой проективно приведенной системой системы (12.8) является проективно особая система 
\\[2.2ex]
\mbox{}\hfill                      %(12.9)
$
\dfrac{d\xi}{d\tau}=\theta (1+\xi^2),
\qquad
\dfrac{d\theta}{d\tau}=1+\xi^2-\theta^2+\xi\theta^2, 
\quad \ \theta d\tau=dt,
$
\hfill (12.9)
\\[2.5ex]
у которой нет 
\vspace{0.5ex}
линейных и лежащих в конечной части  проективной фазовой плоскости $\R\P(\xi,\theta)$ циклов. 
\vspace{0.5ex}
Гипербола $\theta^2-\xi^2=1$ --- предельный разомкнутый цикл системы (12.9). Других предельных циклов в 
проективной фазовой плоскости $\R\P(\xi,\theta)$ 
\vspace{0.35ex}
у системы (12.9) нет. 
Трансцендентная функция 
\\[2ex]
\mbox{}\hfill
$
F^{(1)}_{\phantom 1}\colon (\xi,\theta) \to\ 
\dfrac{1+\xi^2}{\theta^2-\xi^2-1}\ \exp({}-2\arctg\xi)
\quad 
\forall (\xi,\theta)\in G^{(1)}_{\phantom 1}
\hfill
$ 
\\[1.75ex]
является 
\vspace{0.35ex}
общим автономным интегралом системы (12.9) на любой области из множества 
$G^{(1)}_{\phantom 1}=\{(\xi,\theta)\colon \theta^2-\xi^2\ne 1\}.$
\\[3.5ex]
\mbox{}\hfill
{\unitlength=1mm
\begin{picture}(42,42)
\put(0,0){\includegraphics[width=42mm,height=42mm]{r12-2a.eps}}
\put(18,41){\makebox(0,0)[cc]{ $y$}}
\put(40.2,18.2){\makebox(0,0)[cc]{ $x$}}
\end{picture}}
\quad
{\unitlength=1mm
\begin{picture}(42,42)
\put(0,0){\includegraphics[width=42mm,height=42mm]{r12-2b.eps}}
\put(18,41){\makebox(0,0)[cc]{ $\theta$}}
\put(40.2,17.8){\makebox(0,0)[cc]{ $\xi$}}
\put(21,-7){\makebox(0,0)[cc]{Рис. 12.2}}
\end{picture}}
\quad
{\unitlength=1mm
\begin{picture}(42,42)
\put(0,0){\includegraphics[width=42mm,height=42mm]{r12-2c.eps}}
\put(18,41){\makebox(0,0)[cc]{ $\zeta$}}
\put(40.2,18){\makebox(0,0)[cc]{ $\eta$}}
\end{picture}}
\hfill\mbox{}
\\[7.5ex]
\indent
Второй проективно приведенной системой системы (12.8) является проективно особая система 
\\[2.2ex]
\mbox{}\hfill                      %(12.10)
$
\dfrac{d\eta}{d\nu}=1-\eta^2+\zeta^2-\eta^2\zeta, 
\qquad
\dfrac{d\zeta}{d\nu}={}-\eta(1+\zeta^2), 
\quad \ \ \eta d\nu=dt,
$
\hfill (12.10)
\\[2.5ex]
у которой гипербола $\eta^2-\zeta^2=1$ 
\vspace{0.25ex}
является предельным разомкнутым циклом.
Других циклов (в том числе и предельных циклов) 
\vspace{0.25ex}
в проективной фазовой плоскости $\R\P(\eta,\zeta)$ у системы (12.10) нет.
Функция 
\\[1.75ex]
\mbox{}\hfill
$
F^{(2)}_{\phantom 1}\colon (\eta,\zeta) \to\ 
\dfrac{1+\zeta^2}{\eta^2-\zeta^2-1}\ 
\exp\Bigl({}-2\arctg\dfrac{1}{\zeta}\,\Bigr)
\quad
\forall (\eta,\zeta)\in G^{(2)}_{\phantom 1},
\hfill
$ 
\\[1.5ex]
является 
\vspace{0.35ex}
общим автономным интегралом системы (12.10) на любой области из множества 
$G^{(2)}_{\phantom 1}=\{(\eta,\zeta)\colon \zeta\ne 0\,\ \&\,\ \eta^2-\zeta^2\ne 1\}.$
\vspace{0.5ex}

Проективный атлас траекторий системы (12.8) построен на рис. 12.2. 
\vspace{0.5ex}

Проективно особая система
\\[2ex]
\mbox{}\hfill                         %(12.11)
$
\dfrac{dx}{dt}={}-2x-y+3x^2+y^2 -x(x^2+y^2),
\quad \
\dfrac{dy}{dt}={}-1+x+2xy -y(x^2+y^2)
$
\hfill (12.11)
\\[2.25ex]
получена из системы Дарбу (12.8) 
\vspace{0.15ex}
параллельным переносом начала координат фазовой плоскости $Oxy$
в точку $A(1,0).$ 
\vspace{0.25ex}
Поэтому в конечной части  проективной фазовой плоскости $\R\P(x,y)$ у системы (12.11) 
\vspace{0.25ex}
одно состояние равновесия $A(1,0),$  
являющееся неустойчивым фокусом, и один предельный цикл $(x-1)^2+y^2=1.$
Функция 
\\[2ex]
\mbox{}\hfill
$
F\colon (x,y) \to\ 
\dfrac{(x-1)^2+y^2}{1-(x-1)^2-y^2}\ \exp\Bigl({}-2\arctg\dfrac{y}{x-1}\Bigr)
\quad
\forall (x,y)\in G
\hfill
$
\\[1.75ex]
является общим автономным интегралом системы (12.11)
\vspace{0.35ex}
на любой области из множества 
$G=\{(x,y)\colon x \ne 1\ \, \&\,\ (x-1)^2+y^2\ne 1\}.$
\vspace{0.5ex}

Первой проективно приведенной системой системы (12.11) является проективно особая система 
\\[2ex]
\mbox{}\hfill                      %(12.12)
$
\dfrac{d\xi}{d\tau}={}-\xi +\theta+2\xi\theta - \theta^2-\xi^3+\xi^2\theta,
\quad
\dfrac{d\theta}{d\tau}=1-3\theta+\xi^2+2\theta^2-\xi^2\theta+\xi\theta^2, \ \ \ \theta d\tau=dt,
$
\hfill (12.12)
\\[2.25ex]
у которой парабола $\xi^2-2\theta+1=0$ 
\vspace{0.35ex}
является предельным разомкнутым циклом. Других циклов 
(в том числе и предельных циклов) в проективной фазовой плоскости $\R\P(\xi,\theta)$ у системы (12.12) нет. 
Функция 
\\[2ex]
\mbox{}\hfill
$
F^{(1)}_{\phantom 1}\colon (\xi,\theta) \to\ 
\dfrac{\xi^2+(\theta-1)^2}{2\theta-1-\xi^2}\ 
\exp\Bigl(2\arctg\dfrac{\xi}{\theta-1}\Bigr)
\quad
\forall (\xi,\theta)\in G^{(1)}_{\phantom 1}
\hfill
$ 
\\[1.75ex]
является общим автономным интегралом системы (12.12)
\vspace{0.35ex}
на любой области из множества 
$G^{(1)}_{\phantom 1}=\{(\xi,\theta)\colon \theta\ne 1\ \,\&\,\ 1-2\theta+\xi^2\ne 0\}.$ 
\vspace{0.5ex}

Второй проективно приведенной системой системы (12.11) является проективно особая система 
\\[2.2ex]
\mbox{}\hfill                      %(12.13)
$
\dfrac{d\eta}{d\nu}=1-2\eta\zeta+\zeta^2+\eta^3-\eta^2\zeta\,, 
\quad
\dfrac{d\zeta}{d\nu}=1-\eta-2\eta\zeta +\zeta^2+ \eta^2\zeta-\eta\zeta^2, \ \ \ \eta d\nu=dt,
$
\hfill (12.13)
\\[2.5ex]
у которой гипербола $\zeta^2-2\eta\zeta+1=0$ 
\vspace{0.35ex}
является предельным  разомкнутым циклом.
Других циклов (в том числе и предельных циклов) 
\vspace{0.25ex}
в проективной фазовой плоскости $\R\P(\eta,\zeta)$ 
у системы (12.13) нет. 
Функция 
\\[2ex]
\mbox{}\hfill
$
F^{(2)}_{\phantom 1}\colon (\eta,\zeta) \to \
\dfrac{1+(\eta-\zeta)^2}{2\eta\zeta-\zeta^2-1}\ 
\exp\Bigl(2\arctg\dfrac{1}{\eta-\zeta}\,\Bigr)
\quad
\forall (\eta,\zeta)\in G^{(2)}_{\phantom 1}
\hfill
$ 
\\[1.75ex]
является общим автономным интегралом системы (12.13) 
\vspace{0.35ex}
на любой области из множества 
$G^{(2)}_{\phantom 1}=\{(\eta,\zeta)\colon \eta\ne \zeta\,\ \&\,\ 2\eta\zeta-\zeta^2- 1\ne0\}.$
\vspace{0.5ex}

Проективный атлас траекторий системы (12.11) построен на рис. 12.3. 
\\[4.25ex]
\mbox{}\hfill
{\unitlength=1mm
\begin{picture}(42,42)
\put(0,0){\includegraphics[width=42mm,height=42mm]{r12-3a.eps}}
\put(18,41){\makebox(0,0)[cc]{ $y$}}
\put(40.2,18.2){\makebox(0,0)[cc]{ $x$}}
\end{picture}}
\quad
{\unitlength=1mm
\begin{picture}(42,42)
\put(0,0){\includegraphics[width=42mm,height=42mm]{r12-3b.eps}}
\put(18,41){\makebox(0,0)[cc]{ $\theta$}}
\put(40.2,17.8){\makebox(0,0)[cc]{ $\xi$}}
\put(21,-7){\makebox(0,0)[cc]{Рис. 12.3}}
\end{picture}}
\quad
{\unitlength=1mm
\begin{picture}(42,42)
\put(0,0){\includegraphics[width=42mm,height=42mm]{r12-3c.eps}}
\put(18,41){\makebox(0,0)[cc]{ $\zeta$}}
\put(40.2,18){\makebox(0,0)[cc]{ $\eta$}}
\end{picture}}
\hfill\mbox{}
\\[7.75ex]
\indent
Проективно особая система
\\[2.25ex]
\mbox{}\hfill                                       %(12.14)
$
\displaystyle
\dfrac{dx}{dt}=1-(15+9a-8b)x+y-(5+a)x^2-cxy-x^2(x-y),
\hfill
$
\\[0.5ex]
\mbox{}\hfill (12.14)
\\[0.25ex]
\mbox{}\hfill
$
\displaystyle
\dfrac{dy}{dt}=10y-(5+a)xy-cy^2-xy(x-y),
\hfill
$
\\[2.75ex]
где 
\vspace{0.35ex}
$
a={}-10^{-13},\, 
b={}-10^{-52},\, 
c={}-10^{-200}, 
$
не имеет предельных циклов, лежащих в конечной части проективной фазовой плоскости $\R\P(x,y),$ так как 
\vspace{0.25ex}
у нее одно состояние равновесия, которое является седлом.
\vspace{0.35ex}

Второй проективно приведенной системой системы (12.14) является система 
\\[2.2ex]
\mbox{}\hfill                           %(12.15)
$
\displaystyle
\dfrac{d\eta}{d\nu}=c\eta-\zeta -10\eta^2+(5+a)\eta \zeta +\zeta^2, 
\ \ \
\dfrac{d\zeta}{d\nu}=\eta(1+\eta-(25+9a-8b)\zeta),
\ \ 
\eta\, d\nu =dt,
\hfill
$
\mbox{}\hfill (12.15)
\\[2.75ex]
где 
\vspace{0.5ex}
$
a={}-10^{-13},\, 
b={}-10^{-52},\, 
c={}-10^{-200}, 
$
у которой в конечной части проективной фазовой плоскости $\R\P(\eta,\zeta)$
\vspace{0.25ex}
два  состояния равновесия $O^{(2)}(0,0)$ и $A^{(2)}(0,1),$ являющиеся фокусами.
\vspace{0.35ex}

Система (12.15) в конечной части проективной фазовой плоскости $\R\P(\eta,\zeta)$
\vspace{0.25ex}
имеет не менее четырех предельных циклов [16]. 
\vspace{0.25ex}
Система (12.15) --- проективно неособая и  имеет бесконечно удаленное 
 состояние равновесия (седло [17, c. 182]). 
\vspace{0.25ex}
У системы (12.15) нет линейных и разомкнутых циклов (в том числе и предельных).
\vspace{0.25ex}
Используя фазовый портрет поведения траекторий системы (12.15) на проективном круге $\P\K(\eta,\zeta),$
\vspace{0.25ex}
построенный на рис. 7 в [17, c. 182] (на случай четырех предельных циклов), а также, 
\vspace{0.25ex}
приведенное на рис.~5.1, правило отображения проективных кругов 
\vspace{0.25ex}
$\P\K(x,y),\ \P\K(\xi,\theta),\ \P\K(\eta,\zeta),$ на рис.~12.4 построен проективный атлас траекторий системы (12.15).
\vspace{0.35ex}

Таким образом, система (12.14) не имеет 
\vspace{0.25ex}
линейных и лежащих в конечной части проективной  фазовой плоскости $\R\P(x,y)$ 
\vspace{0.25ex}
предельных циклов, а имеет не менее четырех 
разомкнутых предельных циклов.
\vspace{0.25ex}

Первой проективно приведенной системой системы (12.14) является проективно особая система 
\\[2.5ex]
\mbox{}\hfill                      %(12.16)
$
\dfrac{d\xi}{d\tau}=\xi\theta(25+9a -8b-\xi -\theta),
\hfill
$
\\[0.75ex]
\mbox{}\hfill(12.16)
\\[0.5ex]
\mbox{}\hfill
$
\dfrac{d\theta}{d\tau}=1-\xi +(5+a)\theta+c\xi \theta+(15+9a-8b) \theta^2-\xi\theta^2-\theta^3, \ \ \ \theta d\tau=dt,
\hfill
$
\\[2.75ex]
где 
\vspace{0.5ex}
$
a={}-10^{-13},\, 
b={}-10^{-52},\, 
c={}-10^{-200}. 
$
У системы (12.16) нет линейных циклов, но есть разомкнутые и лежащие в конечной части 
\vspace{0.25ex}
проективной  фазовой плоскости $\R\P(\xi,\theta)$ предельные циклы,  которых в сумме не менее четырех. 
\\[5.5ex]
\mbox{}\hfill
{\unitlength=1mm
\begin{picture}(42,42)
\put(0,0){\includegraphics[width=42mm,height=42mm]{r12-4a.eps}}
\put(18,41){\makebox(0,0)[cc]{ $y$}}
\put(40.2,18.5){\makebox(0,0)[cc]{ $x$}}
%\put(21,-3){\makebox(0,0)[cc]{ $Oxy$}}
%\put(22.5,-6){\makebox(0,0)[cc]{Рис. 1}}
\end{picture}}
\qquad
{\unitlength=1mm
\begin{picture}(42,42)
\put(0,0){\includegraphics[width=42mm,height=42mm]{r12-4b.eps}}
\put(18,41){\makebox(0,0)[cc]{ $\theta$}}
\put(40.2,18){\makebox(0,0)[cc]{ $\xi$}}
%\put(21,-3){\makebox(0,0)[cc]{ $O^{{}^{(1)}}uz$}}
\put(21,-10){\makebox(0,0)[cc]{Рис. 12.4}}
\end{picture}}
\qquad
{\unitlength=1mm
\begin{picture}(42,42)
\put(0,0){\includegraphics[width=42mm,height=42mm]{r12-4c.eps}}
\put(18,41){\makebox(0,0)[cc]{ $\zeta$}}
\put(40.2,18){\makebox(0,0)[cc]{ $\eta$}}
%\put(21,-3){\makebox(0,0)[cc]{ $O^{{}^{(2)}}zv$}}
%\put(22.5,-6){\makebox(0,0)[cc]{Рис. 1}}
\end{picture}}
\hfill\mbox{}
\\[12ex]
\centerline{
{\bf  13. 
Симметpичность фазового поля направлений
}
}
\\[1.75ex]
\indent
Через регулярную точку $M(x,y)$ 
\vspace{0.5ex}
фазовой плоскости $Oxy$ проведем отрезок прямой с направляющим вектором 
$\vec{a}(x,y)=(X(x,y), Y(x,y)).$
\vspace{0.5ex}
Будем считать, что отрезок --- ненаправленный единичной длины
с серединой в точке $M.$ 
\vspace{0.25ex}
В середине отрезка траектория системы (D) касается его (т.е. $M$ --- контактная точка отрезка).
\vspace{0.25ex}
Множество таких отрезков, построенных в каждой регулярной точке системы (D), 
\vspace{0.25ex}
назавем 
{\it фазовым полем направлений} системы (D). 
В состоянии равновесия системы (D) ее фазовое поле направлений не определено.
\vspace{0.35ex}

Из геометрических 
\vspace{0.15ex}
соображений получаем критерии симметричности фазового поля 
направлений системы (D).
\vspace{1ex}

{\bf Свойство 13.1} [18, с. 4 --- 5].
\vspace{0.5ex}
%\marginpar{[2150]}
{\it 
Фазовое поле направлений системы {\rm (D)} симметpично относительно пpямой
$Ax + By + C = 0 \ (|A| + |B| \ne 0)$ тогда и только тогда, когда
\\[2.5ex]
\mbox{}\hfill
$
R\Bigl(x - \dfrac{2A(Ax + By + C)}{A^{2} + B^{2}}\,, \, y
- \dfrac{2B(Ax + By + C)}{A^{2} + B^{2}} \Bigr)  =
\dfrac{2AB - \bigl(A^2 - B^2\,\bigr)R(x,y)}{A^{2} - B^{2} +2AB\,R(x,y)}
\ \
\forall (x,y)\! \in\! \Omega,
\hfill
$
\\[2.25ex]
где рациональная функция} 
\\[2.25ex]
\mbox{}\hfill
$
R\colon (x,y)\to\ \dfrac{Y(x,y)}{X(x,y)}
\quad 
\forall (x,y) \in \Omega=\{(x,y)\colon X(x,y)\ne 0\}.
\hfill
$
\\[2.25ex]
\indent
У системы (D) 
\vspace{0.15ex}
с фазовым полем направлений, симметричным относительно некоторой прямой $l,$ 
\vspace{0.15ex}
для каждой траектории существует траектория, симметричная ей относительно прямой $l.$ 
\vspace{0.15ex}
При этом не исключается возможность наличия траекторий, симметричных относительно прямой $l.$

Например, симметричность фазового поля направлений системы (D) 
\vspace{0.25ex}
относительно биссектрис координатных углов 
\vspace{0.15ex}
фазовой плоскости $Oxy$ устанавливается с помощью следующих свойств. 
\vspace{0.75ex}

{\bf Свойство 13.2.} 
\vspace{0.35ex}
{\it 
Фазовое поле направлений системы {\rm (D)} симметpично
относительно пpямой $y = x$ тогда и только тогда, когда
\\[1.5ex]
\mbox{}\hfill
$
X(x,y)\;\!X(y,x) \;\! -\;\!  Y(x,y)\;\!Y(y,x)  =  0
\quad 
\forall (x,y)\in \R^2.
\hfill
$
}
\\[2ex]
\indent
{\bf Свойство 13.3.} 
\vspace{0.35ex}
{\it 
Фазовое поле 
\vspace{0.1ex}
направлений системы {\rm (D)} симметpично
относительно пpямой $y ={} - x$ тогда и только тогда, когда
\\[1.5ex]
\mbox{}\hfill
$
X({}-x,{}-y)\;\!X(y,x)  \;\!-\;\! Y({}-x,{}-y)\;\!Y(y,x)  =  0
\quad
\forall (x,y)\in \R^2.
\hfill
$
}
\\[2ex]
\indent
Аналитическими критериями симметричности 
\vspace{0.15ex}
фазового поля направлений дифференциальной системы (D) 
\vspace{0.25ex}
относительно координатных осей и начала координат 
фазовой плоскости $Oxy$ являются 
\vspace{0.75ex}

{\bf Свойство 13.4} [12, с. 50 --- 51].
\vspace{0.35ex}
%\marginpar{[296]}
{\it 
Фазовое поле направлений системы {\rm (D)} 
симметpично относительно оси $Ox$ тогда и только тогда, когда
\\[1.5ex]
\mbox{}\hfill
$
X(x,y)\;\!Y(x,{}-y)  \;\! +\;\!  X(x,{}-y)\;\!Y(x,y)  =  0
\quad
\forall (x,y)\in \R^2.
\hfill
$
}
\\[2ex]
\indent
{\bf Свойство 13.5.} 
\vspace{0.35ex}
{\it 
Фазовое поле направлений  системы {\rm (D)} 
симметpично относительно оси $Oy$ тогда и только тогда, когда
\\[1.5ex]
\mbox{}\hfill
$
X(x,y)\;\!Y({}-x,y)\;\! +\;\! X({}-x,y)\;\!Y(x,y)  =  0
\quad
\forall (x,y)\in\R^2.
\hfill
$
}
\\[2ex]
\indent
{\bf Свойство 13.6.} 
\vspace{0.35ex}
{\it 
Фазовое поле направлений системы {\rm (D)} симметpично
относительно начала кооpдинат $O(0,0),$ 
если и только если
\\[1.5ex]
\mbox{}\hfill
$
X(x,y)\;\!Y({}-x,{}-y)\;\!  -\;\! X({}-x,{}-y)\;\!Y(x,y)  = 0
\quad
\forall (x,y)\in \R^2.
\hfill
$
}
\\[2ex]
\indent
Симметричность фазовых полей направлений первой и второй проективно приведенных систем относительно 
координатных осей 
\vspace{0.15ex}
и начала координат устанавливается с помощью свойств 13. 4 --- 13.6, 
\vspace{0.25ex}
а также может 
быть установлено на основании связей, содержащихся в свойствах 13.7 --- 13.9.
\vspace{1ex}

{\bf Свойство 13.7.} 
{\it 
Равносильными являются следующие утверждения\;\!{\rm:}
\vspace{0.5ex}

{\rm 1.} 
Фазовое поле направлений системы {\rm (D)} симметpично
\vspace{0.25ex}
относительно начала кооpдинат фазовой плоскости $Oxy;$ 
\vspace{0.75ex}

{\rm 2.} 
\vspace{0.75ex}
Выполняется тождество
$
X(x,y)Y({}-x,{}-y)  - X({}-x,{}-y)Y(x,y)  = 0
\;\;
\forall (x,y)\in \R^2;
$

{\rm 3.} 
Фазовое поле направлений первой проективно 
\vspace{0.35ex}
приведенной системы {\rm (6.3)} симметpично
относительно кооpдинатной оси 
$O_{\phantom1}^{(1)}\xi;$
 \vspace{1ex}

{\rm 4.} 
Выполняется тождество
\vspace{0.35ex}
$
\Xi(\xi,\theta)\;\!\Theta(\xi,{}-\theta) \;\!+\;\! \Xi(\xi,{}-\theta)\;\!\Theta(\xi,\theta)  = 0
\;\;
\forall (\xi,\theta)\in \R^2;
$

\newpage

{\rm 5.} 
Фазовое поле направлений второй 
\vspace{0.35ex}
проективно приведенной системы {\rm (6.6)} симметpично
относительно кооpдинатной оси $O_{\phantom1}^{(2)}\zeta;$
\vspace{0.75ex}
 
{\rm 6.} 
Выполняется тождество
$
H(\eta,\zeta)\;\! Z({}-\eta,\zeta)\;\! +\;\! H({}-\eta,\zeta)\;\! Z(\eta,\zeta)  = 0
\;\;
\forall (\eta,\zeta)\in \R^2.
$
}
\vspace{1.5ex}

{\bf Свойство 13.8.} 
{\it 
Равносильными являются следующие утверждения\;\!{\rm:}
\vspace{0.5ex}

{\rm 1.} 
Фазовое поле 
\vspace{0.25ex}
направлений системы {\rm (D)} симметpично
относительно координатной оси $Ox;$ 
\vspace{0.75ex}

{\rm 2.} 
Выполняется тождество
\vspace{0.75ex}
$
X(x,y)\;\!Y(x,{}-y)  + X(x,{}-y)\;\!Y(x,y)  = 0
\;\;
\forall (x,y)\in \R^2;
$

{\rm 3.} 
Фазовое поле направлений 
\vspace{0.35ex}
первой проективно приведенной системы {\rm (6.3)} симметpично
относительно кооpдинатной оси $O_{\phantom1}^{(1)}\theta;$
\vspace{0.75ex}
 
{\rm 4.} 
Выполняется тождество
\vspace{1ex}
$
\Xi(\xi,\theta)\;\!\Theta({}-\xi,\theta)\;\! +\;\! \Xi({}-\xi,\theta)\;\!\Theta(\xi,\theta)  = 0
\;\;
\forall (\xi,\theta)\in \R^2;
$

{\rm 5.} 
Фазовое поле направлений 
\vspace{0.35ex}
второй проективно приведенной системы {\rm (6.6)} симметpично
относительно начала координат фазовой плоскости $O_{\phantom1}^{(2)}\eta\zeta;$
\vspace{1ex}

{\rm 6.} 
Выполняется тождество
\vspace{1.5ex}
$
H(\eta,\zeta)\;\! Z({}-\eta,{}-\zeta) - H({}-\eta,{}-\zeta)\;\! Z(\eta,\zeta)  = 0
\;\;
\forall (\eta,\zeta)\in \R^2.
$
}

{\bf Свойство 13.9.} 
{\it 
Равносильными являются следующие утверждения\;\!{\rm:}
\vspace{0.5ex}

{\rm 1.} 
Фазовое поле 
\vspace{0.25ex}
направлений системы {\rm (D)} симметpично
относительно координатной оси $Oy;$ 
\vspace{0.75ex}

{\rm 2.} 
Выполняется тождество
\vspace{1ex}
$
X(x,y)Y({}-x,y)  + X({}-x,y)Y(x,y)  = 0
\;\;
\forall (x,y)\in \R^2;
$

{\rm 3.} 
Фазовое поле 
\vspace{0.35ex}
направлений первой проективно приведенной системы {\rm (6.3)} симметpично
относительно начала координат фазовой плоскости $O_{\phantom1}^{(1)}\xi\theta;$
\vspace{1ex}

{\rm 4.} 
Выполняется тождество
\vspace{1ex}
$
\Xi(\xi,\theta)\Theta({}-\xi,{}-\theta) - \Xi({}-\xi,{}-\theta)\Theta(\xi,\theta)  = 0
\;\;
\forall (\xi,\theta)\in \R^2;
$

{\rm 5.} 
Фазовое поле 
\vspace{0.35ex}
направлений второй проективно приведенной системы {\rm (6.6)} симметpично
относительно координатной оси $O_{\phantom1}^{(2)}\eta;$
 \vspace{1ex}

{\rm 6.} 
Выполняется тождество
$
H(\eta,\zeta) Z(\eta,{}-\zeta) + H(\eta,{}-\zeta) Z(\eta,\zeta)  = 0
\;\;
\forall (\eta,\zeta)\in \R^2.
$
}
\\[5ex]
\centerline{
{\bf  14. Множества проективно неособых и проективно особых систем
}
}
\\[1.5ex]
\indent
Пусть $l$--- некоторая прямая плоскости $Oxy.$ 
\vspace{0.25ex}
Тогда существует линейное невырожденное преобразование плоскости $Oxy$ такое, 
\vspace{0.25ex}
что в новой системе координат прямая $l$ будет осью ординат (абсцисс).

Учитывая это обстоятельство, на основании свойств 10.3 и 11.3 получаем
\vspace{0.5ex}

{\bf Свойство 14.1.}
\vspace{0.25ex}
{\it 
Если  некоторая прямая фазовой плоскости $Oxy$ состоит из траекторий системы {\rm (D)}, 
\vspace{0.25ex}
то существует линейное невырожденное преобразование, с помощью которого и 
первого {\rm (}второго{\rm )}  преобразования Пуанкаре система {\rm (D)} 
\vspace{0.25ex}
приводится к проективно неособой системе}.
\vspace{0.5ex}

Все возможные прямые фазовой плоскости системы (D) 
\vspace{0.15ex}
не могут одновременно состоять из ее траекторий.
\vspace{0.15ex}
Поэтому всегда можно указать линейное невырожденное преобразование фазовой плоскости $Oxy,$
\vspace{0.15ex}
с помощью которого система (D) приводится к такой системе, что в новой  системе координат ось ординат (абсцисс) 
\vspace{0.15ex}
не будет состоять из ее траекторий. На основании свойств 10.3, 10.4 и 11.3, 11.4 получаем 
\vspace{0.5ex}

{\bf Свойство 14.2.}
\vspace{0.15ex}
{\it 
Для системы {\rm (D)} всегда можно указать линейное невырожденное преобразование, переводящее ее
\vspace{0.15ex}
в такую систему, у которой первая {\rm (}вторая{\rm )} проективно приведенная система будет проективно особой}. 
\vspace{0.5ex}

Иначе говоря, любая система (D)  
\vspace{0.15ex}
с помощью линейного невырожденного преобразования и первого (второго) преобразования Пуанкаре 
приводится к проективно особой системе.

Пусть $D$ --- 
множество всех обыкновенных автономных дифференциальных систем второго порядка 
с полиномиальными правыми частями любой степени;
\vspace{0.15ex}

$A\subset D$ --- множество всех систем, которые являются проективно особыми или с помощью линейного невырожденного преобразования и первого или второго преобразования Пуанкаре приводятся к проективно особым системам;
\vspace{0.15ex}

$B\subset D$ --- множество всех систем, которые являются проективно неособыми или с помощью линейного невырожденного преобразования и первого или второго преобразования Пуанкаре приводятся к проективно неособым системам;
\vspace{0.15ex}

$C\subset D$ --- множество всех систем, которые являются проективно особыми или с помощью линейного невырожденного преобразования и преобразований Пуанкаре не приводятся к проективно неособым системам.
\vspace{0.15ex}

Тогда $B\cap C=\O,\ B\cup C=D.$ 
\vspace{0.15ex}
По свойству 14.2, $A=D,$ а значит, $B\subset A,\ C\subset A,$ при этом дополнением множества $B$ до множества 
\vspace{0.15ex}
$D=A$ является множество $C,$ дизъюнктивное с множеством $B.$

Докажем существование систем, принадлежащих множеству $C.$ При этом будем использовать
\vspace{0.25ex}

{\bf Свойство 14.3.}
\vspace{0.25ex}
{\it 
Множеству $C$ принадлежит та и только та система, у которой в проективной фазовой плоскости нет прямых,
состоящих из ее траекторий}. 
\vspace{0.35ex}

{\sl Действительно}, 
\vspace{0.15ex}
по свойству 9.2, бесконечно удаленная прямая не состоит из траекторий системы (D) 
тогда и только тогда, когда система (D) проективно особая. 
\vspace{0.15ex}
Поэтому каждая система множества $C$ 
\vspace{0.15ex}
такая, что бесконечно удаленная прямая не состоит из ее траекторий.
В противном случае, т.е. когда бесконечно удаленная прямая состоит из траекторий, 
система является проективно неособой (свойство 9.1) 
\vspace{0.15ex}
и, следовательно, множеству $C$ не принадлежит.
\vspace{0.25ex}

Согласно свойству 10.4 (свойству 11.4) 
первая (вторая) приведенная система является проективно особой тогда и только тогда, когда ось ординат (абсцисс) не состоит из траекторий.
Если на фазовой плоскости $Oxy$ 
\vspace{0.15ex}
существует прямая, состоящая из траекторий системы (D),  
\vspace{0.15ex}
то, по свойству 14.1, система (D) приводится к проективно неособой системе. 
Поэтому каждая система множества $C$ 
\vspace{0.15ex}
такая, что на фазовой плоскости нет прямых, 
которые состояли бы из ее траекторий. 
\vspace{0.15ex}
Наличие же прямых, состоящих из траекторий системы (D), означает, 
\vspace{0.25ex}
что система (D)  приводится к проективно неособой системе, а значит, система ${\rm (D)}\not\in C.\ \k$
\vspace{0.75ex}

У каждой проективно особой 
\vspace{0.25ex}
дифференциальной системы (12.8) --- (12.13) на проективной фазовой плоскости нет прямых, состоящих из ее траекторий. 
\vspace{0.35ex}
По свойству 14.3, системы (12.8) --- (12.13)
принадлежат множеству $C,$ а поэтому $C\ne\O.$
\vspace{1ex}

{\bf Теорема 14.1.}
{\it 
Система} (D) {\it при $n=0,1,2$ не принадлежит множеству} $C.$
\vspace{0.75ex}

{\sl Доказательство}.
\vspace{0.35ex}
Система (D) при $n=0$ имеет вид (7.8)
и является проективно неособой (пример~7.1), а значит, она не принадлежит множеству $C.$
\vspace{0.5ex}

При $n=1$ 
\vspace{0.35ex}
система (D) является линейной стационарной системой (7.11), которая будет проективно особой системой при 
\vspace{0.75ex}
$a_2^{}=b_1^{}=0,\, b_2^{}=a_1^{}\ne 0$ (пример 7.2). 
Однако, в этом случае любая прямая семейства
\vspace{0.75ex}
$C_1^{}(a_1^{}y+b_{_0}^{})+C_2^{}(a_1^{}x+a_{_0}^{})=0$
состоит из траекторий системы (7.11). По свойству 14.1, система (7.11) 
\vspace{0.5ex}
не принадлежит классу $C.$

Если $n=2,$ 
\vspace{0.5ex}
то система (D) имеет вид (7.16) и является проективно особой при 
$a_5^{}=b_3^{}=0,\ b_4^{}=a_3^{},\ b_5^{}=a_4^{},\ |a_3^{}|+|a_4^{}|\ne 0$
\vspace{0.75ex}
(пример 7.3). Уравнением траекторий проективно особой системы (7.16)
является уравнение Якоби, 
\vspace{0.35ex}
которое всегда имеет прямую, состоящую из его траекторий
\vspace{0.35ex}
[10, с. 15 -- 16]. По свойству 14.1, система (7.16) не принадлежит множеству $C.$\k

\newpage

\mbox{}
\\[-2ex]
\centerline{
{\bf  15. Топологическая эквивалентность дифференциальных систем }
}
\\[0.1ex]
\centerline{
{\bf на проективном круге и на проективной сфере 
}
}
\\[1.5ex]
\indent
В зависимости от того рассматривается ли поведение траекторий на проективном круге или на 
проективной сфере  будем различать их поведение с точностью до топологической эквивалентности.
\vspace{0.35ex}

{\bf Определение 15.1} [11, с. 34].
{\it 
Две автономные полиномиальные дифференциальные системы второго порядка 
\textit{\textbf{ топологически эквивалентны на проективном круге}}, если существует гомеоморфизм 
проективных кругов, переводящий траектории одной системы в траектории другой системы}.
\vspace{0.35ex}

Поскольку линейное невырожденное преобразование плоскости является гомеоморфизмом, то имеет место
\vspace{0.35ex}

{\bf Свойство 15.1.}
{\it 
При линейном невырожденном преобразовании фазовой плоскости $(x,y)$ сохраняется{\rm:}
{\rm а)} топологическая эквивалентность системы {\rm (D);} 
{\rm б)} степень $n$ системы {\rm (D);} 
{\rm в)} тип системы  {\rm (D)} на проективной фазовой плоскости.
}
\vspace{0.5ex}

Например, топологически эквивалентными на проективном круге являются системы (7.8) и
\\[1.5ex]
\mbox{}\hfill
$
\dfrac{dx}{dt}=\widetilde{a}_{_0},
\qquad 
\dfrac{dy}{dt}=\widetilde{b}_{_0},
\qquad
|\widetilde{a}_{_0}|+|\widetilde{b}_{_0}|\ne0.
\hfill
$
\\[2.25ex]
\indent
Действительно, параллельным переносом $v=x+\widetilde{a}_{_0}-a_{_0},\ w=y+\widetilde{b}_{_0}-b_{_0}$ получаем:
\\[2.25ex]
\mbox{}\hfill
$
\dfrac{dv}{dt}=\dfrac{dx}{dt}+\widetilde{a}_{_0}-a_{_0}=\widetilde{a}_{_0},
\qquad 
\dfrac{dw}{dt}=\dfrac{dy}{dt}+\widetilde{b}_{_0}-b_{_0}=\widetilde{b}_{_0}.
\hfill
$
\\[2.25ex]
\indent
Поэтому фазовый портрет траекторий системы (7.8) с точностью до топологической эквивалентности на 
проективном круге построен как на  круге $\P\K(x,y)$ из рис. 8.1, так и на  круге $\P\K(x,y)$ из рис. 8.2.
\vspace{0.35ex}

{\bf Определение 15.2.}
{\it 
Две автономные полиномиальные дифференциальные системы второго порядка 
\textit{\textbf{ топологически эквивалентны на проективной сфере}}, если существует гомеоморфизм 
проективных сфер, переводящий траектории одной системы в траектории другой системы}.
\vspace{0.35ex}

Топологическая эквивалентность поведения траекторий систем  на проективной сфере  
сохраняется при суперпозиции линейного невырожденного преобразования и преобразований Пуанкаре. 
Поэтому справедливы следующие два свойства.
\vspace{0.35ex}

{\bf Свойство 15.2.}
{\it 
Системы {\rm (D), (6.3)} и {\rm (6.6)} 
топологически эквивалентны на проективной сфере.
}
\vspace{0.35ex}

{\bf Свойство 15.3.}
{\it 
У топологически эквивалентных на проективной сфере  систем сумма линейных, 
разомкнутых и обычных предельных циклов одинакова.
}
\vspace{0.35ex}

Например, в соответствии со свойством 15.2 топологически эквивалентными 
на проективной сфере являются система (8.1) с проективным атласом на рис. 8.3 
и система (8.2) с проективным атласом на рис. 8.4. 
При этом система (8.2) является первой проективно приведенной системой системы (8.1), 
а система (8.1) является второй проективно приведенной системой системы (8.2).
\vspace{0.5ex}

Если ${\rm deg}\,({\rm D})=n,$ то систему (D) будем обозначать $({\rm D}^n).$
\vspace{0.5ex}

{\bf Свойство 15.4.}
\vspace{0.25ex}
{\it 
Проективно неособая системы $({\rm D}^n)$ топологически эквивалентна на проективной сфере проективно особой системе 
$({\rm D}^{n+1}),$ причем система $({\rm D}^{n+1})$ имеет прямую, состоящую из ее траекторий.
}
\vspace{0.35ex}

{\sl Доказательство} является следствием свойств 7.3, 10.4, 11.3, 15.1 и теоремы 4.1. \k
\vspace{0.5ex}

Это свойство позволяет вместо качественного исследования в целом проективно особой системы при наличии прямой, 
состоящей из ее траекторий, выполнить качественное исследование в целом соответствующей проективно неособой системы. 

Для систем класса $A$ такая возможность перехода к исследованию проективно неособых систем исключается.

Например, 
\vspace{0.25ex}
система Якоби (7.21)  имеет прямую, состоящую из ее траекторий [10, 
\linebreak
c. 15 -- 16]. Тогда согласно свойству 15.4 справедливо
\vspace{0.5ex}

{\bf Свойство 15.5.}
{\it 
Cистема Якоби {\rm(7.21)} топологически эквивалентна на проективной сфере линейной стационарной системе 
{\rm(7.11)}.
}
\vspace{0.5ex}

В этой связи систему Якоби назовем {\it проективно линейной дифференциальной системой}.

В соответствии со свойством 7.4 система $({\rm D}^{2})$ 
\vspace{0.15ex}
будет проективно особой,
если и только если она является системой Якоби.
Поэтому на основании свойства 15.5 можем утверждать следующее.
\vspace{0.5ex}

{\bf Свойство 15.6.}
{\it 
Проективно особая система $({\rm D}^{2})$  топологически эквивалентна на проективной сфере 
 системе $({\rm D}^{1}).$
}
\vspace{0.75ex}

{\bf Свойство 15.7.}
{\it 
Проективно особая система $({\rm D}^{2})$ является 
проективно линейной дифференциальной системой.
}
\vspace{0.5ex}

При переходе к топологически эквивалентным на проективной сфере системам 
могут быть использованы дробно-линейные преобразования
\\[2.25ex]
\mbox{}\hfill               %(15.1)
$
\xi=\dfrac{{}-\beta x+\alpha y}{\alpha x+\beta y+\gamma}\,,
\qquad 
\theta=\dfrac{\sqrt{\alpha^2+\beta^2}}{\alpha x+\beta y+\gamma}
$
\hfill (15.1)
\\[1.5ex]
и 
\\[1.5ex]
\mbox{}\hfill               %(15.2)
$
\eta=\dfrac{\sqrt{\alpha^2+\beta^2}}{\alpha x+\beta y+\gamma}\,,
\qquad 
\zeta= \dfrac{\beta x-\alpha y}{\alpha x+\beta y+\gamma}\,,
$
\hfill (15.2)
\\[2.5ex]
где $\alpha^2+\beta^2\ne 0,$ 
\vspace{0.35ex}
с помощью которых прямую $\alpha x+\beta y+\gamma=0$ переводим 
в бесконечно удаленные прямые проективных фазовых плоскостей 
\vspace{0.35ex}
$\!\R\P(\xi,\theta)\!\!$ и $\!\R\P(\eta,\zeta)\!$ соответственно.

В частности, 
\vspace{0.25ex}
если прямая $\alpha x+\beta y+\gamma=0$ состоит из траекторий системы (D), 
то систему (D) преобразованиями (15.1)  и (15.2) приводим к проективно неособым системам.
\\[4.25ex]
\centerline{
{\bf  16. Примеры глобального качественного исследования траекторий
}
}
\\[0.35ex]
\centerline{
{\bf
дифференциальных систем на проективной фазовой плоскости}
}
\\[1.75ex]
\indent
{\bf 16.1.
Траектории дифференциальной системы Дарбу} [14, с. 109 -- 111; 10]
\\[2ex]
\mbox{}\hfill        % (16.1)
$
\dfrac{dx}{dt}={}-y+x^3\equiv X(x,y),
\qquad 
\dfrac{dy}{dt}=x(1+xy)\equiv Y(x,y).
$
\hfill(16.1)
\\[2.25ex]
\indent
{\it Интегральный базис системы} (16.1). 
На любой области множества $\{(t,x,y)\colon x\ne 0\}$
функционально независимые первые интегралы
\\[1.75ex]
\mbox{}\hfill
$
F_1^{}\colon (t,x,y)\to\ 
{}-t+\arctg \dfrac{y}{x}
\quad
\forall (t,x,y) \in \{(t,x,y)\colon x\ne 0\}
\hfill
$
\\[1.5ex]
и 
\\[1.5ex]
\mbox{}\hfill
$
F_2^{}\colon (t,x,y)\to\
 t+ \dfrac{1+xy}{x^2+y^2}
\quad
\forall (t,x,y) \in \{(t,x,y)\colon |x|+|y|\ne 0\}
\hfill
$
\\[1.75ex]
образуют интегральный базис [19] системы (16.1).
\vspace{0.5ex}

{\it Общий автономный интеграл системы} (16.1).  
\vspace{0.35ex}
На любой области из множества $\{(x,y)\colon x\ne 0\}$
трансцендентная функция
\\[1.75ex]
\mbox{}\hfill
$
F\colon (x,y)\to\
\dfrac{1+xy}{x^2+y^2}+\arctg \dfrac{y}{x}
\quad
\forall (x,y)\in\{(x,y)\colon x\ne 0\}
\hfill
$
\\[1.75ex]
является общим автономным интегралом [20, c. 112 -- 114] системы (16.1).

{\it Состояния равновесия системы} (16.1) 
\vspace{0.25ex}
{\it  в конечной части проективной фазовой плоскости} $\R\P(x,y).$ 
Кубическая парабола $y=x^3$ 
\vspace{0.25ex}
не имеет общих точек с гиперболой $xy+1=0,$
\vspace{0.25ex}
а с прямой $x=0$ пересекается в точке $O(0,0).$
В конечной части проективной фазовой плоскости $\R\P(x,y)$ 
\vspace{0.35ex}
у системы (16.1) одно состояние равновесия $O(0,0).$

На основании характеристического уравнения $\lambda^2+1=0$ 
\vspace{0.25ex}
устанавливаем, что 
$O(0,0)$ является центром или фокусом [3, c. 139 -- 145].
\vspace{0.25ex}

Семейство траекторий системы (16.1) в совмещенной полярной системе координат 
$O\rho\varphi$ задается уравнением 
\\[1.5ex]
\mbox{}\hfill
$
\rho^{{}-2}=C-\varphi-\dfrac{1}{2}\, \sin 2\varphi,
\hfill
$
\\[1.5ex]
а значит,  $O(0,0)$ --- неустойчивый фокус.
\vspace{0.5ex}

{\it Движение радиуса-вектора образующей точки вдоль траекторий системы} (16.1).
Функция 
\\[1.5ex]
\mbox{}\hfill
$
W\colon (x,y)\to\  x\;\!Y(x,y)-y\;\!X(x,y)=x^2+y^2>0
\quad  
\forall (x,y)\in \R^2\backslash\{(0,0)\}.
\hfill
$
\\[1.5ex]
\indent
При движении вдоль траекторий системы (16.1) угол между 
\vspace{0.15ex}
радиусом-век\-то\-ром образующей точки и положительным направлением оси $Ox$ возрастает.
\vspace{0.5ex}

{\it Симметричность фазового поля направлений системы}  (16.1).
Фазовое поле направлений системы (16.1) симметрично 
относительно начала координат
фазовой плоскости $Oxy$ (свойство 13.6).

Для каждой траектории системы (16.1) 
\vspace{0.15ex}
существует траектория, которая симметрична ей относительно 
начала координат фазовой плоскости $Oxy.$ 
\vspace{0.5ex}

{\it Нулевые и ортогональные изоклины системы}  (16.1).\!
\vspace{0.25ex}
Разрешив уравнение $\!Y(x,y)\!=\!0$ при $X(x,y)\ne 0,$ находим, что нулевыми изоклинами системы (16.1) 
\vspace{0.35ex}
являются гипербола $xy+1=0$ и прямая $x=0,$ из которой удалена точка $O(0,0).$
\vspace{0.25ex}

Касательная к траектории системы (16.1) в
\vspace{0.25ex}
каждой точке гиперболы $xy+1=0$ и   в каждой точке прямой $x=0,$ 
\vspace{0.25ex}
отличной от начала координат фазовой плоскости $Oxy,$ параллельна оси $Ox.$
\vspace{0.25ex}

Из уравнения $X(x,y)=0$ при $Y(x,y)\ne 0$ 
\vspace{0.25ex}
находим, что ортогональной изоклиной системы (16.1) 
является кубическая парабола $y=x^3,$ 
\vspace{0.25ex}
из которой удалена точка $O(0,0).$

Касательная к траектории системы (16.1) 
\vspace{0.25ex}
в каждой точке кубической параболы 
\linebreak
$y=x^3,$
отличной от начала координат фазовой плоскости $Oxy,$ 
\vspace{0.5ex}
параллельна оси  $Oy.$

{\it Области знакоопределенности фазового поля направлений  системы}  (16.1).
\vspace{0.25ex}
Из неравенства $X(x,y)\;\!Y(x,y)>0$
\vspace{0.25ex}
находим, что областями положительности фазового поля направлений системы (16.1) являются области 
\\[2ex]
\mbox{}\hfill
$
\Omega_{1}^{+}=\Bigl\{(x,y)\colon x>0\ \&\ {}-\dfrac1{x}<y<x^3\Bigr\}
$ 
\ \ и \ \  
$
\Omega_{2}^{+}=\Bigl\{(x,y)\colon x<0\ \&\ x^3<y<{}-\dfrac1{x}\Bigr\}. 
\hfill
$
\\[2ex]
\indent
Касательная к траектории системы (16.1) в каждой точке множества 
$
\Omega_{}^{+}=\Omega_{1}^{+}\sqcup\Omega_{2}^{+}
$
образует острый угол с положительным направлением оси $Ox.$
\vspace{0.35ex}

Разрешив неравенство  $X(x,y)\;\!Y(x,y)<0,$ 
\vspace{0.25ex}
получаем, что областями
отрицательности фазового поля направлений  системы (16.1) являются области
\\[2ex]
\mbox{}\hfill
$
\Omega_{1}^{-}=\{(x,y)\colon x>0\ \&\ y>x^3\},
\qquad
\Omega_{2}^{-}=\Bigl\{(x,y)\colon x<0\ \&\ y>{}-\dfrac1{x}\Bigr\},
\hfill
$
\\[2.5ex]
\mbox{}\hfill
$ 
\Omega_{3}^{-}=\{(x,y)\colon x<0\ \&\ y<x^3\}, 
\qquad 
\Omega_{4}^{-}=\Bigl\{(x,y)\colon x>0\ \&\ y<{}-\dfrac1{x}\Bigr\}.
\hfill
$
\\[2ex]
\indent
В каждой точке множества 
\vspace{0.35ex}
$
\Omega_{}^{-}=\Omega_{1}^{-}\sqcup\Omega_{2}^{-}\sqcup\Omega_{3}^{-}\sqcup\Omega_{4}^{-}
$
касательная к траектории системы (16.1) образует  тупой угол   с положительным направлением оси $Ox.$
\vspace{0.5ex}

{\it  Контактные точки системы} (16.1) 
\vspace{0.5ex}
{\it координатных осей фазовой плоскости} $Oxy.$ 
Нулевые изоклины 
$\{(x,y)\colon x=0\ \&\ y\ne 0\}$ и $\{ (x,y) \colon  xy+1=0\}$
\vspace{0.5ex}
не имеют общих точек с осью $Ox.$
На оси $Ox$ нет контактных точек системы (16.1).
\vspace{0.5ex}

Ортогональная изоклина 
\vspace{0.5ex}
$\{ (x,y) \colon  y=x^3\ \&\ x\ne 0\}$ не имеет общих точек с осью $Oy.$ 
На оси $Oy$  нет контактных точек  системы (16.1).
\vspace{0.75ex}

{\it  Проективный тип системы} (16.1).  
Функция 
\\[2ex]
\mbox{}\hfill
$
W_3^{}\colon  (x,y)\to\  
x\;\!Y_3^{}(x,y)-y\;\!X_3^{}(x,y)=0
\quad
\forall (x,y)\in \R^2.
\hfill
$
\\[2ex]
\indent
Система (16.1) --- проективно особая.

Бесконечно удаленная прямая проективной фазовой плоскости $\R\P(x,y)$
\vspace{0.15ex}
не состоит из  траекторий системы (16.1) (свойство 9.2). 
\vspace{0.75ex}

{\it Первая проективно приведенная система системы} (16.1). 
\vspace{0.5ex}
Проективно особую систему (16.1) первым преобразованием 
Пуанкаре $x=\dfrac{1}{\theta}\,,\, y=\dfrac{\xi}{\theta}$
\vspace{0.5ex}
приводим к первой проективно приведенной  системе (свойство 9.2)
\\[2.5ex]
\mbox{}\hfill        % (16.2)
$
\dfrac{d\xi}{d\tau}= \theta+\xi^2\;\!\theta\equiv \Xi(\xi,\theta),
\qquad 
\dfrac{d\theta}{d\tau}={}-1+\xi\;\!\theta^2\equiv \Theta(\xi,\theta),
\qquad
\theta\,d\tau=dt.
$
\hfill (16.2)
\\[3ex]
\indent
{\it Вторая проективно приведенная система системы} (16.1). 
\vspace{0.5ex}
Проективно особую систему (16.1) вторым преобразованием 
Пуанкаре  $x=\dfrac{\zeta}{\eta}\,,\, y=\dfrac{1}{\eta}$
\vspace{0.5ex}
приводим к второй проективно приведенной  системе (свойство 9.2)
\\[2.5ex]
\mbox{}\hfill        % (16.3)
$
\dfrac{d\eta}{d\nu}= {}-\zeta^2-\eta^2\zeta\equiv H(\eta,\zeta),
\qquad 
\dfrac{d\zeta}{d\nu}={}-\eta-\eta\zeta^2 \equiv Z(\eta,\zeta),
\qquad
\eta\,d\nu=dt.
$
\hfill (16.3)
\\[2.5ex]
\indent
{\it Состояния равновесия системы} (16.2) 
\vspace{0.35ex}
{\it на координатной оси $O^{(1)}_{\phantom1} \xi$
в конечной части проективной фазовой плоскости} $\R\P(\xi,\theta).$ 
Так как 
\\[2ex]
\mbox{}\hfill
$
\Theta(\xi,0)={}-1\ne 0  
\quad 
\forall  \xi \in \R,
\hfill
$
\\[1.75ex]
то на координатной оси $O^{(1)}_{\phantom1} \xi$
\vspace{0.5ex}
в конечной части проективной фазовой плоскости $\R\P(\xi,\theta)$ 
у системы (16.2) нет  состояний равновесия. 
\vspace{0.75ex}

{\it Состояние равновесия системы} (16.3) 
\vspace{0.25ex}
{\it в начале координат  фазовой плоскости} $O^{(2)}_{\phantom1} \eta\zeta.$
\vspace{0.5ex}
Начало координат фазовой плоскости $O^{(2)}_{\phantom1} \eta\zeta$ является сложным 
 состоянием равновесия системы (16.3) с характеристическим уравнением $\lambda^2=0.$ 
\vspace{0.35ex}
По теореме 6.2.1 из [21, c. 128 -- 129], состояние равновесия 
\vspace{0.25ex}
$O^{(2)}_{\phantom1} (0,0)$ системы (16.3) --- двухсепаратрисное седло, 
сепаратрисы которого  примыкают в направлении координатной оси 
$O^{(2)}_{\phantom1} \zeta.$
\vspace{0.75ex}

{\it Состояния равновесия системы} (16.1) 
\vspace{0.25ex}
{\it в проективной фазовой плоскости} $\R\P(x,y).$ 
У дифференциальной системы (16.1) в проективной фазовой плоскости $\R\P(x,y)$ 
\vspace{0.25ex}
два состояния равновесия: 
неустойчивый фокус $O(0,0)$ 
\vspace{0.25ex}
и двухсепаратрисное седло, лежащее на <<концах>>  прямой $x=0.$
\vspace{0.5ex}

{\it Экваториально контактные точки  системы} (16.1). 
Уравнение 
\\[1.75ex]
\mbox{}\hfill
$
X_3^{}(1,\xi)=0
\hfill
$ 
\\[1.75ex]
не имеет корней. 
\vspace{0.25ex}
У системы (16.1) нет  экваториально контактных точек на <<концах>>  прямых $y=ax,\ a\in\R.$
\vspace{0.25ex}
На <<концах>> оси $Oy$ 
лежит двухсепаратрисное седло, сепаратрисы которого ортогональны оси $Oy.$
\vspace{0.15ex}
У проективно особой системы (16.1) нет  экваториально контактных точек.

\newpage

{\it Предельные циклы системы} (16.1) {\it в проективной фазовой плоскости $\R\P(x,y).$}
У векторного поля 
\\[1.75ex]
\mbox{}\hfill
$
\vec{a}\colon (x,y)\to\ 
({}-y+x^3,x+x^2y)
\quad
\forall (x,y)\in\R^2
\hfill
$
\\[1.5ex]
расходимость 
\\[1.25ex]
\mbox{}\hfill
$
{\rm div}\ \vec{a}(x,y)=4x^2
\quad
\forall (x,y)\in\R^2
\hfill
$
\\[2ex]
является положительной на односвязной области $\R^2\backslash\{(0,0)\}.$
\vspace{0.35ex}

По признаку Бендиксона [4, c. 120], система (16.1) 
\vspace{0.25ex}
не имеет предельных циклов в конечной части проективной фазовой плоскости $\R\P(x,y).$
\vspace{0.25ex}

Среди траекторий проективно особой системы (16.1) 
\vspace{0.25ex}
на проективной фазовой плоскости $\R\P(x,y)$ нет прямых.
Система (16.1)  не имеет предельных линейных  циклов. 
\vspace{0.25ex}

Проективно особая система (16.1) имеет одно бесконечно удаленное состояние равновесия, которое является 
двухсепаратрисным седлом, и не имеет экваториально контактных точек.
Предельный цикл не может окружать двухсепаратрисное седло.
У системы (16.1)  нет предельных разомкнутых циклов.
\vspace{0.5ex}

{\it Общий автономный интеграл системы} (16.2). 
Функция
\\[2ex]
\mbox{}\hfill
$
F^{(1)}_{\phantom1} \colon (\xi,\theta)\to\
\dfrac{\xi +\theta^2}{1+\xi^2}+\arctg \xi 
\quad 
\forall (\xi,\theta)\in \R^2
\hfill
$
\\[2ex]
является общим автономным интегралом системы (16.2) на $\R^2.$
\vspace{0.75ex}

{\it Состояния равновесия системы} (16.2) 
\vspace{0.25ex}
{\it в проективной фазовой плоскости} $\R\P(\xi,\theta).$ 
У системы (16.2) в проективной фазовой плоскости $\R\P(\xi,\theta)$ 
\vspace{0.25ex}
 два состояния равновесия:
неустойчивый фокус на <<концах>>  оси $O^{(1)}_{\phantom1}\theta$ и 
\vspace{0.25ex}
двухсепаратрисное седло  на <<концах>>  оси $O^{(1)}_{\phantom1}\xi,$ 
сепаратрисы которого примыкают в направлении оси $O^{(1)}_{\phantom1}\xi.$ 
\vspace{0.75ex}

{\it Движение радиуса-вектора образующей точки вдоль траекторий системы} (16.2).
Функция 
\\[1.25ex]
\mbox{}\hfill
$
W^{(1)}_{\phantom1}\colon (\xi,\theta)\to\ 
\xi \,\Theta(\xi,\theta)-\theta \,\Xi(\xi,\theta)={}-\xi +\theta^2
\quad 
\forall (\xi,\theta)\in \R^2
\hfill
$
\\[2ex]
является отрицательной, когда $\xi>\theta^2,$ а положительной, когда $\xi<\theta^2.$
\vspace{0.5ex}

При движении вдоль частей траекторий системы (16.2), 
\vspace{0.5ex}
расположенных в области: 
\linebreak
а) $\{(\xi,\theta)\colon \xi>\theta^2\};$  \
б) $\{(\xi,\theta)\colon \xi<\theta^2\},$ 
\vspace{0.35ex}
угол между радиусом-век\-то\-ром образующей точки и положительным направлением оси 
$O^{(1)}_{\phantom1}\xi\colon\!\!$  а) убывает; \  б) возрастает. 
\vspace{0.5ex}

Через каждую точку, лежащую на параболе 
\vspace{0.35ex}
$\xi=\theta^2,$ траектория системы (16.2) проходит в направлении 
радиуса-вектора этой точки.
\vspace{0.75ex}

{\it Симметричность фазового поля направлений системы}  (16.2).
\vspace{0.15ex}
Фазовое поле направлений системы (16.2) симметрично 
относительно  координатной оси 
\vspace{0.15ex}
$O^{(1)}_{\phantom1} \xi$ (свойство 13.6).
Для каждой траектории системы (16.2) 
\vspace{0.25ex}
существует траектория, которая симметрична ей относительно 
координатной оси $O^{(1)}_{\phantom1} \xi.$
\vspace{0.15ex}
Каждая траектория системы (16.2), пересекающая координатную  ось $O^{(1)}_{\phantom1} \xi,$ симметрична 
относительно этой координатной оси.
\vspace{0.75ex}

{\it Нулевые и ортогональные изоклины системы}  (16.2).
\vspace{0.5ex}
Разрешая уравнение $\Theta(\xi,\theta)=0$ относительно $\xi,$ получаем, 
\vspace{0.5ex}
что нулевой изоклиной системы (16.2) является квадратичная гипербола
$\xi=\dfrac{1}{\theta^2}\,.$ 
\vspace{0.5ex}
Касательная к траектории системы (16.2) в каждой точке квадратичной гиперболы 
$\xi=\dfrac{1}{\theta^2}$ параллельна оси $O_{\phantom1}^{(1)}\xi.$
\vspace{0.75ex}

Из уравнения $\Xi(\xi,\theta)=0$ 
\vspace{0.35ex}
находим ортогональную изоклину $\theta=0$ системы (16.2).
Траектории системы (16.2) пересекают координатную ось  $O_{\phantom1}^{(1)}\xi$
под прямым углом.
\vspace{0.5ex}

{\it Области знакоопределенности фазового поля направлений  системы}  (16.2).
\vspace{0.35ex}
Из неравенства $\Xi(\xi,\theta)\;\!\Theta(\xi,\theta)>0$
\vspace{0.25ex}
находим, что областями положительности фазового поля нап\-рав\-ле\-ний системы (16.2) являются области 
\\[1.75ex]
\mbox{}\hfill
$
\Omega_{1}^{{}+(1)}=\Bigl\{(\xi,\theta)\colon \xi>\dfrac1{\theta^2}\ \ \& \ \theta>0\Bigr\}
$ 
\ \ 
и 
\ \ 
$
\Omega_{2}^{{}+(1)}=\Bigl\{(\xi,\theta)\colon \xi<\dfrac1{\theta^2}\ \ \& \ \theta<0\Bigr\}.
\hfill
$
\\[2ex]
\indent
В каждой точке множества 
\vspace{0.25ex}
$
\Omega_{\phantom1}^{{}+(1)}=\Omega_{1}^{{}+(1)}\sqcup\Omega_{2}^{{}+(1)}
$
касательная к траектории системы (16.2)
образует острый угол с положительным направлением оси $O_{\phantom1}^{(1)}\xi.$
\vspace{0.75ex}

Разрешая неравенство  $\Xi(\xi,\theta)\;\!\Theta(\xi,\theta)<0,$ 
\vspace{0.35ex}
получаем, что областями
отрицательности фазового поля направлений  системы (16.2) являются области
\\[1.75ex]
\mbox{}\hfill
$
\Omega_{1}^{{}-(1)}=\Bigl\{(\xi,\theta)\colon \xi<\dfrac1{\theta^2}\ \ \& \ \theta>0\Bigr\}
$
\ \ и \ \ 
$ 
\Omega_{2}^{{}-(1)}=\Bigl\{(\xi,\theta)\colon \xi>\dfrac1{\theta^2}\ \ \& \ \theta<0\Bigr\}.
\hfill
$
\\[2ex]
\indent
В каждой точке множества 
\vspace{0.25ex}
$
\Omega_{\phantom1}^{{}-(1)}=\Omega_{1}^{{}-(1)}\sqcup\Omega_{2}^{{}-(1)}$
касательная к траектории системы (16.2) 
образует  тупой угол   с положительным направлением оси $O_{\phantom1}^{(1)}\xi.$
\vspace{0.5ex}

{\it Контактные точки системы} (16.2) 
\vspace{0.35ex}
{\it координатных осей фазовой плоскости} $O^{(1)}_{\phantom1} \xi\theta.$
Нулевая изоклина $\xi=\dfrac{1}{\theta^2}$ 
\vspace{0.25ex}
не пересекает ось $O_{\phantom1}^{(1)}\xi.$
На оси $O_{\phantom1}^{(1)}\xi$ нет контактных точек системы (16.2).
\vspace{0.3ex}

Ортогональная изоклина $\theta=0$ пересекает ось  $O_{\phantom1}^{(1)}\theta$ в одной 
точке $O_{\phantom1}^{(1)}(0,0).$ 
Начало координат $O_{\phantom1}^{(1)}(0,0)$  является единственной контактной точкой оси
$O_{\phantom1}^{(1)}\theta$  системы (16.2).
В достаточно малой окрестности точки $O_{\phantom1}^{(1)}(0,0)$ контактная 
\vspace{0.35ex}
$O_{\phantom1}^{(1)}$\!-траектория системы (16.2) лежит в полуплоскости 
$\xi\leq 0,$ так как $\Theta(0,0)\,\partial_{\theta}^{}\,\Xi(0,0)={}-1<0.$
\vspace{0.5ex}

{\it Проективный тип системы} (16.2). 
Прямая $x=0$ не состоит из траекторий системы (16.1). 
Система  (16.2) --- проективно особая (свойство 10.4). 

Бесконечно удаленная прямая проективной фазовой плоскости $\R\P(\xi,\theta)$ 
 не состоит из траекторий системы (16.2) (свойство 9.2).
\vspace{0.5ex}

{\it Экваториально контактные точки  системы} (16.2). 
На оси $Oy$ у системы (16.1) нет контактных точек.
\vspace{0.25ex}
На  <<концах>> оси $Oy$ лежит состояние равновесия (двухсепаратрисное седло) системы (16.1).
\vspace{0.5ex}
У системы (16.2) нет экваториально контактных точек. 

{\it Предельные циклы системы} (16.2) {\it в проективной фазовой плоскости $\R\P(\xi,\theta).$} 
Cистема (16.1) не имеет линейных,  разомкнутых и лежащих в конечной части проективной 
фазовой плоскости $\R\P(x,y)$ предельных циклов. 

У системы (16.2) линейных, 
\vspace{0.15ex}
разомкнутых и лежащих в конечной части проективной 
фазовой плоскости $\R\P(\xi,\theta)$ предельных циклов нет.
\vspace{0.5ex}

{\it Общий автономный интеграл системы} (16.3). 
Функция
\\[1.5ex]
\mbox{}\hfill
$
F^{(2)}_{\phantom1} \colon (\eta,\zeta)\to\
\dfrac{\zeta+\eta^2}{1+\zeta^2} + \arctg\dfrac{1}{\zeta} 
\quad
\forall 
(\eta,\zeta)\in G^{(2)}_{\phantom1},
\quad   
G^{(2)}_{\phantom1}=\{(\eta,\zeta)\colon \zeta\ne0\},
\hfill
$
\\[1.25ex]
есть общий автономный интеграл системы (16.3) 
\vspace{0.5ex}
на  любой области из множества $G^{(2)}_{\phantom1}.$

{\it Состояния равновесия системы} (16.3) 
\vspace{0.25ex}
{\it в проективной фазовой плоскости} $\R\P(\eta,\zeta).$ 
У системы (16.3) в проективной фазовой плоскости $\R\P(\eta,\zeta)$  
\vspace{0.35ex}
два состояния равновесия:
двухсепаратрисное седло $O^{(2)}_{\phantom1} (0,0),$ 
\vspace{0.35ex}
сепаратрисы которого примыкают в направлении оси 
$O^{(2)}_{\phantom1}\zeta,$ и
неустойчивый фокус, лежащий на <<концах>>  оси $O^{(2)}_{\phantom1}\eta.$ 
\vspace{0.5ex}

{\it Движение радиуса-вектора образующей точки вдоль траекторий системы} (16.3).
Функция 
\\[0.75ex]
\mbox{}\hfill
$
W^{(2)}_{\phantom1}\colon (\eta,\zeta)\to\ 
\eta \,Z(\eta,\zeta)-\zeta \,H(\eta,\zeta)={}-\eta^2 +\zeta^3
\quad  
\forall (\eta,\zeta)\in \R^2
\hfill
$
\\[1.75ex]
является отрицательной, когда $\zeta<\sqrt[\scriptstyle3]{\eta^2}\,,$ а положительной, когда 
$\zeta>\sqrt[\scriptstyle3]{\eta^2}\,.$

При движении вдоль 
\vspace{0.75ex}
частей траекторий системы (16.3), расположенных в области: 
\linebreak
а) $\{(\eta,\zeta)\colon \zeta<\sqrt[\scriptstyle3]{\eta^2}\,\};$  
\vspace{0.5ex}
б) $\{(\eta,\zeta)\colon \zeta>\sqrt[\scriptstyle3]{\eta^2}\,\},$ 
угол между радиусом-век\-то\-ром образующей точки и положительным направлением оси 
$O^{(2)}_{\phantom1}\eta\colon\!\!$  а) убывает; б) возрастает. 
\vspace{0.75ex}

Через каждую точку, лежащую на кривой $\zeta=\sqrt[\scriptstyle3]{\eta^2}\,,$ 
\vspace{0.35ex}
траектория системы (16.3) проходит в направлении 
радиуса-вектора этой точки.
\vspace{0.5ex}

{\it Симметричность фазового поля направлений системы}  (16.3).
Фазовое поле направлений системы (16.3) симметрично 
\vspace{0.25ex}
относительно  координатной оси $O^{(2)}_{\phantom1} \zeta$ (свойство~13.6).
\vspace{0.25ex}
Для каждой траектории системы (16.3) существует траектория, которая симметрична ей относительно 
координатной оси $O^{(2)}_{\phantom1} \zeta.$
\vspace{0.15ex}
Каждая траектория системы (16.3), пересекающая ось $O^{(2)}_{\phantom1} \zeta$ 
\vspace{0.15ex}
в точке, отличной от начала координат фазовой плоскости $O^{(2)}_{\phantom1} \eta\zeta,$
симметрична  относительно этой координатной оси.
\vspace{0.75ex}

{\it Нулевые и ортогональные изоклины системы}  (16.3).
\vspace{0.35ex}
Из уравнения $Z(\eta,\zeta)=0$ находим, что прямая $\eta=0,$ из  которой удалена точка 
$O_{\phantom1}^{(2)}(0,0),$ является нулевой изоклиной системы (16.3).
\vspace{0.35ex}
Траектории системы (16.3) пересекают координатную ось 
$O_{\phantom1}^{(2)}\zeta$ под прямым углом в каждой точке $(0,\zeta)$ при $\zeta\ne0.$
\vspace{0.5ex}

Из уравнения $H(\eta,\zeta)=0$ находим ортогональные изоклины $\zeta=0$ 
при $\eta\ne0$ и $\zeta={}-\eta^2$ при $\eta\ne0$ системы (16.3).
\vspace{0.35ex}
Траектории системы (16.3) пересекают ось  $O_{\phantom1}^{(2)}\eta$
под прямым углом в каждой точке $(\eta,0)$ при $\eta\ne0.$
\vspace{0.35ex}
Касательная к траектории системы (16.3) в каждой точке параболы 
\vspace{0.25ex}
$\zeta={}-\eta^2,$ отличной от начала  координат фазовой плоскости 
$O_{\phantom1}^{(2)}\eta\zeta,$ параллельна оси $O_{\phantom1}^{(2)}\zeta.$
\vspace{0.5ex}

{\it Области знакоопределенности фазового поля направлений  системы}  (16.3).
\vspace{0.35ex}
Из неравенства $H(\eta,\zeta)\;\!Z(\eta,\zeta)>0$
\vspace{0.35ex}
находим, что областями положительности фазового поля направлений системы (16.3) являются области 
\\[1.75ex]
\mbox{}\hfill
$
\Omega_{1}^{{}+(2)}=\{(\eta,\zeta)\colon \eta>0\ \ \& \ \zeta>0\},
\qquad  
\Omega_{2}^{{}+(2)}=\{(\eta,\zeta)\colon \eta<0\ \ \& \ {}-\eta^2<\zeta<0\},
\hfill
$
\\[2.5ex]
\mbox{}\hfill
$  
\Omega_{3}^{{}+(2)}=\{(\eta,\zeta)\colon \eta>0\ \ \& \ \zeta<{}-\eta^2\}
\hfill
$
\\[2.5ex]
\indent
В каждой точке множества 
\vspace{0.35ex}
$\Omega_{\phantom1}^{{}+(2)}=\Omega_{1}^{{}+(2)}\sqcup\Omega_{2}^{{}+(2)}\sqcup\Omega_{3}^{{}+(2)}$
касательная к траектории системы (16.3)
образует острый угол с положительным направлением оси $O_{\phantom1}^{(2)}\eta.$
\vspace{0.75ex}

Разрешая неравенство $H(\eta,\zeta)\;\!Z(\eta,\zeta)<0,$ 
\vspace{0.35ex}
получаем, что областями
отрицательности фазового поля направлений  системы (16.3) являются области
\\[1.75ex]
\mbox{}\hfill
$
\Omega_{1}^{{}-(2)}=\{(\eta,\zeta)\colon \eta<0\ \ \& \ \zeta>0\},
\qquad 
\Omega_{2}^{{}-(2)}=\{(\eta,\zeta)\colon \eta<0\ \ \& \ \zeta<{}-\eta^2\},
\hfill
$
\\[2.5ex]
\mbox{}\hfill
$  
\Omega_{3}^{{}-(2)}=\{(\eta,\zeta)\colon \eta>0\ \ \& \ {}-\eta^2<\zeta<0\}.
\hfill
$
\\[2.5ex]
\indent
В каждой точке  множества 
\vspace{0.35ex}
$
\Omega_{\phantom1}^{{}-(2)}=\Omega_{1}^{{}-(2)}\sqcup\Omega_{2}^{{}-(2)}\sqcup\Omega_{3}^{{}-(2)}$
касательная к траектории системы (16.3)
образует  тупой угол   с положительным направлением оси $O_{\phantom1}^{(2)}\eta.$
\vspace{0.5ex}

{\it Контактные точки системы} (16.3) 
\vspace{0.35ex}
{\it координатных осей фазовой плоскости} $O^{(2)}_{\phantom1} \eta\zeta.$
Нулевая изоклина 
\vspace{0.5ex}
$\{(\eta,\zeta)\colon \eta=0\ \,\&\,\ \zeta\ne 0\}$ не пересекает ось $O_{\phantom1}^{(2)}\eta.$
На оси $O_{\phantom1}^{(2)}\eta$ нет контактных точек системы (16.3).
\vspace{0.5ex}

Ортогональные изоклины 
\vspace{0.75ex}
$\{(\eta,\zeta)\colon \zeta=0\ \,\&\,\ \eta\ne 0\}$
и
$\{(\eta,\zeta)\colon \zeta={}-\eta^2\ \,\&\,\ \eta\ne 0\}$
не пересекают ось  $O_{\phantom1}^{(2)}\zeta.$
На оси $O_{\phantom1}^{(2)}\zeta$ нет контактных точек системы (16.3).
\vspace{0.75ex}

{\it Проективный тип системы} (16.3). 
\vspace{0.35ex}
Прямая $y=0$ не состоит из траекторий системы (16.1). 
Система  (16.3) --- проективно особая (свойство 10.4). 

\newpage

Бесконечно удаленная прямая проективной фазовой плоскости $\R\P(\eta,\zeta)$ 
\vspace{0.15ex}
 не состоит из траекторий системы (16.3) (свойство 9.2).
\vspace{0.5ex}

{\it Экваториально контактные точки  системы} (16.3). 
\vspace{0.15ex}
На оси $Ox$ у системы (16.1) нет контактных точек.
\vspace{0.15ex}
Система (16.3) не имеет экваториально контактных точек на  <<концах>> 
прямых $\zeta=a\;\!\eta$ при любом вещественном коэффициенте  $a.$
\vspace{0.35ex}

Точка $O_{\phantom1}^{(1)}$ --- контактная точка прямой $\xi=0$ 
\vspace{0.15ex}
системы (16.2), причем в достаточно малой окрестности 
\vspace{0.25ex}
точки $O_{\phantom1}^{(1)}$ контактная $O_{\phantom1}^{(1)}\!$-траектория лежит в полуплоскости $\xi\leq 0.$ 
На  <<концах>> оси $O_{\phantom1}^{(2)}\zeta$ лежит 
\vspace{0.25ex}
экваториально контактная точка  $O_{\phantom1}^{(1)}$ системы (16.3),  
а экваториально контактная 
\vspace{0.25ex}
$O_{\phantom1}^{(1)}\!$-траектория системы (16.3) в достаточно малой окрестности 
бесконечно удаленной прямой проективной фазовой плоскости 
\vspace{0.15ex}
$\R\P(\eta,\zeta)$ лежит в полуплоскости $\zeta<0.$
\vspace{0.5ex}

{\it Предельные циклы системы} (16.3) 
\vspace{0.15ex}
{\it в проективной фазовой плоскости $\R\P(\eta,\zeta).$} 
Cистема (16.1) не имеет линейных,  разомкнутых и лежащих в конечной части проективной 
фазовой плоскости $\R\P(x,y)$ предельных циклов. 
\vspace{0.15ex}

У системы (16.3) 
\vspace{0.25ex}
линейных, разомкнутых и лежащих в конечной части проективной 
фазовой плоскости $\R\P(\eta,\zeta)$ предельных циклов нет.
\vspace{0.5ex}

{\it Пересечение траекторий системы} (16.1) 
\vspace{0.15ex}
{\it бесконечно удаленной прямой проективной фазовой плоскости 
$\R\P(x,y).$}
\vspace{0.25ex}
На <<концах>> оси $Oy$ лежит двухсепаратрисное седло, сепаратрисы которого ортогональны оси $Oy.$ 
\vspace{0.15ex}
Прямая $\theta=0$ является ортогональной изоклиной системы (16.2). 
\vspace{0.15ex}
В каждой  точке граничной окружности проективного круга $\P\K(x,y),$ не лежащей на оси $Oy,$ 
\vspace{0.15ex}
траектория системы (16.1) ортогонально пересекает граничную окружность.
\vspace{0.5ex}

{\it Пересечение траекторий системы} (16.2) 
\vspace{0.15ex}
{\it бесконечно удаленной прямой проективной фазовой плоскости 
$\R\P(\xi,\theta).$}
\vspace{0.5ex}
На <<концах>>  координатных осей $O_{\phantom1}^{(1)}\xi$ и $O_{\phantom1}^{(1)}\theta$
лежат состояния равновесия.  
\vspace{0.35ex}
Прямая $x=0,$ из которой удалена точка $O(0,0),$ является нулевой изоклиной системы (16.1). 
\vspace{0.15ex}
В каждой  точке граничной окружности проективного круга $\P\K(\xi,\theta),$ не лежащей на 
координатных осях  $O_{\phantom1}^{(1)}\xi$ и $O_{\phantom1}^{(1)}\theta,$
\vspace{0.35ex}
траектория системы (16.2) ортогонально пересекает граничную окружность.
\vspace{0.5ex}

{\it Проективный атлас траекторий системы} (16.1) построен на рис. 16.1.
\\[3.5ex]
\mbox{}\hfill
{\unitlength=1mm
\begin{picture}(42,42)
\put(0,0){\includegraphics[width=42mm,height=42mm]{r16-1a.eps}}
\put(18,41){\makebox(0,0)[cc]{ $y$}}
\put(40.2,18.2){\makebox(0,0)[cc]{ $x$}}
\end{picture}}
\qquad
{\unitlength=1mm
\begin{picture}(42,42)
\put(0,0){\includegraphics[width=42mm,height=42mm]{r16-1b.eps}}
\put(18,41){\makebox(0,0)[cc]{ $\theta$}}
\put(40.2,17.8){\makebox(0,0)[cc]{ $\xi$}}
\put(21,-7.75){\makebox(0,0)[cc]{Рис. 16.1}}
\end{picture}}
\qquad
{\unitlength=1mm
\begin{picture}(42,42)
\put(0,0){\includegraphics[width=42mm,height=42mm]{r16-1c.eps}}
\put(18,41){\makebox(0,0)[cc]{ $\zeta$}}
\put(40.2,18){\makebox(0,0)[cc]{ $\eta$}}
\end{picture}}
\hfill\mbox{}
\\[9.5ex]
\indent
{\bf 16.2.
Траектории дифференциальной системы} [1, с. 84 -- 85]
\\[2ex]
\mbox{}\hfill        % (16.4)
$
\dfrac{dx}{dt}=1-x^2-y^2\equiv X(x,y),
\qquad 
\dfrac{dy}{dt}=xy-1\equiv Y(x,y).
$
\hfill (16.4)
\\[2.5ex]
\indent
В конечной части проективной фазовой плоскости $\R\P(x,y)$ 
\vspace{0.35ex}
у системы (16.4) нет состояний равновесия, так как
окружность $x^2+y^2=1$ и гипербола $xy=1$ не пересекаются.

Фазовое поле направлений системы (16.4) симметрично 
относительно начала координат
фазовой плоскости $Oxy.$ 
Для каждой траектории системы (16.4) существует траектория, которая симметрична ей относительно 
начала координатной фазовой плоскости $Oxy.$ Траектория системы (16.4), проходящая через начало 
координат фазовой плоскости $Oxy,$  симметрична относительно 
начала координатной этой фазовой плоскости.

Нулевой изоклиной системы (16.4) является гипербола $xy=1.$ Касательная к траектории системы (16.4) в
каждой точке гиперболы $xy=1$ параллельна оси $Ox.$

Ортогональной изоклиной системы (16.4) является окружность $x^2+y^2=1.$ Касательная к траектории системы (16.4) 
в каждой точке окружности $x^2+y^2=1$ параллельна оси $Oy.$

В каждой точке области положительности 
\\[1.75ex]
\mbox{}\hfill
$
\Omega_{\phantom1}^{+}=\{(x,y)\colon  x^2+y^2>1\ \&\ xy<1\}
\hfill
$ 
\\[1.75ex]
фазового поля направлений касательная к траектории системы (16.4) 
образует острый угол с положительным направлением оси $Ox.$

Области
\\[1.5ex]
\mbox{}\hfill
$
\Omega_{1}^{-}=\{(x,y)\colon x^2+y^2<1\},
\ \ 
\Omega_{2}^{-}=\{(x,y)\colon x<0\ \&\ xy>1\},
\ \
\Omega_{3}^{-}=\{(x,y)\colon x>0\ \&\ xy>1\}
\hfill
$
\\[1.75ex]
являются областями отрицательности фазового поля направлений  системы (16.4).
\vspace{0.15ex}
Касательная к траектории системы (16.4) в каждой точке множества 
$\Omega_{\phantom1}^{-}=\Omega_1^-\sqcup \Omega_2^-\sqcup \Omega_3^-$
образует  тупой угол   с положительным направлением  оси $Ox.$
\vspace{0.15ex}

Нулевая изоклина $xy=1$ 
\vspace{0.15ex}
не имеет общих точек с осью $Ox.$ На оси $Ox$ нет контактных точек  системы  (16.4).
\vspace{0.15ex}

Ортогональная изоклина $x^2+y^2=1$ 
\vspace{0.25ex}
пересекает ось $Oy$ в точках $(0,{}\pm1).$ 
На оси $Oy$ у системы (16.4) две контактные точки $A_1^{}(0,{}-1)$ и $A_2^{}(0,1).$ 
\vspace{0.35ex}

Произведение $Y(0,y)\;\!\partial_y X(0,y)= 2y$ при $y={}-1$
\vspace{0.35ex}
отрицательно, а при $y=1$ положительно.
\vspace{0.25ex}
Контактная $A_1^{}$-траектория в достаточно малой окрестности точки $A_1^{}$ лежит в полуплоскости $x\leq 0,$ 
\vspace{0.25ex}
контактная $A_2^{}$-траектория в достаточно малой окрестности точки $A_2^{}$ лежит в полуплоскости $x\geq 0.$ 
\vspace{0.25ex}

Функция 
\\[1ex]
\mbox{}\hfill
$
W_2^{}\colon (x,y)\to\  
y(2x^2+y^2)
\quad
\forall (x,y)\in \R^2
\hfill
$
\\[2ex]
тождественно не равна  нулю на  $\R^2.$
Система  (16.4) --- проективно неособая. Бесконечно удаленная прямая проективной фазовой плоскости 
$\R\P(x,y)$ состоит из  траекторий системы (16.4). 

Первой и второй проективно приведенными системами системы (16.4)
соответственно являются  системы  
\\[2ex]
\mbox{}\hfill        % (16.5)
$
\dfrac{d\xi}{d\tau}= 2\xi-\theta^2+\xi\bigl(\xi^2-\theta^2\bigr)\equiv \Xi(\xi,\theta),
\quad 
\dfrac{d\theta}{d\tau}=\theta\bigl(1+\xi^2-\theta^2\bigr)\equiv \Theta(\xi,\theta),
\quad
\theta\, d\tau=dt,
$
\hfill (16.5)
\\[2ex]
и
\\[1.5ex]
\mbox{}\hfill        % (16.6)
$
\dfrac{d\eta}{d\nu}= {}-\eta\bigl(\zeta-\eta^2\bigr)\equiv H(\eta,\zeta),
\quad 
\dfrac{d\zeta}{d\nu}={}-1+\eta^2-2\zeta^2 + \eta^2\zeta\equiv Z(\eta,\zeta),
\quad
\eta\,d\nu=dt.
$
\hfill (16.6)
\\[3ex]
\indent
Система уравнений 
\\[1ex]
\mbox{}\hfill
$\Xi(\xi,0)=0,
\quad 
\Theta(\xi,0)=0
\hfill
$ 
\\[1.5ex]
имеет одно решение $\xi=0.$
\vspace{0.15ex}
На координатной оси $O^{(1)}_{\phantom1} \xi$ в конечной части проективной фазовой плоскости $\R\P(\xi,\theta)$ 
\vspace{0.25ex}
у системы (16.5) одно состояние равновесия $O^{(1)}_{\phantom1}(0,0),$
которое является  неустойчивым узлом 
\vspace{0.35ex}
c характеристическим уравнением $\lambda^2-3\lambda+2=0.$

Начало координат 
\vspace{0.35ex}
фазовой плоскости $O^{(2)}_{\phantom1} \eta\zeta$
не является состоянием равновесия системы (16.6),
так как  $Z(0,0)={}-1\ne 0.$

В проективной фазовой плоскости $\R\P(x,y)$ 
\vspace{0.15ex}
у системы (16.4) одно  состояние равновесия $O_{\phantom1}^{(1)},$ лежащее на <<концах>> прямой $y=0$
\vspace{0.25ex}
и  являющееся   неустойчивым узлом.

В конечной части  проективной фазовой плоскости $\R\P(x,y)$ 
система (16.4) не имеет состояний равновесия и, следовательно, предельных циклов.
\vspace{0.15ex}

Система  (16.4) --- проективно неособая и имеет бесконечно удаленное состояние равновесия. 
Линейных и разомкнутых предельных циклов у системы (16.4) нет.
\vspace{0.25ex}

В проективной фазовой плоскости $\R\P(\xi,\theta)$
\vspace{0.15ex}
у  системы (16.5)  одно  состояние равновесия $O_{\phantom1}^{(1)}(0,0),$ 
которое является  неустойчивым узлом.
\vspace{0.25ex}

Линейных, разомкнутых 
\vspace{0.25ex}
и лежащих в конечной части проективной фазовой плоскости $\R\P(\xi,\theta)$
циклов (в том числе и предельных) у системы (16.5) нет.

Фазовое поле направлений системы (16.5) симметрично 
относительно  оси $O^{(1)}_{\phantom1} \xi.$ 
\vspace{0.15ex}
Для каждой траектории системы (16.5) 
\vspace{0.15ex}
существует траектория, которая симметрична ей относительно 
оси $O^{(1)}_{\phantom1} \xi.$
\vspace{0.25ex}

Уравнение $\Theta(\xi,\theta)=0$ 
\vspace{0.35ex}
распадается на два уравнения $\theta=0$ и $\theta^2-\xi^2=1.$
Нулевыми изоклинами системы (16.5) являются  
\vspace{0.25ex}
гипербола $\theta^2-\xi^2=1$ и прямая $\theta=0,$
из которой удалена точка $O_{\phantom1}^{(1)}(0,0).$
\vspace{0.25ex}
Прямая $\theta=0$  состоит из траекторий системы (16.5).
Касательная к траектории системы (16.5) 
\vspace{0.15ex}
в каждой точке  гиперболы $\theta^2-\xi^2=1$ параллельна оси 
$O_{\phantom1}^{(1)}\xi.$ 

Уравнение $\Xi(0,\theta)=0$ 
\vspace{0.25ex}
имеет один корень $\theta=0,$ но $\Theta(0,0)=0.$ 
На оси $O^{(1)}_{\phantom1} \theta$
 нет контактных точек системы (16.5).
\vspace{0.15ex}

Прямая $x=0$ не состоит из траекторий системы (16.4). 
\vspace{0.15ex}
Система  (16.5) --- проективно особая. 
\vspace{0.15ex}
Бесконечно удаленная прямая проективной фазовой плоскости $\R\P(\xi,\theta)$ 
 не состоит из траекторий системы (16.5).
\vspace{0.15ex}

Контактным точкам $A_1^{}(0,{}-1)$ и $A_2^{}(0,1)$ оси $Oy$ системы (16.4) 
\vspace{0.25ex}
соответствуют экваториально контактные точки системы (16.5), лежащие 
\vspace{0.25ex}
соответственно на <<концах>>  прямых $\theta={}-\xi$ и  $\theta=\xi.$ 
\vspace{0.25ex}
Экваториально контактные  $A_1^{}\!$- и $A_2^{}\!$-траектории  системы (16.5) 
\vspace{0.25ex}
в достаточно малой окрестности бесконечно удаленной прямой 
проективной фазовой плоскости $\R\P(\xi,\theta)$ 
лежат в полуплоскости $\xi>0.$
\vspace{0.35ex}

В проективной фазовой плоскости $\R\P(\eta,\zeta)$
\vspace{0.15ex}
система (16.6) имеет  одно  состояние равновесия $O_{\phantom1}^{(1)},$ 
лежащее  на <<концах>>  прямой $\eta=0,$ которое является  неустойчивым узлом.
%\vspace{0.15ex}

Линейных, 
\vspace{0.25ex}
разомкнутых и лежащих в конечной части проективной фазовой плоскости
$\R\P(\eta,\zeta)$
циклов (в том числе и предельных циклов) у системы (16.6) нет.
\vspace{0.15ex}

Фазовое поле направлений системы (16.6) симметрично 
\vspace{0.15ex}
относительно оси $O^{(2)}_{\phantom1} \zeta.$ 
Для каждой траектории системы (16.6) 
\vspace{0.15ex}
существует траектория, которая симметрична ей относительно 
оси $O^{(2)}_{\phantom1} \zeta.$
\vspace{0.15ex}
Прямая-траектория $\eta=0$ системы (16.6) симметрична относительно 
начала координат фазовой плоскости $O^{(2)}_{\phantom1} \eta\zeta.$
\vspace{0.5ex}

Уравнение 
\vspace{0.25ex}
$H(\eta,\zeta)=0$ распадается на два уравнения $\eta=0$ и $\zeta=\eta^2.$
Ортогональными изоклинами системы (16.6) являются прямая 
\vspace{0.25ex}
$\eta=0$ и  парабола $\zeta=\eta^2.$
Ортогональная изоклина $\eta=0$ системы (16.6) является ее траекторией. 
\vspace{0.15ex}
Касательная   к траектории системы (16.6) в каждой точке параболы 
$\zeta=\eta^2$ параллельна оси $O^{(2)}_{\phantom1} \zeta.$
\vspace{0.5ex}

Уравнение 
$Z(\eta,0)=0$
имеет два корня $\eta_1^{}={}-1$ и $\eta_2^{}=1,$ причем 
\\[1.5ex]
\mbox{}\hfill
$
H({}\pm1,0)=\pm1\ne 0.
\hfill
$
\\[1.5ex]
На оси $O^{(2)}_{\phantom1} \eta$  две контактные точки 
\vspace{0.5ex}
$A_1^{(2)}({}-1,0)$ и $A_2^{(2)}(1,0)$ системы (16.6).
Произведение 
\vspace{0.25ex}
$H(\eta,0)\;\!\partial_\eta^{} Z(\eta,0)=2\eta^4$  положительно как при $\eta={}-1,$ так и при $\eta=1.$
Траектории системы (16.6), проходящие через контактные точки $A_1^{(2)}$ и $A_2^{(2)},$ 
\vspace{0.5ex}
в достаточно малой окрестности каждой из этих точек лежат в полуплоскости $\zeta\geq 0.$ 

Прямая $y=0$ не состоит из траекторий системы (16.4). 
\vspace{0.15ex}
Система  (16.6) --- проективно особая. 
Бесконечно удаленная прямая проективной фазовой плоскости $\R\P(\eta,\zeta)$ 
 не состоит из траекторий системы (16.6).
\vspace{0.15ex}

На оси $Ox$ контактных точек дифференциальная система (16.4)  не имеет, а на 
\vspace{0.15ex}
<<концах>> прямой $\!y\!=\!0$ лежит состояние равновесия дифференциальной системы (16.4).
У  системы (16.6) нет экваториально контактных точек.
\vspace{0.15ex}

Направление движения вдоль траекторий систем (16.4) --- (16.6) определяется неустойчивостью узла  
$O^{(1)}_{\phantom1}.$
\vspace{0.15ex}

Проективный атлас траекторий системы (16.4) построен на рис. 16.2.
\\[3.75ex]
\mbox{}\hfill
{\unitlength=1mm
\begin{picture}(42,42)
\put(0,0){\includegraphics[width=42mm,height=42mm]{r16-2a.eps}}
\put(18,41){\makebox(0,0)[cc]{ $y$}}
\put(40.2,18.5){\makebox(0,0)[cc]{ $x$}}
%\put(21,-3){\makebox(0,0)[cc]{ $Oxy$}}
%\put(22.5,-6){\makebox(0,0)[cc]{Рис. 1}}
\end{picture}}
\qquad
{\unitlength=1mm
\begin{picture}(42,42)
\put(0,0){\includegraphics[width=42mm,height=42mm]{r16-2b.eps}}
\put(18,41){\makebox(0,0)[cc]{ $\theta$}}
\put(40.2,17.8){\makebox(0,0)[cc]{ $\xi$}}
%\put(21,-3){\makebox(0,0)[cc]{ $O^{{}^{(1)}}uz$}}
\put(21,-7){\makebox(0,0)[cc]{Рис. 16.2}}
\end{picture}}
\qquad
{\unitlength=1mm
\begin{picture}(42,42)
\put(0,0){\includegraphics[width=42mm,height=42mm]{r16-2c.eps}}
\put(18,41){\makebox(0,0)[cc]{ $\zeta$}}
\put(40.2,18){\makebox(0,0)[cc]{ $\eta$}}
%\put(21,-3){\makebox(0,0)[cc]{ $O^{{}^{(2)}}zv$}}
%\put(22.5,-6){\makebox(0,0)[cc]{Рис. 1}}
\end{picture}}
\hfill\mbox{}
\\[9ex]
\indent
{\bf 16.3.
Траектории дифференциальной системы} [1, c. 85 --- 87; 2, c. 209 --- 212] 
\\[2ex]
\mbox{}\hfill        % (16.7)
$
\dfrac{dx}{dt}={}-1+x^2+y^2\equiv X(x,y),
\quad \ 
\dfrac{dy}{dt}={}-5+5xy\equiv Y(x,y).
$
\hfill (16.7)
\\[2.5ex]
\indent
В конечной части проективной фазовой плоскости $\R\P(x,y)$ 
у системы (16.7) нет состояний равновесия.

Фазовое поле направлений системы (16.7) симметрично 
относительно начала координат фазовой плоскости $Oxy.$  
Для каждой траектории системы (16.7) существует траектория, которая симметрична ей относительно 
начала координат фазовой плоскости $Oxy.$ Траектория системы (16.7), 
\vspace{0.15ex}
проходящая через начало 
координат фазовой плоскости $Oxy,$  симметрична относительно 
начала координат этой фазовой плоскости.
\vspace{0.15ex}

Нулевой изоклиной системы (16.7) является гипербола $xy=1.$ Касательная к траектории системы (16.7) в
каждой точке гиперболы $xy=1$ параллельна оси $Ox.$
\vspace{0.15ex}

Ортогональной изоклиной системы (16.7) 
\vspace{0.15ex}
является окружность $x^2+y^2=1.$ Касательная к траектории системы (16.7) 
в каждой точке окружности $x^2+y^2=1$ параллельна оси $Oy.$

Области 
\\[1.5ex]
\mbox{}\hfill
$
\Omega_{1}^{+}=\{(x,y)\colon  x^2+y^2<1\},
\ \
\Omega_{2}^{+}=\{(x,y)\colon  x<0\ \&\ xy>1\},
\ \
\Omega_{3}^{+}=\{(x,y)\colon  x>0\ \&\ xy>1\}
\hfill
$
\\[1.75ex]
являются областями положительности 
\vspace{0.25ex}
фазового поля направлений системы (16.7).
Касательная к траектории системы (16.7) в каждой точке множества 
\vspace{0.15ex}
$
\Omega_{}^{+}=
\Omega_{1}^{+}\sqcup
\Omega_{2}^{+}\sqcup
\Omega_{3}^{+}
$
образует острый угол с положительным направлением оси $Ox.$
\vspace{0.15ex}

Касательная к траектории системы (16.7) в каждой точке области отрицательности
фазового поля направлений
\\[1.5ex]
\mbox{}\hfill
$
\Omega_{}^{-}=\{(x,y)\colon x^2+y^2>1\ \&\ xy<1\}
\hfill
$
\\[1.5ex]
образует  тупой угол   с положительным направлением оси $Ox.$
\vspace{0.15ex}

Нулевая изоклина $xy=1$ 
\vspace{0.15ex}
не имеет общих точек с осью $Ox.$ На оси $Ox$ нет контактных точек  системы  (16.7).

Ортогональная изоклина $x^2+y^2=1$ 
\vspace{0.35ex}
пересекает ось $Oy$ в точках $(0,{}\pm1).$ 
На оси $Oy$ две контактные точки $A_1^{}(0,{}-1)$ и $A_2^{}(0,1)$ системы (16.7). 
\vspace{0.35ex}
Произведение $Y(0,y)\;\!\partial_y X(0,y)= {}-10y$ при $y={}-1$
\vspace{0.35ex}
положительно, а при $y=1$ отрицательно.
Контактная $A_1^{}$\!-траектория системы (16.7) 
\vspace{0.25ex}
в достаточно малой окрестности точки $A_1^{}$ лежит в полуплоскости $x\geq 0,$ а 
контактная $A_2^{}$\!-траектория 
\vspace{0.25ex}
в достаточно малой окрестности точки $A_2^{}$ лежит в полуплоскости $x\leq 0.$ 
\vspace{0.25ex}

Функция 
\\[1ex]
\mbox{}\hfill
$
W_2^{}\colon (x,y)\to\  
y(4x^2-y^2)
\quad
\forall (x,y)\in \R^2
\hfill
$
\\[2ex]
не является тождественным  нулем на $\R^2.$ Система  (16.7) --- проективно неособая. 
Бесконечно удаленная прямая проективной фазовой плоскости 
\vspace{0.15ex}
$\R\P(x,y)$ состоит из  траекторий системы (16.7). 

Первой и второй проективно приведенными системами системы (16.7)
соответственно являются  системы  
\\[2.5ex]
\mbox{}\hfill        % (16.8)
$
\dfrac{d\xi}{d\tau}= 4\xi-5\theta^2-\xi\bigl(\xi^2-\theta^2\bigr)\equiv \Xi(\xi,\theta),
\quad 
\dfrac{d\theta}{d\tau}={}-\theta\bigl(1+\xi^2-\theta^2\bigr)\equiv \Theta(\xi,\theta),
\quad
\theta\, d\tau=dt,
$
\hfill (16.8)
\\[2.25ex]
и
\\[1.75ex]
\mbox{}\hfill        % (16.9)
$
\dfrac{d\eta}{d\nu}= {}-5\eta\bigl(\zeta-\eta^2\bigr)\equiv H(\eta,\zeta),
\quad 
\dfrac{d\zeta}{d\nu}=1-\eta^2-4\zeta^2 + 5\eta^2\zeta\equiv Z(\eta,\zeta),
\quad
\eta\,d\nu=dt.
$
\hfill (16.9)
\\[3ex]
\indent
На координатной оси $O^{(1)}_{\phantom1} \xi$
в конечной части проективной фазовой плоскости $\R\P(\xi,\theta)$ 
у системы (16.8) три состояния равновесия:
\vspace{0.35ex}
седло  $O^{(1)}_{\phantom1}(0,0);$
устойчивые узлы $A^{(1)}_{\phantom1}({}-2,0)$ и $B^{(1)}_{\phantom1}(2,0).$
\vspace{0.75ex}
Характеристическим уравнением этих состояний равновесия является
$\lambda^2+(4\xi^2-3)\lambda+(\xi^2+1)(3\xi^2-4)=0.$
\vspace{0.35ex}

Начало координат  фазовой плоскости 
\vspace{0.25ex}
$O^{(2)}_{\phantom1} \eta\zeta$
не является состоянием равновесия системы (16.9), 
так как  $Z(0,0)=1\ne 0.$
\vspace{0.35ex}

В проективной фазовой плоскости $\R\P(x,y)$  у системы (16.7) 
\vspace{0.15ex}
три состояния равновесия: седло
$O_{\phantom1}^{(1)},$ лежащее на <<концах>> прямой $y=0,$ 
\vspace{0.15ex}
и два устойчивых узла, одно из которых $A^{(1)}_{\phantom1}$  
лежит на <<концах>> прямой $y={}-2x,$  
\vspace{0.15ex}
а другое $B^{(1)}_{\phantom1}$  лежит на <<концах>> прямой $y=2x.$  
\vspace{0.15ex}

В конечной части  проективной фазовой плоскости $\R\P(x,y)$ 
\vspace{0.15ex}
система (16.7) не имеет состояний равновесия и, следовательно, предельных циклов.

Система  (16.7) --- проективно неособая и имеет бесконечно удаленные состояния равновесия. 
Линейных и разомкнутых предельных циклов у системы (16.7) нет.
\vspace{0.25ex}

В проективной фазовой плоскости $\R\P(\xi,\theta)$ у
\vspace{0.35ex}
системы (16.8) три состояния равновесия: $O_{\phantom1}^{(1)}(0,0)$ --- седло,
$A^{(1)}_{\phantom1}({}-2, 0)$  и $B^{(1)}_{\phantom1}(2, 0)$  --- устойчивые узлы.
\vspace{0.35ex}

Линейных, 
\vspace{0.15ex}
разомкнутых и лежащих в конечной части проективной фазовой плоскости $\R\P(\xi,\theta)$
предельных циклов  у системы (16.8) нет.

Фазовое поле направлений системы (16.8) симметрично 
относительно  оси $O^{(1)}_{\phantom1} \xi.$ 
Для каждой траектории системы (16.8) существует траектория, которая симметрична ей относительно 
оси $O^{(1)}_{\phantom1} \xi.$
\vspace{0.25ex}

Из уравнения $\Theta(\xi,\theta)=0$ 
\vspace{0.15ex}
находим, что система (16.8) имеет одну нулевую изоклину 
$\theta=0,$ которая состоит из ее траекторий.
\vspace{0.35ex}

Уравнение $\Xi(0,\theta)=0$ 
\vspace{0.15ex}
имеет один корень   $\theta=0,$
но $\Theta(0,0)=0.$ 
На координатной оси $O^{(1)}_{\phantom1} \theta$
 нет контактных точек системы (16.8).
\vspace{0.25ex}

Прямая $x=0$ не состоит из траекторий системы (16.7). 
\vspace{0.15ex}
Система  (16.8) --- проективно особая. 
Бесконечно удаленная прямая проективной фазовой плоскости $\R\P(\xi,\theta)$ 
 не состоит из траекторий системы (16.8).

Контактным точкам $A_1^{}(0,{}-1)$ и $A_2^{}(0,1)$ оси $Oy$ системы (16.7) 
\vspace{0.35ex}
соответствуют экваториально контактные точки системы (16.8), лежащие 
\vspace{0.15ex}
соответственно на <<концах>>  прямых $\theta={}-\xi$ и  $\theta=\xi.$ 
\vspace{0.15ex}
Экваториально контактные траектории системы (16.8) 
в достаточно малой окрестности бесконечно удаленной прямой 
\vspace{0.35ex}
проективной фазовой плоскости 
$\R\P(\xi,\theta)$ лежат в полуплоскости $\xi<0.$ 
\vspace{0.75ex}

В проективной фазовой плоскости $\R\P(\eta,\zeta)$
у системы (16.9) три  состояния равновесия: 
седло $O_{\phantom1}^{(1)}$ на <<концах>>  прямой $\eta=0,$
неустойчивый узел $A_{\phantom1}^{(2)}\Bigl(0, {}-\dfrac{1}{2}\Bigr)$ и
устойчивый узел $B_{\phantom1}^{(2)}\Bigl(0, \dfrac{1}{2}\Bigr).$
\vspace{0.5ex}

Линейных, 
\vspace{0.35ex}
разомкнутых и лежащих в конечной части проективной фазовой плоскости
$\R\P(\eta,\zeta)$ предельных циклов у системы (16.9) нет.
\vspace{0.25ex}

Фазовое поле направлений системы (16.9) симметрично 
\vspace{0.25ex}
относительно оси $O^{(2)}_{\phantom1} \zeta.$ 
Для каждой траектории системы (16.9) 
\vspace{0.15ex}
существует траектория, которая симметрична ей относительно 
оси $O^{(2)}_{\phantom1} \zeta.$
\vspace{0.5ex}

Уравнение 
\vspace{0.25ex}
$H(\eta,\zeta)=0$ распадается на два уравнения $\eta=0$ и $\zeta=\eta^2.$
Ортогональными изоклинами системы (16.9) являются прямая $\eta=0$ и  парабола $\zeta=\eta^2.$
Ортогональная изоклина $\eta=0$ системы (16.9) состоит из ее траекторий. 
Касательная   к траектории системы (16.9) в каждой точке параболы 
\vspace{0.75ex}
$\zeta=\eta^2$ параллельна оси $O^{(2)}_{\phantom1} \zeta.$

Уравнение 
$Z(\eta,0)=0$
имеет два корня $\eta_1^{}={}-1$ и $\eta_2^{}=1,$ причем 
\\[2ex]
\mbox{}\hfill
$
H({}\pm1,0)=\pm1\ne 0.
\hfill
$
\\[1.75ex]
На координатной оси $O^{(2)}_{\phantom1} \eta$  две контактные точки 
\vspace{0.75ex}
$A_1^{(2)}({}-1,0)$ и $A_2^{(2)}(1,0)$ системы (16.9).
Произведение 
$H(\eta,0)\;\!\partial_\eta^{} Z(\eta,0)={}-10\eta^4$
отрицательно как при $\eta={}-1,$ так и при $\eta=1.$
Траектории системы (16.9), проходящие через контактные точки 
\vspace{0.5ex}
$A_1^{(2)}$ и $A_2^{(2)},$ в достаточно малой окрестности каждой из этих точек лежат 
\vspace{0.5ex}
в полуплоскости $\zeta\leq 0.$ 

Прямая $y=0$ не состоит из траекторий системы (16.7). 
\vspace{0.25ex}
Система  (16.9) --- проективно особая. 
\vspace{0.15ex}
Бесконечно удаленная прямая проективной фазовой плоскости $\R\P(\eta,\zeta)$ 
не состоит из траекторий системы (16.9).
\vspace{0.25ex}

На оси $Ox$ контактных точек  система (16.7)  не имеет, а на 
\vspace{0.15ex}
<<концах>> прямой $y=0$ лежит состояние равновесия системы (16.7).
У  системы (16.9) нет экваториально контактных точек.
\\[4.25ex]
\mbox{}\hfill
{\unitlength=1mm
\begin{picture}(42,42)
\put(0,0){\includegraphics[width=42mm,height=42mm]{r16-3a.eps}}
\put(18,41){\makebox(0,0)[cc]{ $y$}}
\put(40.2,18.5){\makebox(0,0)[cc]{ $x$}}
%\put(21,-3){\makebox(0,0)[cc]{ $Oxy$}}
%\put(22.5,-6){\makebox(0,0)[cc]{Рис. 1}}
\end{picture}}
\qquad
{\unitlength=1mm
\begin{picture}(42,42)
\put(0,0){\includegraphics[width=42mm,height=42mm]{r16-3b.eps}}
\put(18,41){\makebox(0,0)[cc]{ $\theta$}}
\put(40.2,17.8){\makebox(0,0)[cc]{ $\xi$}}
%\put(21,-3){\makebox(0,0)[cc]{ $O^{{}^{(1)}}uz$}}
\put(21,-7){\makebox(0,0)[cc]{Рис. 16.3}}
\end{picture}}
\qquad
{\unitlength=1mm
\begin{picture}(42,42)
\put(0,0){\includegraphics[width=42mm,height=42mm]{r16-3c.eps}}
\put(18,41){\makebox(0,0)[cc]{ $\zeta$}}
\put(40.2,18){\makebox(0,0)[cc]{ $\eta$}}
%\put(21,-3){\makebox(0,0)[cc]{ $O^{{}^{(2)}}zv$}}
%\put(22.5,-7){\makebox(0,0)[cc]{Рис. 1}}
\end{picture}}
\hfill\mbox{}
\\[9ex]
\indent
Проективный атлас траекторий системы (16.7) построен на рис. 16.3.
Для этого целесообразно сначала установить 
\vspace{0.15ex}
поведение траекторий проективно особой системы (16.8) на 
проективном круге $\P\K(\xi, \theta),$ 
\vspace{0.35ex}
а затем использовать отображения проективных кругов 
$\P\K(x, y),\ \P\K(\xi, \theta),\ \P\K(\eta, \zeta).$
\vspace{0.35ex}
Поведение траекторий на проективном круге $\P\K(\xi, \theta)$ системы (16.8) определяется однозначно из того, что 
две сепаратрисы седла 
\vspace{0.25ex}
$O^{(1)}_{\phantom1}$ лежат на оси $O^{(1)}_{\phantom1}\xi,$ 
\vspace{0.15ex}
а две другие сепаратрисы расположены симметрично оси $O^{(1)}_{\phantom1}\xi$ и лежат в полуплоскости $\xi>0.$ 
\vspace{0.35ex}
Последнее следует из того, что на оси $O^{(1)}_{\phantom1}\theta$ нет контактных точек системы (16.8), 
\vspace{0.35ex}
а в полуплоскости $\xi<0$ лежат экваториально контактные траектории системы (16.8). 
\vspace{0.25ex}

Направление движения вдоль траекторий систем (16.7) --- (16.9)  
\vspace{0.25ex}
определяется устойчивостью узлов  $A^{(1)}_{\phantom1}$ и $B^{(1)}_{\phantom1}.$
\\[2ex]
\indent
{\bf 16.4. Траектории дифференциальной системы} [1, с. 88]
\\[2ex]
\mbox{}\hfill        % (16.10)
$
\begin{array}{l}
\dfrac{dx}{dt}=x(x^2+y^2-1)-y(x^2+y^2+1)
\equiv X(x,y),
\\[3.5ex] 
\dfrac{dy}{dt}=y(x^2+y^2-1)+x(x^2+y^2+1)
\equiv Y(x,y).
\end{array}
$
\hfill (16.10)
\\[3ex]
\indent
На любой области из множества 
\vspace{0.25ex}
$D=\{(t,x,y)\colon x\ne 0\}$ интегральный базис системы (16.10)
образуют первые интегралы
\\[2ex]
\mbox{}\hfill
$
F_1^{}\colon (t,x,y)\to\ 
(x^2+y^2)\exp \Bigl(4t-2\arctg \dfrac{y}{x}\Bigr)
\quad
\forall (t,x,y) \in D
\hfill
$
\\[1.5ex]
и 
\\[1.5ex]
\mbox{}\hfill
$
F_2^{}\colon (t,x,y)\to\ 
(x^2+y^2-1)\exp \Bigl(2t-2\arctg \dfrac{y}{x}\Bigr)
\quad
\forall (t,x,y) \in D.
\hfill
$
\\[2ex]
\indent
Общим автономным интегралом системы (16.10) 
\vspace{0.35ex}
на любой области из множества 
$G=\{(x,y)\colon x^2+y^2\ne 1\ \, \& \ \, x\ne 0\}$ является функция 
\\[2.25ex]
\mbox{}\hfill
$
F\colon (x,y)\to\ 
\dfrac{x^2+y^2}{(x^2+y^2-1)^2}\,\exp \Bigl(2\arctg \dfrac{y}{x}\Bigr)
\quad
\forall (x,y) \in G.
\hfill
$
\\[2ex]
\indent
Система уравнений 
\\[1ex]
\mbox{}\hfill
$
X(x,y)=0,
\quad 
Y(x,y)=0
\hfill
$ 
\\[1.75ex]
имеет одно решение $x=0,\ y=0.$ 
\vspace{0.25ex}
В конечной части проективной фазовой плоскости $\R\P(x,y)$ 
у системы (16.10) одно состояние равновесия $O(0,0),$ 
\vspace{0.35ex}
которое является устойчивым фокусом
с характеристическим уравнением $\lambda^2+2\lambda+2=0.$
\vspace{0.25ex}

Функция 
\\[1.5ex]
\mbox{}\hfill
$
W\colon (x,y)\to\ 
x Y(x,y)-yX(x,y)=(x^2+y^2)(x^2+y^2+1)>0 
\quad 
\forall (x,y)\in\R^2\backslash\{(0,0)\}.
\hfill
$
\\[1.75ex]
\indent
При движении вдоль траекторий системы (16.10) угол между 
радиусом-век\-то\-ром образующей точки и положительным направлением оси $Ox$ возрастает.
\vspace{0.15ex}

Фазовое поле направлений дифференциальной системы (16.10) симметрично 
относительно начала координат фазовой плоскости $Oxy.$  
Для каждой траектории дифференциальной системы (16.10) 
\vspace{0.15ex}
существует траектория, которая симметрична ей относительно 
начала координат фазовой плоскости $Oxy.$ 
\vspace{0.5ex}

Уравнение $Y(x,0)=0$ имеет один корень $x=0,$ но $X(0,0)=0.$
\vspace{0.5ex}

Уравнение $X(0,y)=0$ имеет один корень $y=0,$ но $Y(0,0)=0.$ 
\vspace{0.5ex}

На осях $Ox$ и $Oy$ у системы (16.10) нет контактных точек.
\vspace{0.25ex}

Функция 
\\[1ex]
\mbox{}\hfill
$
W_3^{}\colon (x,y)\to\ 
(x^2+y^2)^2
\quad
\forall (x,y)\in \R^2
\hfill
$ 
\\[1.75ex]
не является тождественным  нулем на $\R^2.$ Система  (16.10) --- проективно неособая.

\newpage

Первым преобразованием Пуанкаре $x=\dfrac{1}{\theta}\,,\ y=\dfrac{\xi}{\theta}$
\vspace{0.5ex}
проективно неособую систему (16.10)
приводим к первой проективно приведенной  системе 
\\[2ex]
\mbox{}\hfill        % (16.11)
$
\begin{array}{c}
\dfrac{d\xi}{d\tau}= (1+\xi^2)(1+\xi^2+\theta^2)\equiv \Xi(\xi,\theta),
\\[4ex] 
\dfrac{d\theta}{d\tau}={}-\theta(1+\xi^2-\theta^2-\xi(1+\xi^2+\theta^2)) \equiv \Theta(\xi,\theta),
\quad \
\theta^2\, d\tau=dt.
\end{array}
$
\hfill (16.11)
\\[2.25ex]
\indent
Вторым преобразованием Пуанкаре $x=\dfrac{\zeta}{\eta}\,,\ y=\dfrac{1}{\eta}$
\vspace{0.5ex}
проективно неособую систему (16.10) приводим к второй проективно приведенной  системе 
\\[2ex]
\mbox{}\hfill        % (16.12)
$
\begin{array}{c}
\dfrac{d\eta}{d\nu}= {}-\eta(1-\eta^2+\zeta^2+\zeta(1+\eta^2+\zeta^2))\equiv H(\eta,\zeta),
\\[4ex]  
\dfrac{d\zeta}{d\nu}={}-(1+\zeta^2)(1+\eta^2+\zeta^2) \equiv Z(\eta,\zeta),
\quad \
\eta^2\,d\nu=dt.
\end{array}
$
\hfill (16.12)
\\[2.5ex]
\indent
В конечной части проективной фазовой плоскости $\R\P(\xi,\theta)$ 
\vspace{0.35ex}
у системы (16.11) нет  состояний равновесия, так как 
$\Xi(\xi,\theta)>0 \;\; \forall  (\xi,\theta)\in \R^2.$
\vspace{0.5ex}

В конечной части проективной фазовой плоскости $\R\P(\eta,\zeta)$ 
\vspace{0.35ex}
у системы (16.12) нет  состояний равновесия, так как 
$Z(\eta,\zeta)<0 \;\;\forall (\eta,\zeta)\in \R^2.$
\vspace{0.5ex}

В проективной фазовой плоскости $\R\P(x,y)$ 
\vspace{0.25ex}
у системы (16.10) одно состояние равновесия
$O(0,0)$ --- устойчивый фокус.
\vspace{0.35ex}

В [1, c. 88] доказано, что окружность $x^2+y^2=1$
\vspace{0.35ex}
является предельным циклом системы (16.10).
\vspace{0.35ex}
Из уравнения $F(x,y)=C$ семейства траекторий системы (16.10)
находим, что окружность 
\vspace{0.35ex}
$x^2+y^2=1$ --- единственный предельный цикл системы (16.10)
в  конечной части  проективной фазовой плоскости $\R\P(x,y).$
\vspace{0.25ex}

Система  (16.10) --- проективно неособая 
\vspace{0.15ex}
и не имеет состояний равновесия на бесконечно удаленной прямой 
проективной фазовой плоскости $\R\P(x,y).$  
\vspace{0.15ex}
Разомкнутых предельных циклов у системы (16.10) нет.
Единственным линейным предельным циклом системы (16.10) является предельный $\ell_\infty^{}$-цикл.
\vspace{0.75ex}

В проективной фазовой плоскости $\R\P(\xi,\theta)$ у
\vspace{0.35ex}
системы (16.11) одно состояние равновесия 
$O,$ лежащее на <<концах>>  прямой $\xi=0$
\vspace{0.75ex}
и являющееся устойчивым фокусом.

В проективной фазовой плоскости $\R\P(\xi,\theta)$ 
\vspace{0.35ex}
система (16.11) имеет два предельных цикла:
линейный предельный цикл $\theta=0$
\vspace{0.5ex}
и разомкнутый предельный цикл $\theta^2-\xi^2=1.$

Фазовое поле направлений системы (16.11) симметрично 
\vspace{0.15ex}
относительно  координатной оси $O^{(1)}_{\phantom1}\xi.$ 
\vspace{0.15ex}
Для каждой траектории системы (16.11) существует траектория, которая симметрична ей относительно 
координатной оси $O^{(1)}_{\phantom1} \xi.$
\vspace{0.15ex}
Через начало координат фазовой плоскости $O^{(1)}_{\phantom1} \xi\theta$
\vspace{0.15ex}
проходит прямая-траектория $\theta=0,$ совпадающая с 
координатной осью $O^{(1)}_{\phantom1} \xi.$
\vspace{0.35ex}

Прямая $\theta=0$  
\vspace{0.35ex}
является траекторией системы (16.11). Уравнение $\Xi(0,\theta)=0$ не имеет решений. 
На координатной оси $O^{(1)}_{\phantom1} \theta$  нет контактных точек системы (16.11).
\vspace{0.5ex}

Прямая $x=0$ 
\vspace{0.25ex}
не состоит из траекторий системы (16.10). Система  (16.11) --- проективно особая. Бесконечно удаленная прямая \vspace{0.15ex}
проективной фазовой плоскости $\R\P(\xi,\theta)$  не состоит из траекторий системы (16.11).
\vspace{0.15ex}

На оси $Oy$ контактных точек система (16.10) не имеет. 
\vspace{0.15ex}
Прямая $\theta=0$ является траекторией системы (16.11). 
У проективно особой системы (16.11) 
эк\-ва\-то\-ри\-аль\-но контактных точек  нет.

\newpage

В проективной фазовой плоскости $\R\P(\eta,\zeta)$
\vspace{0.35ex}
у системы (16.12) одно состояние равновесия 
$O,$ лежащее на <<концах>>  прямой $\eta=0$
\vspace{0.5ex}
и являющееся устойчивым фокусом.

В проективной фазовой плоскости $\R\P(\eta,\zeta)$ 
\vspace{0.35ex}
система (16.12) имеет два предельных цикла:
линейный предельный цикл $\eta=0$
\vspace{0.5ex}
и разомкнутый предельный цикл $\eta^2-\zeta^2=1.$

Фазовое поле направлений системы (16.12) симметрично 
\vspace{0.15ex}
относительно  координатной оси $O^{(2)}_{\phantom1} \zeta.$ 
\vspace{0.15ex}
Для каждой траектории системы (16.12) существует траектория, которая симметрична ей относительно 
координатной оси $O^{(2)}_{\phantom1} \zeta.$
\vspace{0.25ex}
Через начало координат фазовой плоскости $O^{(2)}_{\phantom1} \eta\zeta$
\vspace{0.25ex}
проходит прямая-траектория $\eta=0,$ совпадающая с 
координатной осью $O^{(2)}_{\phantom1} \zeta.$
\vspace{0.15ex}

Уравнение $Z(\eta,0)=0$ не имеет решений. 
\vspace{0.35ex}
На координатной оси $O^{(2)}_{\phantom1} \eta$  нет контактных точек системы (16.12).
Прямая $\eta=0$  является траекторией системы (16.12).
\vspace{0.35ex}

Прямая $y=0$ не состоит из траекторий системы (16.10). 
\vspace{0.15ex}
Система  (16.12) --- проективно особая. 
Бесконечно удаленная прямая 
\vspace{0.15ex}
проективной фазовой плоскости $\R\P(\eta,\zeta)$  
не состоит из траекторий системы (16.12).
\vspace{0.15ex}

На оси $Ox$ контактных точек  система (16.10) не имеет.
\vspace{0.25ex}
Прямая $\eta=0$ является траекторией системы (16.12). 
У проективно особой системы (16.12) эк\-ва\-то\-ри\-аль\-но контактных точек нет.
\vspace{0.15ex}

Проективный атлас траекторий системы (16.10) построен на рис. 16.4.
\\[3.75ex]
\mbox{}\hfill
{\unitlength=1mm
\begin{picture}(42,42)
\put(0,0){\includegraphics[width=42mm,height=42mm]{r16-4a.eps}}
\put(18,41){\makebox(0,0)[cc]{ $y$}}
\put(40.2,18.2){\makebox(0,0)[cc]{ $x$}}
\end{picture}}
\qquad
{\unitlength=1mm
\begin{picture}(42,42)
\put(0,0){\includegraphics[width=42mm,height=42mm]{r16-4b.eps}}
\put(18,41){\makebox(0,0)[cc]{ $\theta$}}
\put(40.2,17.8){\makebox(0,0)[cc]{ $\xi$}}
\put(21,-7){\makebox(0,0)[cc]{Рис. 16.4}}
\end{picture}}
\qquad
{\unitlength=1mm
\begin{picture}(42,42)
\put(0,0){\includegraphics[width=42mm,height=42mm]{r16-4c.eps}}
\put(18,41){\makebox(0,0)[cc]{ $\zeta$}}
\put(40.2,18){\makebox(0,0)[cc]{ $\eta$}}
\end{picture}}
\hfill\mbox{}
\\[9ex]
\indent
{\bf 16.5. Траектории дифференциальной системы} [1, с. 88 -- 90]
\\[2ex]
\mbox{}\hfill                                       %(16.13)
$
\begin{array}{l}
\dfrac{dx}{dt}=x(x^2+y^2-1)(x^2+y^2-9)-y(x^2+y^2-2x-8)
\equiv X(x,y),
\\[3.5ex]
\dfrac{dy}{dt}=y(x^2+y^2-1)(x^2+y^2-9)+x(x^2+y^2-2x-8)
\equiv Y(x,y).
\end{array}
$
\hfill (16.13)
\\[2.5ex]
\indent
Из системы уравнений 
\\[1.25ex]
\mbox{}\hfill
$
X(x,y)=0,
\quad 
Y(x,y)=0
\hfill
$ 
\\[1.75ex]
находим, что в конечной части проективной фазовой плоскости $\R\P(x,y)$ 
\vspace{0.5ex}
у системы (16.13) три  состояния равновесия $O(0,0),\, 
A\biggl(\dfrac12\,,\dfrac{\sqrt{35}}{2}\,\biggr)$ и
$B\biggl(\dfrac12\,,{}-\dfrac{\sqrt{35}}{2}\,\biggr).$
\vspace{0.75ex}

Cостояния равновесия  $A$ и $B$ лежат на окружности $x^2+y^2=9,$
\vspace{0.35ex}
которая состоит из траекторий системы (16.13).
\vspace{0.35ex}

Cостояние равновесия  $O(0,0)$ имеет 
характеристическое уравнение 
\\[1.25ex]
\mbox{}\hfill
$
\lambda^2-18\lambda+145=0
\hfill
$
\\[1.25ex]
и является неустойчивым фокусом.

Cостояние равновесия  $A\biggl(\dfrac12\,,\dfrac{\sqrt{35}}{2}\,\biggr)$ с 
характеристическим уравнением 
\\[1.5ex]
\mbox{}\hfill
$
\lambda^2-\bigl( 144+\sqrt{35}\,\bigr)\lambda+25+144\;\!\sqrt{35}=0
\hfill
$
\\[1.25ex]
является неустойчивым узлом.
\vspace{0.35ex}

Cостояние равновесия  $B\biggl(\dfrac12\,,{}-\dfrac{\sqrt{35}}{2}\,\biggr)$
\vspace{0.35ex}
является седлом, так как у его 
характеристического уравнения свободный член 
$\Delta=25-144\;\!\sqrt{35}<0.$
\vspace{0.35ex}

Функция 
\\[1.5ex]
\mbox{}\hfill
$
W\colon (x,y)\to\ 
x Y(x,y)-yX(x,y)=(x^2+y^2)(x^2+y^2-2x-8)
\quad
\forall (x,y)\in \R^2
\hfill
$
\\[1.75ex]
является положительной при $(x-1)^2+y^2>9,$
\vspace{0.5ex}
а отрицательной, когда $(x-1)^2+y^2<9$ и $x^2+y^2\ne 0.$
\vspace{0.5ex}
При движении вдоль частей траекторий системы (16.13),
расположенных в множестве: 
\vspace{0.5ex}
a)  $\!\!\bigl\{ (x,y)\colon (x-1)^2+y^2<9 \ \& \ x^2+y^2\ne 0\bigr\};$
б)  $\!\!\bigl\{ (x,y)\colon (x-1)^2+y^2>9\bigr\},$
угол между радиусом-век\-то\-ром 
\vspace{0.35ex}
образующей точки и положительным направлением оси $Ox\colon$ 
a) убывает; б) возрастает.
\vspace{0.35ex}
Через каждую точку, лежащую на окружности $\!(x-1)^2\!+y^2\!=\!9,$
траектория системы (16.13) проходит в направлении радиуса-вектора этой точки.
\vspace{0.35ex}

Уравнение $Y(x,0)=0$ имеет три корня $x=0,\ x={}-2,\ x=4$ таких, что 
\\[1.5ex]
\mbox{}\hfill
$
X(0,0)=0,\ \ 
X({}-2,0)=30\ne 0,\ \ 
X(4,0)=420\ne0.
\hfill
$
\\[1.5ex]
\indent
У системы (16.13) на оси $Ox$ две контактные точки
$C_1^{}({}-2,0)$ и $C_2^{}(4,0).$  
Поскольку 
\\[1.75ex]
\mbox{}\hfill
$
X({}-2,0)\;\!\partial_x^{} Y({}-2,0)=360>0,
\quad 
X(4,0)\;\! \partial_x^{} Y(4,0)=420\cdot 24>0,
\hfill
$
\\[1.5ex]
то контактные $C_1^{}$\!- и $C_2^{}$\!-траектории системы (16.13)
\vspace{0.35ex}
в достаточно малых окрестностях точек 
$C_1^{}$ и $C_2^{}$ лежат в полуплоскости  $y\geq 0.$  
\vspace{0.5ex}

Уравнение $X(0,y)=0$ имеет три корня $y=0$ и $y={}\pm 2\;\!\sqrt{2}$ таких, что 
\\[1.75ex]
\mbox{}\hfill
$
Y(0,0)=0,
\quad
Y(0,{}-2\;\!\sqrt{2}\,)=14\;\!\sqrt{2}\,\ne 0,
\quad 
Y(0,2\;\!\sqrt{2}\,)={}-14\;\!\sqrt{2}\,\ne 0.
\hfill
$
\\[1.75ex]
\indent
У системы (16.13) на оси $Oy$ две контактные точки
$D_1^{}(0,{}-2\;\!\sqrt{2}\,)$ и $D_2^{}(0,2\;\!\sqrt{2}\,).$  
\vspace{0.15ex}

Произведение 
\\[1.15ex]
\mbox{}\hfill
$
Y(0,{}-2\;\!\sqrt{2}\,)\;\! \partial_y^{} X(0,{}-2\;\!\sqrt{2}\,)=14\;\!\sqrt{2}\,\cdot ({}-16)<0.
\hfill
$
\\[1.75ex]
\indent
Контактная $D_1^{}$\!-траектория системы (16.13)
\vspace{0.25ex}
в достаточно малой окрестности точки $D_1^{}$
лежит в полуплоскости  $x\leq 0.$
\vspace{0.15ex}

Произведение 
\\[1.15ex]
\mbox{}\hfill
$
Y(0,2\;\!\sqrt{2}\,)\;\! \partial_y^{} X(0,2\;\!\sqrt{2}\,)={}-14\;\!\sqrt{2}\,\cdot ({}-16)>0.
\hfill
$
\\[1.75ex]
\indent
Контактная $D_2^{}$\!-траектория системы (16.13)
\vspace{0.25ex}
в достаточно малой окрестности точки $D_2^{}$
лежит в полуплоскости  $x\geq 0.$
\vspace{0.15ex}

Функция
\\[1ex]
\mbox{}\hfill
$
W_5^{}\colon (x,y)\to\
x\;\!Y_5^{}(x,y)-y\;\!X_5^{}(x,y)=0
\quad
\forall (x,y)\in \R^2.
\hfill
$
\\[1.5ex]
\indent
Система  (16.13) --- проективно особая. 
\vspace{0.15ex}
Бесконечно удаленная прямая проективной фазовой плоскости $\R\P(x,y)$
не состоит из траекторий системы (16.13).
\vspace{0.15ex}

Первой и второй проективно приведенными системами системы (16.13)
соответственно являются  системы  
\\[2.15ex]
\mbox{}\hfill        % (16.14)
$
\dfrac{d\xi}{d\tau}= \theta\bigl(1+\xi^2\bigr)\bigl(1-2\theta +\xi^2-8\theta^2\bigr)\equiv \Xi(\xi,\theta),
\hfill
$
\\
\mbox{}\hfill (16.14)      
\\[0.25ex]
\mbox{}
$
\dfrac{d\theta}{d\tau}={}-1-2\xi^2+10\theta^2+\xi \theta^2-\xi^4+10\xi^2\;\!\theta^2-
2\xi \theta^3-9\theta^4+\xi\;\!\theta^2\bigl(\xi^2-8\theta^2\bigr)
 \equiv \Theta(\xi,\theta),
\quad
\theta^3\, d\tau=dt,
\hfill
$
\\[1ex]
и
\\[2ex]
\mbox{}\       % (16.15)
$
\dfrac{d\eta}{d\nu}= {}-1+10\eta^2-2\zeta^2-\eta^2\;\! \zeta-9\eta^4+10\eta^2\;\!\zeta^2
-\zeta^4+\eta^2\;\!\zeta\bigl( 8\eta^2+2\eta\;\!\zeta-\zeta^2\bigr)
\equiv H(\eta,\zeta),
\hfill
$
\\[0.75ex]
\mbox{}\hfill (16.15)      
\\[0.25ex]
\mbox{}\hfill
$
\dfrac{d\zeta}{d\nu}=\eta\bigl(1+\zeta^2\bigr)\bigl({}-1+8\eta^2+2\eta\zeta-\zeta^2\bigr) \equiv Z(\eta,\zeta),
\quad
\eta^3\,d\nu=dt.
\hfill 
$
\\[2.75ex]
\indent
На координатной оси $O^{(1)}_{\phantom1} \xi$
\vspace{0.25ex}
в конечной части проективной фазовой плоскости $\R\P(\xi,\theta)$ 
у системы (16.14) нет  состояний равновесия,
так как 
\\[1.5ex]
\mbox{}\hfill
$
\Theta(\xi,0)={}-1-2\xi^2-\xi^4\ne 0
\quad 
\forall  \xi \in \R.
\hfill
$
\\[1.75ex]
\indent
Начало координат  фазовой плоскости $O^{(2)}_{\phantom1} \eta\zeta$
\vspace{0.35ex}
не является  состоянием равновесия системы (16.15), так как  $H(0,0)={}-1\ne 0.$
\vspace{0.5ex}
 
В проективной фазовой плоскости $\R\P(x,y)$ 
\vspace{0.75ex}
у системы (16.13) три состояния равновесия 
$O(0,0),\ A\biggl(\dfrac12\,,\dfrac{\sqrt{35}}{2}\,\biggr)$ и 
$B\biggl(\dfrac12\,,{}-\dfrac{\sqrt{35}}{2}\,\biggr).$
\vspace{0.75ex}

У системы (16.13) нет  экваториально контактных точек, так как  
\\[1.5ex]
\mbox{}\hfill
$
X_5^{}(1,\xi)=(1+\xi^2)^2\ne 0 
\quad
\forall \xi \in \R
$
\ \ и \ \ 
$
Y_5^{}(\zeta,1)= (1+\zeta^2)^2\ne 0 
\quad
\forall \zeta \in \R.
\hfill
$
\\[1.75ex]
\indent
В [1, c. 88 -- 90] доказано, 
\vspace{0.25ex}
что окружность $x^2+ y^2=1$ --- единственный предельный цикл системы (16.13) в 
конечной части проективной фазовой плоскости $\R\P(x,y).$
\vspace{0.25ex}

Среди тракторий проективно особой системы  (16.13) 
\vspace{0.15ex}
в проективной фазовой плоскости $\R\P(x,y)$ нет прямых. 
\vspace{0.15ex}
Система (16.13) не имеет линейных предельных циклов.

Проективно особая система  (16.13) не имеет бесконечно удаленных состояний равновесия и
экваториально контактных точек.
У системы (16.13) нет  разомкнутых предельных циклов.

В проективной фазовой плоскости $\R\P(\xi,\theta)$ у
\vspace{0.35ex}
системы (16.14) три состояния равновесия:
\vspace{0.5ex}
неустойчивый узел $A^{(1)}_{\phantom1} (\sqrt{35}\,, 2),$
седло  $B^{(1)}_{\phantom1} ({}-\sqrt{35}\,, 2)$
и неустойчивый  фокус, лежащий на <<концах>>  прямой $\xi=0.$
\vspace{0.35ex}

Уравнение $\Xi(\xi,\theta)=0$
\vspace{0.5ex}
распадается на два уравнения $\theta=0$
и $8\theta^2+2\theta-\xi^2-1=0.$
\vspace{0.35ex}
Прямая $\theta=0$
и гипербола $8\theta^2+2\theta-\xi^2-1=0$
являются ортогональными изоклинами системы (16.14).
\vspace{0.15ex}
Касательная к траектории системы (16.14) в каждой точке этих прямой и гиперболы 
параллельна координатной оси $O^{(1)}_{\phantom1} \theta.$
\vspace{0.25ex}

Прямая $\theta=0$  
\vspace{0.25ex}
является ортогональной изоклиной системы (16.14).
На координатной оси $O^{(1)}_{\phantom1} \xi$  у системы (16.14)  нет контактных точек.
\vspace{0.5ex}
Ортогональные изоклины
$\theta=0$
и $8\theta^2+2\theta-\xi^2-1=0$
пересекают прямую $\xi=0$ в трех точках 
$O^{(1)}_{\phantom1} (0,0)\,,\ 
C^{(1)}_{1} \Bigl(0,{}-\dfrac12\Bigr)$
и $C^{(1)}_{2} \Bigl(0,\dfrac14\Bigr),$
которые и будут контактными точками координатной оси  $O^{(1)}_{\phantom1} \theta.$
\vspace{0.5ex}

Произведение $\Theta(0,0)\;\! \partial_\theta^{} \Xi(0,0)={}-1<0.$
\vspace{0.25ex}
Контактная $O^{(1)}_{\phantom1}\!$-траектория системы (16.14) в достаточно малой окрестности точки 
\vspace{0.25ex}
$O^{(1)}_{\phantom1}$ лежит в полуплоскости  $\xi\leq 0.$
Учитывая расположение контактных
\vspace{0.25ex}
$C_1^{}$\!- и  $C_2^{}$\!-траекторий системы (16.13), получаем, что 
 контактная $C_1^{(1)}$\!-траектория системы (16.14) 
\vspace{0.25ex}
в достаточно малой окрестности точки $C^{(1)}_{1}$ лежит в полуплоскости  $\xi\leq 0;$
\vspace{0.35ex}
контактная $C_2^{(1)}$\!-траектория системы (16.14) 
в достаточно малой окрестности точки $C^{(1)}_{2}$ лежит в полуплоскости  $\xi\geq 0.$
\vspace{0.35ex}

Прямая $x=0$ не состоит из траекторий системы (16.13). 
\vspace{0.25ex}
Система  (16.14) --- проективно особая. 
\vspace{0.15ex}
Бесконечно удаленная прямая проективной фазовой плоскости $\R\P(\xi,\theta)$ 
 не состоит из траекторий системы (16.14).

\newpage

На <<концах>> прямой  $\xi=0$ лежит фокус.
\vspace{0.25ex}
Бесконечно удаленная точка  
на <<концах>> оси $O^{(1)}_{\phantom1} \theta$
не является экваториально контактной точкой системы (16.14). 
\vspace{0.15ex}

Уравнение 
\\[1ex]
\mbox{}\hfill
$
\Xi_5^{}(1,a)=0,
\hfill
$ 
\\[1ex]
где 
\\[1ex]
\mbox{}\hfill
$
\Xi_5^{}(\xi,\theta)=\xi^2\;\!\theta(\xi^2-8\theta^2)
\quad 
\forall (\xi,\theta)\in\R^2,
\hfill
$
\\[1.75ex]
имеет три корня $a=0$ и $a={}\pm \dfrac{\sqrt{2}}{4}\,,$
при которых 
\\[1.75ex]
\mbox{}\hfill
$
W^{(1)}_4(1,a)=a\, \Xi_4^{}(1,a)- \Theta_4^{}(1,a)\ne 0.
\hfill
$
\\[1.5ex]
\indent
У системы (16.14) три экваториально контактные точки 
$O^{(2)}_{\phantom1} ,\, D^{(1)}_{1}$ и $D^{(1)}_{2},$
лежащие на <<концах>> прямых $\theta=0,\,\theta={}-\dfrac{\sqrt{2}}{4}\ \xi$ и 
$\theta=\dfrac{\sqrt{2}}{4}\ \xi$ соответственно.
\vspace{0.25ex}

Учитывая расположение контактных траекторий оси $Oy$ системы (16.13),
получаем, что экваториально контактные траектории системы (16.14) 
в достаточно малой окрестности бесконечно удаленной прямой
проективной фазовой плоскости  $\R\P(\xi,\theta)$ 
лежат в  полуплоскости  $\xi> 0.$
\vspace{0.25ex}

У системы (16.14) нет 
\vspace{0.25ex}
линейных и лежащих в конечной части проективной фазовой плоскости 
$\R\P(\xi,\theta)$   предельных циклов.
\vspace{0.15ex}
Гипербола $\theta^2-\xi^2=1$ ---  разомкнутый предельный цикл системы (16.14).
\vspace{0.25ex}

В проективной фазовой плоскости $\R\P(\eta,\zeta)$ у
\vspace{0.15ex}
системы (16.15) три состояния равновесия:
\vspace{0.5ex}
неустойчивый фокус, лежащий на <<концах>>  прямой $\zeta=0,$
неустойчивый узел $A^{(2)}_{\phantom1}\biggl(\dfrac{2\sqrt{35}}{35}\,,\dfrac{\sqrt{35}}{35}\,\biggl)$ и 
седло 
$B^{(2)}_{\phantom1}\biggl({}-\dfrac{2\sqrt{35}}{35}\,,{}-\dfrac{\sqrt{35}}{35}\,\biggl).$  
\vspace{1ex}

Уравнение $Z(\eta,\zeta)=0$
\vspace{0.35ex}
распадается на два уравнения $\eta=0$
и $8\eta^2+2\eta\zeta-\zeta^2-1=0.$
Прямая $\eta=0$
и гипербола $8\eta^2+2\eta\zeta-\zeta^2-1=0$
\vspace{0.35ex}
являются нулевыми изоклинами системы (16.15).
\vspace{0.25ex}
Касательная к траектории системы (16.15) в каждой точке этих прямой и гиперболы 
параллельна  оси $O^{(2)}_{\phantom1} \eta.$

Прямая $\eta=0$  является нулевой изоклиной системы (16.15). 
\vspace{0.25ex}
На  оси $O^{(2)}_{\phantom1} \zeta$ у системы (16.15)  нет контактных точек.
\vspace{0.35ex}

Нулевые изоклины 
\vspace{0.5ex}
$\eta=0$ и $8\eta^2+2\eta\zeta-\zeta^2-1=0$
пересекают прямую $\zeta=0$  в трех точках 
\vspace{0.25ex}
$O^{(2)}_{\phantom1} (0,0)\,,\, D^{(2)}_{1} \Bigl({}-\dfrac{\sqrt{2}}{4}\,,0\Bigr)$
и $D^{(2)}_{2} \Bigl(\dfrac{\sqrt{2}}{4}\,,0\Bigr),$
которые и будут контактными точками оси  $O^{(2)}_{\phantom1} \eta.$
\vspace{1ex}

Произведение 
\\[1ex]
\mbox{}\hfill
$
H(0,0)\;\! \partial_\eta^{} Z(0,0)=1>0.
\hfill
$
\\[1.5ex]
\indent
Контактная 
\vspace{0.35ex}
$O^{(2)}_{\phantom1}\!$-траектория системы (16.13) в достаточно малой окрестности точки 
$O^{(2)}_{\phantom1}$ лежит в полуплоскости  $\zeta\geq 0.$
\vspace{0.35ex}

Учитывая расположение контактных
\vspace{0.35ex}
$D_1^{}$\!- и  $D_2^{}$\!-траекторий системы (16.13), получаем, что 
 контактные $D_1^{(2)}$\!- и  $D_2^{(2)}$\!-траектории 
\vspace{0.35ex}
системы (16.15) в достаточно малых окрестностях 
точек $D_1^{(2)}$ и  $D_2^{(2)}$  лежат в полуплоскости  $\zeta\geq 0.$
\vspace{0.5ex}

Прямая $y=0$ не состоит из траекторий системы (16.13). 
\vspace{0.25ex}
Система  (16.15) --- проективно особая. 
Бесконечно удаленная прямая проективной фазовой плоскости $\R\P(\eta,\zeta)$ 
 не состоит из траекторий системы (16.15).
\vspace{0.15ex}

На <<концах>> прямой  $\zeta=0$ лежит фокус.
\vspace{0.25ex}
Бесконечно удаленная точка  
на <<концах>> оси $O^{(2)}_{\phantom1} \eta$
не является экваториально контактной точкой системы (16.15). 

Уравнение 
\\[1ex]
\mbox{}\hfill
$
Z_5^{}(a,1)=0,
\hfill
$ 
\\[1ex]
где 
\\[1ex]
\mbox{}\hfill
$
Z_5^{}(\eta,\zeta)=\eta \zeta^2\;\!(8\eta^2+2\eta\zeta-\zeta^2)
\quad 
\forall (\eta,\zeta)\in\R^2,
\hfill
$
\\[2ex]
имеет три корня $a=0,\ a={}-\dfrac{1}{2}$ и $a=\dfrac{1}{4}\,,$
при которых 
\\[2ex]
\mbox{}\hfill
$
W^{(2)}_4(a,1)=H_4^{}(a,1)-aZ_4^{}(a,1)\ne 0.
\hfill
$
\\[1.75ex]
\indent
У системы (16.15) три экваториально контактные точки 
\vspace{0.35ex}
$O^{(1)}_{\phantom1} ,\, C^{(2)}_{1}$ и $C^{(2)}_{2},$
лежащие на <<концах>> прямых $\eta=0,\,\zeta={}-2\eta$ и 
$\zeta=4\eta$ соответственно.
\vspace{0.35ex}

Учитывая расположение контактных траекторий оси $Ox$ системы (16.13),
\vspace{0.15ex}
получаем, что экваториально контактная $O^{(1)}_{\phantom1}\!$-траектория системы (16.15) 
\vspace{0.15ex}
в достаточно малой окрестности бесконечно удаленной точки 
\vspace{0.15ex}
$O^{(1)}_{\phantom1}$ лежат в  полуплоскости  $\zeta< 0,$
а экваториально контактные $C^{(2)}_{1}\!$- и  $C^{(2)}_{2}\!$-траектории системы (16.15) 
\vspace{0.35ex}
в достаточно малой окрестности бесконечно удаленной прямой
\vspace{0.25ex}
проективной фазовой плоскости  $\R\P(\eta,\zeta)$ 
лежат в  полуплоскости  $\eta> 0.$
\vspace{0.35ex}

У системы (16.15) нет  линейных  и 
\vspace{0.35ex}
лежащих в конечной части проективной фазовой плоскости  $\R\P(\eta,\zeta)$ предельных циклов.
\vspace{0.25ex}
Гипербола $\eta^2-\zeta^2=1$ ---  разомкнутый предельный цикл системы (16.15).
\vspace{0.35ex}

Прямая $\theta=0$ является ортогональной изоклиной системы (16.14), а 
\vspace{0.15ex}
прямая $\eta=0$ является нулевой изоклиной системы (16.15).
\vspace{0.15ex}
Каждая траектория системы (16.13), проходящая через бесконечно удаленную прямую 
\vspace{0.25ex}
проективной фазовой плоскости 
$\R\P(x,y),$ ортогональна граничной окружности проективного круга $\R\K(x,y).$
\vspace{0.5ex}

Проективный атлас траекторий системы (16.13) построен на рис. 16.5.
\\[3.75ex]
\mbox{}\hfill
{\unitlength=1mm
\begin{picture}(42,42)
\put(0,0){\includegraphics[width=42mm,height=42mm]{r16-5a.eps}}
\put(18,41){\makebox(0,0)[cc]{ $y$}}
\put(40.2,18.2){\makebox(0,0)[cc]{ $x$}}
\end{picture}}
\qquad
{\unitlength=1mm
\begin{picture}(42,42)
\put(0,0){\includegraphics[width=42mm,height=42mm]{r16-5b.eps}}
\put(18,41){\makebox(0,0)[cc]{ $\theta$}}
\put(40.2,17.8){\makebox(0,0)[cc]{ $\xi$}}
\put(21,-7){\makebox(0,0)[cc]{Рис. 16.5}}
\end{picture}}
\qquad
{\unitlength=1mm
\begin{picture}(42,42)
\put(0,0){\includegraphics[width=42mm,height=42mm]{r16-5c.eps}}
\put(18,41){\makebox(0,0)[cc]{ $\zeta$}}
\put(40.2,18){\makebox(0,0)[cc]{ $\eta$}}
\end{picture}}
\hfill\mbox{}
\\[10ex]
\indent
{\bf 16.6. Траектории дифференциальной системы}  [1, с. 90 -- 91]
\\[2ex]
\mbox{}\hfill                                       %(16.16)
$
\!\!\!
\begin{array}{l}
\dfrac{dx}{dt}=x(2x^2+2y^2+1)\Bigl((x^2+y^2)^2+x^2-y^2+\dfrac{1}{10}\Bigr)-y(2x^2+2y^2-1)
\equiv X(x,y),
\\[3.75ex]
\dfrac{dy}{dt}=y(2x^2+2y^2-1)\Bigl((x^2+y^2)^2+x^2-y^2+\dfrac{1}{10}\Bigr)+x(2x^2+2y^2+1)
\equiv Y(x,y).
\end{array}
$
\hfill (16.16)
\\[3ex]
\indent
В конечной части проективной фазовой плоскости $\R\P(x,y)$ 
\vspace{0.5ex}
у системы (16.16) три  состояния равновесия $O(0,0),\, 
A\biggl(0,\,\dfrac{\sqrt{2}}{2}\,\biggr)$ и
$B\biggl(0,{}-\dfrac{\sqrt{2}}{2}\,\biggr).$
\vspace{0.5ex}

\newpage

Cостояние равновесия  $O(0,0)$ --- седло, так как у его 
\vspace{0.5ex}
характеристического уравнения свободный член 
$\Delta(0,0)={}-0{,}99<0.$
\vspace{0.5ex}

Cостояния равновесия  
\vspace{0.75ex}
$A\biggl(0,\,\dfrac{\sqrt{2}}{2}\,\biggr)$ и
$B\biggl(0,{}-\dfrac{\sqrt{2}}{2}\,\biggr)$
--- устойчивые фокусы с характеристическим уравнением 
$\lambda^2+\dfrac{3}{5}\,\lambda+\dfrac{409}{100}=0.$
\vspace{0.75ex}

Фазовое поле направлений системы (16.16) симметрично 
\vspace{0.15ex}
относительно начала координат фазовой плоскости $Oxy.$
\vspace{0.15ex}
Для каждой траектории системы (16.16) существует траектория, которая симметрична ей относительно 
\vspace{0.15ex}
начала координатной фазовой плоскости $Oxy.$ 
\vspace{0.35ex}

Уравнение $Y(x,0)=0$ не имеет  корней.  На оси $Ox$ у системы (16.16) нет контактных точек. 
Уравнение $X(0,y)=0$ имеет три корня $y=0$ и $y={}\pm \dfrac{\sqrt{2}}{2}\,,$ но 
$Y(0,0)=0$ и $Y\biggl(0, {}\pm \dfrac{\sqrt{2}}{2}\,\biggr)=0.$
На оси $Oy$ у системы (16.16) нет контактных точек. 
\vspace{0.75ex}

Функция 
\\[1ex]
\mbox{}\hfill
$
W_7^{}\colon (x,y)\to\ 
x\;\!Y_7^{}(x,y)-y\;\!X_7^{}(x,y)=0
\quad
\forall (x,y)\in \R^2.
\hfill
$
\\[1.75ex]
\indent
Система  (16.16) --- проективно особая. 
\vspace{0.15ex}
Бесконечно удаленная прямая проективной фазовой плоскости $\R\P(x,y)$
не состоит из траекторий системы (16.16).
\vspace{0.25ex}

Первой и второй проективно приведенными системами системы (16.16)
соответственно являются  системы  
\\[3ex]
\mbox{}\       % (16.17)
$
\dfrac{d\xi}{d\tau}= \theta\Bigl({}-2\xi +2\theta^2 -4\xi^3 -2\xi\theta^2+ 4\xi^2\theta^2+\theta^4 -2\xi^5+2\xi^3\theta^2-
\dfrac{1}{5}\, \xi\theta^4+\xi^2\theta^2\bigl(2\xi^2-\theta^2\bigr)\Bigr)\equiv \Xi(\xi,\theta),
\hfill
$
\\[0.75ex]
\mbox{}\hfill (16.17)      
\\
\mbox{}\hfill
$
\dfrac{d\theta}{d\tau}={}-2-6\xi^2-3\theta^2-6\xi^4-2\xi^2\;\!\theta^2-
\dfrac{6}{5}\, \theta^4+2\xi\theta^4\ - 
\hfill
$
\\[2.5ex]
\mbox{}\hfill
$
-\ 2\xi^6+\xi^4\theta^2+
\dfrac{4}{5}\,\xi^2\theta^4-\dfrac{1}{10}\, \theta^6+\xi\theta^4\bigl(2\xi^2-\theta^2\bigr)
 \equiv \Theta(\xi,\theta),
\qquad 
\theta^5\, d\tau=dt,
\hfill
$
\\[1.5ex]
и
\\[1.5ex]
\mbox{}\hfill       % (16.18)
$
\dfrac{d\eta}{d\nu}= 
{}-2+3\eta^2-6\zeta^2-\dfrac{6}{5}\,\eta^4+2\eta^2\;\!\zeta^2
-6\zeta^4-2\eta^4\;\!\zeta\ +
\hfill
$
\\[2.5ex]
\mbox{}\hfill
$
+\ \dfrac{1}{10}\,\eta^6+ 
\dfrac{4}{5}\,\eta^4\;\!\zeta^2-
\eta^2\;\!\zeta^4-2\zeta^6-\eta^4\;\!\zeta\bigl(\eta^2+2\zeta^2\bigr)
\equiv H(\eta,\zeta),
\hfill
$
\\[0.5ex]
\mbox{}\hfill (16.18)      
\\[0.5ex]
\mbox{}\hfill
$
\dfrac{d\zeta}{d\nu}=\eta\;\!\Bigl(2\zeta-2\eta^2-2\eta^2\zeta+4\zeta^3+\eta^4
-4\eta^2\;\!\zeta^2\ +
\hfill
$
\\[2.5ex]
\mbox{}\hfill
$
+\ \dfrac{1}{5}\,\eta^4\;\!\zeta+ 
2\eta^2\;\!\zeta^3+2\zeta^5-\eta^2\;\!\zeta^2\bigl(\eta^2+2\zeta^2\bigr)\Bigr)
\equiv Z(\eta,\zeta),
\qquad
\eta^5\,d\nu=dt.
\hfill 
$
\\[2.5ex]
\indent
На координатной оси $O^{(1)}_{\phantom1} \xi$
\vspace{0.35ex}
в конечной части проективной фазовой плоскости $\R\P(\xi,\theta)$ 
у системы (16.17) нет  состояний равновесия, так как  
\\[1.5ex]
\mbox{}\hfill
$
\Theta(\xi,0)={}-2-6\xi^2-6\xi^4-2\xi^6\ne 0  
\quad 
\forall  \xi \in \R.
\hfill
$
\\[1.75ex]
\indent
Начало координат  фазовой плоскости $O^{(2)}_{\phantom1} \eta\zeta$
\vspace{0.35ex}
не является  состоянием равновесия системы (16.18),
так как 
$
H(0,0)={}-2\ne 0.$
\vspace{0.5ex}

В проективной фазовой плоскости $\R\P(x,y)$  
\vspace{0.35ex}
у системы (16.16) три состояния равновесия 
$O(0,0),\, A\biggl(0,\,\dfrac{\sqrt{2}}{2}\,\biggr)$ и $B\biggl(0,{}-\dfrac{\sqrt{2}}{2}\,\biggr).$

У проективно особой системы (16.16) нет  экваториально контактных точек, так как 
\\[1.5ex]
\mbox{}\hfill
$
X_7^{}(1,\xi)=2(1+\xi^2)^3\ne 0 
\quad
\forall \xi \in \R
$
\ \ и \ \ 
$
Y_7^{}(\zeta,1)=2(1+\zeta^2)^3\ne 0 
\quad
\forall \zeta \in \R.
\hfill
$
\\[1.75ex]
\indent
В [1, c. 90 -- 91] доказано, что 
\vspace{0.35ex}
предельными циклами системы (16.16) в 
конечной части проективной фазовой плоскости $\R\P(x,y)$
\vspace{0.35ex}
являются две непересекающиеся замкнутые кривые, заданные уравнением 
$(x^2+y^2)^2+x^2-y^2+\dfrac{1}{10}=0.$ 
\vspace{0.75ex}

Среди тракторий проективно особой системы  (16.16) 
\vspace{0.15ex}
на проективной фазовой плоскости $\R\P(x,y)$ нет прямых. 
\vspace{0.25ex}
Система (16.16) не имеет линейных предельных циклов.

Проективно особая система  (16.16) не имеет бесконечно удаленных состояний равновесия и экваториально контактных точек. У системы (16.16) нет  разомкнутых предельных циклов.

В проективной фазовой плоскости $\R\P(\xi,\theta)$ у
\vspace{0.25ex}
системы (16.17) три состояния равновесия:
\vspace{0.5ex}
устойчивые фокусы $A^{(1)}_{\phantom1}$ и 
$B^{(1)}_{\phantom1},$ лежащие соответственно  на <<концах>>  прямых $\xi=\dfrac{\sqrt2}{2}\,\theta$ и 
$\xi={}-\dfrac{\sqrt2}{2}\,\theta,$  а также
седло, лежащее на <<концах>>  прямой $\xi=0.$
\vspace{0.75ex}

У системы (16.17) 
\vspace{0.15ex}
нет линейных и лежащих в конечной части проективной фазовой плоскости
$\R\P(\xi,\theta)$  предельных циклов.

Кривые, заданные уравнением 
\\[1.25ex]
\mbox{}\hfill
$
10(\xi^2+1)^2-10\theta^2(\xi^2-1)^2+\theta^4=0,
\hfill
$ 
\\[1.25ex]
образуют два разомкнутых предельных цикла системы (16.17).
\vspace{0.25ex}

Фазовое поле направлений системы (16.17) симметрично 
\vspace{0.15ex}
относительно   оси $O^{(1)}_{\phantom1} \xi.$
Для каждой траектории системы (16.17) существует траектория, которая симметрична ей относительно 
оси $O^{(1)}_{\phantom1} \xi.$
Каждая траектория системы (16.17), пересекающая  ось $O^{(1)}_{\phantom1} \xi,$ 
симметрична относительно этой координатной оси. 

Прямая $\theta=0$   является ортогональной изоклиной системы (16.17).
\vspace{0.25ex}
На оси $O^{(1)}_{\phantom1} \xi$  у системы (16.17)  нет контактных точек.
\vspace{0.35ex}

Уравнение $\Xi(0,\theta)=0$ имеет один корень $\theta=0,$ a $\Theta(0,0)={}-2\ne0.$
\vspace{0.35ex}

На оси $O^{(1)}_{\phantom1} \theta$ у системы (16.17) одна контактная точка $O^{(1)}_{\phantom1} (0,0).$
В достаточно малой окрестности точки $O^{(1)}_{\phantom1}(0,0)$ контактная  
\vspace{0.25ex}
$O^{(1)}_{\phantom1}\!$-траектория системы (16.17) лежит в  полуплоскости $\theta\geq 0.$
\vspace{0.35ex}

Прямая $x=0$ не состоит из траекторий системы (16.16). 
\vspace{0.25ex}
Система  (16.17) --- проективно особая. 
\vspace{0.15ex}
Бесконечно удаленная прямая проективной фазовой плоскости $\R\P(\xi,\theta)$ 
 не состоит из траекторий системы (16.17).
\vspace{0.25ex}

На  оси  $Oy$ у системы (16.16) нет контактных точек. 
У системы (16.17) нет экваториально контактных точек на <<концах>>  
\vspace{0.25ex}
прямых $\xi=a\theta$ при любом вещественном коэффициенте $a.$ 
\vspace{0.35ex}
Так как $\Xi_7^{}(1,0)=0,\,
W_6^{(1)}(1,0)={}-2\ne 0,$ то на <<концах>> оси $O^{(1)}_{\phantom1}\xi$ лежит экваториально контактная  точка 
системы (16.17), а экваториально контактная траектория 
\vspace{0.25ex}
в достаточно малой окрестности бесконечно удаленной прямой проективной фазовой плоскости $\R\P(\xi,\theta)$  лежит в полуплоскости $\xi<0.$
\vspace{0.5ex}

В проективной фазовой плоскости $\R\P(\eta,\zeta)$ у
\vspace{0.75ex}
системы (16.18) три состояния равновесия:
устойчивые фокусы $A^{(2)}_{\phantom1}(\sqrt2\,,0)$ и 
$B^{(2)}_{\phantom1}({}-\sqrt2\,,0);$ 
седло, лежащее на <<концах>>  прямой $\zeta=0.$
\vspace{0.25ex}

У системы (16.18) нет  
\vspace{0.15ex}
линейных и лежащих в конечной части проективной фазовой плоскости
$\R\P(\eta,\zeta)$  предельных циклов.

Кривые, заданные уравнением 
\\[1.15ex]
\mbox{}\hfill
$
10(\zeta^2+1)^2+10\eta^2(\zeta^2-1)^2+\eta^4=0,
\hfill
$ 
\\[1ex]
образуют два разомкнутых предельных цикла системы (16.18).

Фазовое поле направлений системы (16.18) симметрично 
\vspace{0.15ex}
относительно  оси $O^{(2)}_{\phantom1} \zeta.$
Для каждой траектории системы (16.18) 
\vspace{0.15ex}
существует траектория, которая симметрична ей относительно 
оси $O^{(2)}_{\phantom1} \zeta.$
\vspace{0.15ex}
Каждая траектория системы (16.18), пересекающая ось $O^{(2)}_{\phantom1} \zeta,$ 
симметрична относительно этой координатной оси. 
\vspace{0.25ex}

Прямая $\eta=0$  является нулевой изоклиной системы (16.18). 
\vspace{0.25ex}
На оси $O^{(2)}_{\phantom1} \zeta$ у системы (16.18)  нет контактных точек.
\vspace{0.5ex}

Уравнение $Z(\eta,0)=0$ имеет три корня $\eta=0$ и  $\eta={}\pm \sqrt2\,,$ а 
\\[1.25ex]
\mbox{}\hfill
$
H(0,0)={}-2,
\quad
H({}\pm \sqrt2\,,0)=0.
\hfill
$ 
\\[1.5ex]
\indent
На оси $O^{(2)}_{\phantom1} \eta$ у системы (16.18) одна контактная точка $O^{(2)}_{\phantom1} (0,0).$
\vspace{0.25ex}
В достаточно малой окрестности точки $O^{(2)}_{\phantom1}(0,0)$ контактная  
\vspace{0.25ex}
$O^{(2)}_{\phantom1}\!$-траектория системы (16.18) лежит в  полуплоскости $\zeta\leq 0.$
\vspace{0.35ex}

Прямая $y=0$ не состоит из траекторий системы (16.16). 
\vspace{0.15ex}
Система  (16.18) --- проективно особая. 
\vspace{0.15ex}
Бесконечно удаленная прямая проективной фазовой плоскости $\R\P(\eta,\zeta)$ 
 не состоит из траекторий системы (16.18).
\vspace{0.25ex}

На оси $Ox$ у системы (16.16) нет контактных точек.
\vspace{0.15ex}
Система (16.18) не имеет экваториально контактных точек на  <<концах>> 
\vspace{0.25ex}
прямых $\eta=a\zeta$ при любом вещественном коэффициенте  $a.$
Точка $O_{\phantom1}^{(1)}$ --- контактная точка прямой $\xi=0$ системы (16.17), 
\vspace{0.35ex}
причем в достаточно малой окрестности точки $O_{\phantom1}^{(1)}$ 
\vspace{0.35ex}
контактная $O_{\phantom1}^{(1)}\!$-траектория лежит в полуплоскости $\xi\geq 0.$ 
\vspace{0.25ex}
На  <<концах>> оси $O_{\phantom1}^{(2)}\zeta$ лежит 
экваториально контактная точка  $O_{\phantom1}^{(1)}$ системы (16.18),  
а экваториально контактная $O_{\phantom1}^{(1)}\!$-траектория системы (16.18) 
\vspace{0.25ex}
в достаточно малой окрестности 
бесконечно удаленной прямой 
\vspace{0.25ex}
проективной фазовой плоскости $\R\P(\eta,\zeta)$ лежит в полуплоскости $\zeta>0.$
\vspace{0.5ex}

Прямая $\theta=0$ является ортогональной изоклиной системы (16.17), а 
\vspace{0.15ex}
прямая $\eta=0$ является нулевой изоклиной системы (16.18).
\vspace{0.15ex}
Каждая траектория системы (16.16), проходящая через бесконечно удаленную прямую 
\vspace{0.15ex}
проективной фазовой плоскости 
$\R\P(x,y),$ ортогональна граничной окружности проективного круга $\R\K(x,y).$
\vspace{0.35ex}

Проективный атлас траекторий системы (16.16) построен на рис. 16.6.
\\[3.75ex]
\mbox{}\hfill
{\unitlength=1mm
\begin{picture}(42,42)
\put(0,0){\includegraphics[width=42mm,height=42mm]{r16-6a.eps}}
\put(18,41){\makebox(0,0)[cc]{ $y$}}
\put(40.2,18.2){\makebox(0,0)[cc]{ $x$}}
\end{picture}}
\qquad
{\unitlength=1mm
\begin{picture}(42,42)
\put(0,0){\includegraphics[width=42mm,height=42mm]{r16-6b.eps}}
\put(18,41){\makebox(0,0)[cc]{ $\theta$}}
\put(40.2,17.8){\makebox(0,0)[cc]{ $\xi$}}
\put(21,-7){\makebox(0,0)[cc]{Рис. 16.6}}
\end{picture}}
\qquad
{\unitlength=1mm
\begin{picture}(42,42)
\put(0,0){\includegraphics[width=42mm,height=42mm]{r16-6c.eps}}
\put(18,41){\makebox(0,0)[cc]{ $\zeta$}}
\put(40.2,18){\makebox(0,0)[cc]{ $\eta$}}
\end{picture}}
\hfill\mbox{}
\\[11.75ex]
\indent
{\sl Замечание.}
\vspace{0.25ex}
Основные результаты данной статьи опубликованы автором  в 
работах 
\linebreak
\text{[5; 22 -- 24]}.

\newpage

\mbox{}
\\[-2.25ex]

%\begin{thebibliography}{99}
{\Large\bf Список литературы}
\vspace{1.75ex}

1.
{\it Пуанкаре А.}
О кривых, определяемых дифференциальными уравнениями. -- М.;Л.: ГИТТЛ, 1947. -- 392 с.
\vspace{0.75ex}

2.
{\it Лефшец С.} Геометрическая теория дифференциальных уравнений. -- М.: ИЛ, 1961. -- 387 с.
\vspace{0.75ex}

3.
{\it Андронов А.А., Леонтович Е.А., Гордон И.И., Майер А.Г.} 
Качественная теория динамических систем второго порядка.  -- М.: Наука, 1966.  -- 568 с.
\vspace{0.75ex}

4.
{\it Баутин Н.Н., Леонтович Е.А.} 
Методы и приемы качественного исследования динамических систем на плоскости.  --
М.:  Наука, 1976.  -- 496 с.
\vspace{0.75ex}

5.
{\it Горбузов В.Н., Королько И.В.} 
Траектории полиномиальной дифференциальной системы на сфере Пуанкаре // Дифференц. уравнения.
-- 2002. -- Т. 38, \No\ 6. -- 
\linebreak
С. 845 -- 846.
\vspace{0.75ex}

6.
{\it Математическая} энциклопедия. Т. 4. -- М.: Сов. энцикл., 1984. -- 1216 с.
\vspace{0.75ex}

7.
{\it Мищенко А.С., Фоменко А.Т.}
Курс дифференциальной геометрии и топологии. -- М.: Изд-во Моск. ун-та, 1980. -- 439 с.
\vspace{0.75ex}

8.
{\it Дубровин Б.А., Новиков С.П., Фоменко А.Т.}
Современная геометрия: Методы и приложения. -- 2-е изд., перераб. 
-- М.: Наука, 1986. -- 760 с.
\vspace{0.75ex}

9.
{\it Ефимов Н.В.} Высшая геометрия. -- М.: Наука, 1978. -- 576 с.
\vspace{0.75ex}

10.
{\it Горбузов В.Н., Самодуров А.А.} 
Уравнение Дарбу и  его аналоги.  -- Гродно:  ГрГУ, 1985.  -- 94 с.
\vspace{0.75ex}

11.
{\it Сибирский К.С.}
Введение в алгебраическую теорию инвариантов дифференциальных
уравнений. -- Кишинёв: Штиинца, 1982. -- 168 с.
\vspace{0.75ex}

12.
{\it Сибирский К.С.} 
Алгебраические инварианты
дифференциальных уравнений и матриц.  -- Кишинёв:  Штиинца,
1976.  -- 268 с.
\vspace{0.75ex}

13.
{\it Дюлак Г.} О предельных циклах. -- М.: Наука, 1980. -- 160 с.
\vspace{0.75ex}

14.
{\it Матвеев Н.М.} Методы интегрирования обыкновенных
дифференциальных уравнений. -- СПб.: Изд-во <<Лань>>,
2003.  -- 832 с.
\vspace{0.75ex}

15. 
{\it Вулпе Н.И., Косташ С.И.}
Центроаффинно-инвариантные условия топологического различения
дифференциальной системы с кубическими нелинейностями вида Дарбу. --
Кишинев: Акад. наук Молдавской ССР, 1989. -- 55 c.
\vspace{0.75ex}

16.
{\it Songling Shi.}
A concrete example of the existence of four  limit cycles
for plane quadratic systems // Scintia Sinica. -- 1980. -- Vol. 23,  No.  2. -- P. 153 -- 158.
\vspace{0.75ex}

17.
{\it Амелькин В.В., Лукашевич Н.А., Садовский А.П.}
Нелинейные колебания в системах второго порядка.  -- Минск:
Изд-во БГУ, 1982.  -- 210 с.
\vspace{0.75ex}

18. 
{\it Горбузов В.Н., Тыщенко В.Ю.} 
Симметричность траекторий квадратичных систем второго порядка. В 2-х частях.
Ч. 1.  -- Гродно: ГрГУ, 1992. -- 96 с. 
\vspace{0.75ex}

19.
{\it Gorbuzov V.N.} 
Integral equivalence of multidimensional differential systems // Mathematics.
Dynamical Systems (arXiv: 0909.3220v1 [math.DS]. Cornell Univ., Ithaca, New York). --- 2009. --- P. 1 --- 45.
\vspace{0.75ex}

20.
{\it Горбузов В.Н.}
Интегралы дифференциальных систем: монография. -- Гродно: ГрГУ, 2006. -- 447 с. 
\vspace{0.75ex}

21.
{\it Андреев А.Ф.} 
Особые точки дифференциальных
уравнений. -- Минск: Вышэйшая школа, 1979. -- 136 с.
\vspace{0.75ex}

22.
{\it Горбузов В.Н.} 
Проективный атлас траекторий дифференциальных систем второго порядка // 
Веснiк Гродзенскага дзяржа\u{y}нага \u{y}нiверсiтэта. Сер. 2. -- 2011. -- \No\, 2(111). -- С. 15 -- 26.
\vspace{0.75ex}

23.
{\it Горбузов В.Н.} 
Траектории проективно приведенных дифференциальных систем // 
Веснiк Гродзенскага дзяржа\u{y}нага \u{y}нiверсiтэта. Сер. 2. -- 2012. -- \No\, 1(126). -- 
\linebreak
С. 39 -- 52.
\vspace{0.75ex}

24.
{\it Горбузов В.Н., Павлючик П.Б.} 
Линейные и разомкнутые предельные циклы дифференциальных систем  // 
Веснiк Гродзенскага дзяржа\u{y}нага \u{y}нiверсiтэта. Сер. 2. -- 2013. -- \No\, 3(159). -- С.  23 -- 32.

}
\end{document}